\documentclass[reqno, 10pt]{amsart}

\usepackage[usenames,dvipsnames]{color}
\usepackage[pdfauthor={Andrea Appel, Valerio Toledano Laredo}, pdftitle={Coxeter categories and quantum groups}]{hyperref}
\hypersetup{
    colorlinks=true,       			% false: boxed links; true: colored links
    linkcolor=blue,          			% color of internal links
    citecolor=blue,		 		% color of links to bibliography
    filecolor=blue,      				% color of file links
    urlcolor=blue           			% color of external links
}
\usepackage{amsmath, amssymb, empheq, stmaryrd, enumitem}
\usepackage{tikz}
\usetikzlibrary{shapes,backgrounds}
\usetikzlibrary{arrows,decorations.pathmorphing,positioning,fit,petri}
\usepackage[all,cmtip,knot]{xy}
\xyoption{arc}
\usepackage{multirow}

\newcommand{\Omit}[1]{}

\newtheorem*{theorem}{Theorem}
\newtheorem*{proposition}{Proposition}
\newtheorem*{corollary}{Corollary}
\newtheorem*{lemma}{Lemma}

\theoremstyle{definition}
\newtheorem*{definition}{Definition}

\numberwithin{equation}{section}

\newenvironment{pf}{\paragraph{{\sc Proof}}}{\qed\par\medskip}

\newcommand{\remark}{{\bf Remark.\ }}
\newcommand{\remarks}{{\bf Remarks.\ }}

\newcommand {\IN}{\mathbb{N}}

\newcommand {\IQ}{\mathbb{Q}}
\newcommand {\IZ}{\mathbb{Z}}

\newcommand {\bfI}{\mathbf I}

\newcommand {\B}{\mathcal B}
\newcommand {\C}{\mathcal C}
\newcommand {\D}{\mathcal D}
\newcommand {\E}{\mathcal E}
\newcommand {\F}{{\mathcal F}}
\newcommand {\G}{{\mathcal G}}
\renewcommand {\H}{{\mathcal H}}

\newcommand {\K}{\mathcal K}

\newcommand {\M}{\mathcal M}
\newcommand {\N}{\mathcal N}
\renewcommand {\O}{\mathcal O}
\renewcommand {\P}{\mathcal P}
\newcommand {\Q}{\mathcal Q}
\newcommand {\sQ}{\sfs\Q}

\newcommand {\T}{\mathcal T}
\newcommand {\U}{\mathcal U}
\newcommand {\V}{\mathcal V}

\newcommand {\W}{\mathcal W}

\renewcommand {\a}{\mathfrak a}
\renewcommand {\b}{\mathfrak b}
\renewcommand {\c}{\mathfrak c}
\renewcommand {\d}{\mathfrak d}

\newcommand {\g}{\mathfrak{g}}
\newcommand {\h}{\mathfrak h}
\newcommand {\olh}{\ol{\h}}
\newcommand {\wth}{\wt{\h}}
\newcommand {\wtg}{\wt{\g}}
\newcommand {\wtI}{\wt{I}}
\renewcommand {\ll}{\mathfrak l}

\newcommand {\n}{\mathfrak n}
\newcommand {\p}{\mathfrak p}

\renewcommand {\SS}{\mathfrak{S}}
\newcommand{\fkA}{\mathfrak{A}}
\newcommand{\fkB}{\mathfrak{B}}

\newcommand{\hA}{\fkA}
\newcommand{\hC}{\fkB}

\newcommand{\Cox}[2]{\mathbb{DY}_{#1}^{#2}} %Drinfeld Yetter
\newcommand{\OCox}[2]{\mathbb{O}_{#1}^{#2}} %category O
\newcommand{\Hcox}[1]{\mathbb{H}_{#1}} %Coxeter equivalence
\newcommand{\pCox}[1]{\mathfrak{C}_{#1}} %symbol for a PROPic pre-Coxeter structure 

\newcommand {\sfA}{\mathsf{A}}
\newcommand {\sfB}{\mathsf{B}}

\newcommand {\sfD}{\mathsf{D}}

\newcommand {\ie}{{\it i.e., }}
\newcommand {\eg}{{\it e.g.}, }

\newcommand {\fd}{finite--dimensional }
\newcommand {\fg}{finitely--generated }
\newcommand {\lhs}{left--hand side }
\newcommand {\rhs}{right--hand side }

\newcommand {\wrt}{with respect to }

\newcommand {\ol}{\overline}
\newcommand {\wt}{\widetilde}
\newcommand {\wh}{\widehat}

\newcommand {\mns}{maximal nested set }
\newcommand {\mnss}{maximal nested sets }

\newcommand {\DCP}{De Concini--Procesi }
\newcommand {\DrJ}{Drinfeld--Jimbo }

\newcommand {\nEK}{Etingof--Kazhdan }

\newcommand {\KM}{Kac--Moody }
\newcommand {\KMA}{Kac--Moody algebra }
\newcommand {\KMAs}{Kac--Moody algebras }

\newcommand {\qW}{quantum Weyl group }

\def\iip#1#2{\langle{#1},{#2}\rangle}

\def\slantfrac#1#2{\kern.2em\raise.2em\hbox{$#1$}\kern-.2em\left/\lower.4em\hbox{$#2$}\kern-.2em\right.}
\def\backslantfrac#1#2{\kern.2em\lower.4em\hbox{$#1$}\kern-.3em\left\backslash\raise.4em\hbox{$#2$}\right.}

\DeclareMathOperator{\Hom}{Hom}
\DeclareMathOperator{\End}{End}
\DeclareMathOperator{\Ker}{Ker}

\DeclareMathOperator{\Res}{Res}
\DeclareMathOperator{\id}{id}

\DeclareMathOperator{\rank}{rank}

\DeclareMathOperator{\Rep}{Rep}

\DeclareMathOperator{\vect}{Vect}
\DeclareMathOperator{\ad}{ad}

\newcommand {\fml}{[\negthinspace[\hbar]\negthinspace]}

\newcommand {\cor}[1]{\alpha_{#1}^{\vee}}	
\newcommand {\hcor}[1]{\alpha^\vee_{#1}}

\newcommand {\cow}[1]{\lambda^{\vee}_{#1}}
\newcommand {\odots}[1]{#1\cdots #1}

\newcommand{\ten}{\otimes}

\newcommand {\ul}[1]{\underline{#1}}

\newcommand{\sint}{{\scs\operatorname{int}}}

\newcommand {\Uhg}{U_{\hbar}\g}

\newcommand {\reg}{_{\operatorname{reg}}}
\newcommand {\hreg}{\h\reg}

\newcommand {\aand}{\qquad\text{and}\qquad}

\newcommand{\DCPA}[2]{\Upsilon_{#1#2}}
\newcommand{\oDCPA}[2]{\ol{\Upsilon}_{#1#2}}

\newcommand{\Mns}[1]{\sfMns(#1)}
\newcommand{\Mnsr}[2]{\sfMns_{#2}(#1)}
\newcommand{\sfMns}{\mathsf{Mns}}

\renewcommand{\1}{\mathbf{1}}

\renewcommand{\DJ}{U_{\hbar}}

\newcommand{\Eq}[1]{\E_{#1}}%equicontinuous
\newcommand {\Db}{\mathfrak{D}}
\newcommand{\DrY}[1]{\mathsf{DY}_{#1}}
\newcommand{\hDrY}[2]{\mathsf{DY}_{#1}^{#2}}

\newcommand{\ff}{\mathsf{f}}

\newcommand{\gb}{\g_{\b}}
\newcommand{\ga}{\g_{\a}}

\newcommand {\<}{\langle}
\renewcommand {\>}{\rangle}

\newcommand {\sfk}{\mathsf{k}}

\newcommand {\kvect}{\vect_{\sfk}}

\newcommand{\PROP}{{\sf PROP}}

\newcommand{\preLA}{{\sf LA}}

\newcommand{\HAk}{\sf{HA}(\sfk)}
\newcommand{\sHAk}{\sf{sHA}(\sfk)}
\newcommand{\LBA}{\ul{\sf LBA}}

\newcommand{\DY}{\ul{\sf DY}}

\newcommand{\sfEnd}[1]{\mathsf{End}\left(#1\right)}

\newcommand{\sfAd}[1]{\mathsf{Ad}(#1)}
\newcommand{\sfAut}[1]{\mathsf{Aut}(#1)}

\newcommand{\Cat}{\mathsf{Cat}}

\newcommand{\sfad}{\mathsf{ad}}

\newcommand{\ID}{\mathfrak{P}(D)}

\newcommand{\NsID}{\mathfrak{Ns}(D)}
\newcommand{\PD}{\P(D)}
\newcommand{\PPD}{\sfP_D}

\newcommand{\GCM}[1]{\mathsf{#1}}

\newcommand{\fcw}{\lambda^{\vee}}

\renewcommand{\LBA}{\ul{\sf LBA}}
\renewcommand{\PROP}{\sf PROP}

\renewcommand{\DY}[1]{\ul{\sf DY}_{#1}}

\newcommand{\VDY}[1]{\ul{\mathsf{V}}_{#1}} %MODULES IN DY
\newcommand{\VCDY}[1]{\VDY{#1}} %MODULES IN CDY
\newcommand{\ACDY}[1]{[#1]} % CYCLIC OBJ IN CDY

\newcommand{\pEnd}[2]{{\mathsf{End}}_{#1}\left(#2\right)}

\newcommand {\sfQ}{\mathsf Q}
\renewcommand{\int}{^{\scriptscriptstyle{\operatorname{int}}}}

\newcommand {\sfa}{\mathsf a}

\newcommand {\sfP}{{\mathsf P}}

\newcommand{\sfi}{\mathsf{i}}

\newcommand{\gup}{\b_+}
\newcommand{\gum}{\b_-}
\newcommand{\gupm}{\b_{\pm}}

\newcommand{\gdpm}{\a_{\pm}}

\newcommand{\EK}{U_{\hbar}}
\newcommand{\EKeq}[1]{{H}_{#1}}%{\wt{F}\EEK_{#1}} %tilda removed and renamed H

\newcommand{\fun}[1]{F_{#1}}

\newcommand{\qcc}{Coxeter category }

\newcommand{\XX}{\mathfrak{X}}
\newcommand{\FPC}[3]{\wt{\P}}
\newcommand{\OFPC}[3]{\P}

\newcommand{\Ns}[1]{\mathsf{Ns}(#1)}
\newcommand{\Nsr}[2]{\mathsf{Ns}_{#2}(#1)}

\newcommand{\solved}[1]{}

\newcommand{\cc}[1]{\mathsf{conn}(#1)}

\newcommand{\scs}{\scriptscriptstyle}
\newcommand{\scsop}[1]{\scriptscriptstyle\operatorname{#1}}

\newcommand{\sd}[2]{{#2}{#1}}

\newcommand{\ulm}{\ul{m}}
\newcommand{\BDm}{\B_D^{\ulm}}
\newcommand{\BBm}{\B_B^{\ulm}}
\newcommand{\BBpm}{\B_{B'}^{\ulm}}
\newcommand{\sfS}{\mathsf{S}}

\newcommand{\drc}[1]{\delta_{#1}}
\newcommand{\ekm}[1]{\ol{#1}}
\newcommand{\olg}{\ol{\g}}
\newcommand{\olb}{\ol{\b}}

\newcommand{\bb}[1]{\ekm{\b}_{#1}}
\newcommand{\sfR}{\mathsf{R}}

\newcommand{\hext}[1]{{#1}{\fml}}

\newcommand{\sfK}{\mathsf{K}}

\newcommand{\DYt}{Drinfeld--Yetter }

\newcommand{\QQUE}{\mathsf{QUE}}

\newcommand{\adm}{\scs\operatorname{adm}}
\newcommand{\aDrY}[1]{\mathsf{DY}^{\scs\operatorname{adm}}_{#1}}

\newcommand {\res}{^{\scs\operatorname{res}}}
\newcommand{\grb}{\g_{\b}\res}
\newcommand{\resDb}[1]{(\Db{#1})^{\scs\operatorname{res}}}
\newcommand{\resqD}[1]{(D#1)^{\scs\operatorname{res}}}
\newcommand{\Kar}[1]{\mathsf{Kar}(#1)}
\newcommand{\cKar}[1]{\ul{#1}}%{\wh{\ul{#1}}}%{{#1}^{\scriptscriptstyle{\mathsf{kar,ind}}}}

\newcommand{\preLBA}{\mathsf{LBA}}

\newcommand{\mLBA}[1]{{\LBA}_{#1}}
\newcommand{\dLBA}[1]{{\LBA}_{#1}}

\newcommand{\kLBA}{\preLBA(\sfk)}
\newcommand{\ksLBA}{\mathsf{sLBA}(\sfk)}
\newcommand{\ksMT}{\mathsf{sMT}(\sfk)}

\newcommand{\DYrho}[2]{\rho^{#2}_{#1}} % realisation map for DY univ algebra

\newcommand{\MDY}[2]{\ul{\mathsf{DY}}_{#1}^{#2}}
\newcommand{\MUA}[2]{\mathfrak{U}_{#1}^{#2}}
\newcommand{\CMUA}[2]{\wh{\mathfrak{U}}_{#1}^{#2}}

\newcommand{\Diag}{\operatorname{Diag}}

\newcommand{\dmp}[1]{\theta_{#1}}

\newcommand{\RDY}[2]{\ul{\mathsf{DY}}^{#2}_{#1}}
\newcommand{\RDYUA}[2]{\mathfrak{U}_{#1}^{#2}}
\newcommand{\CRDYUA}[2]{\wh{\mathfrak{U}}_{#1}^{#2}}

\newcommand{\vrtx}[1]{\{i\}}

\newcommand {\khvect}{\mathsf{Vect}_{\hext{\sfk}}}

\newcommand{\DYA}[2]{\U_{#1}^{#2}}  %completion wrt DY --> Vect
\newcommand{\hDYA}[2]{\wh{\U}_{#1}^{#2}}

\newcommand{\Fun}{\mathsf{Fun}}
\newcommand{\bfb}{\mathbf{b}}

\newcommand{\gaugebis}[2]{(#1)^{-1}_{12}\cdot#2\cdot(#1)_1\cdot(#1)_2}
\newcommand{\gaugetris}[2]{(#1)_{12}\cdot#2\cdot(#1)_1^{-1}\cdot(#1)_2^{-1}}

\newcommand{\twistp}[2]{(#1)_{23}^{-1}\cdot(#1)^{-1}_{1,23}\cdot #2\cdot(#1)_{12,3}\cdot(#1)_{12}}

\newcommand{\resped}{_{\a,\b}}
\newcommand{\Uha}{\Q(\a)}
\newcommand{\Uhb}{\Q(\b)}
\newcommand{\Uhc}{\Q(\c)}

\newcommand {\diagr}[2]{\mathsf{\mathsf{Diagr}}_{#1}(#2)}

\newcommand {\sff}{\mathsf{f}}
\newcommand {\sfs}{\mathsf{s}}

\newcommand{\dsg}{\text{D}}

\newcommand{\KQUE}{\mathsf{QUE}(\sfK)}

\newcommand{\KsQUE}{\mathsf{sQUE}(\sfK)}

\newcommand {\SC}[2]{{S}_{#1}^{#2}}

\newcommand{\norm}[1]{#1}

\newcommand{\defDY}[1]{\hDrY{#1}{\hbar}}
\newcommand{\sLBA}{\mathsf{sLBA}}

\newcommand{\asso}[3]{a_{#1#3}^{{\color{white}#1}#2}}
\newcommand{\redasso}[2]{\mathsf{a}_{#2}^{#1}} %reduced vertical joins
\newcommand{\elemasso}[1]{\mathsf{b}_{#1}} %elementary vertical joins
\newcommand{\alpho}[3]{\alpha_{#3#1}^{#2}}
\newcommand{\oalpho}[3]{{\sfa}_{#3#2#1}}

\newcommand {\Co}{\C}%{\C^{(1)}}

\newcommand{\cdg}{Cartan }

\newcommand{\CB}{{\mathcal C}_B}
\newcommand{\CBp}{{\mathcal C}_{B'}}
\newcommand{\CBpp}{{\mathcal C}_{B''}}

\newcommand{\HomA}{\Hom_\sfA}

\newcommand{\Lambdav}{\Lambda^\vee}

\newcommand{\olPi}{\ol{\Pi}}
\newcommand{\olPiv}{\ol{\Pi}^\vee}
\newcommand{\wtV}{\wt{V}}
\newcommand{\Piv}{\Pi^\vee}
\newcommand{\real}{(V,\Pi,\Piv)}

\newcommand {\two}{^{(2)}}
\newcommand {\gtwo}{\g\two}
\newcommand {\btwo}{\b\two}
\newcommand {\olgtwo}{\olg\two}
\newcommand {\olbtwo}{\olb\two}

\newcommand{\zolh}{\olh^{c}}

\newcommand {\zee}{\mathfrak z}

\newcommand {\LQ}{\operatorname{\mathsf{Lie}}_\sfQ}

\newcommand{\BonC}{\triangleright}
\newcommand{\ConB}{\triangleleft}
\newcommand{\dcp}[2]{#1\negthinspace\BonC\negthinspace\negthinspace\ConB\negthinspace#2}

\newcommand{\bfC}{\mathbf{C}}
\newcommand{\Ch}[1]{\mathsf{Ch}(#1)}

\newcommand{\hextsub}[1]{{#1}_{\hbar}}
\newcommand{\hextsup}[1]{{#1}^{\hbar}}

\newcommand{\sQUE}{\mathsf{sQUE}}

\newcommand{\DUA}[2]{\MUA{}{#2}}%{\mathbb{U}_{#1}^{#2}} %
\newcommand{\CDUA}[2]{\CMUA{}{#2}}%{\wh{\mathbb{U}}_{#1}^{#2}} 

\newcommand{\bbT}{\mathbb{T}}
\newcommand{\bfg}{\mathbf{g}}

\newcommand{\ub}{[\b]}
\newcommand{\ua}{[\a]}
\newcommand{\uaa}[1]{[\a_{#1}]}
\newcommand{\ubb}[1]{[\b_{#1}]}
\newcommand{\ucc}[1]{[\c_{#1}]}

\newcommand{\pair}[2]{#1\hookrightarrow #2}

\newcommand{\univb}{[\b]}%{\b_{\scsop{univ}}}
\newcommand{\gstr}{{\scsop{\Upsilon\text{-}str}}}
\newcommand{\astr}{{\scsop{\sfa\text{-}str}}}

\newcommand{\LBAo}{\LBA_{\mathsf{o}}}

\newcommand{\vertice}[1]{\stackrel{#1}{\circ}}

\newcommand{\twF}[1]{K_{#1}}

\newcommand{\sfp}{\mathsf{p}}
\newcommand{\intDY}[1]{\hDrY{#1}{\hbar,\mathsf{int}}}

\title%[Coxeter categories and quantum groups]
{Coxeter categories and quantum groups}
\author[A. Appel]{Andrea Appel}
\address{School of Mathematics,
University of Edinburgh,
James Clerk Maxwell Building, 
Peter Guthrie Tait Road,
Edinburgh, EH9 3FD, UK}
\email{andrea.appel@ed.ac.uk}
\author[V. Toledano Laredo]{Valerio Toledano Laredo}
\address{Department of Mathematics,
Northeastern University,
360 Huntington Avenue,
Boston MA 02115}
\email{V.ToledanoLaredo@neu.edu}
\thanks{The first author was supported in part through the
ERC STG Grant 637618 and the NSF grant DMS--1255334, 
the second author through the NSF grants DMS--1505305
and DMS--1802412.}

 %quantum groups and deformations; Lie bialgebras; quantum groups and related algebraic methods

\begin{document}

\begin{abstract}
We define the notion of braided Coxeter category, which is informally a
monoidal category carrying compatible, commuting actions of a generalised
braid group $B_W$ and Artin's braid groups $B_n$ on the tensor powers
of its objects. The data which defines the action of $B_W$ bears a formal
similarity to the associativity constraints in a monoidal category, but is related
to the coherence of a family of fiber functors. We show that the quantum
Weyl group operators of a quantised Kac--Moody algebra $\DJ{\g}$,
together with the universal $R$--matrices of its Levi subalgebras, give
rise to a braided \qcc structure on integrable, category $\O$--modules
for $\DJ{\g}$. By relying on the 2--categorical extension of \nEK quantisation
obtained in \cite{ATL1-1}, we then prove that this structure can be transferred
to integrable, category $\O$--representations of $\g$. These results are
used in \cite{ATL3} to give a monodromic description of the quantum
Weyl group operators of $\DJ{\g}$, which extends the one obtained
by the second author for a semisimple Lie algebra.
\end{abstract}
  
\maketitle
\setcounter{tocdepth}{1}
\tableofcontents

\newpage

\section{Introduction}
%===============

\subsection{}%problem and result
%--------------

This is the first of a series of three papers the aim of which is to extend
the description of the monodromy of the rational Casimir connection of
a complex semisimple Lie algebra in terms of quantum Weyl groups
obtained in \cite{vtl-2,vtl-3,vtl-4,vtl-6} to the case of an arbitrary symmetrisable
\KM algebra $\g$.

The method we follow is close, in spirit at least, to that of \cite{vtl-4}. It
relies on the notion of {\it braided Coxeter category}, the definition of
which is the first main contribution of the present article. Informally, 
such a category is a monoidal category carrying compatible, commuting actions
of a given generalised braid group and Artin's braid groups on the tensor
products of its objects. This structure arises for example on the category
$\O\int_{\DJ{\g}}$ of integrable, highest weight representations of the
quantum group $\DJ{\g}$, from the quantum Weyl group operators of
$\DJ{\g}$ and the $R$--matrices of its Levi subalgebras.

A rigidity result, proved in the second paper of this series \cite{ATL2},
shows that there is at most one braided Coxeter structure with prescribed
restriction functors, $R$--matrices and local monodromies on the category
$\O\int_\g$ of integrable, highest weight representations of $\g$. 
It follows that the generalised braid group actions arising from quantum
Weyl groups and the monodromy of the Casimir connection \cite{ATL3}
are equivalent, provided the braided Coxeter structure underlying the
former can be transferred from $\O\int_{\DJ{\g}}$ to $\O\int_\g$.
This result is the second main contribution of this article.

\subsection{} % outline of intro

In the rest of the introduction, we outline the definition of a Coxeter
and of a braided Coxeter category. We then focus on two main sources
of examples. The first arises from {\it diagrammatic} Lie bialgebras,
and generalises category $\O$ for a symmetrisable \KMA $\g$. The
second arises from diagrammatic Hopf algebras, and generalises
category $\O$ for the quantum group $\Uhg$. Finally, we explain
how the Etingof--Kazhdan quantisation of Lie bialgebras \cite{ek-1,ek-2},
and its 2--categorical extension recently obtained in \cite{ATL1-1}
give rise to an equivalence between a canonical deformation of
the first class of examples and the second, thus yielding the transfer
theorem alluded to above.

\subsection{}\label{ss:qc vs br} % Quasi-Coxeter vs Braided
%--------------

The definition of a Coxeter category bears some formal similarity to that
of a braided monoidal category, with Artin's braid groups $\{\sfB_n\}_{n\geqslant 2}$
replaced by a given generalised braid group $\sfB_{W}$ of Coxeter type
$W$. If $\C$ is braided monoidal then, for any object $V\in\C$ and
$n\geqslant 2$, there is an action
\[\rho_{b}:\sfB_n\to\sfAut{V^{\ten n}_{b}}\]
for any bracketing $b$ on the non--associative monomial $x_1\cdots x_n
$.\footnote{
The notation $V^{\ten n}_{b}$ indicates that $n$ copies of $V$ have been
tensored together according to $b$. For example, if $b=(x_1x_2)x_3$,
$V_b^{\ten 3}=((V\ten V)\ten V)$.} The choice of $b$ is in a sense immaterial
since, for any two bracketings $b,b'$, the associativity constraint $\Phi_{b'b}:
V_{b}^{\ten n}\to V_{b'}^{\ten n}$ of $\C$ intertwines the corresponding
actions of $\sfB_n$. Similarly, if $V$ is an object in a Coxeter category 
$\Q$, there is an action
\[
\lambda_{\F}:\sfB_W\to\sfAut{V_{\F}}
\]
which depends on a discrete choice $\F$.
%the choice of a \emph{$W$--bracketing} $\F$.
Moreover, for any two such choices $\F,\G$, there is an isomorphism $\DCPA
{\G}{\F}:V_{\F}\to V_{\G}$ which intertwines the actions of $\sfB_W$.

\subsection{} %Nested Sets
%--------------

The relevant discrete choice is that of a {\it maximal nested set} $\F$
on the Dynkin diagram $D$ of $W$, a combinatorial notion introduced
by De Concini--Procesi \cite{DCP2} which generalises that of a bracketing
on $x_1\cdots x_n$ when $W$ is the symmetric group $\SS_n$ with
diagram $\sfA_{n-1}$. Specifically, to a pair of parentheses $x_1\cdots
(x_i\cdots x_j)\cdots x_n$, one can  associate the connected subdiagram
of $\sfA_{n-1}$ with nodes $\{i,\dots, j-1\}$. Under this identification, a
(complete) bracketing on $x_1\cdots x_n$ corresponds to a (maximal)
%nested set, \ie a (maximal)
collection $\F=\{B\}$ of connected subdiagrams
of $\sfA_{n-1}$ which are pairwise {\it compatible}, \ie such that for any
$B,B'\in\F$, one has
\[B\subseteq B',
\qquad
B'\subseteq B
\qquad\text{or}\qquad
B\perp B'\]
where the latter condition means that $B$ and $B'$ have no vertices in
common, and that no edge in $\sfA_{n-1}$ connects a vertex in $B$ to
one in $B'$. Such a collection is called a nested set on $\sfA_{n-1}$,
and may be defined for any Coxeter group, and in fact any diagram $D
$.\footnote{We use the term diagram to denote an undirected graph,
with mo mutiple edges or loops.}

\subsection{}\label{ss:ext int} 
%--------------

Despite the above formal similarities, there is one significant difference
between braided monoidal categories and Coxeter categories. 
In a Coxeter category $\C$, the braid group $\sfB_W$ does not act by
morphism in $\C$.
For example, the quantum Weyl group operators do not commute with
the action of $\DJ\g$. Thus, $\sfB_W$ does not act through morphism
of $\C=\Rep\DJ\g$, but rather automorphisms of the forgetful functor
$F:\Rep\DJ\g\to\vect$.
This is a general feature: in a Coxeter category $\C$, the braid group
$\sfB_W$ acts by automorphisms of a fiber functor from $\C$ to a base
category $\C_\emptyset$. In fact, $\C$ is endowed with a {\it collection}
of such functors $F_{\F}:\C\to\C_\emptyset$, labelled by the \mnss on
$D$. For any such $\F$, and object $V\in\C$, there is a homomorphism
\[
\lambda_{\F}:\sfB_W\to\mathsf{Aut}_{\C_\emptyset}(V_{\F})
\]
where $V_{\F}=F_{\F}(V)$. Further, for any $\F,\G$, there is an isomorphism
of fiber functors $\Upsilon_{\G\F}:F_\F\Rightarrow F_\G$ which give rise
to an identification of $\sfB_W$--modules $V_{\F}\to V_{\G}$. 

\subsection{} % intermediate categories and relative mns
%--------------

In a (braided) Coxeter category, the fiber functors $F_\F$ are additionally
required to factorise vertically in the following sense. For any subdiagram
$B\subseteq D$, one is given a (braided monoidal) category $\C_B$. In
the case of quantum groups, $\C_B$ consists of representations of the
subalgebra $\Uhg_B$ of $\Uhg$ with generators labelled by the vertices
of $B$. Moreover, for any pair of subdiagrams $B'\subseteq B$, there is a
family of (monoidal) functors $F_\F:\C_B\to\C_{B'}$, which can be thought
of as restriction functors. These are labelled by \mnss on $B$ relative to $B'$,
that is nested sets on $B$ whose elements are compatible with, but not strictly
contained in $B'$.\footnote{If $D=\sfA_{n-1}$, $B=D$ and $B'$ corresponds to
the pair of parentheses $x_1\cdots x_{i-1}\cdot (x_i\cdots x_j)\cdot x_{j+1}
\cdots x_n$, a \mns on $B$ relative to $B'$ consists of a complete
bracketing of the monomial $x_1\cdots x_{i-1}\cdot x_{ij} \cdot x_{j+1}\cdots x_n$.}
As in the absolute case $B'=\emptyset\subset D=B$ discussed in \ref{ss:ext int},
the functors $F_\F$ are related by a transitive family of isomorphisms
$\Upsilon_{\G\F}:F_\F\Rightarrow F_\G$. Finally, for any triple of subdiagrams $B''\subseteq B'\subseteq B$, 
a \mns $\F$ on $B$ relative to $B'$ and a \mns $\F'$ on $B'$ relative to $B''$,
the composition $F_{\F'}\circ F_\F:\C_B\to\C_{B''}$ is isomorphic to 
$F_{\F'\cup\F''}$ via a coherent isomorphism.

\subsection{}\label{ss:cox cat} %Definition of Coxeter category
%--------------

Let now $\left(D,\{m_{ij}\}\right)$ be a labelled diagram with set of vertices $\bfI$, and $\sfB_{D}$
the generalised braid group corresponding to $D$\footnote{A labelling on $D$
is the additional data of integers $m_{ij}\in\{2,\dots,\infty\}$ for any two $i\neq j
\in\bfI$ such that $m_{ij}=m_{ji}$ and $m_{ij}=2$ if $i\perp j$.} \ie
\[
\sfB_{ D}=
\<S_i\>_{i\in\bfI}/
_{\underbrace{S_i S_j\ S_i\cdots}_{m_{ij}}
=
\underbrace{S_j S_i S_j\cdots}_{m_{ij}}}
\]
For any pair of subdiagrams $B'\subseteq B$ of $D$, we denote by $\Mns{B,B'}$
the collection of maximal nested sets on $B$ relative to $B'$. 

A braided Coxeter category of type $ D$ consists of the following five pieces
of data.
\begin{enumerate}
%%%
\item \emph{Diagrammatic categories.}
For any subdiagram $B\subseteq D$, a braided monoidal category $\C_B$. 
%%%
\item \emph{Restriction functors.} 
For any pair of subdiagrams $B'\subseteq B$, and maximal nested
set $\F$ on $B$ relative to $B'$, a monoidal functor $F_\F:\C_B\to\C_{B'}$.
\item \emph{De Concini--Procesi associators.} 
For any $B'\subseteq B$, and \mnss $\F,\G$ on $B$ relative to $B'$, an
isomorphism of monoidal functors
\[
\DCPA{\G}{\F}:F_{\F}\Rightarrow F_{\G}
\] 
such that $\Upsilon_{\H\G}\cdot\Upsilon_{\G\F}=\Upsilon_{\H\F}$ for any
$\F,\G,\H\in\Mns{B,B'}$. 
\item \emph{Vertical joins.} 
For any triple of subdiagrams $B''\subseteq B'\subseteq B$, and \mnss
$\F\in\Mns{B,B'}$, $\F'\in\Mns{B',B''}$, an isomorphism of monoidal
functors $\redasso{\F}{\F'}:F_{\F'}\circ F_{\F}\Rightarrow F_{\F'\cup\F}$,
such that the following properties hold 
\begin{itemize}
\item[(a)] \emph{Vertical Factorisation}. For any subdiagrams 
$B''\subseteq B'\subseteq B$, and \mnss $\F,\G\in\Mns{B,B'}$
and $\F',\G'\in\Mns{B',B''}$
\[
\xy
(0,50)*++{\C_B}="C1";
(0,25)*++{\C_{B'}}="C2";
(0,0)*++{\C_{B''}}="C3";
%
%FUNCTORS
{\ar@/_2pc/@{->}_{\fun{\F}} "C1";"C2"};
{\ar@/^2pc/@{->}^{\fun{\G}} "C1";"C2"};
{\ar@/_2pc/@{->}_{\fun{\F'}} "C2";"C3"};
{\ar@/^2pc/@{->}^{\fun{\G'}} "C2";"C3"};
{\ar@/_6pc/@{->}_{\fun{\F'\cup\F}} "C1";"C3"};
{\ar@/^6pc/@{->}^{\fun{\G'\cup\G}} "C1";"C3"};
%LABELS
(-10,12.5)*+++{}="F";(10,12.5)*+++{}="G";
(-10,37.5)*+++{}="F'";(10,37.5)*+++{}="G'";
(-25,25)*+++{}="F+F";(25,25)*+++{}="G+G";
%DCPA
{\ar@{<=}|{\DCPA{\F'}{\G'}}"F";"G"};
{\ar@{<=}|{\DCPA{\F}{\G}}"F'";"G'"};
{\ar@/^1pc/@{<:}|(.8){\DCPA{\F'\cup\F\,}{\G'\cup\G}}"F+F";"G+G"};
\endxy
\]
where the triangular faces are given by $\redasso{\F}{\F'}$ and $\redasso{\G}{\G'}$.
\item[(b)] \emph{Vertical associativity.}
For any $B'''\subseteq B''\subseteq B'\subseteq B$, and \mnss $\F\in
\Mns{B,B'},\F'\in\Mns{B',B''},\F''\in\Mns{B'',B'''}$, the following equality
holds
\[\redasso{\F\cup\F'}{\F''}\circ\redasso{\F}{\F'}=
\redasso{\F}{\F'\cup\F''}\circ\redasso{\F'}{\F''}\]
as natural transformations $F_{\F''}\circ F_{\F'}\circ F_{\F}\Rightarrow
F_{\F\cup\F'\cup\F''}$. 
\end{itemize}
%%%
\item \emph{Local monodromies.}
For any vertex $i\in D$, an element $S_i^{\C}\in\mathsf{Aut}(F_{\{i\}})$,
where $\{i\}$ is the unique element in $\Mns{i,\emptyset}$, satisfying
\begin{itemize}
% braid relations
\item[(a)] \emph{Braid relations.}
For any $B\subseteq D$, $i\neq j\in B$ and maximal nested sets $\F,\G$ on $B$ 
with $i\in\F,j\in\G$,  the following holds in $\sfAut{F_\G}$
\[
\underbrace{\mathsf{Ad}\left(\DCPA{\G}{\F}\right)(S_i^\C)\cdot S_j^\C
\cdot \mathsf{Ad}\left(\DCPA{\G}{\F}\right)(S_i^\C)\cdots}_{m_{ij}}=
\underbrace{S_j^\C\cdot\mathsf{Ad}\left(\DCPA{\G}{\F}\right)(S_i^\C)\cdot S_j^\C\cdots}_{m_{ij}}
\]
where $S_i^\C$ is regarded as an automorphism of $F_\F$ via the factorisation
$\redasso{\F\setminus\{i\}}{\{i\}}:F_{\{i\}}\circ F_{\F\setminus \{i\}}\Rightarrow
F_\F$, and $S_j^\C$ is similarly regarded as an automorphism of $F_\G$.
% coproduct identity
\item[(b)] \emph{Coproduct identity.} 
For any vertex $i\in D$, and $V,W \in\C_i$ the following diagram in $\C_
\emptyset$ is commutative 
\begin{equation}\label{eq:copro diagr}
\xymatrix{
F_{\{i\}}(V)\ten F_{\{i\}}(W) \ar[d]_{J_{\{i\}}^{V,W}} \ar[r]^{S^{\C}_i\ten S^{\C}_i} & F_{\{i\}}(V)\ten F_{\{i\}}(W)  \ar[r]^{c_{\emptyset}} & F_{\{i\}}(W)\ten F_{\{i\}}(V)
\ar[d]^{J_{\{i\}}^{W,V}}\\
F_{\{i\}}(V\ten W) \ar[r]_{S_i^{\C}} & F_{\{i\}}(V\ten W)  \ar[r]_{F_{\{i\}}(c_i)} & F_{\{i\}}(W\ten V)  
}
\end{equation}
where $J_{\{i\}}$ is the tensor structure on $F_{\{i\}}:\C_i\to\C_\emptyset$,
and $c_i, c_{\emptyset}$ are the opposite braidings in $\C_i$ and $\C_
{\emptyset}$, respectively.\footnote{In a braided monoidal category with
braiding $\beta$, the opposite braiding is $\beta^{\scs\operatorname{op}}
_{X,Y}:=\beta_{Y,X}^{-1}$.}
\end{itemize}
\end{enumerate}

\noindent \remarks
\begin{enumerate}
\item The diagram \eqref{eq:copro diagr} codifies the coproduct identity
$\Delta(S_i)=R_i^{21}\cdot S_i\otimes S_i$ satisfied by quantum Weyl group
elements \cite[Prop. 5.3.4]{l}. It relates the failure of $F_{\{i\}}$ to be a braided
monoidal functor and that of $S^\C_i$ to be a monoidal isomorphism.
That is, if \eqref{eq:copro diagr} is commutative, then $S^\C_i$ is
monoidal if and only if $F_{\{i\}}$ is braided.
\item As mentioned in \ref{ss:qc vs br}, the definition of a Coxeter category
$\C$ is tailored to produce a family of equivalent representations of $\sfB_
{ D}$. Specifically, there is a collection of homomorphisms $\lambda_{\F}:
\sfB_{ D}\to\sfAut{F_\F}$, labelled by the maximal nested sets on $ D$,
which is uniquely determined by
\begin{itemize}
\item $\lambda_\F(S_i)=S_i^\C$ if $i\in\F$.
\item $\lambda_\G=\sfAd{\DCPA{\G}{\F}}\circ\lambda_\F$, for any $\F,\G\in\Mns{D}$.
\end{itemize}
\end{enumerate}

\subsection{}\label{ss:ex lbas}%LBA and DY
%--------------

An important class of braided {\it pre--}Coxeter categories, that is
structures satisfying the axioms (1)--(4) of \ref{ss:cox cat} but not
necessarily endowed with local monodromies, arises from split {\it
diagrammatic Lie bialgebras}. Recall first that a Lie bialgebra is a
triple $(\b,[\cdot,\cdot]_{\b},\delta_{\b})$, where $(\b,[\cdot,\cdot]_
{\b})$ is a Lie algebra and $(\b,\delta_{\b})$ a Lie coalgebra such
that the cobracket $\delta_{\b}$ and the bracket $[\cdot,\cdot]_{\b}$
satisfy an appropriate compatibility condition. 

A natural class of representations over a Lie bialgebra $\b$
is that of \emph{\DYt modules} \cite{ek-2}. Such a module is a triple $(V,\pi,\pi^*)$ 
such that $\pi:\b\otimes V\to V$ gives $V$ the structure of a left 
$\b$--module, $\pi^*:V\to\b\otimes V$ that of a right $\b$--comodule, 
and $\pi,\pi^*$ satisfy a compatibility condition. The latter is designed
so as to give rise to a representation of the Drinfeld double $\gb
=\b\oplus\b^*$ of $\b$, with $\phi\in\b^*$ acting on $V$ by $\phi
\otimes\id_V\circ\pi^*$.

If $\b$ is finite--dimensional, the symmetric monoidal category
$\hDrY{\b}{}$ of such modules coincides in fact with that of $\gb
$--modules, with the coaction of $\b$ on $V\in\Rep(\gb)$ given
by $\pi^*(v)=\sum_i b_i\otimes b^i v$, where $\{b_i\},\{b^i\}$ are
dual bases of $\b$ and $\b^*$. For an arbitrary $\b$, $\hDrY{\b}{}$
is isomorphic to the category $\E_{\gb}$ of {\it equicontinuous}
modules over $\gb$ \cite{ek-1}, that is those for which $\b^*$
acts locally finitely. This makes $\hDrY{\b}{}$ more convenient
to study than $\E_{\gb}$.

If $\a\stackrel{i}{\to}\b\stackrel{p}{\to}\a$ is a split embedding of
Lie bialgebras, there is a tensor restriction functor $\Res_{\a,\b}:
\DrY{\b}\to\DrY{\a}$ defined by
\[\Res_{\a,\b}(V,\pi_V,\pi^*_V)=(V,\pi_V\circ i\ten\id_V,p\ten\id_V\circ\pi^*_V)\]
Moreover, if $\a\hookrightarrow\b\hookrightarrow\c$ is a chain of split
embeddings, then $\Res_{\a,\b}\circ\Res_{\b,\c}=\Res_{\a,\c}$.
In terms of Drinfeld doubles, a split embedding gives rise to an isometric embedding
of Lie algebras $j=i\oplus p^t:\ga\to\gb$, and the functor $\Res_{\a,\b}$
corresponds to the pull--back functor $j^*$ from (equicontinuous) modules
over $\gb$ to those over $\ga$.

\subsection{}\label{intro:diag-lba}%diagrammatic LBA and pre-Coxeter 
%--------------

A {\it (split) diagrammatic} Lie (bi)algebra $\b$ over a diagram $D$
is the datum of a family of Lie (bi)algebras $\{\b_B\}_{B\subseteq D}$
labelled by the subdiagrams of $D$, together with (split) morphisms
$\b_{B'}\to\b_B$ for any $B'\subseteq B$. These are assumed to be
transitive under compositions $B''\subseteq B'\subseteq B$, and such
that if $B',B''\subset D$ are orthogonal subdiagrams, $\b_{B'\sqcup B''}$
is isomorphic to $\b_{B'}\oplus\b_{B''}$ as Lie (bi)algebras. 

If $\b$ is a split diagrammatic Lie bialgebra, there is a symmetric
pre--Coxeter category $\Cox{\b}{}$ defined as follows
\begin{enumerate}
	\item For any $B\subseteq D$, $\Cox{\b,B}{}$ is the symmetric monoidal category $\DrY{\b_B}$
	\item For any $B'\subseteq B$ and \mns $\F$ on $B$ relative to $B'$,
	the restriction functor $F_{\F}:\Cox{\b,B}{}\to\Cox{\b,B'}{}$ is given by $\Res_{\b_{B'},\b_{B}}$
	\item For any $B'\subseteq B$ and \mnss $\F,\G$ on $B$ relative to $B'$,
	the associator $\DCPA{\G}{\F}:F_\F\Rightarrow F_\G$ is the identity on $\Res_{\b_{B'},\b_{B}}$
	\item For any $B''\subseteq B'\subseteq B$, and \mnss $\F\in\Mns{B,B'}$,
	$\F'\in\Mns{B',B''}$, the vertical join $\redasso{\F}{\F'}:F_{\F'}\circ F_{\F}\Rightarrow F_{\F'\cup\F}$
	is the equality $\Res_{\b_{B''},\b_{B'}}\circ \Res_{\b_{B'},\b_{B}}=\Res_{\b_{B''},\b_{B}}$.
\end{enumerate}

A deformation of $\Cox{\b}{}$, where the restriction functors $F_\F$
and associators $\Phi_{\G\F}$ genuinely depend on the choice of
\mnss will be outlined in \ref{ss:quant diagr}.

\subsection{}\label{ss:ss g}%semisimple Lie algebras 
%--------------

Semisimple Lie algebras are basic examples of diagrammatic Lie
bialgebras. Specifically, let $\g$ be a complex semisimple Lie algebra,
with opposite Borel subalgebras $\b_\pm\subset\g$, Dynkin diagram
$D$, Serre generators $\{e_i,f_i,h_i\}_{i\in D}$, and standard Lie
bialgebra structure determined by $\b_\pm$ and an invariant inner
product on $\g$ (see \ref{ss:bil on g}). Then, $\g$ is a diagrammatic
Lie bialgebra where, for any $B\subseteq D$, $\g_B\subseteq\g$
is the subalgebra generated by $\{e_i,f_i,h_i\}_{i\in B}$.

The diagrammatic structure on $\g$ determines a split
diagrammatic one on $\b_\pm$ as follows. For any $B
\subseteq D$, let $\b_{\pm,B}=\b_\pm\cap\g_B$ be the
subalgebras generated by $\{h_i,e_i\}_{i\in B}$ and $
\{h_i,f_i\}_{i\in B}$ respectively. If $B'\subseteq B$, let $i_{\pm,BB'}:\b_{\pm,B'}\to\b_
{\pm,B}$ be the embedding, and regard its transpose
$i_{\pm,BB'}^t$ as a map $p_{\mp,B'B}:\b_{\mp,B}\to
\b_ {\mp,B'}$ via the identifications $\b_{\mp,C}\cong
\b_{\pm,C}^*$ given by the inner product. Then, $\{
i_{\pm,BB'},p_{\pm,B'B}\}$ give the required splitting
of $\b_\pm$.

Taking \DYt modules gives rise to a symmetric pre--Coxeter
category $\Cox{\b_\pm}{}$, as explained in \ref{intro:diag-lba}.
Moreover, the realisation of each $\g_B$ as a quotient of the
Drinfeld double of $\b_{\pm,B}$ gives rise to an embedding
of the pre--Coxeter category of $\g$--modules with standard
restriction functors to $\Cox{\b_\pm}{}$.

\subsection{}\label{ss:ext kma}%extended KM
%--------------

The above example does not immediately extend to the case
of a symmetrisable \KMA $\g$, however, since $\g$ need not
be diagrammatic (Sect. \ref{se:diagr kmas}). To remedy this,
we introduce the notion of an {\it extended} \KM algebra.

Fix an $|\bfI|\times|\bfI|$ matrix $\sfA$ with entries in a field $\sfk$.
The extended \KMA $\olg(\sfA)$ corresponding to $\sfA$ is the
quotient of the Lie algebra generated by $\{e_i,f_i,\hcor{i},\cow
{i}\}_{i\in\bfI}$, with relations $[\hcor{i},\hcor{j}]=0$, $[\fcw_i,\fcw
_j]=0$, $[\hcor{i},\fcw_j]=0$, 
\[[\hcor{i},e_j]=a_{ji}e_j,\qquad
[\hcor{i},f_j]=-a_{ji}f_j,\qquad
[\fcw_i,e_j]=\drc{ij}e_j ,
\qquad [\fcw_i,f_j]=-\drc{ij}f_j\]
and $[e_i,f_j]=\drc{ij}h_i$, for any $i,j\in\bfI$, by the maximal
ideal intersecting the span of $\{\hcor{i},\cow{i}\}_{i\in\bfI}$ trivially.

$\olg=\olg(\sfA)$ is non--canonically a split central extension of the
\KMA $\g=\g(\sfA)$ corresponding to $\sfA$ (\ref{ss:h''}). Unlike $\g$,
however, the Lie algebra $\olg$ always possesses a diagrammatic
structure over the Dynkin diagram $D$ of $\sfA$, which is given by
associating to any $B\subseteq D$ the subalgebra $\olg_B\subseteq
\olg$ generated by $\{e_i,f_i,\hcor{i},\fcw_i\}_{i\in B}$. In particular,
$\olg_B$ is the extended \KMA corresponding to $\sfA_B$.

If $\sfA$ is symmetrisable, the Borel subalgebras $\olb_+,\olb_-$
generated by $\{e_i,\hcor{i},\fcw_{i}\}_{i\in\bfI}$ and $\{f_i,\hcor{i}
\fcw_{i}\}_{i\in\bfI}$ respectively, are split diagrammatic Lie bialgebras.
Each gives rise to a symmetric pre--Coxeter category
$\Cox{\olb_\pm}{}$ with diagrammatic categories $\DrY{\olb_
{\pm,B}}$, $B\subseteq D$, and, as in \ref{ss:ss g} there
is a canonical embedding of the pre--Coxeter category of
$\olg$--modules with a locally finite $\b_\mp$--action to
$\Cox{\olb_\pm}{}$.

\subsection{}\label{ss:hopf ex}%quantum groups and admissible DY 
%---------------

A quantum analogue of the symmetric pre--Coxeter category $\Cox{\b}{}$
can be obtained along similar lines from split diagrammatic Hopf algebras.
A \DYt module over a Hopf algebra $\hC$ is a triple $(\V,\rho,\rho^*)$, where
$\rho:\hC\otimes \V\to \V$ is a left $\hC$--module, $\rho^*:\V\to \hC \otimes
\V$ a right $\hC$--comodule, and $\rho,\rho^*$ satisfy an appropriate
compatibility \cite{Y,ek-2}. Such modules form a braided monoidal
category $\hDrY{\hC}{}$, with commutativity constraints
$\beta_{\U,\V}:\U\otimes \V\to \V\otimes\U$ given by
\[\beta_{\U,\V}=
(1\,2)\circ \rho_\U\otimes \id_\V \circ\,(1\,2)\circ \id_\U\otimes\rho_\V^*.\]

%If $\hC$ is finite--dimensional, the category $\hDrY{\hC}{}$ coincides
%with that of representations of the quantum double of $\hC$ \cite{drin-2}.

If $\hC$ is finite--dimensional, the category $\hDrY{\hC}{}$ coincides
with that of representations of the quantum double of $D\hC$ of $\hC$
\cite{drin-2}. As a coalgebra, $D\hC$ is the tensor product $\hC\otimes
\hC^\circ$, where $\hC^\circ$ is the dual Hopf algebra $\hC^*$ endowed
with the opposite coproduct. Moreover, $D\hC$ is endowed with a unique
product such that $\hC,\hC^\circ$ are subalgebras, and $D\hC$
is a quasitriangular Hopf algebra, with $R$--matrix given by the
canonical element in $\hC\otimes\hC^\circ$. The isomorphism
$\hDrY{\hC}{}\cong\Rep(D\hC)$ is obtained by letting $\phi\in\hC
^\circ$ act on $\V\in\hDrY{\hC}{}$ by $\phi\otimes\id_\V\circ\rho^*$,
and conversely defining the coaction of $\hC$ on $\V\in\Rep(D\hC)$
by $\rho^*v=R\;1\otimes v$.

A similar equivalence holds if $\hC$ is a quantised universal enveloping 
algebra (QUE), that is a topological Hopf algebra over $\sfk{\fml}$ such
that $\hC/\hbar \hC$ is a universal enveloping algebra $U\b$.
If $\b$ is finite--dimensional, one can consider the quantised formal group
$\hC'\subset\hC$ corresponding to $\hC$ defined in \cite{drin-2,gav},
define the dual QUE $\hC^\vee$ as $(\hC')^*$, and the quantum double
of $\hC$ as the double crossed product $D\hC=\dcp{\hC}{\; \hC^\vee}$
introduced in \cite{majid-2}. The latter is a quasitriangular QUE, which
quantises the Drinfeld double of $\b$, with $R$--matrix given by the
canonical element in $\hC'\otimes\hC^\vee\subset D\hC^{\otimes 2}$.
The representations of $D\hC$ then coincide, as a braided monoidal
category, with the category $\aDrY{\hC}$ of {\it admissible} \DYt modules
over $\hC$, that is are those for which the coaction $\rho^*:\V\to\hC
\otimes\V$ factors through $\hC'\otimes\V$.

More generally, let $\hC=\bigoplus_{n\geqslant 0}\hC_n$ be an $\IN
$--graded QUE such that $\hC_0$ deforms $U\b$ with $\dim\b<\infty$,
and each $\hC_n$ is \fg over $\hC_0$. Then, admissible \DYt modules
over $\hC$ coincide with modules $\V$ over the quantum double of
$\hC$ such that the action of $\hC^\vee=\bigoplus_{n\geqslant 0}
(\hC'\cap \hC_n)^*$ is locally finite, \ie such that for any $v\in\V$,
$\hC^\vee_n v=0$ for $n$ large enough (Sect. \ref{ss:adm QUE}).

\subsection{}\label{ss:hopf split}%quantum groups and admissible DY 
%--------------

As in the case of Lie bialgebras, a split pair $\hA\stackrel{i}{\to}\hC
\stackrel{p}{\to}\hA$ of Hopf algebras gives rise to a monoidal restriction
functor $\Res_{\hC,\hA}:\DrY{\hC}\to\DrY{\hA}$ defined by
\[\Res_{\hA,\hC}(\V,\rho_\V,\rho^*_\V)=(\V,\rho_\V\circ i\ten\id_\V,p\ten\id_\V\circ\rho^*_\V)\]
If both $\hA,\hC$ are finite--dimensional, $\Res_{\hA,\hC}$ corresponds
to the pullback functor $(i\otimes p^t)^*:\Rep(D\hC)\to\Rep(D\hA)$.  %Moreover, i
If both $\hA,\hC$ are QUEs, $\Res_{\hA,\hC}$ restricts to a functor $\aDrY{\hC}\to\aDrY{\hA}$.
It follows that if $\hC$ is a diagrammatic QUE, there is a braided pre--Coxeter
category $\Cox{\hC}{\adm}$ with diagrammatic categories $\aDrY{\hC_B}$, $B\subseteq D$,
restriction functors $\Res_{\hC_{B'},\hC_{B}}$, $B'\subseteq B$, and trivial associators
and vertical joins.

Such an example arises from a quantised extended \KMA algebra
$\DJ{\olg}$, specifically from the split diagrammatic structure on its
quantum Borel subalgebras $\DJ{\olb_\pm}$. Moreover, the
realisation of $\DJ{\olg}$ as a central quotient of the %(restricted)
quantum double of $\DJ{\olb_\pm}$ yields an embedding of the
pre--Coxeter category of $\DJ{\olg}$--modules with a locally finite
action of $\DJ{\olb_\mp}$ into $\Cox{\olb_\pm}{\adm}$. Moreover,
once attention is restricted to integrable modules, Lusztig's quantum
Weyl group elements extend the structure to that of a braided Coxeter
category.

\subsection{}\label{ss:ek expo} %EK quantisation
%---------------

We now explain how the 2--categorical extension of \nEK quantisation
obtained in \cite{ATL1-1} yields an equivalence between a deformation
of the braided pre--Coxeter category $\Cox{\b}{}$ described in \ref
{intro:diag-lba}, and its quantum counterpart described in \ref{ss:hopf split}.

In \cite{ek-1,ek-2}, Etingof and Kazhdan construct a quantisation 
functor $\Q$ from the category of Lie bialgebras over a field $\sfk$
of characteristic zero to that of QUEs. $\Q$ depends on the choice
of an associator $\Phi$, and is compatible with taking \DYt modules.
Specifically, it is endowed with a braided tensor equivalence
\begin{equation}\label{eq:TEK}
\EKeq\b: \hDrY{\b}{\Phi}\longrightarrow\aDrY{\Q(\b)}
\end{equation}
where $\hDrY{\b}{\Phi}$ is the category of deformation \DYt modules
over the Lie bialgebra $\b$, with deformed associativity constraints
given by $\Phi$ \cite{ek-6, ATL1-1}. If $\g$ is a symmetrisable
(extended) \KM algebra with negative Borel subalgebra $\b$, this
implies in particular the existence of an equivalence between category
$\O$ representations of $\g$ and those of the quantum group
$U_\hbar\g$.

\subsection{} \label{ss:quant diagr}% upgrade of EK
%--------------

Assume now that $\b$ is a split diagrammatic Lie bialgebra. By
functoriality, its quantisation $\Q(\b)$ is a split diagrammatic QUE
and, by \ref {ss:hopf split}, there is a braided pre--Coxeter category
$\Cox{\Q(\b)}{\adm}$ with diagrammatic categories $\aDrY{\Q(\b_B)}$,
$B\subseteq D$. The equivalence \eqref{eq:TEK} then raises the
following question: is there a braided pre--Coxeter category $\Cox
{\b}{\Phi}$ with diagrammatic categories $\hDrY{\b_B}{\Phi}$, such
that the \nEK equivalences $\{H_{\b_B}\}_{B\subseteq D}$ fit within
an equivalence $H:\Cox{\b}{\Phi}\to\Cox{\Q(\b)}{\adm}$?

This involves in particular constructing, for any $B'\subseteq B$
and \mns $\F\in\Mns{B,B'}$, a monoidal restriction functor $F_\F
:\hDrY{\b_B}{\Phi}\to\hDrY{\b_{B'}}{\Phi}$, and a natural isomorphism
$v_\F$ making the following diagram commutative
\begin{equation}\label{eq:compatibility}
\xymatrix@R=1cm@C=2cm{
\hDrY{\b_B}{\Phi}\ar[r]^{H_{\b_B}}\ar[d]_{F_\F} & \aDrY{\Q(\b_B)}\ar[d]^{\Res_{\Q(\b_{B'}),\Q(\b_B)}}\ar@{<=}[dl]_{v_\F}\\
\hDrY{\b_B}{\Phi}\ar[r]_{H_{\b_{B'}}}			& \aDrY{\Q(\b_{B'})}
}
\end{equation}
In order for the pre--Coxeter category $\Cox{\b}{\Phi}$ to fall within
the scope of the rigidity theorem proved in \cite{ATL2}, we require 
further that the non--monoidal functor underlying $F_\F$ be equal
to $\Res_{\b_{B'},\b_B}$, which renders the problem non--trivial.\footnote
{Equivalently, we require that the composition $H_{\b_{B'}}^{-1}\circ
\Res_{\Q(\b_{B'}),\Q(\b_B)}\circ H_{\b_B}$ be isomorphic, as a
non--monoidal functor, to $\Res_{\b_{B'},\b_B}$.} 

%\subsection{}\label{ss:quant diagr} %quantisation of diagrammatic LBA

We answer this question in the affermative by relying on the compatibility
of \nEK quantisation \wrt restrictions  proved in \cite{ATL1-1}, and get the
following (Thms. \ref{ss:main-thm-1} and \ref{ss:main-thm-2})

\begin{theorem}\label{th:intro transfer}
Let $\b$ be a split diagrammatic Lie bialgebra, and $\Phi$ a Lie associator.
\begin{enumerate}
\item There is a canonical braided pre--Coxeter category $\Cox{\b}{\Phi}$
with the following properties.
\begin{itemize}
\item For any $B\subseteq D$, the diagrammatic category $\Cox{\b,B}{\Phi}$
is given by $\hDrY{\b_B}{\Phi}$.
\item For any $B'\subseteq B$, and $\F\in\Mns{B,B'}$, the functor
\[F_\F:\hDrY{\b_B}{\Phi}\to\hDrY{\b_{B'}}{\Phi}\]
is of the form $(\Res_{\b_{B'},\b_B},J_\F)$ for some tensor structure $J_\F$.
\item For any $B''\subseteq B'\subseteq B$, $\F\in\Mns{B,B'}$ and $\F\in\Mns
{B',B''}$, the composition $F_{\F'}\circ F_{\F}$ is equal to $F_{\F\cup\F'}$ as
functors $\hDrY{\b_B}{\Phi}\to \hDrY{\b_{B''}}{\Phi}$, and the %corresponding
vertical join $\redasso{\F'}{\F}:F_{\F'}\circ F_{\F}\Rightarrow F_{\F\cup\F'}$ is
the identity.
%The vertical joins of $\Cox{\b}{\Phi}$ are trivial.
\item $\Cox{\b}{\Phi}$ reduces to $\Cox{\b}{}$ mod $\hbar$.
\end{itemize}
\item The \nEK equivalences $H_{\b_B}:\hDrY{\b_B}{\Phi}\to\aDrY{\Q(\b_B)}$, %$B\subseteq D$,
fit within an equivalence of braided pre--Coxeter categories $\Hcox{\b}:\Cox{\b}{\Phi}\to\Cox{\Q(\b)}{\adm}$.
\end{enumerate}
\end{theorem}

\subsection{} %PROP EK and ATL1.1
%---------------

Recall that the functoriality of \nEK quantisation is a direct consequence of its
realisation in the context of $\PROP$s \cite{ek-2}. Roughly, this consists in
obtaining formulae which define a Hopf algebra $\Q(\univb)$ which quantises
the {\it universal Lie bialgebra} $\univb$ over $\sfk$. By definition, the latter is
the generating object of a $\sfk$--linear, symmetric monoidal category $\LBA$
endowed with a morphism $\univb\otimes\univb\to\univb$, which is antisymmetric
and satisfies the Jacobi identity. The definition of $\LBA$ implies that the category
of Lie bialgebras over $\sfk$ is equivalent to that of monoidal functors $F:\LBA
\to\vect_\sfk$, via the functor mapping $F$ to $F([\b])$. As a consequence, a
quantisation of $\univb$ in $\LBA$ can be applied to any Lie bialgebra $\b$,
and gives rise to a quantisation functor $\b\mapsto\Q(\b)$.

An extension of the $\PROP$ic definition of \nEK quantisation plays an even
greater role in proving the compatibility of the equivalences $H_{\b}$ with the
restriction functors (cf. \eqref{eq:compatibility}), as well as proving that the
functor $H_\b$ is an equivalence \cite{ATL1-1} .

\subsection{} %PROP pre-Coxeter
%--------------

In a similar vein, the braided pre--Coxeter category $\Cox{\b}{\Phi}$ of
Theorem \ref{th:intro transfer} is constructed through suitable $\PROP$s.
To this end, we introduce a universal split diagrammatic Lie bialgebra
$\ub$ by extending the category $\LBA$ by a family of idempotents
$\theta_B\in\End(\ub)$ labelled by the subdiagrams of $D$, which
satisfy $\theta_D=\id$, 
\[\theta_B\circ\theta_{B'}=\theta_{B'}=\theta_{B'}\circ\theta_B
\aand
\theta_{B'\sqcup B''}=\theta_{B'}+\theta_{B''}\]
whenever $B'\subseteq B$ and $B'\perp B''$ respectively. By relying on
\cite{ATL1-1}, we then construct a braided pre--Coxeter structure $\Cox
{\ub}{\Phi}$ on \DYt modules over $\ub$. This structure gives rise to a
braided pre--Coxeter category $\Cox{\b}{\Phi}$ for any split diagrammatic
Lie bialgebra $\b$.

Other than its economy, the use of $\Cox{\ub}{\Phi}$ shows that the
structure constants of each $\Cox{\b}{\Phi}$ are {\it universal}, that is 
admit a lift to the algebras of endomorphisms of tensor products of
\DYt modules over $\ub$. This feature is a crucial requirement of
the rigidity result obtained in \cite{ATL2}.

\subsection{}\label{ss:main res}%diagrammatic/extended KM
%--------------

Finally, we apply these results to an extended symmetrisable \KMA $\olg$,
with negative Borel subalgebra $\olb$ and Dynkin diagram $D$.

The \DrJ quantum group $\DJ\olb$ is a split diagrammatic QUE. As such,
it gives rise to a braided pre--Coxeter category $\Cox{\DJ{\olb}}{\adm}$.
Consider the subcategories defined as follows.
\begin{itemize}
\item $\Cox{\DJ{\olb}}{\adm,\sint }\subset\Cox{\DJ{\olb}}{\adm}$. The
diagrammatic category corresponding to $B\subseteq D$ consists of
admissible \DYt modules over $\DJ\olb_B$ which arise from integrable
$\DJ{\olg_B}$--modules. Specifically, since $\DJ{\olg_B}$ is a quotient
of the quantum double of $\DJ\olb_B$, we require that the action of
$D(\DJ\olb_B)$ factor through an integrable action of $\DJ\olg_B$.
\item $\OCox{\infty,\DJ{\olg}}{\sint}\subset\Cox{\DJ{\olb}}{\adm,\sint }$.
The corresponding diagrammatic categories consist of integrable $\DJ
\olg_B$--modules in category $\O_\infty$.\footnote{The symbol $\infty$
refers to the fact that we allow infinite--dimensional weight spaces. This
is required by the fact that the restriction corresponding to $\olg_{B'}\subset
\olg_B$ or $\DJ{\olg_{B'}}\subset\DJ{\olg_B}$ does not preserve the
finite--dimensionality of weight spaces if $B'\subsetneq B$.}
\end{itemize}
The quantum Weyl group operators of $\DJ\olg$ \cite{l} endow $\Cox
{\DJ{\olb}}{\adm,\sint }$, and therefore $\OCox{\infty,\DJ{\olg}}{\sint}$,
with the structure of a braided Coxeter category.

The combination of Theorem \ref{th:intro transfer} and the isomorphism
of diagrammatic Hopf algebras $\DJ{\olb}\simeq\Q(\olb)$ yields our main
result.

\begin{theorem}
Let $\olg$ be an extended symmetrisable Kac--Moody algebra. 
\begin{enumerate}
\item There is a universal braided Coxeter category $\OCox{\infty,\olg}
{\Phi,\sint}$ such that
\begin{itemize}
\item The diagrammatic category corresponding to $B\subseteq D$
is $\O_{\infty,\olg_B}^{\Phi,\sint}$.
\item The functor
$F_\F:\O_{\infty,\olg_B}^{\Phi,\sint}\to\O_{\infty,\olg_{B'}}^{\Phi,\sint}$
corresponding to $B'\subseteq B$ and $\F\in\Mns{B,B'}$ is the standard
restriction functor endowed with an appropriate tensor structure.
\item The vertical joins $\redasso{\F}{\F'}:F_{\F'}\circ F_\F\Rightarrow
F_{\F\cup\F'}$ are trivial.
\item The underlying braided pre--Coxeter structure is $\PROP$ic, and
trivial modulo $\hbar$.
\end{itemize}
\item The \nEK equivalences $H_{\olb_B}$ restrict to equivalences
$\O_{\infty,\olg_B}^{\Phi,\sint}\to\O_{\infty,\DJ\olg_B}^{\sint}$, and 
fit within an equivalence of braided pre--Coxeter categories
$\OCox{\infty,\olg}{\Phi,\sint}\to\OCox{\infty,\DJ\olg}{\sint}$.
\end{enumerate}
\end{theorem}

%---------------------------,
%OUTLINE
%---------------------------

\subsection{Outline of the paper}
%---------------------------------------
%Section 2
We begin in Section \ref{s:diagrams} by reviewing a number of combinatorial
notions related to diagrams. 
%Sections 3-4
We lay out the definition of a Coxeter object in an arbitrary $2$--category in
Section \ref{s:dcat}, and of a braided Coxeter category in Section \ref{s:qcqtqba}.
%Sections 5-6
In Sections \ref{s:lba} and \ref{s:quant-lba} we produce examples of braided Coxeter
categories through Drinfeld--Yetter modules over diagrammatic Lie bialgebras 
and their quantisations. 
%Section 7
In Section \ref{s:diag-prop}, we introduce a $\PROP$ which describes diagrammatic
Lie bialgebras.
%Sections 8-9
In Sections \ref{s:univ alg} and \ref{s:univ pC}, we describe in terms of $\PROP$s 
a {universal} braided pre--Coxeter structure on the category of Drinfeld--Yetter
modules over a diagrammatic Lie bialgebra.
%Section 10
In Section \ref{s:qcstructure}, we apply the results from \cite{ATL1-1} 
to the case of a diagrammatic Lie bialgebra $\b$. We show in particular that 
the braided pre--Coxeter structure of the Etingof--Kazhdan quantisation 
$\Q(\b)$ is equivalent to a universal braided pre--Coxeter structure on the 
category of Drinfeld--Yetter modules over $\b$.
%Section 11
In Section \ref{se:kmas}, we review the definition and basic properties
of the \KMA associated to an $n\times n$ matrix $\sfA$.
%Section 12
In Section \ref{se:diagr kmas}, we define extended \KM algebras which
are associated to a (non--minimal) realisation of $\sfA$ of dimension
$2n$, and show that they are naturally endowed with a structure of
diagrammatic Lie bialgebras.
%Section 13
In Section \ref{s:qkm}, we show that integrable \DYt modules over an
extended quantum group has a natural structure of braided Coxeter
category.  We then apply the results from Section \ref{s:qcstructure},
and obtain the desired  transport of the braided Coxeter structure of
the quantum group $\DJ\g$ to the category of integrable \DYt modules
for $\g$.
%Appendix
Finally, in Appendix \ref{app-graph-calc}, we provide an alternative 
description of the axioms of a Coxeter object in terms of the standard
graphical calculus for $2$--categories.

\subsection{}
%---------------

The main results of this paper first appeared in more condensed form
in the preprint \cite{ATL0}. The latter is superseded by the present paper,
and its companion \cite{ATL1-1}. 

%--------------------------------------------------------------------
% DIAGRAMS
%---------------------------------------------------------------------

\section{Diagrams and nested sets}\label{s:diagrams}
%--------------------------------------------

We review in this section a number of combinatorial notions associated
to a diagram $D$, in particular the definition of nested sets on $D$
following \cite{DCP2}, and \cite[Section 2]{vtl-4}.

\subsection{Nested sets on diagrams} 
%------------------------------------------------

A {\it diagram} is an undirected graph $D$ with no multiple edges or
loops. A {\it subdiagram} $B\subseteq D$ is a full subgraph of $D$,
that is, a graph consisting of a (possibly empty) subset of vertices
of $D$, together with all edges of $D$ joining any two elements of it.

Two subdiagrams $B_1,B_2\subseteq D$ are {\it orthogonal} if they
have no vertices in common, and no two vertices $i\in B_1$, $j\in B
_2$ are joined by an edge in $D$. We denote by $B_1\sqcup B_2$
the disjoint union of orthogonal subdiagrams. Two subdiagrams $B_1,
B_2\subseteq D$ are {\it compatible} if either one contains the other
or they are orthogonal.

A {\it nested set} on $D$ is a collection $\H$ of pairwise compatible,
connected subdiagrams of $D$ which contains the empty subdiagram
and $\cc{D}$, where $\cc{D}$ denotes the set of connected components 
of $D$. It is easy to see that the cardinality of any maximal nested set
on $D$ is equal to $|D|+1$. 

Let $\Ns{D}$ be the set of nested sets on $D$, and $\Mns{D}$ that
of maximal nested sets. Every (maximal) nested set $\H$ on $D$
is uniquely determined by a collection $\{\H_i\}_{i=1}^r$ of (maximal)
nested sets on the connected components $D_i$ of $D$. We therefore
obtain canonical identifications
\[\Ns{D}=\prod_{i=1}^r \Ns{D_i}\qquad\text{and}\qquad\Mns{D}=\prod_{i=1}^r\Mns{D_i}.\]

\subsection{Relative nested sets}\label{ss:rel-ns}
%-----------------------------------------

If $B'\subseteq B\subseteq D$ are two subdiagrams of $D$, a nested
set on $B$ {\it relative} to $B'$ is a collection of subdiagrams of $B$
which contains $\cc{B}$ and $\cc{B'}$, and in which every element is
compatible with, but not properly contained in any of the connected
components of $B'$.  We denote by $\Ns{B,B'}$ and $\Mns{B,B'}$
the collections of nested sets and maximal nested sets on $B$ relative
to $B'$. In particular, $\Ns{B}=\Ns{B,\emptyset}$ and $\Mns{B}=\Mns
{B,\emptyset}$.

Relative nested sets are endowed with the following operations, which
preserve maximal nested sets.
\begin{enumerate}
\item {\bf Vertical union.}
For any $B''\subseteq B' \subseteq B$, there is an embedding
\[%\begin{equation}\label{eq:union-nes-set}
\cup:\Ns{B,B'}\times\Ns{B',B''}\to\Ns{B,B''},
\]%\end{equation}
given by the union of nested sets. Its image is the collection $\Nsr{B,B''}
{B'}\subseteq\Ns{B,B''}$ of relative nested sets which contain $\cc{B'}$.
\item {\bf Orthogonal union.}
For any $B=B_1\sqcup B_2$ and $B'=B_1'\sqcup B_2'$
with $B_1'\subseteq B_1$, $B_2'\subseteq B_2$, there is 
a bijection
\[\Ns{B_1,B_1'}\times\Ns{B_2,B_2'}\to\Ns{B,B'},\]
mapping $(\H_1,\H_2)\mapsto\H_1\cup\H_2$.
\end{enumerate}

\subsection{Nested sets and chains of subdiagrams}\label{ss:chain-ns}
%-------------------------------------------------------------------

\begin{definition}
A \emph{chain} from $B\subseteq D$ to $B'\subseteq B$ is a sequence
of subdiagrams
\[
\mathbf{C}: B'=B_0\subsetneq B_1\subsetneq\cdots\subsetneq B_m=B.
\] 
A chain is called \emph{maximal} if $|B_k\setminus B_{k-1}|=1$ for every
$k$. The sets of chains and maximal chains from $B$ to $B'$ are denoted
$\mathsf{Ch}(B, B')$ and $\mathsf{MCh}(B, B')$, respectively. 
\end{definition}
\noindent
Note that, unlike the notion of nested set, that of chain is independent of 
the connectivity of the graph and only depends on the underlying set of 
vertices. The following is clear.

\begin{lemma}
There is a surjective map $p:\mathsf{Ch}(B, B')\to\mathsf{Ns}(B,B')$ given by
\[
p(B'=B_0\subsetneq B_1\subsetneq\cdots\subsetneq B_m=B)=\bigcup_{k=0}^{m}\cc{B_k},
\]
%where $\cc{B_k}$ denotes the connected components of $B_k$.
Moreover, $p$ restricts to a surjection $p:\mathsf{MCh}(B, B')\to
\Mns{B,B'}$.
\end{lemma}

\noindent
The operations defined in \ref{ss:rel-ns} naturally extend to chains, and it is
easy to check that the maps $p$ preserve these operations. In particular, 
\begin{itemize}
\item {\bf Vertical union.} For any $B''\subseteq B' \subseteq B$, $\bfC\in\Ch
{B,B'}$, and $\bfC'\in\Ch{B',B''}$, we denote by $\bfC\cup\bfC'\in\Ch{B,B''}$
the chain obtained by vertical composition
\item {\bf Orthogonal union.} For any $B=B_1\sqcup B_2$ and $B'=B_1'\sqcup
B_2'$ with $B_1'\subseteq B_1$, $B_2'\subseteq B_2$, $\bfC\in\Ch{B_1\sqcup
B_2, B'_1\sqcup B'_2}$, we denote by $\bfC_{B_k}\in\Ch{B_k,B'_k}$, $k=1,2$,
the chains determined by $\bfC$ on $B_1$ and $B_2$.
\end{itemize}

Two chains give rise to the same nested set if they differ only at the level of
orthogonal subdiagrams. Specifically, if $B_1'\subsetneq B_1\perp B_2\supsetneq B_2'$,
the chains
\begin{eqnarray}\label{eq:min-chain}
&&\bfC_1: B'_1\sqcup B'_2 \subset B_1\sqcup B'_2 \subset B_1\sqcup B_2\\[1.1ex]
&&\bfC_{\color{white}0}: B'_1\sqcup B'_2\subset B_1\sqcup B_2\nonumber\\[1.1ex]
&&\bfC_2: B'_1\sqcup B'_2 \subset B'_1\sqcup B_2 \subset B_1\sqcup B_2\nonumber
\end{eqnarray}
give rise to the same nested set in $\Ns{B_1\sqcup B_2,B'_1\sqcup B'_2}$.
More generally, for any $B'\subseteq B$, we denote by $\mathbf{G}_{B,B'}$
the graph having $\Ch{B,B'}$ as set of vertices, and an edge between $\bfC
^1$ and $\bfC^2$ if and only if their difference is limited to a subchain of the
form \eqref{eq:min-chain} for the same subdiagrams $B'_1, B'_2, B_1, B_2$. 
More precisely,  $\bfC^1$ and $\bfC^2$ are connected by an edge if and only if
$\bfC^1\neq\bfC^2$ and the following holds
\begin{itemize}
\item $\bfC^1\not\subset\bfC^2$ and $\bfC^2\not\subset\bfC^1$, $\bfC^1,
\bfC^2$ are of the same length, there is an index $i$ such that $B^1_j=B^2
_j$, for $j\neq i$, and subdiagrams $B_1'\subsetneq B_1\perp B_2\supsetneq
B_2'$ such that 
\begin{gather*}
B^1_{i+1}=B_1\sqcup B_2=B^2_{i+1}\\[1.1ex]
B^1_i= B_1\sqcup B'_2\quad\quad B'_1\sqcup B_2=B^2_i\\[1.1ex]
B^1_{i-1}=B'_1\sqcup B'_2=B^2_{i-1}
\end{gather*}
\item $\bfC^1\subset\bfC^2$, there is an index $i$ such that $B^1_j=B^2_{j}$
if $j<i$ and $B^1_j=B^2_{j+1}$ if $j>i$, $B_1'\subsetneq B_1\perp B_2\supsetneq
B_2'$ such that
\begin{align*}
B^1_i=\,	&B_1\sqcup B_2=B^2_{i+1}\\
		&B_1\sqcup B'_2=B^2_i\\
B^1_{i-1}=\,	&B'_1\sqcup B'_2=B^2_{i-1}
\end{align*}
(and similarly for $\bfC^2\subset\bfC^1$).
\end{itemize}
The following is straightforward.

\begin{proposition}
The map $p:\Ch{B,B'}\to\Ns{B,B'}$ descends to a bijection
\[p:\slantfrac{\Ch{B,B'}}{\sim}\to\Ns{B,B'}\]
where $\sim$ is the equivalence relation defined by the graph 
$\mathbf{G}_{B,B'}$, \ie $\bfC\sim\bfC'$ if and only if they are 
connected in $\mathbf{G}_{B,B'}$. 
\end{proposition}

\noindent\remark
The map $p$ admits a canonical section 
$s:\Ns{B,B'}\to\Ch{B,B'}$ which assigns to a nested set $\H$ the
chain $s(\H)$ defined recursively as follows
\begin{itemize}
\item $s(\H)_{\operatorname{top}}=B$
\item $s(\H)_{k-1}$ is the union of the elements of $\H$ which are properly 
contained and maximal in $s(\H)_{k}$
\end{itemize}
Clearly, $p(s(\H))=\H$. Note, however, that $s$ does not preserve the 
vertical union of nested set. Namely, if $\H\in\Ns{B,B'}$ and $\H'\in\Ns{B',B''}$,
then in general $s(\H)\cup s(\H')\neq s(\H\cup\H')$.  Also, $s$ does not map maximal 
nested set to maximal chains. Indeed, if $\F\in\Mns{B,B'}$, $|s(\F)_k\setminus s(\F)_{k-1}|=
|\cc{s(\F)_k}|\geqslant1$.

%%%%%%%%%%%%%%%%%%%
% PRE COXETER CATEGORIES	%
%%%%%%%%%%%%%%%%%%%

\section{Coxeter objects}\label{s:dcat}
%=================

In this section, we define Coxeter objects in an arbitrary 2--category
$\XX$.

\subsection{2--Categories}
%--------------------------------

By definition, a $2$--category is a category enriched over $\mathsf{Cat}$,
the category of categories, functors and natural transformations \cite
{ehresmann,lack}. In particular, a $2$--category is a special example
of a bicategory \cite{benabou}. The difference between the two notions
lies in the composition of $1$--morphisms, which is required to be
associative up to a prescribed isomorphism in a bicategory, and strictly
associative in a $2$--category. In particular, a 2--category with one
object is a strict (small) monoidal category.

For simplicity, in this section we work with a fixed $2$--category $\XX$,
though our definitions easily carry over to a bicategory.

\subsection{The diagrammatic 2--category ${\diagr{}{\XX}}$}\label{ss:Dcatmor}
%-----------------------------------------------------------------------------

Let $B'\subseteq B$ be two diagrams. If $\K\in\Ns{B,B'}$ is a relative
nested set, we denote by $\Mnsr{B,B'}{\K}$ the collection of relative
\mnss on $B$ which contain $\K$. If $C_1,\ldots,C_m\subseteq B$
are compatible diagrams such that $\K=\cc{C_1}\cup\cdots\cup\cc{
C_m}$ is a relative nested set in $\Ns{B,B'}$, we abbreviate $\Mnsr
{B,B'}{\K}$ to $\Mnsr{B,B'}{\{C_1,\ldots,C_m\}}$.

\begin{definition}
The diagrammatic category $\diagr{}{\XX}$ is the following 2--category
\begin{enumerate}
%
% B-objects
\item If $B$ is a diagram, a \emph{$B$--object} is an object $\C_B$ in
$\XX$ labelled by $B$. 
%
% diagrammatic 1-morphisms
\item If $B'\subseteq B$ are diagrams, $\CB$ a $B$--object, $\CBp$
a $B'$--object, and $\K\in\Ns{B,B'}$, a {\it diagrammatic $1$--morphism
$\C\to\C'$ of degree $\K$} is the datum of
\begin{itemize}
\item for any $\F\in\Mnsr{B,B'}{\K}$, a $1$--morphism $F_{\F}:
\CB\to\CBp$
\item for any $\F,\G\in\Mnsr{B,B'}{\K}$, a $2$--isomorphism
$\DCPA{\G}{\F}:F_{\F}\Rightarrow F_{\G}$
\end{itemize}
such that the morphisms $\Upsilon$ are transitive, \ie 
for any $\F,\G,\H\in\Mnsr{B,B'}{\K}$, 
\[\Upsilon_{\H\G}\circ\Upsilon_{\G\F}=\Upsilon_{\H\F}\]
This implies in particular that $\DCPA{\F}{\F}=\id_{F_{\F}}$, and that
$\DCPA{\G}{\F}=\DCPA{\F}{\G}^{-1}$ for any $\F,\G\in\Mnsr{B,B'}{\K}$.
% 1-Homs
We denote the collection of 1--morphisms $\CB\to\CBp$ of degree
$\K$ by $\Hom(\CB,\CBp)[\K]$, and set\footnote{Note that if $\K_1
\subseteq\K_2\in\Ns{B,B'}$ then $\Mnsr{B,B'}{\K_1}\supseteq \Mnsr
{B,B'}{\K_2}$, and there is a forgetful map $\Hom(\C,\C')[\K_1]\to
\Hom(\C,\C')[\K_2]$}
\[\diagr{}{\XX}(\CB,\CBp)=\bigsqcup_{{\K}\in\Ns{B,B'}}
\mathsf{Hom}(\CB,\CBp)[{\K}]\]
% composition of 1--morphisms
\item 
If $B''\subseteq B'\subseteq B$ are encased diagrams, $\CB,\CBp,
\CBpp$ are $B,B'$, and $B''$--objects, $\K\in\Ns{B,B'}$ and $\K'\in
\Ns{B',B''}$, the composition of 1--morphisms
\[F:\CB\to\CBp\aand F':\CBp\to\CBpp\]
of degrees $\K$ and $\K'$ is a 1--morphism $F'\circ F:\CB\to\CBpp$
of degree $\K\cup\K'\in\Ns{B,B''}$.
Specifically, if $\F,\G\in\Mnsr{B,B''}{\K\cup\K'}$, the 1-- and 2--morphisms
\[F_{\F}:\CB\to\CBpp\aand 
\DCPA{\G}{\F}:
F_{\G}\Rightarrow F_{\F}\]
corresponding to $F'\circ F$ are given by the composition $F'_{\F_{B''B'}}
\circ F_{\F_{B'B}}$ and the vertical composition 
$\begin{array}{l}\DCPA{\G_{B'B}}{\F_{B'B}}\\[1.1ex]\DCPA{\G_{B''B'}}{\F_{B''B'}}'\end{array}$
respectively.\footnote{Note that the composition $F'\circ F$ forgets
some of the data of $F$, %in a way similar to \eqref{eq:drop degree},
namely the 1--morphisms $F_{\F}$ and 2--morphisms $\DCPA{\F}{\G}$ 
corresponding to $\F,\G\in\Mnsr{B,B''}{\K}\setminus
\Mnsr{B,B''}{\K\cup\K'}$.}

%
% 2-morphisms
\item
If $F^1,F^2:\CB\to\CBp$ are $1$--morphisms of degrees ${\K}_1,{\K}_2\in
\Ns{B,B'}$ respectively, a {\it diagrammatic $2$--morphism} $u:F^1\Rightarrow
F^2$ is the datum, for any $\F_1\in\Mnsr{B,B'}{{\K}_1}$ and $\F_2\in\Mnsr
{B,B'}{{\K}_2}$, of a $2$--morphism $u_{\F_2\F_1}:F^1_{\F_1}\Rightarrow
F^2_{\F_2}$ in $\XX$ such that, for any $\F_1,\G_1\in\Mnsr{B}{{\K}_1}$
and $\F_2,\G_2\in\Mnsr{B}{{\K}_2}$, 
\begin{equation}\label{eq:cube}
u_{\G_2\G_1}\circ \Upsilon_{\G_1\F_1}^1=
\Upsilon_{\G_2\F_2}^2\circ u_{\F_2\F_1}
\end{equation}
as 2--morphisms $F^1_{\F_1}\Rightarrow F^2_{\G_2}$. This amounts
to the commutativity of %the diagram
\begin{equation}\label{eq:cube diagr}\xy
(0,-5)*{
\xy
(0,50)*++{\CB}="C1";
(0,10)*++{\CBp}="C2";
%%%%%
%FUNCTORS
{\ar@/_5pc/@{->}_{F^2_{\G_2}} "C1";"C2"};
%%%%%
{\ar@/_1pc/@{->}|(.7){F^2_{\F_2}} "C1";"C2"};
%%%%%
{\ar@/^1pc/@{..>}|(.2){F^1_{\G_1}} "C1";"C2"};
%%%%
{\ar@/^5pc/^{F^1_{\F_1}} "C1";"C2"};
%%%%%
(-20,30)*++{}="F2";
(-3,22)*+++{}="G2";
(20,30)*++{}="G1";
(3,42)*+++{}="F1";
%NAT TRANS
{\ar@/_0.8pc/@{<=}|(.45){\DCPA{\G_2}{\F_2}^2}"F2";"G2"};
{\ar@/_0.8pc/@{<=}|{u_{\F_2\F_1}}"G2";"G1"};
{\ar@/^0.8pc/@{<:}|{\DCPA{\G_1}{\F_1}^1}"F1";"G1"};
{\ar@/^0.8pc/@{<:}|(.4){u_{\G_2\G_1}}"F2";"F1"};
\endxy
};
\endxy\end{equation}
%
% Diagr_D
\end{enumerate}
If $D$ is a fixed diagram, we denote by $\diagr{D}{\XX}\subset
\diagr{}{\XX}$ the full 2--subcategory of $B$--objects, where
$B\subseteq D$. 
\end{definition}

\subsection{Pre--Coxeter objects}\label{ss:pre-cox}
%------------------------------------------

Let $D$ be a diagram.

\begin{definition}
A \emph{pre--Coxeter object} of type $D$ in $\XX$ is the datum of 
\begin{itemize}
\item for any $B\subseteq D$, a $B$--object $\C_B$ %in $\XX$
\item for any $B'\subseteq B$, a diagrammatic $1$--morphism
$F_{B'B}:\C_{B}\to\C_{B'}$ of minimal degree $\K=\cc{B}\cup\cc{B'}$
\item for any $B''\subseteq B'\subseteq B$, a diagrammatic
$2$--isomorphism
\begin{equation}\label{eq:alpha BBB}
\xymatrix@R=.8cm@C=1.2cm{
& \C_{B'} \ar[dl]_{F_{B''B'}} \ar|(.55){\alpho{B}{B'}{B''}}@{=>}[d]& \\
\C_{B''}  & & \ar[ll]^{F_{B''B}}\ar[ul]_{F_{B'B}}\C_{B}
}\end{equation}
\end{itemize}
such that
\begin{itemize}
\item \norm
{for any $B'\subseteq B$, $F_{BB}=\id
_{\C_B}$ and $\alpho{B'}{B'}{B}=\id_{F_{B'B}}=\alpho{B'}{B}{B}$.}
\item the 2--morphisms $\alpha$ are associative, \ie for any $B'''\subseteq
B''\subseteq B'\subseteq B$, the following tetrahedron in $\diagr{D}{\XX}$
with 2--faces given by the morphisms $\alpha$ is commutative
\begin{equation}\label{eq:tetra}
\xymatrix{
&\C_{B'''}&&\ar[ll]_{F_{B'''B''}}\C_{B''}&\\
&&&&\C_{B'}\ar|(.55){F_{B'''B'}}@{..>}[ulll]  \ar[ul]_{F_{B''B'}} \\
&&\C_{B}\ar[uul]^{F_{B'''B}}\ar[urr]_{F_{B'B}}\ar|(.5){F_{B''B}}[uur]&&
}
\end{equation}
In other words, the following equality holds
\[\alpho{B}{B''}{B'''}\circ\alpho{B}{B'}{B''}=
\alpho{B}{B'}{B'''}\circ\alpho{B'}{B''}{B'''}\]
as 2--isomorphisms $F_{B'''B''}\circ F_{B''B'}\circ F_{B'B}\Rightarrow F_{B'''B}$.
\end{itemize}
\end{definition}

\subsection{Unfolding the definition}\label{ss:alpha-reduction}
%---------------------------------------------
	
	We give below a more hands--on description of a pre--Coxeter object,
	which will be used throughout this paper to construct examples.

	\begin{proposition}\label{prop:alpha-reduction}
		A pre--Coxeter object of type $D$ in $\XX$ is equivalently described by the
		datum of 
		\begin{itemize}
			\item for any $B\subseteq D$, an object $\C_B\in\XX$
			\item for any $B'\subseteq B$ and $\F\in\Mns{B,B'}$, a 
			$1$--morphism $F_{\F}:\C_{B}\to\C_{B'}$ 
			\item for any $B'\subseteq B$ and $\F,\G\in\Mns{B,B'}$,
			a $2$--isomorphism $\DCPA{\G}{\F}:F_{\F}\Rightarrow F_{\G}$
			\item for any $B''\subseteq B'\subseteq B$, $\F\in\Mns{B,B'}$ and 
			$\F'\in\Mns{B',B''}$, a $2$--isomorphism 
			$\redasso{\F}{\F'}: F_{\F'}\circ F_{\F}\Rightarrow F_{\F'\cup\F}$.
		\end{itemize}
		such that 
		\begin{enumerate}
			\item if $\F$ and $\F'$ are the unique elements in $\Mns{B,B}$ and $\Mns{B',B'}$, respectively, and $\G\in\Mns{B,B'}$, then $F_{\F}=\id_{\C_B}$, $F_{\F'}=\id_{\C_{B'}}$ 
			and $\redasso{\F}{\G}=\id_{F_{\G}}=\redasso{\G}{\F'}$
			\item the $2$--isomorphisms $\DCPA{}{}$ are transitive, \ie for any $\F,\G,\H\in\Mns{B,B'}$,
			$\DCPA{\F}{\G}\circ \DCPA{\G}{\H}=\DCPA{\F}{\H}$
			\item the $2$--isomorphism $\redasso{}{}$ are associative, \ie for any $B'''\subseteq
			B''\subseteq B'\subseteq B$, and \mnss $\F\in\Mns {B,B'},\F'\in\Mns{B',B''},\F''\in
			\Mns{B'',B'''}$, the following holds
			\[\redasso{\F\cup\F'}{\F''}\circ\redasso{\F}{\F'}=
\redasso{\F}{\F'\cup\F''}\circ\redasso{\F'}{\F''}\]
as 2--morphisms $F_{\F''}\circ F_{\F'}\circ F_{\F}\Rightarrow
F_{\F\cup\F'\cup\F''}$. 
			\item for any $B''\subseteq B'\subseteq B$, $\F,\G\in\Mns{B,B'}$,
			and $\F',\G'\in\Mns{B',B''}$, 
			\begin{equation*}%\label{eq:vert-fact-cat}
			\DCPA{\F'\cup\F, }{\G'\cup\G}\circ\redasso{\G}{\G'}=\redasso{\F}{\F'}\circ \begin{array}{l}\DCPA{\F}{\G}\\[1.1ex]\DCPA{\F'}{\G'}\end{array}
			\end{equation*}
		\end{enumerate}
	\end{proposition}
	
	\begin{pf}
		First, we show that any pre--Coxeter object $(\C, F,\alpha)$ 
		gives rise to the datum described above. 
		
		By definition, $F_{B'B}:\C_B\to\C_{B'}$ is a diagrammatic 1--morphism, \ie
		it amounts to a collection of 1--morphisms $F_{\F}:\C_B\to\C_{B'}$ 
		and 2--isomorphisms $\DCPA{\G}{\F}:F_{\F}\to F_{\G}$, labeled by
		$\F,\G\in\Mns{B,B'}$ and satisfying $(1)$.
		
		The diagrammatic 2--isomorphism $\alpho{B}{B'}{B''}$ 
		amounts to a collection of $2$--isomorphisms
		\[
		\xymatrix@R=.35cm@C=1.5cm{
			\C_{B} 
			\ar[dr]^{F_{\G'}} \ar[dd]_{F_{\F}}&\\
			& \C_{B'}\ar@{=>}|(.6){\asso{\F}{\G'}{\G''}}[l]\ar[dl]^{F_{\G''}}\\
			\C_{B''} 
		}
		\]
		labelled by $\F\in\Mnsr{B,B''}{}$, $\G'\in\Mnsr{B,B'}{}$ and $\G''\in\Mnsr{B',B''}
		{}$, satisfying the compatibility condition \eqref{eq:cube} and \eqref{eq:cube diagr}. 
		We set $\redasso{\F}{\F'}\coloneqq\asso{\F'\cup\F, }{\F}{\F'}$.
		Then, the condition \eqref{eq:tetra}, encoding the associativity of the morphisms
		$\alpho{}{}{}$ \eqref{eq:tetra}, clearly implies $(2)$. 
	    Then, it follows from  \eqref{eq:cube diagr} (with $\F_1=\G'\cup\G=\G_1$ and 
	    $\F_2=\G'\cup\G$, $\G_2=\F$) that $\asso{\F}{\G}{\G'}=\DCPA{\F, }{\G'\cup\G}\circ\redasso{\G}{\G'}$. 
		This implies that $\alpha$ is completely determined by the 2--isomorphisms 
		$\{\redasso{\G}{\G'}\}_{\G,\G'}$. Finally, the condition $(3)$ follows directly 
		from \eqref{eq:cube diagr}, by choosing $\wt{\F}\in\Mnsr{B,B''}{B'}$ and setting 
		$\F=\wt{\F}_{B'B}\in\Mns{B,B'}$ and $\F'=\wt{\F}_{B''B'}\in\Mns{B',B''}$. 
		The converse is proved similarly.
	\end{pf}

\subsection{The 2--categories $\ID$ and $\NsID$}\label{ss:diag-cat-model}
%----------------------------------------------------------------

We give below a succinct definition of a pre--Coxeter
object as a 2--functor to the diagrammatic category $\diagr{D}{\XX}$.

% \ID
Let $\ID$ be the $2$--category where
\begin{itemize}
\item the objects are the subdiagrams of $D$
\item the $1$--morphisms $B\to B'$ are the inclusions $B'\subseteq B$
\item the $2$--morphisms are equalities
\end{itemize}

% \NsID
Consider also the 2--category $\NsID$ where
\begin{itemize}
\item the objects are the subdiagrams of $D$
\item the $1$--morphisms $B\to B'$ are the relative nested sets $\K\in\Ns{B,B'}$,
with composition given by union
\item for any $\K_1,\K_2\in\Ns{B,B'}$, there is a unique 
$2$--isomorphism $\K_1\to\K_2$
\end{itemize}

There is a forgetful $2$--functor $\sff_D:\NsID\to\ID$, which is the identity on objects,
maps all $1$--morphisms in $\Ns{B,B'}$ to the inclusion $B'\subseteq B$, and the $2
$--morphisms to the identity. $\sff_D$ has a canonical section $\sfs_{D}:\ID\to\NsID$,
which maps the inclusion $B'\subseteq B$ to $\K_{\operatorname{min}}=
\cc{B}\cup
\cc{B'}\in\Ns{B,B'}$.% and satisfies $\sff_D\circ\iota_D=\id_{\ID}$. 
\footnote{Note that $\sfs_{D}$ is technically a pseudo $2$--functor, since it preserves
the composition only up to a coherent $2$--isomorphism. Namely, for any
$B''\subseteq B'\subseteq B$, set $\K=\cc{B}\cup\cc{B'}$, $\K'=\cc{B'}\cup\cc{B''}$
and $\K''=\cc{B}\cup\cc{B''}$. Then, the $2$--isomorphism $\K'\cup\K\to\K''$ in
$\NsID$ gives an identification$%({\K'}\to{\K})\colon
\sfs_D(B'\to B'')\circ\sfs_D(B\to B')\to\sfs_D(B\to B'')$.}
%COLON gives the right spacing
%We denote by $\diagr{D}{\XX}$ the $2$--category of $B$--objects where $B\subseteq D$. 

Consider now the $2$--functor $\sff_{D,\XX}\colon\diagr{D}{\XX}\to\NsID$, which
maps a $B$--object to the underlying diagram $B\subseteq D$, and a $1$--morphism
$\CB\to\CBp$ to its degree in $\Ns{B,B'}$. 
%In particular, $\ID$ and $\diagr{D}{\XX}$ are \emph{graded} over the $2$--category $\NsID$.
Then, a pre--Coxeter object in $\XX$
is a (pseudo) $2$--functor $\C:\ID\to\diagr{D}{\XX}$ such that $\sff_{D,\XX}\circ\C=\sfs_D$,
that is
\begin{equation}\label{eq:grading diagr}
\xymatrix@C=0.4in@R=0.3in{
\ID \ar[r]^(.4){\C} \ar[dr]_{\sfs_D} & \diagr{D}{\XX} \ar[d]^{\sff_{D,\XX}}\\ & \NsID
}
\end{equation}

\subsection{Morphisms}\label{ss:mor-pre-cox}
%----------------------------------------------------------
A {\it 1--morphism $\C\to\C'$} of pre--Coxeter objects in $\XX$ is a
natural transformation of the corresponding functors $\ID\to\diagr{D}
{\XX}$, which is compatible with \eqref{eq:grading diagr}. %preserves the $\NsID$--grading. 
Concretely, this consists of the datum of
\begin{itemize}
\item for any $B\subseteq D$, a diagrammatic 1--morphism $H_B:\C_B\to\C'_B$
\item for any $B'\subseteq B$, a diagrammatic 2--isomorphism 
\[\xymatrix@C=0.5in@R=0.5in{
\C_{B} \ar[r]^{H_{B}} \ar[d]_{F_{B'B}} & \C'_{B}\ar[d]^{F'_{B'B}} \ar|(.5){\gamma_{B'B}}@{=>}[dl]\\%_{\gamma_{B'B}} \\ 
\C_{B'}\ar[r]_{H_{B'}} & \C'_{B'}}\]
\end{itemize}
such that the morphisms $\gamma$ factorise vertically, \ie for any
$B''\subseteq B'\subseteq B$, the following prism in $\diagr{D}{\XX}$
is commutative
\[\xymatrix@C=.8cm@R=.5cm{
\C_{B}\ar[rr]^{H_{B}} \ar[dd]_{F_{\sd{B}{B''}}}  \ar|{F_{\sd{B}{B'}}}[dr]& & \C'_{B} \ar|(.3){F_{\sd{B}{B''}}}@{..>}[dd] \ar[dr]^{F_{\sd{B}{B'}}}& \\
& \C_{B'}\ar|(.35){H_{B'}}[rr]  \ar|{F_{\sd{B'}{B''}}}[dl] && \C'_{B'} \ar[dl]^{F_{\sd{B'}{B''}}}\\
\C_{B''}\ar[rr]_{H_{B''}}  & & \C'_{B''} & 
}\]
where the rectangular 2--faces are the morphisms $\gamma$, and the
triangular ones the morphisms $\alpha,\alpha'$.\\

	\noindent\remark
	In view of \ref{ss:alpha-reduction}, a 1--morphism of pre--Coxeter object 
	$\C\to\C'$ is equivalently described as the datum of
	\begin{itemize}
		\item for any $B\subseteq D$, a 1--morphism $H_B:\C_B\to\C'_B$
		\item for any $B'\subseteq B$ and $\F\in\Mns{B,B'}$, a 2--isomorphism
		$\gamma_{\F}: F'_{\F}\circ H_B\Rightarrow H_{B'}\circ F_{\F}$
	\end{itemize} 
	such that, for any $B''\subseteq B'\subseteq B$ and $\F'\in\Mns{B,B'}$, $\F''\in\Mns{B',B''}$
	and $\F\in\Mns{B,B''}$, the following prism 
	\[\xymatrix@C=.8cm@R=.5cm{
		\C_{B}\ar[rr]^{H_{B}} \ar[dd]_{F_{\F}}  \ar|{F_{\F'}}[dr]& & \C'_{B} \ar|(.3){F'_{\F}}@{..>}[dd] \ar[dr]^{F'_{\F'}}& \\
		& \C_{B'}\ar|(.35){H_{B'}}[rr]  \ar|{F_{\F''}}[dl] && \C'_{B'} \ar[dl]^{F'_{\F''}}\\
		\C_{B''}\ar[rr]_{H_{B''}}  & & \C'_{B''} & 
	}\]
	where the rectangular 2--faces are the morphisms $\gamma$, and the
	triangular ones the morphisms $a, a'$. Note also that, if $B'=B''$, for 
	$\F,\G\in\Mns{B,B'}$, one just get 
	$\DCPA{\F}{\G}\circ\gamma_{\G}=\gamma_{\F}\circ\DCPA{\F}{\G}^\prime$.\\

If $H^1,H^2:\C\to\C'$ are 1--morphisms of pre--Coxeter objects in $\XX$,
a {\it 2--morphism} $u:H^1\Rightarrow H^2$ is likewise a morphism of the
natural transformations of the corresponding functors $\ID\to\diagr{D}
{\XX}$. Specifically, $u$ consists of the datum of a diagrammatic 2--morphism
$u_B:H^1_B\to H^2_B$ for any $B\subseteq D$ such that, for any
$B'\subseteq B$, the following cylinder in $\diagr{D}{\XX}$ is commutative
\[\xymatrix@R=0.4cm@C=2cm{
&\C'_{B}\ar[dd]^{F'_{\sd{B}{B'}}}\\
\C_{B}\ar[dd]_{F_{\sd{B}{B'}}} \ar@/_15pt/|{H^2_{B}}[ur]  \ar@/^15pt/|{H^1_{B}}[ur]&\\
&\C'_{B'}\\
\C_{B'}\ar@/_15pt/|{H^2_{B'}}[ur]  \ar@/^15pt/|{H^1_{B'}}@{..>}[ur]&
}\]
where the rectangular 2--faces are the morphisms $\gamma,\gamma'$
and the circular ones the morphisms $u_B,u_{B'}$.\\

	\noindent\remark
	In view of \ref{ss:alpha-reduction}, a 2--morphism 
	$u:H^1\to H^2$ is equivalently described as a collection
	of $2$--morphisms $u_B:H^1_B\to H^2_B$, indexed by 
	$B\subseteq D$, such that, for any $B'\subseteq B$ and 
	$\F\in\Mns{B,B'}$, it holds $\gamma^2_{\F}\circ u_B=u_{B'}\circ\gamma^1_{\F}$.

\subsection{$\DCPA{}{}$--strict pre--Coxeter objects} \label{ss:dcp-str-cox-cat}
%-------------------------------------------------

A pre--Coxeter object $\C$ in $\XX$ is 
$\DCPA{}{}$--\emph{strict} if,
for any $B'\subseteq B$, and $\F,\G\in\Mns{B,B'}$ the following holds
\[F_\F=F_\G\aand\Upsilon_{\F\G}=\id_{F_\G}\]
We denote the common value of $\{F_\F\}_{\F\in\Mns{B,B'}}$ by ${F}_{B'B}:\C_B\to
\C_{B'}$. It follows from condition (4) in Definition \ref{ss:alpha-reduction} that,
for any $B''\subseteq B'\subseteq B$, $\F,\G\in\Mns{B,B'}$, $\F',\G'\in\Mns{B',B''}$, 
$\redasso{\F}{\F'}=\redasso{\G}{\G'}$. We denote the common value of $\{\redasso{\F}{\F'}\}_{\F,\F'}$ 
by $\oalpho{B}{B'}{B''}:{F}_{B''B'}\circ{F}_{B'B}\Rightarrow{F}_{B''B}$.

\begin{proposition}\label{pr:U strict}\hfill
\begin{enumerate}
\item
A $\DCPA{}{}$--strict pre--Coxeter object of type $D$ in $\XX$ is equivalently described by the
datum of 
\begin{itemize}
	\item for any $B\subseteq D$, an object $\C_B\in\XX$
	\item for any $B'\subseteq B$, a $1$--morphism ${F}_{B'B}:\C_{B}\to\C_{B'}$ 
	\item for any $B''\subseteq B'\subseteq B$, a $2$--isomorphism 
	\[\oalpho{B}{B'}{B''}: {F}_{B''B'}\circ{F}_{B'B}\Rightarrow{F}_{B''B}\]
\end{itemize}
such that 
\begin{itemize}
	\item for any $B'\subseteq B$, ${F}_{BB}=\id_{\C_B}$ and 
	$\oalpho{B}{B}{B'}=\id_{{F}_{B'B}}=\oalpho{B}{B'}{B'}$
	\item the $2$--isomorphism $\oalpho{}{}{}$ are associative, \ie for any $B'''\subseteq
	B''\subseteq B'\subseteq B$, 
	\[\oalpho{B}{B''}{B'''}\circ\oalpho{B}{B'}{B''}=\oalpho{B}{B''}{B'''}\circ\oalpho{B'}{B''}{B'''}\]
	as 2--morphisms ${F}_{B'''B''}\circ{F}_{B''B'}\circ{F}_{B'B}\Rightarrow{F}_{B'''B}$
\end{itemize}
\item
Every pre--Coxeter object $\C$ in $\XX$ is equivalent to a $\DCPA{}{}$--{strict} 
pre--Coxeter object in $\XX$.
\end{enumerate}
\end{proposition}
\begin{pf}
(1) is clear. 
(2) For any $B'\subseteq B$, choose a maximal nested set $\E(B,B')\in\Mns{B,B'}$.
We denote by $\ol{\C}$ the $\DCPA{}{}$--{strict} pre--Coxeter object with $\ol{\C}
_B\coloneqq\C_B$, ${F}_{B'B}\coloneqq F_{\E(B,B')}$, and 
\[\oalpho{B}{B'}{B''}\coloneqq
\DCPA{\E(B,B'')}{, \E(B,B')\cup\E(B',B'')}\circ\redasso{\E(B,B')}{\E(B',B'')}\]
Then, there is a canonical equivalence of pre--Coxeter objects
$\C\to\ol{\C}$ with $H_{B}\coloneqq\id_{\C_B}$ and 
$\gamma_{\F}\coloneqq\DCPA{\E(B,B')}{, \F}$ for $\F\in\Mns{B,B'}$.
\end{pf}

\noindent\remark We show in Sections \ref{ss:cat-O-infty} and \ref{ss:qC-on-qG} 
that Kac--Moody algebras and their quantum groups naturally give rise to 
$\DCPA{}{}$--{strict}  pre--Coxeter objects in $\Cat^{\ten}$. On the other
hand, we prove in \cite{ATL3} that the monodromy of the Casimir connection 
of a symmetrisable \KMA naturally gives rise to a pre--Coxeter structure
which is not $\DCPA{}{}$--{strict}. The latter, however, is $\sfa$--strict in the
following sense.

\subsection{$\redasso{}{}$--strict pre--Coxeter objects} \label{ss:a-str-cox-cat}
%---------------------------------------------------------------------

A pre--Coxeter object $\C$ in $\XX$ is $\redasso{}{}$--\emph{strict} if, for
any $B''\subseteq B'\subseteq B$, $\F\in\Mns{B,B'},\F'\in\Mns{B',B''}$, and
\[F_{\F'\cup\F}=F_{\F'}\circ F_{\F}\aand\redasso{\F}{\F'}=\id_{F_{\F'}\circ F_{\F}}\]
In contrast with Proposition \ref{pr:U strict}, not every pre--Coxeter object $\C$
is equivalent to an $\sfa$--strict one. We give, however, a sufficient condition
for that to be the case below. Let $B_1'\subset B_1\perp B_2\supset B'_2$,
with $|B_k\setminus B_k'|=1$, and denote by $\F_k$ the unique element in 
$\Mns{B_k, B'_k}$. Consider the diagram 
\begin{equation}\label{eq:rhombus-cat}
\xymatrix{
	& \C_{B_1\sqcup B_2} \ar[dl]_{F_{(\F_1, B_2)}} \ar|{F_{(\F_1,\F_2)}}[dd]\ar[dr]^{F_{(B_1, \F_2)}}  & \\
		\C_{B'_1\sqcup B_2} \ar[dr]_{F_{(B'_1, \F_2)}} & & \C_{B_1\sqcup B'_2} \ar[dl]^{F_{(\F_1, B'_2)}} \\
		& \C_{B'_1\sqcup B'_2} & 
	}
\end{equation}
where the triangular $2$--faces are given by the vertical joins $\redasso{(\F_1, B_2)}
{(B'_1,\F_2)}$ and $\redasso{(B_1,\F_2)}{(\F_1, B'_2)}$ respectively. We say that
\eqref{eq:rhombus-cat} is {\it trivial} if
\begin{gather}
F_{(B'_1,\F_2)}\circ F_{(\F_1,B_2)}=F_{(\F_1, B'_2)}\circ F_{(B_1,\F_2)}
\label{eq:triv one}\\
\redasso{(\F_1, B_2)}{(B'_1,\F_2)}=\redasso{(B_1,\F_2)}{(\F_1,B'_2)}
\label{eq:triv two}
\end{gather}
as 2--morphisms $F_{(B'_1,\F_2)}\circ F_{(\F_1,B_2)}=F_{(\F_1, B'_2)}\circ F_{(B_1,\F_2)}
\Rightarrow F_{(\F_1,\F_2)}$. Note that this is the case if $\C$ is $\sfa$--strict.
 
\begin{proposition}
Let $\C$ be a pre--Coxeter object in $\XX$. If the diagrams \eqref{eq:rhombus-cat}
are trivial, then $\C$ is canonically equivalent to an $\redasso{}{}$--strict pre--Coxeter
object.
\end{proposition}
\begin{pf}
Retain the notation from \ref{ss:chain-ns}. Let $B'\subseteq B$, $\F\in\Mns{B,B'}$
and ${\bf C}: B'=B_0\subsetneq B_1\cdots\subsetneq B_{\ell}=B$ a maximal chain
corresponding to $\F$. Denote by $\ol{F}_{\F}:\C_B\to\C_{B'}$ the composition
$F_{\F_1}\circ\cdots\circ F_{\F_\ell}$, where $\F_k$ is the unique element in
$\Mns{B_k, B_{k-1}}$. By \eqref{eq:triv one}, $\ol{F}_{\F}$ does not depend
upon the choice of ${\bf C}\in p^{-1}(\F)$. Moreover, for any $\F\in\Mns{B,B'}$
and $\F'\in\Mns{B',B''}$, one has $\ol{F}_{\F'}\circ\ol{F}_{\F}=\ol{F}_{\F'\cup\F}$.

For any $\F\in\Mns{B,B'}$, let $u_{\F}:\ol{F}_{\F}\Rightarrow F_{\F}$
be the $2$--morphism obtained as the composition of vertical joins
$\redasso{\F_2\cup\cdots\cup\F_\ell}{\F_1}
\circ\cdots\circ
\redasso{\F_{\ell-1}\cup\F_\ell}{\F_{\ell-2}}
\circ
\redasso{\F_\ell}{\F_{\ell-1}}$. By \eqref{eq:triv two}, $u_{\F}$ is 
independent of the choice of a maximal chain ${\bf C}\in p^{-1}(\F)$.
For any $\F,\G\in\Mns{B,B'}$, set
\[\oDCPA{\F}{\G}\coloneqq u_{\F}^{-1}\circ\DCPA{\F}{\G}\circ u_{\G}:\ol{F}_{\G}\Rightarrow\ol{F}_{\F}\]
Then, the datum of the objects $\ol{\C}_B=\C_B$, 
$1$--morphisms $\ol{F}_{\F}$, and $2$--morphisms $\oDCPA{\F}{\G}$
gives rise to an $\redasso{}{}$--strict pre--Coxeter object $\ol{\C}$. Moreover, there is a canonical equivalence
$\C\to\ol{\C}$ with $H_B=\id_{\C_B}$ and $\gamma_{\F}=u_{\F}^{-1}$.
\end{pf}

\subsection{Generalised braid groups}
%------------------------------------------------

\begin{definition}
A {\it labelling} $\ulm$ of a diagram $D$ is the assignment of an integer
$m_{ij}\in\{2,3,\ldots,\infty\}$ to any pair $i,j$ of distinct vertices of $D$ such
that $m_{ij}=m_{ji}$ and $m_{ij}=2$ if $i$ and $j$ are orthogonal.

The \emph{generalised} braid group corresponding to $D$ and a labelling
$\ulm$ is the group $\BDm$ with generators $\{\sfS_i\}_{i\in D}$ and
relations
\begin{equation}\label{eq:gen-braid}
\underbrace{\sfS_i\cdot \sfS_j\cdot \sfS_i\;\cdots\;}_{m_{ij}}=\underbrace{\sfS_j\cdot \sfS_i\cdot \sfS_j\;\cdots\;}_{m_{ij}}
\end{equation}
If $B\subseteq D$ is a subdiagram, we denote by $\BBm\subseteq\BDm$
the subgroup generated by the elements $\sfS_i$, $i\in B$, which is isomorphic
to the generalised braid group corresponding to $B$ and the labelling $\ulm$ restricted to $B$.
\end{definition}

\subsection{Coxeter objects}\label{ss:cox}
%-----------------------------------

Let $(D,\ulm)$ be a labelled diagram.

\begin{definition}
A {\it Coxeter object} of type $(D,\ulm)$ in $\XX$ is the datum of 
\begin{itemize}
\item a pre--Coxeter object $\left(\C_B,F_{B'B},\alpho{B''}{B'}{B}\right)$ of type $D$ in $\XX$
\item for any $i\in D$, a diagrammatic 2--isomorphism $S_i:F_{\emptyset i}\Rightarrow F_{\emptyset i}$
\end{itemize}
such that for any subdiagram $B\subseteq D$, and $i,j\in B$ with $i\neq j$
\begin{equation}\label{eq:cox-braid-rel}
\underbrace{S^B_i\cdot S^B_j\cdot S^B_i\;\cdots\;}_{m_{ij}}=\underbrace{S^B_j\cdot S^B_i\cdot S^B_j\;\cdots\;}_{m_{ij}}
\end{equation}
where $\SC{i}{B}:F_{\emptyset B}\Rightarrow F_{\emptyset B}$ is the diagrammatic $2$--morphism 
\[
\xymatrix@=1.2cm{
F_{\emptyset B}\ar@{=>}[r]^{(\alpho{\emptyset}{i}{B})^{-1}} & F_{\emptyset i}\circ F_{i B} 
\ar@{=>}[r]^{\SC{i}{}}& 
F_{\emptyset i}\circ F_{i B} \ar@{=>}[r]^{\alpho{\emptyset}{i}{B}} &  F_{\emptyset B}
}
\]
\end{definition} 

	\noindent\remark
	More explicitly, the equation \eqref{eq:cox-braid-rel} reads as follows.
	Let $\F,\G\in\Mns{B}$ be two maximal nested sets on $B$ such that $\{i\}\in\F$,
	$\{j\}\in\G$, so that $\G=\G_j\cup\G'$, with $\G_{j}=\{\emptyset, \{j\}\}$. 
	Let $\xi_{\G}^{j}:\sfEnd{F_{\emptyset j}}\to\sfEnd{F_{\G}}$ be the natural 
	isomorphism induced by the map $\redasso{\G'}{\G_j}:F_{\G_j}\circ F_{\G'}\Rightarrow F_{\G}$,
	and set $\xi^i_{\G\F}\coloneqq\sfAd{\DCPA{\G}{\F}}\circ\xi^i_{\F}$, so that
	$\xi^i_{\G\F}:\sfEnd{F_{\emptyset i}}\to\sfEnd{F_{\G}}$. 
	Then, \eqref{eq:cox-braid-rel} reads
	\begin{equation*}
		\underbrace{\xi_{\G\F}^i({S}_i)\cdot \xi^j_{\G}({S}_j)\cdot\xi_{\G\F}^i({S}_i)\;\cdots\;}_{m_{ij}}=
		\underbrace{\xi^j_{\G}({S}_j)\cdot\xi_{\G\F}^i({S}_i)\cdot 
			\xi^j_{\G}({S}_j)\;\cdots\;}_{m_{ij}}\ ,
	\end{equation*}
	as an identity in $\sfEnd{F_{\G}}$.\\

A $1$--morphism $\C\to\C'$ of Coxeter objects in $\XX$ is one of the 
underlying pre--Coxeter objects, which preserves the braid group operators $S$. 
That is, it consists of a datum $(H_B,\gamma_{B'B})$ defined as in \ref{ss:mor-pre-cox}
such that, for any $i\in D$, 
\[\sfAd{\gamma_{\emptyset i}}(H_{\emptyset}(S_i))=S'_i|_{H_i}\]
in $\diagr{D}{\XX}(F'_{\emptyset i}\circ H_i, F'_{\emptyset i}\circ H_{i})$.
A $2$--morphism is defined as in \ref{ss:mor-pre-cox}.

\subsection{Braid group actions}\label{ss:braid-gp-act}
%----------------------------------------

Let $(D,\ulm)$ be a labelled diagram. , and $H:\C\to\D$ a $1$--isomorphism of Coxeter objects.

\begin{proposition}\hfill
\begin{enumerate}
\item
Let $\C$ be Coxeter object of type $(D,\ulm)$ in $\XX$.
For any subdiagram $B\subseteq D$, there is a unique homomorphism
$\rho_B^{\C}:\BBm\to \diagr{D}{\XX}(F_{\emptyset B},F_{\emptyset B})$, 
such that, for any $i\in B$, $\rho^{\C}_{B}(\SC{i}{})=S^B_i$. Moreover, for any 
$B'\subseteq B$, the following diagram is commutative
\[
\xymatrix{
\BBm \ar[r]^(.3){\rho_B} & \diagr{D}{\XX}(F_{\emptyset B},F_{\emptyset B})\\
\BBpm \ar[r]_(.3){\rho_{B'}} \ar[u] & \diagr{D}{\XX}(F_{\emptyset B'},F_{\emptyset B'}) 
\ar[u]
}
\]
where the vertical right arrow is induced by the $2$--isomorphism
$\alpho{B}{B'}{\emptyset}: 
F_{\emptyset B'}\circ F_{B'B}\Rightarrow F_{\emptyset B}$.
\item Let $\C,\D$ be Coxeter objects of type $(D,\ulm)$ in $\XX$ and
$H:\C\to\D$ a $1$--isomorphism of Coxeter objects. For any subdiagram
$B\subseteq D$, the representations $\rho^{\C}_B$ and $\rho^{\D}_B$
of $\BBm$ are equivalent, \ie the following diagram is commutative
\[
\xymatrix@R=0.5cm{
	& \diagr{D}{\XX}(F^{\C}_{\emptyset B},F^{\C}_{\emptyset B})\\
	\BBm \ar[ur]^(.3){\rho^{\C}_B} \ar[dr]_(.3){\rho^{\D}_{B}} &\\
	& \diagr{D}{\XX}(F^{\D}_{\emptyset B},F^{\D}_{\emptyset B}) \ar[uu]\\
}
\]
where the vertical arrow is induced by the $2$--isomorphism 
$\gamma_B: F^{\D}_{\emptyset B}\circ H_B\Rightarrow F^{\C}_{\emptyset B}$.
\end{enumerate}
\end{proposition}

\begin{pf}
(1) The existence of the homomorphisms $\rho_B$, $B\subseteq D$, 
follows by construction. For the commutativity of the diagram, it is enough to observe that the map 
$\diagr{D}{\XX}(F_{\emptyset B'},F_{\emptyset B'})\to\diagr{D}{\XX}(F_{\emptyset B},F_{\emptyset B})$
sends a $2$--endomorphism $\phi$ to $(\alpho{B}{B'}{\emptyset})\circ\phi|_{F_{B'B}}\circ(\alpho{B}{B'}{\emptyset})^{-1}$. Therefore, for any $i\in B'$, one has
\[\begin{split}
(\alpho{B}{B'}{\emptyset})\circ\SC{i}{B'}\circ(\alpho{B}{B'}{\emptyset})^{-1}
&=
(\alpho{B}{B'}{\emptyset})\circ\left((\alpho{B'}{i}{\emptyset})\circ\SC{i}{}
\circ(\alpho{B'}{i}{\emptyset})^{-1}\right)|_{F_{B'B}}\circ(\alpho{B}{B'}{\emptyset})^{-1}\\
&=(\alpho{B}{i}{\emptyset})\circ\SC{i}{}\circ(\alpho{B}{i}{\emptyset})^{-1}\\
&=\SC{i}{B}
\end{split}\]
where the second equality follows from the associativity of $\alpha$.
(2) follows immediately from the definition of $1$--morphism of Coxeter
objects (cf. \ref{ss:cox}).
\end{pf}
\noindent
\remark
In the $2$--category $\XX$, the representations $\rho_B$ are described as
follows. For any $B\subseteq D$ and $\F\in\Mns{B}$, there is a collection of
homomorphisms $\rho_{\F}:\BBm\to\mathsf{Aut}_{\XX}(F_{\F})$, $\F\in\Mns
{B}$, uniquely determined by the conditions
\begin{itemize}
\item $\rho_{\F}(\sfS_i)=S_i^{\F}$, if $\{i\}\in\F$
\item $\rho_{\G}=\sfAd{\DCPA{\G}{\F}}\circ\rho_{\F}$
\end{itemize}

\subsection{Lax diagrammatic algebras \cite[Sec. 3]{vtl-4}}\label{ss:D-alg}
%---------------------------------------

A {\it lax diagrammatic algebra}\footnote{The terminology adopted here differs
from the one in \cite{vtl-4}, where the adjective lax is not used in particular. In
the present paper, we reserve the term diagrammatic algebra for a lax diagrammatic
algebra such that $m_B\circ i_{BB'}\otimes i_{BB''}:A_{B'}\otimes A_{B''}\to A_
B$ is an isomorphism for any $B=B'\sqcup B''$, which implies in particular that
$A_\emptyset=\sfk$ (see Remark \ref{rk:diagrammatic}).} is the datum of
\begin{itemize}
\item for any $B\subseteq D$, a $\sfk$--algebra $A_B$
\item for any $B'\subseteq B$, a homomorphism $i_{BB'}:A_{B'}\to A_B$
\end{itemize}
such that 
\begin{itemize}
\item for any $B\subseteq D$, $i_{BB}=\id_{A_B}$
\item for any $B''\subseteq B'\subseteq B$, $i_{BB'}\circ i_{B'B''}=i_{BB''}$
\item for any $B=B'\sqcup B''$, with $B'\perp B''$, $m_B\circ i_{BB'}\otimes i_{BB''}$
is a morphism of algebras $A_{B'}\ten A_{B''}\to A_{B}$, where $m_B$ denotes the
multiplication in $A_B$.
\end{itemize}

A morphism of lax diagrammatic algebras $\varphi:A\to A'$ is a collection of
homomorphisms $\varphi_B:A_B\to A'_B$ such that $\varphi_B\circ i_{BB'}
=i'_{BB'}\circ\varphi_{B'}$ for any $B'\subseteq B$.\footnote{In \cite{vtl-4}, a
morphism of lax diagrammatic algebras is referred to as a \emph{strict} morphism.} 

\subsection{Pre--Coxeter categories from lax diagrammatic algebras}\label{ss:Cox D-alg}
%------------------------------------------------------------------------------

A lax diagrammatic algebra $A$ gives rise to an $(\oalpho{}{}{},\DCPA{}{})$--strict pre--Coxeter object $\C=\Rep
(A)$ in $\XX=\Cat$ given by\footnote{Note that the commutativity of $A_{B'}, A_{B''}$
in $A_{B}$, for any $B',B''\subseteq B$ with $B'\perp B''$, has no relevance in the
above construction of pre--Coxeter structure on $\Rep(A)$. On the other hand, this
feature is particularly convenient in the construction of examples arising from 
the quantisation of Lie bialgebras (cf. Section \ref{s:qcstructure}, in particular Lemma
\ref{sss:twist}).}
\begin{itemize}
\item For any $B\subseteq D$, $\C_B=\Rep(A_B)$
\item For any $B'\subseteq B$, 
${F}_{B'B}: \C_B\to\C_{B'}$ is the pullback functor $i^*_{BB'}$ 
\end{itemize}
Moreover, a morphism of lax diagrammatic algebras $\varphi:A\to A'$ gives rise to a 
morphism of pre--Coxeter objects $\Rep(A')\to\Rep(A)$.

If $(D,\ulm)$ is a labelled diagram, the group algebra $\sfk\B_{D}^{\ulm}$ is
naturally endowed with a lax diagrammatic algebra structure. If a lax diagrammatic algebra $A$ is
further endowed with a morphism of lax diagrammatic algebras $\rho_B:\sfk\BBm\to A_B$, $B\subseteq D$,
then the elements 
$\rho(\SC{i}{})\in A_i=\sfEnd{{F}_{\emptyset i}}$ give rise to the structure of
Coxeter object on $\Rep(A)$.

This construction can be generalised by replacing the categories $\Rep(A_B)$
by a collection of subcategories $\C_B\subseteq \Rep(A_B)$ stable under 
restrictions, and $\rho$ by a morphism of lax diagrammatic algebras $\sfk\B_D^{\ulm}\to\sfEnd{F_
{\emptyset D}}=:\wh{A}$. We show in Section \ref{s:qkm} that an example
of such Coxeter objects is provided by quantum Weyl groups of quantised
Kac--Moody algebras.

\subsection{Topological definition}\label{ss:top def}
%------------------------------------------

In \cite{fink-sch}, Finkelberg and Schechtman propose an alternative definition
of a (pre--)Coxeter object in $\Cat$ for Dynkin diagrams of finite type, which is
akin to Deligne's topological definition of a braided monoidal category \cite{Del}.
This is given by a category $\C_B$ for every diagram $B\subseteq D$, together with 
\begin{itemize}
\item for any $B'\subseteq B$, a Weyl group equivariant local system of restriction
functors $\mathfrak{F}_{B'B}:\C_B\to\C_{B'}$, defined over $(\h_{B/B'})_{\operatorname{reg}}$
\footnote{
Here, $\h_B$ is the Cartan subalgebra of $\g_B\subseteq\g_D$, 
$\h_{B/B'}\subseteq\h_{B}$ is the orthogonal complement of $\h_{B'}$, and
$(\h_{B/B'})_{\operatorname{reg}}$ is the complement in $\h_{B/B'}$ to the root
hyperplanes in $\h_B$ not containing $\h_{B/B'}$. 
}
\item for any $B''\subseteq B'\subseteq B$, a suitable analogue of the factorisation
isomorphism $\alpho{B}{B'}{B''}$.
\end{itemize} 

\noindent
This gives rise to a Coxeter object in the sense of \ref{ss:pre-cox}, where,
for each $\F\in\Mns{B,B'}$, the functor $F_{\F}:\C_B\to\C_{B'}$, $\F\in\Mns
{B,B'}$, is the limit of $\mathfrak{F}_{B'B}$ at the point at infinity $p_\F$ in the
De Concini--Procesi compactification of $(\h_{B/B'})_{\operatorname{reg}}$
\cite{DCP2}.

\subsection{Example: rational Cherednik algebras}
%----------------------------------------------------------------

Let $\h$ be a \fd complex vector space, and $W\subset GL(\h)$ a finite
complex reflection group. Let $c$ be a conjugation invariant function on
the set $S$ of reflections in $W$, and  $H_c(W,\h)$ the corresponding
rational Cherednik algebra. Let $\O(W,\h)$ be the category of highest
weight $H_c(W,\h)$--modules, $W'\subset W$ a parabolic subgroup,
$\h'=\h/\h^{W'}$, $c'$ the restriction of $c$ to $S\cap W'$.

In \cite{beze}  Bezrukavnikov and Etingof construct a parabolic restriction
functor 
\[
\Res_b:\O(W,\h)\to\O(W',\h')
\]
where $b\in\hreg^{W'}$. In \cite[Cor. 2.5]{shan}, Shan shows that the composition
of two parabolic restriction functors is isomorphic to a parabolic restriction
functor, compatibly with the parameter $b$. 
If $W$ is a Weyl group with Dynkin diagram $D$, these functors and their
factorisation isomorphisms give rise to topological Coxeter object in $\Cat$, 
in the sense sketched in \ref{ss:top def}.

%%%%%%%%%%%%%%%%%%%%%
% BRAIDED COXETER CATEGORIES %
%%%%%%%%%%%%%%%%%%%%%

\section{Braided Coxeter categories}\label{s:qcqtqba}
%==========================

\subsection{}
\label{ss:braid-qC-cat}
%---------------------------

Denote by $\Cat^{\ten}$ (resp. $\Cat^{\ten,\beta}$) the $2$--category
of monoidal (resp. braided monoidal) categories.

\begin{definition}\label{def:qc-cat}
Let $D$ be a diagram.
\begin{enumerate}
\item A \emph{braided pre--Coxeter category of type $D$} is a tuple $(\C_B,
F_{B'B},\alpho{B''}{B'}{B})$ such that 
\begin{itemize}
\item $\C_B$ is a $B$--object in $\Cat^{\ten,\beta}$
\item $(\C_B, F_{B'B},\alpho{B''}{B'}{B})$ is a pre--Coxeter object in $\Cat^
{\ten}$
\end{itemize}
\item If $\ulm$ is a labelling on $D$, a \emph{braided Coxeter category of
type $(D,\ulm)$} is a tuple $(\C_B,F_{B'B},\alpho{B''}{B'}{B},S_i)$ such that
\begin{itemize}
\item $\C_B$ is a $B$--object in $\Cat^{\ten,\beta}$
\item $(\C_B,F_{B'B},\alpho{B''}{B'}{B})$ is a pre--Coxeter object in $\Cat^{\ten}$
\item $(\C_B,F_{B'B},\alpho{B''}{B'}{B},S_i)$ is a Coxeter object in $\Cat$
\end{itemize}
and, for any $i\in D$, the following holds in $\sfAut{F_i\ten F_i}$
\begin{equation}\label{eq:coprod-id}
J_i^{-1}\circ 
F_i(c_i)\circ\Delta(S_i)\circ J_i=
%J_i\circ
c_{\emptyset}\circ S_i\ten S_i
\end{equation}
where $F_i=F_{\emptyset i}$, $J_i$ is the tensor structure on $F_i$ and $c_i, c_{\emptyset}$ are
the opposite braidings in $\C_i$ and $\C_{\emptyset}$, respectively.\footnote{
In a braided monoidal category with braiding $\beta$, the opposite
braiding is $\beta^{\scs\operatorname{op}}_{X,Y}:=\beta_{Y,X}^{-1}$.
} In other words, the following diagram is commutative for any $V,W
\in\C_i$, 
\[
\xymatrix{
F_i(V)\ten F_i(W) \ar[d]_{J_i^{V,W}} \ar[r]^{S^V_i\ten S^W_i} & F_i(V)\ten F_i(W)  \ar[r]^{c_{\emptyset}} & F_i(W)\ten F_i(V)
\ar[d]^{J_i^{W,V}}\\
F_i(V\ten W) \ar[r]_{S_i^{V\otimes W}} & F_i(V\ten W)  \ar[r]_{F_i(c_i)} & F_i(W\ten V)  
}
\]
\item A functor of braided Coxeter categories $\C\to\C'$ is a tuple 
$(H_B,\gamma_{B'B})$ such that
\begin{itemize}
\item $H_B:\C_B\to\C'_{B}$ is a $1$--morphism of $B$--objects in $\Cat^{\ten,\beta}$;
\item $(H_B,\gamma_{B'B})$ is a $1$--morphism of pre--Coxeter objects in $\Cat^{\ten}$.
\end{itemize}
Finally, a natural transformation $u:H\Rightarrow H'$ is a $2$--morphism of $B$--objects
in $\Cat^{\ten,\beta}$.
\end{enumerate}
\end{definition}

\noindent\remarks\hfill
\begin{itemize}\itemsep0.25cm
\item
The identity \eqref{eq:coprod-id} relates the 
failure of $(F_i,J_i)$ to be a braided monoidal functor and that of $S_i$
to be a monoidal isomorphism. That is, if \eqref{eq:coprod-id} holds, then
$S_i$ is monoidal if and only if $J_i$ is braided. Conversely, if $S_i$ is
monoidal and $J_i$ is braided, then \eqref{eq:coprod-id} automatically
holds. In particular, every Coxeter object in $\Cat^{\ten,\beta}$ is a braided
Coxeter category.
\item
The main examples of braided Coxeter categories arise as representations
of a \emph{quasi--Coxeter quasitriangular quasibialgebra}, as defined in
\cite[Sec. 3]{vtl-4}.
\item 
In \cite{ATL2}, we only consider {$\sfa$--strict} braided Coxeter categories
and, for simplicity, refer to them as braided Coxeter categories.
\end{itemize}

\subsection{Unfolded definition}
In view of \ref{prop:alpha-reduction},	braided Coxeter category of
type $(D,\ulm)$ is equivalently described by the datum of 
\begin{itemize}
	\item for any $B\subseteq D$, a braided monoidal category $\C_B\in\XX$
	\item for any $B'\subseteq B$ and $\F\in\Mns{B,B'}$, a 
	(not necessarily braided) monoidal functor $F_{\F}:\C_{B}\to\C_{B'}$
	\item for any $B'\subseteq B$ and $\F,\G\in\Mns{B,B'}$,
	an isomorphism of monoidal functors $\DCPA{\G}{\F}:F_{\F}\Rightarrow F_{\G}$
	\item for any $B''\subseteq B'\subseteq B$, $\F\in\Mns{B,B'}$ and 
	$\F'\in\Mns{B',B''}$, an isomorphism of monoidal functors
	$\redasso{\F}{\F'}: F_{\F'}\circ F_{\F}\Rightarrow F_{\F'\cup\F}$
	\item for any $i\in D$, an isomorphism of functors $S_i:F_{\emptyset\{i\}}\Rightarrow F_{\emptyset\{i\}}$
	(not necessarily preserving the tensor structure)
\end{itemize}
satisfying the properties listed in \ref{ss:alpha-reduction}, \ref{ss:cox}, and the 
coproduct identity \eqref{eq:coprod-id}.

\subsection{Balanced categories}
%-----------------------------------------

In \cite{fink-sch}, the coproduct identity \eqref{eq:coprod-id} is replaced by
the assumption that the categories $\C_i$ are balanced categories (in fact,
that $\C_B$ is balanced for any $B\subseteq D$). We point out below that,
in general, this assumption is stronger than \eqref{eq:coprod-id}.

Recall that a braided monoidal category $(\C,\ten,b,\Phi)$ is \emph{balanced}
if there is a $\theta\in\sfAut{\id_{\C}}$ such that 
\begin{equation}\label{eq:balance-4}
\theta_{V\ten W}=b_{W,V}\circ b_{V,W}\circ\theta_V\ten\theta_W
\end{equation}
for any $V,W\in\C$.
\begin{proposition}\label{pr:balance}
Let $\C$ be a braided Coxeter category such that
\begin{enumerate}
\item $\C_{\emptyset}$ is symmetric
\item $S_i^2=F_i(\theta_i)$ for some $\theta_i\in\sfAut{\id_{\C_i}}$
\item $F_i:\C_i\to\C_\emptyset$ is faithful
\end{enumerate}
Then $\C_i$ is a balanced monoidal category with balance $\theta_i$.
\end{proposition}
\begin{pf}
Squaring the \rhs of \eqref{eq:coprod-id} yields
\[\left(c_{\emptyset}\circ S_i\ten S_i\right)^2=
c_{\emptyset}^2 \circ S_i^2\ten S_i^2=
F_i(\theta_i)\otimes F_i(\theta_i)\]
where we used the binaturality of $c_\emptyset$ and the assumptions
(1) and (2). On the other hand, the square of the \rhs of \eqref{eq:coprod-id}
is equal to
\[\begin{split}
J_i^{-1}\circ F_i(c_i)\circ\Delta(S_i)\circ F_i(c_i)\circ\Delta(S_i)\circ J_i
&=
J_i^{-1}\circ F_i(c_i^2)\circ\Delta(S_i^2)\circ J_i\\
&=
J_i^{-1}\circ F_i(c_i^2)\circ F_i(\theta_i\circ\otimes)\circ J_i
\end{split}\]
where we used the naturality of $S_i$. Since $J_i\circ F_i(\theta_i)
\otimes F_i(\theta_i)\circ J_i^{-1}=F_i(\theta_i\otimes\theta_i)$ by
naturality of $J_i$, we get
\[ F_i(c_i^2 \circ  \theta_i\circ\otimes)=F_i(\theta_i\otimes\theta_i)\]
hence the required result since $F_i$ is faithful.
\end{pf}

\noindent\remark
\begin{itemize}
\item The converse of Proposition \ref{pr:balance} does not hold in general.
That is, the existence of a balance does not imply \eqref{eq:coprod-id}. 
Instead, the correct categorical interpretation of \eqref{eq:coprod-id} corresponds to
the braided monoidal categories $\C_i$ (with the tensor functors $F_i$) being \emph{half--balanced} 
(cf. \cite[Sec. 4]{tingley-sny}).
\item
Finally, we note that the coproduct identity \eqref{eq:coprod-id} cannot in general
be extended to subdiagrams with more than one vertex. Specifically, in the examples of braided Coxeter 
structures described in Sections \ref{s:qcstructure} and \ref{s:qkm}, the categories $\C_B$, 
with $|B|>1$, do not in general admit a half--balanced structure.
\end{itemize}

%%%%%%%%%%%%%%%%%%%%%
% DIAGRAMMATIC LIE BIALGEBRAS  %
%%%%%%%%%%%%%%%%%%%%%

\section{Diagrammatic Lie bialgebras}\label{s:lba}
%===========================

In this section, we introduce the notion of a diagrammatic Lie bialgebra
$\b$. We then show that \DYt modules over the canonical subalgebras 
of $\b$ give rise to a symmetric  pre--Coxeter category.

\subsection{Lie bialgebras \cite{drin-2}}\label{ss:lba}
%---------------------------------

A Lie bialgebra is a triple $(\b,[\,,\,]_{\b}, \delta_{\b})$ where $(\b,[\,,\,]_\b)$
is a Lie algebra, $(\b,\delta_\b)$ a Lie coalgebra, and the cobracket $\delta_
\b:\b\to\b\otimes\b$ satisfies the cocycle condition
\[\delta_\b([X,Y]_\b)=\ad(X)\,\delta_\b(Y)-\ad(Y)\,\delta_b(X)\]

\subsection{Manin triples \cite{drin-2,ek-1}}\label{ss:lbamt}
%------------------------------------------------------

A Manin triple is the data of a Lie algebra $\g$ with 
\begin{itemize}
\item a nondegenerate invariant symmetric bilinear form $\iip{\cdot}{\cdot}$
\item isotropic Lie subalgebras $\gupm\subset\g$
\end{itemize}
such that 
\begin{itemize}
\item $\g=\gum\oplus\gup$ as vector spaces
\item the inner product defines an isomorphism $\gup\to\gum^*$  
\item the Lie bracket of $\g$ is continuous with respect to the topology 
obtained by putting the discrete and the weak topologies on $\gum$
and $\gup$  respectively. Equivalently, the bracket on $\gup$ is
continuous \wrt the weak topology.
\end{itemize}

Under these assumptions, the commutator on $\gup\simeq\gum^*$
induces a cobracket $\delta:\gum\to\gum\otimes\gum$ which satisfies
the cocycle condition, thus endowing $\gum$ with a Lie bialgebra
structure. In general, however, $\gup$ is only a topological
Lie bialgebra.

One can similarly consider {\it restricted} Manin triples, where
\begin{itemize}
\item $\g$ is $\IZ$--graded as a Lie algebra, with \fd %homogeneous
components $\{\g_n\}_{n\in\IZ}$
\item the inner product satisfies $\iip{\g_n}{\g_m}=0$ unless $n+m=d$,
for a given $d\in\IZ$
\item $\g=\gum\oplus\gup$ as vector spaces, with the isotropic subalgebra
$\b_-$ (resp. $\b_+$) concentrated in non--negative (resp. non--positive)
degrees
\end{itemize}
In this case, the inner product induces an isomorphism $\b_\pm\to\b
_\mp^{\star}$, where $\b_\mp^{\star}=\bigoplus_{n}(\b_{\mp,n})^*$ is
the restricted dual of $\b_\mp$. The joint continuity of the bracket on
$\g$ is automatic, and both $\b_-$ and $\b_+$ are Lie bialgebras with
a cobracket of degree $d$.\footnote{Note that the Lie algebra grading
on $\b_+$ inherited from $\g$ differs from that induced by the identification
$\b_+\cong\b_-^\star$ by a shift since the inner product yields an
isomorphism $(\b_{-,n})^*\cong \b_{+,-n+d}$. Note also that the
isotropy of $\b_\pm$ implies that $\b_{-,n}=0$ if $n\leqslant d-1$ and
$\b_{+,n}=0$ if $n\geqslant d+1$.}

\subsection{Example}
%--------------------------

A \fd Lie algebra $\ll$ with an invariant inner product
$(-,-)$ gives rise to a restricted Manin triple as follows.
\[\g=\ll[t,t^{-1}]
\qquad
\b_-=\g[t]
\qquad
\b_+=t^{-1}\ll[t^{-1}]\]
with the standard grading $\deg(\ll\otimes t^m)=m$, and inner
product given by the residue pairing $\iip{f}{g}=\Res_{t=0}(f(t),
g(t))$, so that $\iip{X\otimes t^m}{Y\otimes t^n}=(X,Y)\delta_{m
+n,-1}$. In this case, $\b_-$ has a degree $d=-1$ cobracket
given by
\[\delta(f)(t,s)
=\left[f(t)\otimes 1+1\otimes f(s),\frac{\Omega}{s-t}\right]\]%\in\ll[t]\otimes\ll[s]\]
%\[\delta(f)=[f(t)\otimes 1+1\otimes f(s),\frac{\Omega}{s-t}]\]
where $\Omega\in(\ll\otimes\ll)^\ll$ corresponds to $\iip{\cdot}{\cdot}$.

The corresponding Manin triple is $\left(\ll((t^{-1})),\ll[t],t^{-1}\ll[[t^{-1}]]\right)$.

%--------------------------------------------------
\subsection{Drinfeld double \cite{drin-2}}\label{ss:drinf-double}
%--------------------------------------------------

The Drinfeld double of a Lie bialgebra $(\b,[\,,\,]_{\b},\delta_{\b})$ is the
Lie algebra $\gb$ defined as follows. As a vector space, $\gb=\b\oplus
\b^*$. The duality pairing $\b^*\ten\b\to\sfk$ extends uniquely to a symmetric,
non--degenerate bilinear form $\iip{\cdot}{\cdot}$ on $\gb$, \wrt which
both $\b$ and $\b^*$ are isotropic subspaces. The Lie bracket on $\gb$ is
defined as the unique bracket which coincides with $[\,,\,]_{\b}$ on $\b$,
with $\delta_{\b}^t$ on $\b^*$, and is compatible with $\iip{\cdot}{\cdot}$,
\ie satisfies $\iip{[x,y]}{z}=\iip{x}{[y,z]}$ for all $x,y,z\in\gb$. The mixed
bracket of $x\in\b$ and $\phi\in\b^*$ is then given by
\[[x,\phi]
=\sfad^*(x)(\phi)+\phi\ten\id_\b\circ\delta(x)
\]
where $\sfad^*$ is the coadjoint action of $\b$ on $\b^*$. $(\gb,\b,\b^*)$
is a Manin triple, and any such triple arises this way.

Similarly, if $\b$ is a Lie bialgebra which is $\IN$--graded with \fd components,
and such that the bracket and cobracket are homogeneous of degrees 0 and
$d\in\IZ$ respectively,\footnote{In the sequel, we shall abusively refer to such
a $\b$ as an $\IN$--graded Lie bialgebra.} the restricted double of $\b$ is
defined as $\grb=\b\oplus\b^{\star}[d]$, where $\b^{\star}[d]_n=(\b_{-n+d})^*$,
and is a restricted Manin triple.

%---------------------------------------------------------------------
\subsection{Drinfeld--Yetter modules \cite{ek-2}}\label{ss:pre-DY}
%---------------------------------------------------------------------

A \DYt module over a Lie bialgebra $\b$ is a triple $(V,\pi_V,\pi_V^*)$, where
$(V,\pi_V)$ is a left $\b$--module, $(V,\pi_V^*)$ a right $\b$--comodule, and
the maps $\pi_V:\b\otimes V\to V$ and $\pi_V^*:V\to\b\otimes V$ satisfy the
following relation in $\End(\b\ten V)$
\[\id_{\b}\ten\pi_V\circ(12)\circ\id_{\b}\ten\pi_V^*-\pi_V^*\circ\pi_V=
-[\cdot,\cdot]_{\b}\ten\id_V\circ\id_{\b}\ten\pi_V^*+\id_{\b}\ten\pi_V\circ\delta_{\b}\ten\id_V\]

The category $\DrY{\b}$ of \DYt modules over $\b$ is a symmetric tensor
category. For any $V,W\in\DrY{\b}$, the action and coaction on the tensor product 
$V\ten W$ are defined, respectively, by
\begin{align*}
\pi_{V\ten W}	&=\pi_{V}\ten\id_W+\id_V\ten\pi_{W}\circ(1\,2)\otimes\id_W\\
\pi^*_{V\ten W}	&=\pi^*_{V}\ten\id_W+(1\,2)\otimes\id_W\circ\id_V\ten\pi^*_{W}
\end{align*}
The associativity constraints are trivial, and the braiding is defined by
$\beta_{VW}=(1\,2)$.

%-----------------------------------------------------------------
\subsection{Representations of the Drinfeld double}\label{ss:drinf-db-rep}
%-----------------------------------------------------------------

The category $\DrY{\b}$ is canonically isomorphic to the category $\Eq
{\gb}$ of {\it equicontinuous} $\gb$--modules \cite{ek-1}, \ie those endowed
with a locally finite $\b^*$--action. This condition yields a functor $E:\Eq{\gb}
\to\DrY{\b}$, which assigns to any $V\in\Eq{\gb}$ the \DYt $\b$--module
$(V,\pi,\pi^*)$, where $\pi$ is the restriction of the action of $\gb$ to $\b$,
and the coaction $\pi^*$ is given by
\[\pi^*(v)=\sum_{i}b_i\ten b^i\,v\in\b\ten V\]
where $\{b_i\},\{b^i\}$ are dual bases of $\b$ and $\b^*$. The inverse
functor is obtained by letting $\phi\in\b^*\subset\gb$ act on $V\in\DrY
{\b}$ by $\phi\otimes\id_V\circ\pi^*$.

If $\b$ is $\IN$--graded with \fd homogeneous components, the formulae
defining $E$ similarly give rise to an isomorphism $E\res$ between the
category $\Eq{\grb}$ of equicontinuous modules over the restricted double
of $\b$ and $\DrY{\b}$. Moreover, the categories $\Eq{\gb}$ and $\Eq{\grb}$
are isomorphic, since any locally finite action of $\b^\star$ extends
uniquely to one of $\b^*$, and the following diagram is commutative
\[\xymatrix@C=1.5cm{
\Eq{\gb}\ar[rr]\ar[rd]_{E}&&\Eq{\grb}\ar[ld]^{E\res}\\
&\DrY{\b}&
}\]

% SPLIT MATERIAL

\subsection{Split pairs of Lie bialgebras \cite{ATL1-1}}\label{ss:split lbas}
%--------------------------------------------------------------------

A {\it split pair} of Lie bialgebras $(\b,\a)$ is the datum of two Lie
bialgebras $\a,\b$, together with Lie bialgebra morphisms $i:\a\to
\b$ and $p:\b\to\a$ such that $p\circ i=\id_{\a}$. 

As mentioned in \ref{ss:drinf-double}, the assignment $\b\mapsto\gb$ 
gives rise to a one--to--one correspondence between Lie bialgebras
and Manin triples. Similarly, there is a one--to--one correspondence
between split pairs of Lie bialgebras and \emph{split morphisms} of
Manin triples. 
A morphism of Manin triples $i:(\ga,\a_-,\a_+)\to(\gb,\b_-,\b_+)$ is a
morphism of Lie algebras $i:\ga\to\gb$ which is continuous, preserves
inner products, and is such that $i(\gdpm)\subset\gupm$.\footnote{Note
that such an $i$ is necessarily an embedding.} Set
\[i_{\pm}=\left.i\right|_{\a_{\pm}}:\a_\pm\to\b_\pm
\aand
p_{\pm}=i_{\mp}^t:\b_{\pm}\to\a_{\pm}\]
$i$ is \emph{split} if the projections $p_{\pm}$ are morphisms of Lie
algebras. The following holds \cite[Prop. 3.3]{ATL1-1}
\begin{itemize}
	\item If $i:(\ga,\a_-,\a_+)\longrightarrow(\gb,\b_-,\b_+)$ is a split inclusion
	of Manin triples, then $(\a_-,\b_-,i_-,p_-)$ is a split pair of Lie bialgebras.
	\item Conversely, if $(\a,\b,i,p)$ is a split pair of Lie bialgebras, then 
	$i\oplus p^t:(\ga, \a, \a^*)\longrightarrow(\gb, \b, \b^*)$ is a split
	inclusion of Manin triples.
\end{itemize}

This correspondence may be reformulated as follows. Let $\ksLBA$ be the
category of split Lie bialgebras. The objects of $\ksLBA$ are the same as
those of $\kLBA$, and the morphisms are given by
\begin{equation}\label{eq:sLBA hom}
\Hom_{\ksLBA}(\a,\b)=\{(i,p)\in\Hom_{\kLBA}(\a,\b)\times\Hom_{\kLBA}(\b,\a)\;|\; p\circ i=\id_{\a}\}
\end{equation}
Let $\ksMT$ be the category of Manin triples and split morphisms. Then, the
assignment $\b\to\gb$, $(i,p)\to i\oplus p^t$ is an isomorphism of categories
$\ksLBA\to\ksMT$.

\subsection{Split pairs and restriction functors  \cite{ATL1-1}}
%-----------------------------------------------------------------------------

For any split pair of Lie bialgebras $(\b,\a)$, there is a monoidal
restriction functor $\Res_{\a,\b}:\DrY{\b}\to\DrY{\a}$ defined by
\[\Res_{\a,\b}(V,\pi_V,\pi^*_V)=
(V,\pi_V\circ i\ten\id_V,p\ten\id_V\circ\pi^*_V)\]
Moreover, if $\a\hookrightarrow\b\hookrightarrow\c$ is a chain of split
embeddings, then $\Res_{\a,\b}\circ\Res_{\b,\c}=\Res_{\a,\c}$.
Under the identification of $\DrY{\b},\DrY{\a}$ with the categories
of equicontinuous modules over the doubles $\gb$ and $\ga$
respectively, $\Res_{\a,\b}$ coincides with the pullback functor
corresponding to the morphism $i\oplus p^t:\ga\to\gb$.

% DIAGRAMMATIC MATERIAL

\subsection{Diagrammatic Lie bialgebras}\label{ss:diag-LBA}
%-----------------------------------------------------

A {\it diagrammatic Lie (bi)algebra} $\b$ is the datum of
\begin{itemize}
\item a diagram $D$
\item for any $B\subseteq D$, a Lie (bi)algebra $\b_{B}$
\item for any $B'\subseteq B$, a Lie (bi)algebra morphism $i_{BB'}
:\b_{B'}\to\b_{B}$
\end{itemize}
such that 
\begin{itemize}
\item for any $B\subseteq D$, $i_{BB}=\id_{\b_B}$
\item for any $B''\subseteq B'\subseteq B$, $i_{BB'}\circ i_{B'B''}=i_{BB''}$
\item for any $B=B'\sqcup B''$ %with $B'\perp B''$,
\[i_{BB'}+i_{BB''}:\b_{B'}\oplus\b_{B''}\to\b_{B}\]
is an isomorphism of Lie (bi)algebras. 
\end{itemize}
The above properties imply in particular that $\b_{\emptyset}=0$,
and that $U\b$ is a diagrammatic algebra, with $(U\b)_B=U\b_B$
(cf.\;\ref{ss:D-alg}).

A morphism $\varphi:\b\to\c$ of diagrammatic Lie (bi)algebras with the
same underlying diagram $D$ is a collection of Lie (bi)algebra morphisms
$\varphi_B:\b_B\to\c_B$ labelled by the subdiagrams $B\subseteq
D$ such that, for any $B'\subseteq B$, $\varphi_B\circ i^\b_{BB'}=
i^\c_{BB'}\circ\varphi_{B'}$.

%In most examples of diagrammatic Lie (bi)algebras 

\subsection{Split diagrammatic Lie bialgebras and Manin triples}\label{ss:diag-sLBA}
%-----------------------------------------------------------------------------------

A diagrammatic Lie (bi)algebra $\b$ is {\it split} if there are Lie (bi)algebra
morphisms $p_{B'B}:\b_{B}\to\b_{B'}$ for any $B'\subseteq B$, such
that $p_{B'B}\circ i_{BB'}=\id_{\b_{B'}}$, and

\begin{itemize}
\item for any $B\subseteq D$, $p_{BB}=\id_{\b_B}$
\item for any $B''\subseteq B'\subseteq B$, $p_{B''B'}\circ p_{B'B}=p_{B''B}$
\item for any $B=B'\sqcup B''$ %with $B'\perp B''$,
\[p_{B'B}\oplus p_{B''B}:\b_{B}\to\b_{B'}\oplus\b_{B''}\]
is an isomorphism of Lie (bi)algebras, and is the inverse of $i_{BB'}+i_{BB''}
$.\footnote{The requirements
on $p_{B'B}$ are formulated so as to mirror those in \ref{ss:diag-LBA}.
Note, however, that 1) $p_{BB}=\id_{\b_B}$ follows from 
$p_{BB}\circ i_{BB}=\id_{b_B}$ and $i_{BB}=\id_{\b_B}$ and 2) the
fact that $p_{B'B}\oplus p_{B''B}$ is the inverse of $i_{BB'}+i_{BB''}$
implies that it is a Lie (bi)algebra morphism. Note also that since $p_{CB}\circ
i_{BC}=\id_{b_C}$ for $C=B',B''$, the requirement that $p_{B'B}\oplus
p_{B''B}=(i_{BB'}+i_{BB''})^{-1}$ is equivalent to $p_{B'B}\circ i_{BB''}
=0$ for any $B'\perp B''$.}
\end{itemize}

A morphism $\varphi:\b\to\c$ of split diagrammatic Lie (bi)algebras with
the same underlying diagram is one of the underlying diagrammatic
Lie (bi)algebras such that, for any $B'\subseteq B$, $p^\c_{B'B}\circ
\varphi_B=\varphi_{B'}\circ p^\b_{B'B}$.

One can define similarly a diagrammatic Manin triple as a diagrammatic
Lie algebra $\g=\{\g_B\}_{B\subseteq D}$, where each $\g_B$ is a Manin
triple, and the maps $i_{BB'}:\g_{B'}\to\g_B$ are split morphisms of Manin
triples (see \ref{ss:split lbas}). The equivalence of categories $\ksLBA\cong\ksMT$
implies that a split diagrammatic Lie bialgebra $\b=\{b_B\}_{B\subseteq D}$
gives rise to a diagrammatic Manin triple $\gb=\{\g_{\b_B}\}_{B\subseteq D}$,
which will be referred to as the double of $\b$, and that any such triple arises
this way.

Similarly, if $\b$ is an $\IN$--graded split diagrammatic Lie bialgebra 
with \fd homogeneous components (\ie for any $B\subseteq D$, $\b_
B$ is $\IN$--graded, with \fd homogeneous components and, for any
$B'\subseteq B$, the morphisms $i_{B'B}$ and $p_{B'B}$ are homogeneous
of degree 0), one can similarly define a diagrammatic Lie bialgebra
$\gb\res$, with $(\gb\res)_B\coloneqq\g_{\b_B}\res$, endowed with
a canonical morphism of diagrammatic Lie bialgebras $\b\to\gb\res$.

%===============
\subsection{Example}\label{ss:ex-simple-LA-diag}
%===============

Let $\g$ be a complex semisimple Lie algebra, with opposite
Borel subalgebras $\b_\pm\subset\g$, Dynkin diagram $D$,
Serre generators $\{e_i,f_i,h_i\}_{i\in D}$, and standard Lie
bialgebra structure determined by $\b_\pm$ and an invariant
inner product on $\g$ (see \ref{ss:bil on g}). Then $\g$ is a
diagrammatic Lie bialgebra where, for any $B\subseteq D$,
$\g_B\subseteq\g$ is the subalgebra generated by $\{e_i,f_i,
h_i\}_{i\in B}$.

The diagrammatic structure on $\g$ determines a split
diagrammatic one on $\b_\pm$ as follows. For any $B
\subseteq D$, let $\b_{\pm,B}=\b_\pm\cap\g_B$ be the
subalgebras generated by $\{h_i,e_i\}_{i\in B}$ and $
\{h_i,f_i\}_{i\in B}$ respectively. If $B'\subseteq B$, let $i_{\pm,BB'}:\b_{\pm,B'}\to\b_
{\pm,B}$ be the standard embedding, and regard $p_{\pm,B'B}
=i_{\mp,BB'}^t$ as a map $\b_{\pm,B}\to\b_ {\pm,B'}$
via the identifications $\b_{\mp,C}^*\cong\b_{\pm,C}$
given by the inner product. Then,
$\ker(p_{\pm, B'B})$ is a Lie subalgebra in $\b_{\pm, B}$, 
and therefore $\{p_{\pm,B'B}\}$ give the required splitting 
of the Lie bialgebra $\b_\pm$.

\subsection{\DYt modules over diagrammatic Lie bialgebras}\label{ss:diag-DYLBA}
%-----------------------------------------------------------------------------

The following is straightforward.
\begin{proposition}\label{pr:precox on bialg}
Let $\b$ be a split diagrammatic Lie bialgebra. Then, $\b$ gives
rise to an $(\oalpho{}{}{},\DCPA{}{})$--strict symmetric pre--Coxeter category
$\Cox{\b}{}$, which is defined as follows.
\begin{itemize}
\item For any $B\subseteq D$, $\Cox{\b,B}{}$ is the symmetric monoidal
category $\DrY{\b_B}$.
\item For any $B'\subseteq B$, the functor $F_{B'B}:\Cox{\b,B}{}\to \Cox{\b,B'}{}$ 
is the restriction functor $\Res_{\b_{B'},\b_{B}}:\DrY{\b_B}\to \DrY{\b_{B'}}$. 
\end{itemize}
\end{proposition}

Note that the orthogonality condition $\b_{B'\sqcup B''}\simeq\b_{B'}\oplus
\b_{B''}$ is not needed to define the pre--Coxeter category $\Cox{\b}{}$.
However, it is convenient to construct its deformations as we explain in
Sections \ref{s:univ pC}--\ref{s:qcstructure}.

%%%%%%%%%%%%%%%%%%%%
\subsection{Partial monoidal categories}
%%%%%%%%%%%%%%%%%%%%

The notion of diagrammatic Lie bialgebra may be reformulated
in terms of monoidal functors between \emph{partial} monoidal categories.
A partial monoidal category generalises a monoidal category, in that the
tensor product is only assumed to be defined on a full subcategory
$\C^{(2)}\subseteq\C\times\C$. A  monoidal functor %ORTHOGONALITY
\[(F,J):(\C,\C^{(2)},\ten_{\C},\Phi_{\C})\to(\D,\D^{(2)},\ten_{\D},\Phi_{\D})\] 
between two such categories is the datum of 
\begin{itemize}
\item a functor $F:\Co\to\D$ which preserves the unit, and is such that
$F\times F$ maps $\C^{(2)}$ to $\D^{(2)}$
\item an isomorphism over $\C^{(2)}$ 
\[J:\ten_{\D}\circ F^2\to F\circ\ten_{\C}\]
which is compatible with the unit and the associativity constraint.
\end{itemize}

%%%%%%%%%%%%%%%%%%%%%%%%%%%%%%%%%
\subsection{Functorial description of diagrammatic Lie bialgebras}\label{ss:funct-diag-LBA}
%%%%%%%%%%%%%%%%%%%%%%%%%%%%%%%%%

Let $\PD$ be the category whose objects are the subdiagrams of $D$,
and the morphisms $B'\to B$ are given by inclusions $B'\subseteq
B$. The union $\sqcup$ of orthogonal diagrams is a (symmetric,
strict) partial tensor product on $\PD$, with $\emptyset$ as unit
object.\footnote{Note that $\PD$ is the opposite category to the
category $\ID$ introduced in \ref{ss:diag-cat-model}.}
Let $(\kLBA, \oplus)$ be the category of Lie bialgebras, with 
monoidal structure given by the direct sum, and $0$ as unit
object.

\begin{proposition}\label{pr:partial ten}
The category of diagrammatic Lie bialgebras is isomorphic to that
of monoidal functors $\PD\to\kLBA$. Specifically,
\begin{enumerate}
\item A monoidal functor
\[(F,J):(\PD, \sqcup)\to(\kLBA,\oplus)\]
gives rise to a diagrammatic Lie bialgebra $\b$ defined as follows
\begin{itemize}
\item for any $B\subseteq D$, $\b_B=F(B)$
\item for any $B'\subseteq B$, $i_{BB'}=F(B'\to B)$
\end{itemize} 
Conversely, any diagrammatic Lie bialgebra arises this way
for a unique monoidal functor $(F,J)$.
\item A natural transformation 
of monoidal functors $(\PD, \sqcup)\to(\kLBA,\oplus)$ gives rise to a morphism of
the corresponding diagrammatic Lie bialgebras, and any such natural
transformation arises this way.
\end{enumerate}
\end{proposition}
\begin{pf}
(1) It is clear that $i_{BB}=\id_{\b_B}$, and that $i_{BB'}\circ i_{B'B''}=
i_{BB''}$ for any $B''\subseteq B'\subseteq B$. The key point is to observe
that the existence of the natural isomorphism $J_{B',B''}:F(B')\oplus F(B'')
\to F(B'\sqcup B'')$ for $B'\perp B''$ is equivalent to the requirement
that $i_{BB'}+i_{BB''}:\b_{B'}\oplus\b_{B''}\to\b_{B'\sqcup B''}$ be an
isomorphism of Lie bialgebras.

To this end, note that the naturality of $J$ implies the commutativity of
the following diagram
\[
\xymatrix@C=2.5cm{
F(B')\oplus F(\emptyset) \ar[d]_{J_{B',\emptyset}} \ar[r]^{F(\id_{B'})\oplus F(\emptyset\to B'')} 
& F(B')\oplus F(B'') \ar|{J_{B',B''}}[d] & 
F(\emptyset)\oplus F(B'') \ar[d]^{J_{\emptyset,B''}} \ar[l]_{F(B'\leftarrow\emptyset)\oplus F(\id_{B''})}\\
F(B') \ar[r]_{F(B'\to B'\sqcup B'')} & F(B'\sqcup B'') & \ar[l]^{F(B'\sqcup B''\leftarrow B'')} F(B'')
}\]
Since $F(\emptyset)=0$, it follows that $F(\emptyset\to B'')=0=F(B'\leftarrow
\emptyset)$. Moreover, the compatibility of $J$ with the unit, that is $J_{C,
\emptyset}=\id_{F(C)}=J_{\emptyset,C}$, implies that the above diagram
reduces to
\[\xymatrix@C=2.5cm{
& F(B')\oplus F(B'') \ar|{J_{B',B''}}[d] & \\
F(B') \ar[ur]^{\id\oplus 0}\ar[r]_{F(B'\to B'\sqcup B'')} & F(B'\sqcup B'') & \ar[l]^{F(B'\sqcup B''\leftarrow B'')}\ar[ul]_{0\oplus\id}F(B'')
}\]
so that $J_{B',B''}=i_{BB'}+i_{BB''}$.

(2) If $(F,J),(G,K)$ are monoidal functors, a natural transformation 
$F\Rightarrow G$ of the underlying functors is clearly the same as a morphism
$\varphi:\b\to\c$ of the corresponding diagrammatic Lie bialgebras. The only
point is to observe that $\varphi$ is automatically compatible with the tensor
structures, which follows from the commutativity of the following diagram for
any $B=B'\sqcup B''$
\[\xymatrix{
\b_{B'}\oplus \b_{B''} \ar[rr]^{i_{BB'}^\b\oplus i_{BB''}^\b} \ar[d]_{\varphi_{B'}\oplus\varphi_{B''}}&& \b_B\oplus\b_B\ar[r]^{+}\ar|{\varphi_B\oplus\varphi_B}[d] &\b_B\ar[d]^{\varphi_B} \\
\c_{B'}\oplus \c_{B''} \ar[rr]_{i_{BB'}^\c\oplus i_{BB''}^\c} && \c_B\oplus\c_B\ar[r]_{+} &\c_B\\
}\]
\end{pf}

\noindent
Split diagrammatic Lie bialgebras can be described in similar terms. Let 
$\ksLBA$ be the category of split Lie bialgebras \eqref{eq:sLBA hom}.
Then, the category of monoidal functors $(F,J):(\PD,\sqcup)\to(\ksLBA,\oplus)$
is canonically isomorphic to that of split diagrammatic Lie bialgebras.
Note also that any such functor is automatically symmetric.\\

\noindent\remark\label{rk:diagrammatic}
In view of Proposition \ref{pr:partial ten}, it is natural to define a diagrammatic
object in a monoidal category $(\C,\ten)$ as a monoidal functor $(\PD,\sqcup)
\to(\C,\ten)$, and a morphism of such objects as a natural transformation of 
the corresponding functors.

%%%%%%%%%%%%%
%HOPF DIAGRAMMATIC
%%%%%%%%%%%%%

\section{Diagrammatic Hopf algebras}\label{s:quant-lba}
%===========================

In this section, we introduce the notion of diagrammatic Hopf algebra
and quantised universal enveloping algebra (QUE). We then point out
that the quantisation $\Q(\b)$ of a diagrammatic Lie bialgebra $\b$ is
a diagrammatic QUE, and that admissible %the category of
\DYt modules over $\Q(\b)$ and its canonical subalgebras give rise
to a braided  pre--Coxeter category.

\subsection{Drinfeld--Yetter modules over a Hopf algebra \cite{ek-2,Y}}\label{ss:DY Hopf}
%--------------------------------------------------------------------------

A \DYt module over a Hopf algebra $\hC$ is a triple $(\V,\pi_{\V},\pi_{\V}^*)
$, where $(\V,\pi_\V)$ is a left $\hC$--module, $(\V,\pi^*_\V)$ a right $\hC$--comodule,
and the maps $\pi_\V:\hC\otimes\V\to\V$ and $\pi_\V^*:\V\to \hC\otimes\V$
satisfy the following compatibility condition in $\End(\hC\ten \V)$
\[\pi_{\V}^*\circ\pi_{\V}=
m^{(3)}\otimes\pi_{\V} \circ (1\,3)(2\,4) \circ
S^{-1}\otimes\id^{\otimes 4}\circ \Delta^{(3)}\otimes\pi_{\V}^*\]
where $m^{(3)}:\hC^{\otimes 3}\to \hC$ and $\Delta^{(3)}:\hC\to \hC^{\otimes 3}$
are the iterated multiplication and comultiplication respectively, and $S:\hC
\to \hC$ is the antipode.

The category $\DrY{\hC}$ of such modules is a braided monoidal
category. For any $\V,\W\in\DrY{\hC}$, the action and coaction on
the tensor product  $\V\ten \W$ are defined by
\[
\pi_{\V\ten \W}=\pi_{\V}\ten\pi_{\W}\circ (2\,3)\circ\Delta\ten\id_{\V\otimes\W}
\quad\mbox{and}\quad
\pi^*_{\V\ten \W}=m^{21}\ten\id_{\V\otimes\W}\circ(2\,3)\circ\pi^*_{\V}\ten\pi^*_{\W}
\]
The associativity constraints are trivial, and the braiding is defined
by $\beta_{\V\W}=(1\,2)\circ R_{\V\W}$, where the $R$--matrix
$R_{\V\W}\in\End(\V\ten \W)$ is defined by
\[R_{\V\W}=\pi_{\V}\ten\id_\W\circ(1\,2)\circ\id_\V\ten\pi_{\W}^*\]
The linear map $R_{\V\W}$ is invertible, with inverse
\[R_{\V\W}^{-1}=
\pi_{\V}\ten\id_\W\circ S\ten \id_{\V\otimes\W}\circ(1\,2)\circ\id_\V\ten\pi_{\W}^*\]
The braiding $\beta_{\V\W}$ is therefore invertible, with inverse $R^{-1}
_{\V\W}\circ(1\,2)$.

\subsection{The finite quantum double \cite{drin-2}}\label{ss:DY-qD}
%----------------------------------------------------------

Let $\hC$ be a finite--dimensional Hopf algebra, and $\hC^{\circ}$ the dual
Hopf algebra $\hC^*$ with opposite coproduct. The quantum double of $\hC$
is the unique quasitriangular Hopf algebra $(D\hC, R)$ such that 1) $D\hC=
\hC\otimes \hC^\circ$ as vector spaces 2) $\hC$ and $\hC^\circ$ are Hopf subalgebras
of $D\hC$ and 3) $R$ is the canonical element in $\hC\ten \hC^{\circ}\subset
D\hC\ten D\hC$. The multiplication in $D\hC$ is given in Sweedler's notation by
\begin{equation}\label{eq:q-double-mult}
b\ten f \cdot b'\ten f'=
\langle S^{-1}(b'_1), f_1\rangle\langle b'_3, f_3\rangle\,b\cdot b'_2
\otimes f_2\cdot f'
\end{equation}
where $b,b'\in \hC$, $f,f'\in \hC^{\circ}$, and $\iip{\cdot}{\cdot}:\hC\otimes \hC^\circ
\to\sfk$ is the duality pairing \cite[Sec. 13]{drin-2}. %It is well--known that
The quantum double can also be realised as the double cross product Hopf
algebra $\dcp{\hC\,}{\,\,\hC^*}$ associated to a matched pair of Hopf algebras.
given by the coadjoint actions of $\hC$ on $\hC^*$ and of $\hC^*$ on $\hC$
\cite{majid-2} (see also \cite[Appendix A]{ATL1-1}).

The category $\Rep D\hC$ is canonically isomorphic, as a braided
monoidal category, to $\DrY{\hC}$. Namely, there are two braided monoidal
functors
\begin{equation}\label{eq:iso-DY-DB}
\xymatrix{
	\DrY{\hC} \ar@<2pt>[r]^(.4){\Xi} \ar@<-2pt>@{<-}[r]_(.4){\Theta}
	&\Rep D\hC
}
\end{equation}
which are defined as follows 
\begin{itemize}
	\item For any $D\hC$--module $(\V,\xi_{\V})$, $\Theta(\V,\xi_{\V})=(\V,\pi_{\V},\pi_{\V}^*)$
	is the Drinfeld--Yetter $\hC$--module whose action $\pi_{\V}$ is given by restricting
	$\xi_{\V}$ to $\hC$, and coaction iby the formula $\pi_{\V}^*(v)=R\,1\ten v$.
	\item  For any Drinfeld--Yetter $\hC$--module $(\V,\pi_\V,\pi^*_\V)$, 
	$\Xi(\V,\pi_\V,\pi^*_\V)=(\V,\xi_\V)$ is the $D\hC$--module such that
	$\hC$ acts by $\pi_\V$, and $\phi\in\hC^\circ$ by $\phi\otimes
	\id_\V\circ\pi^*_\V$.
\end{itemize}
One checks easily that the two functors are well--defined, and
are each other's inverses \cite[Prop. A.4]{ATL1-1}.

\subsection{Quantum double for QUEs}\label{ss:dualityQUE}
%--------------------------------------------------

The construction of the quantum double can be adapted for 
quantised universal enveloping algebras (QUE). Recall that
a QUE is a Hopf algebra $\hC$ over $\sfK=\sfk\fml$ which 
%quantises a Lie bialgebra $\b$, \ie it
reduces modulo $\hbar$ to
an enveloping algebra $U\b$ for some Lie bialgebra $\b$, and 
is such that, for any $x\in\b$,
\[
\delta(x)=\frac{\Delta(\wt{x})-\Delta^{21}(\wt{x})}{\hbar}\,\mod\hbar
\]
where $\wt{x}\in \hC$ is any lift of $x$. A QUE is of finite type if the
underlying Lie bialgebra $\b$ is finite--dimensional. In this case,
the dual $\hC^*=\Hom_{\sfK}(\hC,\sfK)$ is a {quantised formal
series Hopf algebra} (QFSH), \ie a topological Hopf algebra over
$\sfK$ which reduces modulo $\hbar$ to $\wh{S\b}=\prod_nS^n\b$. 
Conversely, the dual of a QFSH of finite type is a QUE (cf. \cite
{drin-2, gav} or \cite[Sec. 2.19]{ATL1-1}).

If $\hC$ is a QUE, set
\[\hC'=\{b\in \hC\;|\; (\id-\iota\circ\varepsilon)^{\ten n}\circ\Delta^{(n)}(b)\in\hbar^n\hC^{\ten n}
\;\text{for any $n\geqslant 0$}\}\]
where $\Delta^{(n)}:\hC\to \hC^{\otimes n}$ is the iterated coproduct.
Then, $\hC'$ is a Hopf subalgebra of $\hC$, and a QFSH \cite
{drin-2, gav}. In particular,
if $\hC$ is of finite type, $\hC^\vee\coloneqq(\hC')^*$ is a QUE. 
As in \ref{ss:DY-qD}, $(\hC, \hC^\vee)$ is a matched pair of Hopf algebras
\cite[A.5]{ATL1-1}.
The double cross product $D\hC=\dcp{\hC}{\; \hC^\vee}$ is a 
quasitriangular QUE, whose $R$--matrix is the canonical element 
$R\in \hC'\ten \hC^\vee$ and underlying Lie bialgebra
is the Drinfeld double $\gb=\b\oplus\b^*$. 

This construction extends to the case of \emph{finitely $\IN$--graded}
QUEs, \ie $\IN$--graded Hopf algebras $\hC=\bigoplus_{n\geqslant 0}
\hC_n$ such that $\hC_0$ is a QUE of finite type, and each $\hC_n$ is a
finitely generated $\hC_0$--module. Note that such a QUE is a quantisation
of an $\IN$--graded Lie bialgebra with \fd components
and cobracket of degree $d=0$ (cf. \ref{ss:drinf-double}).
Moreover, $\hC'=\bigoplus_{n\geqslant 0}(\hC'\cap \hC_n)$ is also graded, 
and its \emph{restricted dual} $\hC^\star\coloneqq\bigoplus_{n\geqslant 0}(\hC'\cap \hC_n)^*$
is a finitely $\IN$--graded QUE quantising the restricted dual Lie bialgebra $\b^\star$.
The double cross product $(D\hC)\res\coloneqq\dcp{\hC\;}{\;\;\hC^\star}$
is called the \emph{restricted} quantum double of $\hC$. $(D\hC)\res$
is a quasitriangular, finitely $\IZ$--graded QUE whose $R$--matrix
is the canonical element  in the graded completion of $\hC'\ten \hC^\star$,
and underlying Lie bialgebra is the restricted Drinfeld double $\gb
\res=\b\oplus\b^\star$. 

\subsection{Admissible \DYt modules over a QUE}\label{ss:adm QUE}
%----------------------------------------------------------------

The isomorphism \eqref{eq:iso-DY-DB} between the categories 
of modules over the quantum double and Drinfeld--Yetter modules 
does not hold as is for a QUE and needs to be corrected.

An {\it admissible} \DYt module over a QUE $\hC$ is a \DYt
module $(\V,\pi_\V,\pi^*_\V)$ for which the coaction $\pi^*_\V:\V\to
\hC\otimes\V$ factors through $\hC'\otimes\V$. We denote the category
of such modules by $\aDrY{\hC}$.\footnote{The notion of admissible
\DYt module is due to P. Etingof (private communication), and is
studied in detail in \cite[2.20--2.22]{ATL1-1}.} We show in \cite[Prop. 2.22]{ATL1-1} 
that $\aDrY{\hC}$ reduces modulo $\hbar$ to $\DrY{\b}$.

The following holds.
\begin{itemize}
	\item If $\hC$ is a QUE of finite type, since $R\in \hC'\ten \hC^\vee$, the functors 
	$\Xi,\Theta$ from \eqref{eq:iso-DY-DB} define an isomorphism of braided monoidal 
	categories between $\aDrY{\hC}$ and $\Rep D\hC$. Moreover, this reduces modulo $\hbar$
	to the isomorphism between $\DrY{\b}$ and $\Rep U\gb$.
	\item If $\hC$ is a finitely $\IN$--graded QUE, since $R$ belongs to the grading completion
	of $\hC'\ten \hC^\star$, the functors $\Xi,\Theta$ define an isomorphism of braided monoidal 
	categories between $\aDrY{\hC}$ and the category of $D\hC$--modules whose action of $\hC^\star$
	is locally finite (\ie for any $v\in\V$, $(\hC'\cap \hC_n)^*v=0$ for $n\gg 0$). Moreover, this reduces
	modulo $\hbar$ to the isomorphism $E\res$ between $\DrY{\b}$ and $\E_{\gb\res}$ (cf. 
	\ref{ss:drinf-db-rep}).
\end{itemize}

\subsection{Diagrammatic Hopf algebras}\label{ss:diag-HA}
%-----------------------------------------------------

Let $D$ be a diagram. A {\it diagrammatic} Hopf algebra with underyling
diagram $D$ is a monoidal functor
\[(F,J):(\PD,\sqcup)\to(\HAk,\otimes)\]
where $\HAk$ is the category of Hopf algebras over $\sfk$ (cf.
Remark \ref{rk:diagrammatic}). Concretely, this consists of the
datum of
\begin{itemize} 
%\item a diagram $D$
\item for any $B\subseteq D$, a Hopf algebra $\hC_B$
\item for any $B'\subseteq B$, a morphism of Hopf algebras
$i_{BB'}:\hC_{B'}\to \hC_B$
\end{itemize}
such that
\begin{itemize}
\item for any $B\subseteq D$, $i_{BB}=\id_{\hC_B}$
\item for any $B''\subseteq B'\subseteq B$, $i_{BB'}\circ i_{B'B''}=i_{BB''}$
\item for any $B=B'\sqcup B''$, % with $B'\perp B''$,
\[m_B\circ i_{BB'}\ten i_{BB''}:\hC_{B'}\ten \hC_{B''}\to \hC_{B}\]
is an isomorphism of Hopf algebras, where $m_B$ is the multiplication
of $\hC_B$.
\end{itemize}
The above properties imply in particular that $\hC_{\emptyset}$ is
equal to $\sfk$. Diagrammatic QUEs are defined similarly.

A morphism $\varphi:\hC\to \hC'$ of diagrammatic Hopf algebras (resp.
QUEs) is a collection of Hopf algebra morphisms $\varphi_B:\hC_B
\to \hC'_B$ labelled by the subdiagrams $B\subseteq D$ such that,
for any $B'\subseteq B$, $\varphi_B\circ i^\hC_{BB'}=i^{\hC'}_{BB'}
\circ\varphi_{B'}$.

\subsection{Split diagrammatic Hopf algebras}\label{ss:diag-sHA}
%-----------------------------------------------------

Recall that a split pair of Hopf algebras is the datum of two Hopf algebras
$\hA,\hC$ together with Hopf algebra morphisms $\hA\xrightarrow{i} \hC\xrightarrow
{p}\hA$ such that $p\circ i=\id_\hA$ \cite[Sec. 4.6]{ATL1-1}. We denote by $(\sHAk,
\ten)$ the monoidal category of split Hopf algebras. The objects in $\sHAk$
are the same as those in $\HAk$, and the morphisms are
\[\Hom_{\sHAk}(\hA,\hC)=\{(i,p)\in\Hom_{\HAk}(\hA,\hC)\times\Hom_{\HAk}(\hA,\hC)\;|\; p\circ i=\id_\hA\}\]

A {\it split diagrammatic} Hopf algebra is a monoidal functor $(\PD,\sqcup)\to
(\sHAk,\otimes)$. Concretely, this consists of a diagrammatic Hopf algebra
$\hC=\{\hC_B\}_{B\subseteq D}$, together with Hopf algebra morphisms
$p_{B'B}:\hC_{B}\to \hC_{B'}$ for any $B'\subseteq B$, such that $p_{B'B}\circ i
_{BB'}=\id_{\hC_{B'}}$ and 
\begin{itemize}
\item for any $B$, $p_{BB}=\id_{\b_B}$
\item for any $B''\subseteq B'\subseteq B$, $p_{B''B'}\circ p_{B'B}=p_{B''B}$
\item for any $B=B'\sqcup B''$, $p_{B'B}\ten p_{B''B}\circ\Delta_B:\hC_{B}\to \hC_{B'}\ten \hC_{B''}$
is a morphism of Hopf algebras, and the inverse of $m_B\circ i_{BB'}\ten i_{BB''}$.
\end{itemize}

Split diagrammatic QUEs are defined similarly. A morphism $\varphi:\hC\to \hC'$
of split diagrammatic Hopf algebras (resp. QUEs) is one of the underlying
diagrammatic Hopf algebras (resp. QUEs) such that, for any $B'\subseteq
B$, $p^\hC_{B'B}\circ\varphi_B=\varphi_{B'}\circ p^{\hC'}_{B'B}$.\\

\noindent\remark 
One can formulate in this context a quantum analogue of the Drinfeld double
of a diagrammatic Lie bialgebra defined in \ref{ss:diag-sLBA}. If $\hC$ is a split
diagrammatic  Hopf algebra, where $\hC_B$ are \fd Hopf algebras (resp. 
finitely $\IN$--graded QUE), there is a diagrammatic Hopf algebra $D\hC$ 
with $(D\hC)_B=D\hC_B$ (resp. $(D\hC)\res$ with $(D\hC)\res_B=(D\hC_
B)\res$), endowed  with a canonical embedding of diagrammatic Hopf algebras
$\hC\to D\hC$ (resp. $\hC\to (D\hC)\res$).

\subsection{Drinfeld--Yetter modules over split diagrammatic Hopf algebras}\label{ss:DYdiag}
%--------------------------------------------------------------------------------------------------

If $\hA\leftrightarrows \hC$ is a split pair of Hopf algebras, there is a
monoidal restriction functor $\Res_{\hA,\hC}:\DrY{\hC}\to\DrY{\hA}$ given by 
\[\Res_{\hA,\hC}(\V,\pi_\V,\pi^*_\V)=
(\V,\pi_\V\circ i\otimes\id_\V,p\otimes\id_\V\circ\pi_\V^*)\]
If $\hA,\hC$ are QUEs, $\Res_{\hA,\hC}$ restricts to a functor $\aDrY{\hC}
\to\aDrY{\hA}$.

\begin{proposition}\label{pr:diagr HA}
Let $\hC$ be a split diagrammatic Hopf algebra. Then, $\hC$ gives rise to 
an $(\oalpho{}{}{},\DCPA{}{})$--strict braided pre--Coxeter category $\Cox{\hC}{}$, which is
defined as follows.
\begin{itemize}
\item For any $B\subseteq D$, $\Cox{\hC, B}{}$ is the braided monoidal
category $\DrY{\hC_B}$.
\item For any $B'\subseteq B$, the functor $F_{B'B}:\Cox{\hC, B}{}\to \Cox{\hC, B'}{}$
is the restriction functor $\Res_{\hC_{B'},\hC_{B}}:\DrY{\hC_B}\to \DrY{\hC_{B'}}$. 
\end{itemize}
\end{proposition}

Similarly, a split diagrammatic QUE $\hC$ gives rise to a braided pre--Coxeter
category $\Cox{\hC}{\adm}$ given by $\Cox{\hC, B}{\adm}=\aDrY{\hC_B}$.

\subsection{Quantisation of diagrammatic Lie bialgebras}\label{ss:ek-quantisation}
%-------------------------------------------------------------------------

In \cite{ek-1,ek-2}, Etingof and Kazhdan construct a quantisation functor
$\Q$ from the category of Lie bialgebras over $\sfk$ to the category of
{quantised universal enveloping algebras} over $\sfK=\hext{\sfk}$. One
checks easily that $\Q$ respects direct sums, \ie for any Lie bialgebras
$\a,\b$, there is an isomorphism of Hopf algebras $J_{\a,\b}:\Q(\a)\ten
\Q(\b)\to\Q(\a\oplus\b)$. In fact, this holds for any quantisation functor.

\begin{proposition}
Every quantisation functor $\Q$ is canonically endowed with a monoidal
structure $(\Q,J):(\kLBA,\oplus)\to(\KQUE,\ten)$.
\end{proposition}

\begin{pf}
The result is an easy consequence of Radford's theorem \cite{rad2}. Namely, 
let $i_{\a}:\a\to\a\oplus\b$ and $p_{\a}:\a\oplus\b\to\a$ be the canonical injection 
of and projection to $\a$ and set $\pi_{\a}=i_{\a}\circ p_{\a}$. 
Then, $\Q(\a\oplus\b)$ projects onto $\Q(\a)$ through $\Q(i_{\a})$ and $\Q(p_{\a})$.
By Radford's theorem, $\Q(\a\oplus\b)$ is canonically isomorphic, as a Hopf 
algebra, to the \emph{Radford product} $\Q(\a)\star L$, where 
$L=\{x\in\Q(\a\oplus\b)\;|\;\Q(\pi_{\a})\ten\id\circ\Delta(x)=1\ten x\}$.
It is easy to show that, in this case, $L=\Q(\b)$ and 
$\Q(\a)\star\Q(\b)=\Q(\a)\ten\Q(\b)$. The isomorphism 
$J_{\a,\b}:\Q(\a)\ten\Q(\b)\to\Q(\a\oplus\b)$
is given by $J_{\a,\b}=m_{\Q(\a\oplus\b)}\circ\Q(i_{\a})\ten\Q(i_{\b})$, it is
natural and defines a monoidal structure on $\Q$.
\end{pf}

The same holds for $\ksLBA$ and $\KsQUE$, since the quantisation 
of a split pair of Lie bialgebras is a split pair of QUEs.

\begin{corollary}
The quantisation of a (split) diagrammatic Lie bialgebra is a (split) diagrammatic QUE.
\end{corollary}

\begin{pf}
A (split) diagrammatic Lie bialgebra is a monoidal functor $(\PD,\sqcup)\to(\mathsf{(s)}\kLBA,\oplus)$.
By composition with the quantisation functor, we obtain a monoidal functor 
$(\PD,\sqcup)\to(\mathsf{(s)}\KQUE,\oplus)$, \ie a (split) diagrammatic QUE.
\end{pf}
%ORTHOGONALITY

\subsection{Drinfeld--Yetter $\Q(\b)$--modules}\label{ss:quantum-DY}
%------------------------------------------------------------

The following is a direct consequence of Propositions \ref{ss:ek-quantisation}
and \ref{pr:diagr HA}.

\begin{corollary}
Let $\Q:\kLBA\to\KQUE$ be a quantisation functor, and $\b$ a split diagrammatic
Lie bialgebra. Then, there is an $(\oalpho{}{}{},\DCPA{}{})$--strict braided pre--Coxeter category 
$\Cox{\Q(\b)}{\adm}$ defined by the following data
\begin{itemize}
\item For any $B\subseteq D$, $\Cox{\Q(\b), B}{\adm}$ is the braided monoidal category $\aDrY{\Q(\b_B)}$
\item For any $B'\subseteq B$ and $\F\in\Mns{B,B'}$, the functor
\[F_{\F}:\Cox{\Q(\b), B}{\adm}\to\Cox{\Q(\b), B'}{\adm}\]
is the restriction functor $\Res_{\Q(\b_{B'}),\Q(\b_B)}$.
\end{itemize}
\end{corollary}

One checks easily that $\Cox{\Q(\b)}{\adm}$ reduces modulo $\hbar$ to the
braided pre--Coxeter category $\Cox{\b}{}$defined in \ref{ss:diag-DYLBA}.
In \ref{ss:main-thm-2}, we construct an equivalence of pre--Coxeter categories 
between $\Cox{\Q(\b)}{\adm}$ and a (non $\oalpho{}{}{}$--strict) deformation of $\Cox{\b}{}$.

%%%%%%%%%%%%%%%%%
\section{Diagrammatic $\PROP$s}\label{s:diag-prop}
%%%%%%%%%%%%%%%%%

We review in this section the definition of $\PROP$s, and introduce
a $\PROP$ which governs split diagrammatic Lie bialgebras.

%-------------------------------------------------------------
\subsection{$\PROP$s \cite{La,mac,ee,ATL1-1}}\label{ss:prop-intro} %
%-------------------------------------------------------------

A $\PROP$ is a $\sfk$--linear, strict, symmetric monoidal category
$\sfP$ whose objects are the non--negative integers, and such that
$[n]\ten[m]=[n+m]$. In particular, $[0]$ is the unit object and $[1]^
{\ten n}=[n]$. A morphism of $\PROP$s is a symmetric monoidal
functor $\G:\sfP\to\sfQ$ which is the identity on objects,
and is endowed with the trivial tensor structure
\[\id:\G([m]_{\sfP})\otimes\G([n]_{\sfP})=
[m]_{\sfQ}\otimes[n]_{\sfQ}=[m+n]_{\sfQ}=\G([m+n]_{\sfP})\]

Fix henceforth a complete bracketing $b_n$ on $n$ letters for
any $n\geqslant 2$, and set $\bfb=\{b_n\}_{n\geqslant 2}$. A {\it module}
over $\sfP$ in a symmetric monoidal category $\N$ is a symmetric
monoidal functor $(\G,J):\sfP\to\N$ such that\footnote{In a monoidal
category $(\C,\otimes)$, $V^{\otimes n}_{b_n}$ denotes the $n$--fold
tensor product of $V\in\C$ bracketed according to $b_n$. For example
$V^{\otimes 3}_{(\bullet\bullet)\bullet}=(V\otimes V)\otimes V$.}
\begin{equation*}
\G([n])=\G([1])^{\otimes n}_{b_n}
\end{equation*}
and the following diagram is commutative
\begin{equation}\label{eq:J Phi}
\xymatrix{
\G([m])\otimes \G([n]) \ar@{=}[d]\ar[rr]^{J_{[m],[n]}} && \G([m+n])\ar@{=}[d]\\
\G([1])^{\otimes m}_{b_m}\otimes \G([1])^{\otimes n}_{b_n}\ar[rr]_{\Phi}&& \G([1])^{\otimes(m+n)}_{b_{m+n}}
}
\end{equation}
where $\Phi$ is the associativity constraint in $\N$.\footnote{Note that
the requirement \eqref{eq:J Phi} determines $J$ uniquely. In fact, given
any functor $\G:\sfP\to\N$ such that $\G([n])=\G([1])^{\otimes n}_{b_n}$,
\eqref{eq:J Phi} defines a family of isomorphisms $J_{m,n}:\G([m])\otimes\G([n])\to
\G([m+n])$, which is easily seen to be compatible with the commutativity
and associativity constraints in $\sfP$ and $\N$. Such a $J$, however,
need not be natural \wrt morphisms in $\sfP$, that is satisfy $\G(f\otimes
g)=J_{m_2,n_2}\cdot G(f)\otimes G(g)\cdot J_{m_1,n_1}^{-1}$ for any
$f\in\sfP([m_1],[m_2])$ and $g\in\sfP([n_1],[n_2])$. For example, if 
$\N$ is strict, then $J=\id$, and $J$ is natural if and only if $\G$ is multiplicative
\wrt tensor products of morphisms.} A {\it morphism} of modules over
$\sfP$ is a natural transformation of functors. The category of $\sfP
$--modules in $\N$ is denoted by $\Fun_{\bfb}^\otimes(\sfP,\N)$.

\subsection{The $\PROP$s $\preLA$, $\mathsf{LCA}$ and $\mathsf{LBA}$}\label{ss:prop-la}
%--------------------------------------------------------------------------------------------------

Let $\preLA$ be the $\PROP$ generated by a morphism $\mu:[2]\to[1]$,
subject to the relations
\[\mu\circ(\id_{[2]}+(1\,2))=0
\aand \mu\circ(\mu\ten\id_{[1]})\circ(\id_{[3]}+(1\,2\,3)+(3\,1\,2))=0\]
as morphisms $[2]\to[1]$ and $[3]\to[1]$ respectively. Then, there is a
canonical isomorphism of categories $\Fun_{\bfb}(\preLA,\vect_{\sfk})
\simeq\preLA(\sfk)$, where $\preLA(\sfk)$ is the category of Lie algebras
over $\sfk$. We denote by $\mathsf{LCA}$ and $\mathsf{LBA}$ the
$\PROP$s corresponding to the notions of Lie coalgebras and Lie
bialgebras. 

\subsection{The Karoubi envelope}\label{ss:karoubi}
%--------------------------------------------

Recall that the Karoubi envelope of a category $\C$ is the category $\Kar
{\C}$ whose objects are pairs $(X,\pi)$, where $X\in\C$ and $\pi:X\to X$
is an idempotent. The morphisms in $\Kar{\C}$ are defined as
\[
{\Kar{\C}}((X,\pi), (Y,\rho))=\{f\in\C(X,Y)\;|\; \rho\circ f=f=f\circ\pi\}
\]
with $\id_{(X,\pi)}=\pi$. In particular, ${\Kar{\C}}((X,\id),(Y,\id))={\C}(X,Y)$, 
so that the functor $\C\to\Kar{\C}$ which maps $X\mapsto(X,\id)$ and
$f\mapsto f$ is fully faithful.

Every idempotent in $\Kar{\C}$ splits canonically. Namely,
if $q\in\Kar{\C}((X,\pi),(X,\pi))$ satisfies $q^2=q$, the maps
\[i=q:(X,q)\to(X,\pi)\aand p=q:(X,\pi)\to(X,q)\]
satisfy $i\circ p=q$ and $p\circ i=\id_{(X,q)}$. 
 
If $\sfP$ is a $\PROP$, we denote by $\cKar{\sfP}$ the closure
under infinite direct sums of the Karoubi completion of $\sfP$. If
$\N$ is a symmetric monoidal category, a {\it module} over $\cKar
{\sfP}$ in $\N$ is a symmetric monoidal functor $\cKar{\sfP}\to\N$
such that the composition $\sfP\to\cKar{\sfP}\to\N$ is a module
over $\sfP$. We denote the category of such modules by $\Fun
_{\bfb}^\otimes(\cKar{\sfP},\N)$. It is clear that, if $\N$ is Karoubi
complete and closed under infinite direct sums, the pull--back functor
\[\Fun_{\bfb}^\otimes(\cKar{\sfP},\N)\to\Fun_{\bfb}^\otimes(\sfP,\N)\]
is an equivalence of categories.

%%%%%%%%%%%%%%%%%%

\subsection{Diagrammatic $\PROP$s}\label{ss:diagr-prop}
%------------------------------------------------

Let $D$ be a non--empty diagram. We denote by $\PPD$ the $\PROP$
generated by an idempotent $\dmp{B}:[1]\to[1]$ for any $B\subseteq D$
subject to the relations
\begin{itemize}
\item $\dmp{D}=\id_{[1]}$
\item for any $B'\subseteq B$, $\dmp{B'}\circ\dmp{B}=\dmp{B'}=\dmp{B}\circ\dmp{B'}$ 
\item  for any $B'\perp B''$, $\dmp{B'\sqcup B''}=\dmp{B'}+\dmp{B''}$. 
\end{itemize}
The above relations imply that $\dmp{\emptyset}=0$, and that
$\dmp{B'}\circ\dmp{B''}=0=\dmp{B''}\circ\dmp{B'}$ for any
$B'\perp B''$ since if $p,q$ are idempotents, $p+q$ is an
idempotent if and only if $pq=0=qp$.\footnote{If $p,q$ are
idempotents, $(p+q)^2=p+q$ is equivalent to $pq=-qp$.
This implies $pq=pq^2=-qpq=q^2p=qp$, and therefore
$pq=0$.}

Let $\sfQ$ be a $\PROP$, and consider the $\PROP$ $\sfQ_D$
generated by the morphisms in $\sfQ$ and $\PPD$ subject to the
relation
\[\dmp{B}^{\ten m}\circ f=f\circ\dmp{B}^{\ten n}\]
for any $f\in\sfQ([n],[m])$ and $B\subseteq D$.

\subsection{The $\PROP$ $\preLBA_D$}
%----------------------------------------------------

By definition, $\mathsf{LBA}_D$ is generated by a Lie bialgebra
object $([1],\mu,\delta)$, and idempotents $\dmp{B}\in\End([1])$,
$B\subseteq D$, which are Lie bialgebra maps.

For any category $\C$, denote by $\sfs\C$ the category with the
same objects as $\C$, and with a morphism $X\to Y$ in $\sfs\C$
given by a pair of morphisms $i:X\to Y$, $p:Y\to X$ in $\C$ such
that $p\circ i=\id_X$.

\begin{proposition}
Let $\N$ be a $\sfk$--linear, symmetric monoidal category, and
$\mathsf{LBA}(\N)$ the category of Lie bialgebras in $\N$. Let
$(\PD,\sqcup)$ be the partial monoidal category of subdiagrams
of $D$ introduced in \ref{ss:funct-diag-LBA}. Then, there is a
canonical isomorphism of categories
\[\Fun_{\bfb}(\LBA_D,\N)\simeq\Fun_{\ten}\left((\PD,\sqcup),(\sLBA(\N),\oplus)\right)\]
In particular,
%$if $\N$ is Karoubi complete and closed under infinite direct sums,
the notions of module over $\LBA_D$ and split diagrammatic Lie
bialgebra in $\N$ coincide.
\end{proposition}
\begin{pf}
Let $\T:\PD\to\ul{\sfs}\LBA_D$ be the functor given by
\begin{itemize}
\item $\T(B)=([1],\dmp{B})$
\item $\T(B'\subseteq B)=
\left(i=\dmp{B'}:([1],\dmp{B'})\to ([1],\dmp{B}),p=\dmp{B'}:([1],\dmp{B})\to ([1],\dmp{B'})\right)$
\end{itemize}
$\T$ is a tensor functor $(\PD,\sqcup)\to(\ul{\sfs}\LBA_D,\oplus)$ with
the (iso)morphism $\T(B')\oplus\T(B'')\to\T(B'\sqcup B'')$ given by
the pair of morphisms
\begin{align*}
i&=\dmp{B'}+\dmp{B''}:([1]\oplus [1],\dmp{B'}\oplus\dmp{B''})\to ([1],\dmp{B'\sqcup B''})\\
p&=\dmp{B'}\oplus\dmp{B''}:([1],\dmp{B'\sqcup B''})\to ([1]\oplus [1],\dmp{B'}\oplus\dmp{B''})
\end{align*}
which are each other's inverses because $\dmp{B'\sqcup B''}=\dmp{B'}+\dmp{B''}$.

The functor $\Fun_{\bfb}(\LBA_D,\N)\to\Fun_{\ten}(\PD,\sLBA(\N))$ is
defined by precomposition with $\T$, and is easily seen to be an isomorphism.
\end{pf}

%%%%%%%%%%%%%%%%%%%%%%%%%%%%
\section{Universal algebras}			\label{s:univ alg}
%%%%%%%%%%%%%%%%%%%%%%%%%%%%

In this section, we define a family of algebras which are universal analogues
of the tensor powers $U\gb^{\otimes n}$ of the enveloping algebra of the
double of a diagrammatic Lie bialgebra.

%----------------------------------------
\subsection{Colored $\PROP$s}
%----------------------------------------

A \emph{colored} $\PROP$ $\mathsf{P}$ is a $\sfk$--linear, strict,
symmetric monoidal category whose objects are finite sequences
over a set $\sfA$, \ie
\[\mathsf{Obj}(\mathsf{P})=\coprod_{n\geqslant0}\sfA^n\]
with tensor product given by the concatenation of sequences, and
tensor unit given by the empty sequence. Modules over a colored
$\PROP$ $\mathsf{P}$ and its closure $\ul{\sfP}$ are defined as in
\ref{ss:prop-intro} and \ref{ss:karoubi}, respectively.

%-------------------------------------------------------------
\subsection{Universal Drinfeld--Yetter modules}\label{ss:dyclba}\label{ss:CLBAcompletions}
%-------------------------------------------------------------

Given a diagram $D$ and $n\geqslant0$, the category $\MDY{\dsg}{n}$
is the colored $\PROP$ generated by $n+1$ objects, $\ACDY{1}$ and
$\{\VCDY{k}\}_{k=1}^n$, and morphisms
\begin{itemize}
\item $\dmp{B}:\ACDY{1}\to\ACDY{1}$, $B\subseteq D$\\
\item $\mu:\ACDY{2}\to\ACDY{1}$, $\delta:\ACDY{1}\to\ACDY{2}$\\[.05ex]
\item $\pi_k:\ACDY{1}\ten \VCDY{k}\to \VCDY{k}$ and $\pi_k^*:\VCDY{k}\to\ACDY{1}\ten \VCDY{k}$
\end{itemize}
such that 
\begin{itemize}
\item $(\ACDY{1},\{\dmp{B}\}_{B\subseteq D},\mu,\delta)$ is an $\dLBA{\dsg}$--module 
in $\MDY{\dsg}{n}$\\
\item every $(\VCDY{k},\pi_k,\pi_k^*)$ is a \DYt module over $\ACDY{1}$
\end{itemize}
In particular, $\MDY{\dsg}{0}=\dLBA{\dsg}$.

\subsection{Modules over $\MDY{\dsg}{n}$}\label{ss:MDY}
%-------------------------------------------------------

If $\N$ is a $k$--linear symmetric monoidal category, %the category of
$\MDY{\dsg}{n}$--modules in $\N$ are isomorphic to the category whose
objects are tuples $(\b;V_1,\ldots,V_n)$ consisting of a diagrammatic
Lie bialgebra $\b$ in $\N$, and $n$ \DYt modules $V_1,\ldots,V_n\in\N$
over $\b_D$. A morphism $(\b;V_1,\ldots,V_n)\mapsto(\c;W_1,\ldots,W_n)$
is a tuple $(\phi;f_1,\ldots,f_n)$, where $\phi:\b\to\c$ is a morphism of diagrammatic
Lie bialgebras, and $f_i:V_i\to W_i$ are such that the following diagrams
are commutative
\[
\xymatrix{
\b_D\otimes V_i\ar[r]^{\pi_{V_i}}\ar[d]_{\phi_D\otimes f_i}	&V_i\ar[d]^{f_i}
&&V_i\ar[r]^{\pi^*_{V_i}}\ar[d]_{f_i}	&\b_D\otimes V_i\ar[d]^{\phi_D\otimes f_i}\\
\c_D\otimes W_i\ar[r]_{\pi_{W_i}}					&W_i
&&W_i\ar[r]_{\pi^*_{W_i}}						&\c_D\otimes W_i\\
}
\]
so that $f_i$ is a morphism of $\b_D$--modules $V_i\to\phi_D^*W_i$
as well as a morphism of $\c_D$--comodules $(\phi_D)_*V_i\to W_i$.

\subsection{Universal algebras}\label{ss:univ-alg}\label{ss:D-univ-algebras-1}
%---------------------------------------

Let $\MUA{\dsg}{n}$ be the algebra defined by
\[\MUA{\dsg}{n}=
\pEnd{\MDY{\dsg}{n}}{\VCDY{1}\ten\cdots\ten\VCDY{n}}\]
Let $\N$ be a symmetric tensor category and $(\b;V_1,\ldots,V_n)$
a $\MDY{\dsg}{n}$--module  in $\N$. The corresponding realisation
functor $\G_{\b;V}:\MDY{\dsg}{n}\to\N$ yields a homomorphism
$\MUA{\dsg}{n}\to\End_\N(V_1\odots{\otimes}V_n)$.
We shall need the following.
\begin{lemma}\label{le:diag trick}
Let $(\b;V_1,\ldots,V_n)$ and $(\c;W_1,\ldots,W_n)$ be two $\MDY{\dsg}
{n}$--modules in $\N$, $\phi:\b\to\c$ a morphism of split diagrammatic
Lie bialgebras, and
\[f:V_1\odots{\otimes}V_n\longrightarrow W_1\odots{\otimes}W_n\]
a morphism which intertwines the action of $\b_D$ and the coaction
of $\c_D$ on each tensor factor.
Then, $f$ intertwines the action of $\MUA{D}{n}$
on $V_1\odots{\otimes}V_n$ and $W_1\odots{\otimes}W_n$.
\end{lemma}
\begin{pf}
Let $\G_{\b;V},\G_{\c;W}:\MDY{\dsg}{n}\to\N$ be the realisation
functors corresponding to $(\b;V_1,\ldots,V_n)$ and $(\c;W_1,\ldots,
W_n)$.

By \ref{ss:MDY}, the result holds if $f$ is of the form $f_1\odots
{\otimes}f_n$, where each $f_k:V_k\to W_k$ intertwines the action
of $\b_D$ and coaction of $\c_D$. Indeed, in that case $(\phi;f_1,
\ldots,f_n)$ gives rise to a morphism $\G_{\b;V}\to\G_{\c;W}$,
whose value on $V_1\odots{\otimes}V_n$ is $f$.

More generally, consider the colored $\PROP$ $\MDY{\dsg}{1,n}$
generated by an $\dLBA{\dsg}$--module $(\ACDY{1},\{\dmp{B}\}_
{B\subseteq D},\mu,\delta)$, together with an object $\VCDY
{}$ endowed with $n$ commuting actions $\pi_k:\ACDY{1}\ten
\VCDY{}\to \VCDY{}$, and $n$ commuting coactions $\pi_k^*:
\VCDY{}\to\ACDY{1}\ten\VCDY{}$ such that $(\VCDY{},\pi_k,
\pi_k^*)$ is a \DYt module over $[1]$ for any $1\leq k\leq n$.
There is a natural tensor functor $\Delta:\MDY{\dsg}{1,n}\to\MDY{\dsg}{n}$
which maps $[1]$ to $[1]$
and $\VCDY{}$ to $\VCDY{1}\odots{\otimes}\VCDY{n}$.

The pair $(\phi;f)$ gives rise to a morphism of functors $\G_{\b;V}
\circ\Delta\to\G_{\c;W}\circ\Delta$, so that $f$ intertwines the
action of $\End_{\MDY{\dsg}{1,n}}(\VCDY{})$ on $V_1\odots
{\otimes}V_n$ and $W_1\odots{\otimes}W_n$.
The result now follows because the functor $\Delta$ is full.
\end{pf}

\subsection{Diagrammatic structure on universal algebras}
\label{ss:D univ}\label{ss:D-univ-algebras-2}
%--------------------------------------------------------------------------

For any $B'\subseteq B$, there is a canonical realisation functor
$\RDY{B'}{n}\to\RDY{B}{n}$ which sends the object $[1]_{B'}$ in
$\RDY{B'}{n}$ to the Lie bialgebra $\theta_{B'}([1]_B)=([1]_B,\dmp
{B'})$ in $\RDY{B}{n}$, and each $(\VCDY{B',k},\pi_{B',k},\pi_{B',
k}^*)$ to
\[\Res_{\theta_{B'}([1]_B),[1]_B}(\VCDY{B,k},\pi_{B,k},\pi_{B,k}^*)=
(\VCDY{B,k},\pi_{B,k}\circ \theta_{B'}\otimes\id,\theta_{B'}\otimes
\id\circ\pi_{B,k}^*)\]
where $\theta_{B'}$ is regarded both as the split injection $([1]_B,
\dmp{B'})\to[1]_B$ and projection $[1]_B\to([1]_B,\dmp{B'})$ (cf.
\ref{ss:karoubi}). The functor induces a homomorphism $\sfi_{BB'}:
\RDYUA{B'}{n}\to\RDYUA{B}{n}$, and it is clear that $\sfi_{BB}=
\id_{\RDYUA{B}{n}}$, and $\sfi_{BB'}\circ\sfi_{B'B''}=\sfi_{BB''}$
for any $B''\subseteq B'\subseteq B$.

\begin{proposition}\label{pr:univ D alg}
The algebras $\{\MUA{B}{n}\}_{B\subseteq D}$ and maps $\{\sfi_{BB'}
\}_{B'\subseteq B}$ give rise to a lax diagrammatic algebra, which we
denote by $\DUA{D}{n}$.
\end{proposition}
\begin{pf}
We need to prove that if $B'\perp B''$, the images of $\sfi_{DB'}$ and
$\sfi_{DB''}$ commute in $\MUA{D}{n}$. This can be proved by a direct
computation \cite[Prop. 10.6]{ATL2}. We give a more conceptual proof
below.

% first reduction
By Lemma \ref{le:diag trick}, it suffices to show that the action of $\MUA
{B''}{n}$ on $\VCDY{1}\ten\cdots\ten\VCDY{n}\in\MDY{\dsg}{n}$ commutes
with the action and coaction of $[1]_{B'}$ on each $\VCDY{k}$.
%
% second reduction
It is easy to check that each of these
%the the action and coaction of $[1]_{B'}$ on any $\VCDY{k}\in\MDY{\dsg}{n}$
commutes with both the action and the coaction of $[1]_{B''}$ on $\VCDY{k}$.
This implies that the maps
\begin{align*}
\pi_{B',k}&:
\VCDY{1}\odots{\otimes}
\left([1]_{B'}\otimes \VCDY{k}\right)\odots{\otimes}\VCDY{n}
\longrightarrow
\VCDY{1}\odots{\otimes}\VCDY{n}\\[.4ex]
\pi_{B',k}^*&:
\VCDY{1}\odots{\otimes}\VCDY{n}
\longrightarrow
\VCDY{1}\odots{\otimes}
\left([1]_{B'}\otimes \VCDY{k}\right)\odots{\otimes}\VCDY{n}
\end{align*}
commute with the action and coaction of $[1]_{B''}$ on each tensor
factor,  where $[1]_{B'}$ is given the structure of trivial \DYt module
over $[1]_{B''}$. 

By Lemma \ref{le:diag trick}, if $x''\in\MUA{B''}{n}$, and $x''_{\VCDY{1},
\ldots,\VCDY{n}}$ (resp. $x''_{\VCDY{1},\ldots,[1]_{B'}\otimes\VCDY{k},
\ldots,\VCDY{n}}$) denote its action on $\VCDY{1}\ten\cdots\ten\VCDY
{n}$ (resp. $\VCDY{1}\odots{\ten}\left([1]_{B'}\otimes\VCDY{k}\right)\odots
{\ten}\VCDY{n}$), then
\begin{align*}
x''_{\VCDY{1},\ldots,\VCDY{n}}\cdot \pi_{B',k}
&=
\pi_{B',k}\cdot x''_{\VCDY{1},\ldots,[1]_{B'}\otimes\VCDY{k},\ldots,\VCDY{n}}\\[.4ex]
\pi_{B',k}^*\cdot x''_{\VCDY{1},\ldots,\VCDY{n}}
&=
x''_{\VCDY{1},\ldots,[1]_{B'}\otimes\VCDY{k},\ldots,\VCDY{n}}
\cdot \pi_{B',k}^*
\end{align*}
The conclusion now follows from the fact that, since $[1]_{B'}$ is regarded
as a trivial \DYt module over $[1]_{B''}$, 
\[x''_{\VCDY{1},\ldots,[1]_{B'}\otimes\VCDY{k},\ldots,\VCDY{n}}=\id_{[1]_{B'}}
\otimes x''_{\VCDY{1},\ldots,\VCDY{n}}\]
under the identification $\VCDY{1}\ten\cdots\ten\left([1]_{B'}\otimes\VCDY{k}\right)\ten
\cdots\ten\VCDY{n}\cong [1]_{B'}\otimes \VCDY{1}\ten\cdots\ten\VCDY{n}$.
\end{pf}

\noindent\remark\label{re:inj D}
We show in \cite{ATL2} that, for any $B'\subseteq B$, the homomorphism
$\sfi_{BB'}:\RDYUA{B'}{n}\to\RDYUA{B}{n}$ is injective. We shall therefore
regard $\RDYUA{B'}{n}$ as a subalgebra of $\RDYUA{B}{n}$ and, for $x\in
\RDYUA{B'}{n}$, write $x\in\RDYUA{B}{n}$ instead of $\sfi_{BB'}(x)\in\RDYUA
{B}{n}$. Moreover, $\{\RDYUA{B}{n}\}_{B\subseteq D}$ is a diagrammatic algebra, 
since multiplication induces an isomorphism $\RDYUA{B_1\sqcup B_2}{n}
\cong\RDYUA{B_1}{n}\ten\RDYUA{B_2}{n}$ \cite[Prop. 10.6 (4)]{ATL2}.

\subsection{Fiber functors and diagrammatic structures}
\label{ss:D fiber}
%-----------------------------------------------------------------------

Let now $\b$ be a split diagrammatic Lie bialgebra. For any $B\subseteq
D$, let
\[\ff_{\b_B}:\DrY{\b_B}\to{\kvect}
\aand
\DYA{\b_B}{n}=\sfEnd{\ff_{\b_B}^{\boxtimes n}}\]
be the forgetful functor and algebra of endomorphisms of $\ff_{\b_B}
^{\boxtimes n}$. By definition, an element of $\DYA{\b_B}{n}$ is a collection 
$x_{V_1,\ldots,V_n}\in\End_\sfk(V_1\otimes\cdots\otimes V_n)$ labelled
by $V_1,\ldots,V_n\in\DrY{\b_B}$ such that if $f_k\in\End_{\DrY{\b_B}}
(V_k,W_k)$, $1\leq k\leq n$, then
\[f_1\otimes\cdots\otimes f_n\circ x_{V_1,\ldots,V_n}=
x_{W_1,\ldots,W_n}\circ f_1\otimes\cdots\otimes f_n\]
The equivalence between $\DrY{\b_B}$ and equicontinuous modules
over $\g_{\b_B}$ (Section \ref {ss:drinf-db-rep}) gives rise to a map $U\g_{\b_
B}^{\ten n}\to\DYA{\b_B}{n}$, which is an isomorphism if $\dim\b_B<
\infty$.

For any $B'\subseteq B$, $(\b_{B'},\b_{B})$ is a split pair of Lie bialgebras.
The corresponding restriction functor $\DrY{\b_{B}}\to\DrY{\b_{B'}}$ induces
a homomorphism $i_{BB'}:\DYA{\b_{B'}}{n}\to\DYA{\b_B}{n}$, which clearly
satisfies $i_{BB}=\id_{\DYA{\b_B}{n}}$ and $i_{BB'}\circ i_{B'B''}=i_{BB''}$ for any $B''\subseteq B' \subseteq B$.

\begin{proposition}
The algebras $\{\DYA{\b_B}{n}\}_{B\subseteq D}$ and maps $\{i_{BB'}\}_
{B'\subseteq B}$
give rise to a lax diagrammatic algebra, which we denote by $\DYA{\b}{n}$.
\end{proposition}
\begin{pf}
We need to prove that if $B'\perp B''$, the images of $i_{DB'}$ and $i_{DB''}$
commute in $\DYA{\b_D}{n}$. It is easy to check that the action and coaction
of $\b_{B'}$ commute with those of $\b_{B''}$ on any $V\in\DrY{\b_D}$. Thus,
$\b_{B'}$ acts and coacts on each tensor factor of $V_1\otimes
\cdots\otimes V_n$, $V_k\in\DrY{\b_D}$, through morphisms in $\DrY{\b_{B''}}$.
By definition of $\DYA{\b_{B''}}{n}$, the action of the latter on $V_1\otimes
\cdots\otimes V_n$ therefore commutes with the action and coaction of $b_{B'}$
on each tensor factor. Thus, $\DYA{\b_{B''}}{n}$ acts by tensor products of
morphisms in $\DrY{\b_{B'}}$ and thefore commutes with $\DYA{\b_{B'}}{n}$.
\end{pf}

\subsection{Universal algebras as endomorphisms of fiber functors}\label{ss:univ fiber} % relating the two D-algebras
%--------------

%Let $\b$ be a split diagrammatic Lie bialgebra.
The following shows that the lax diagrammatic algebra $\DUA{D}{n}$ 
obtained in \ref{ss:D univ} is a universal analogue of the lax diagrammatic algebra 
$\DYA{\b}{n}$ obtained in \ref{ss:D fiber}.

Let $B\subseteq D$. For any $n$--tuple $\{V_k,\pi_k,\pi_
k^*\}_{k=1}^ n$ of \DYt modules over $\b_B$, let
\[\G_{(\b_B;V_1,\dots, V_n)}:\MDY{B}{n}\longrightarrow{\vect_{\sfk}}\]
be the corresponding realisation functor.

\begin{proposition}\label{pr:universal action}\hfill
\begin{enumerate}
\item %For any $B\subseteq D$, t
There is an algebra homomorphism
\[\DYrho{\b_B}{n}:\MUA{B}{n}\to\DYA{\b_B}{n}\]
which assigns to any $T\in\MUA{B}{n}$, and any $V_1,\ldots,V_n\in
\DrY{\b_B}$ the endomorphism $\G_{(\b_B;V_1,\dots, V_n)}(T)\in\End_
\sfk(V_1\otimes\cdots\otimes V_n)$.
\item The collection of homomorphisms $\{\DYrho{\b_B}{n}\}_{B\subseteq D}$
is a morphism of lax diagrammatic algebras $\DYrho{\b}{n}:\DUA{D}{n}\to\DYA{\b}{n}$.
\end{enumerate}
\end{proposition}
\begin{pf}
(1) follows from Lemma \ref{le:diag trick}. (2) is clear.
\end{pf}

%%%%%%%%%%%%%%%%%%

%----------------------------------------------------------------
\subsection{Cosimplicial structure of $\DYA{\b}{\bullet}$}\label{ss:cosim}
%----------------------------------------------------------------

For any $B\subseteq D$, the tensor structure on $\DrY{\b_B}$ endows
the tower $\{\DYA{\b_B}{n}\}_{n \geqslant 0}$, with the structure of a
cosimplicial complex of algebras
\[\xymatrix@C=.5cm{\sfk\ar@<-2pt>[r]\ar@<2pt>[r] & 
\sfEnd{\ff_{\b_B}} \ar@<3pt>[r] \ar@<0pt>[r] \ar@<-3pt>[r] & 
\sfEnd{\ff_{\b_B}^{\boxtimes 2}} \ar@<-6pt>[r]\ar@<-2pt>[r]\ar@<2pt>[r]\ar@<6pt>[r] & \sfEnd{\ff_{\b_B}^{\boxtimes 3}}
\quad\cdots}\]
which is compatible with the cosimplicial structure on $\{U\g_{\b_B}^
{\otimes n}\}_{n\geq 0}$ induced by the coproduct, via the maps $U\g_{\b_B}^{\otimes n}\to \DYA{\b_B}{n}$.

The corresponding face morphisms $d_i^n:\sfEnd{\ff_{\b_B}^{\boxtimes n}}\to
\sfEnd{\ff_{\b_B}^{\boxtimes n+1}}$, $i=0,\dots, n+1$ are given by $(d_0^0
\varphi)_V=(d_1^0 \varphi)_V=\varphi\cdot\id_V$, for $\varphi\in\sfk$ and
$V\in\DrY{\b_B}$, and, for $n\geqslant 1$, $\varphi\in\sfEnd{\ff_{\b_B}^
{\boxtimes n}}$, and $V_i\in\DrY{\b_B}$, $1\leqslant i\leqslant n+1$,
\[
(d_i^n \varphi)_{V_1,\dots, V_{n+1}}=
\left\{
\begin{array}{cc}
\id_{V_1}\ten \varphi_{V_2,\dots, V_{n+1}} & i=0\\[1.1ex]
\varphi_{V_1,\dots, V_i\ten V_{i+1},\dots V_{n+1}} & 1\leqslant i \leqslant n\\[1.1ex]
\varphi_{V_1,\dots, V_n}\ten\id_{V_{n+1}} & i=n+1
\end{array}
\right.
\]
The degeneration homomorphisms $\varepsilon_n^i:\sfEnd{\ff_{\b_B}^{\boxtimes n}}\to
\sfEnd{\ff_{\b_B}^{\boxtimes n-1}}$, $i=1,\dots, n$, are
\[(\varepsilon_n^i \varphi)_{X_1, \dots, X_{n-1}}=\varphi_{X_1,\dots, X_{i-1}, {\bf 1}, X_{i}, \dots, X_{n-1} }\]
where ${\bf 1}$ is the trivial \DYt module. The morphisms $\varepsilon_n^i$, $d_i^n$
satisfy the standard relations 
\[
\begin{array}{lr}
d_{n+1}^jd^i_n=d^i_{n+1}d_n^{j-1} & i<j \\
\varepsilon_n^j\varepsilon_{n+1}^i=\varepsilon_n^i\varepsilon_{n+1}^{j+1} & i\leqslant j
\end{array}
\]
\[
\varepsilon_{n+1}^jd_i^n=\left\{
\begin{array}{cc}
d_{n-1}^i\varepsilon_n^{j-1} & i<j\\
\id & i=j,j+1\\
d_{n-1}^{i-1}\varepsilon_n^j & i>j+1
\end{array}
\right.
\]
and give in particular rise to the Hochschild differential
\[d^n=
\sum_{i=0}^{n+1}(-1)^{i}
d_i^n:\sfEnd{\ff_{\b_B}^{\boxtimes n}}\to\sfEnd{\ff_{\b_B}^{\boxtimes n+1}}\]

\noindent
The cosimplicial structure is compatible with the maps $\{i_{BB'}\}_{B'\subseteq
B\subseteq D}$ and therefore determines a cosimplicial lax diagrammatic algebra
$\DYA{\b}{\bullet}$. 

%--------------------------------------------------------------------------
\subsection{Cosimplicial structure of $\RDYUA{}{\bullet}$}\label{ss:cosimp-DY}
%--------------------------------------------------------------------------

The above construction can be lifted to the $\PROP$s $\MDY{B}{n}$.
For every $B\subseteq D$, $n\geqslant1$ and $i=0,\dots, n+1$, there
are faithful functors
\[
\D_{i}^n:\MDY{B}{n}\to\MDY{B}{n+1}
\]
mapping $[1]$ to $[1]$, and given by
\[\D^n_0(\VDY{k})=\VDY{k+1}
\aand
\D^n_{n+1}(\VDY{k})=\VDY{k}\]
for $1\leqslant k\leqslant n$, and, for $1\leqslant i\leqslant n$,
\[\D^n_i(\VDY{k})=
\left\{\begin{array}{cc}
\VDY{k} 				& 1\leqslant k\leqslant i-1	\\[1.1ex]
\VDY{i}\otimes \VDY{i+1} 	& k=i					\\[1.1ex]
\VDY{k+1}				& i+1\leqslant k\leqslant n	
\end{array}\right.\]
and $\E_n^{(i)}:\MDY{B}{n}\to\MDY{B}{n-1}$
\[
\E_n^{(i)}=\G_{([1], \VDY{1},\dots,\VDY{i-1},\1,\VDY{i+1},\dots, \VDY{n-1})}
\]
where $\1$ is the tensor unit in $\MDY{B}{n}$, regarded as trivial
\DYt module. These induce algebra homomorphisms
\[
\Delta_i^n:\RDYUA{B}{n}\to\RDYUA{B}{n+1}
\]
which are universal analogues of the insertion/coproduct maps on $U\g
_{\b_B}^{\otimes n}$. They endow the tower $\{\RDYUA{B}{n}\}_{n\geqslant
0}$ with the structure of a cosimplicial algebra, with Hochschild differential
$d^n=\sum_{i=0}^{n+1}(-1)^i\Delta_i^n:\RDYUA{B}{n}\to\RDYUA{B}{n+1}$.
This structure is compatible with the maps $\{\sfi_{BB'}\}_{B'\subseteq B
\subseteq D}$ and therefore it extends to a cosimplicial structure on
the lax diagrammatic algebras $\DUA{D}{n}$.

The following is straightforward.

\begin{proposition}
The morphism of lax diagrammatic algebras $\rho^\bullet_{\b}:\DUA{D}{\bullet}
\to\DYA{\b}{\bullet}$ obtained in \ref{ss:univ fiber} is compatible with
the cosimplicial structures.
\end{proposition}

%%%%%%%%%%%%%%%%%%%%%%%%%
%%%%%%%%%%%%%%%%%%%%%%%%%
\section{Universal braided pre--Coxeter structures}\label{s:univ pC}
%%%%%%%%%%%%%%%%%%%%%%%%%
%%%%%%%%%%%%%%%%%%%%%%%%%

We introduce in this section a class of braided pre--Coxeter categories
related to split diagrammatic Lie bialgebras. They are universal, in that
they are arise from the $\PROP$s $\RDY{D}{n}$ defined in Section
\ref{s:univ alg}.

%--------------------------
\subsection{Gradings}\label{ss:grading}
%---------------------------

Let $B\subseteq D$. The $\PROP$ $\MDY{B}{n}$ has a natural $\IN$--bigrading 
given by $\deg(\sigma)=(0,0)=\deg(\dmp{B'})$ for any $\sigma\in\SS_N$ and 
$B'\subseteq B$,
\[\deg(\mu)=(1,0)=\deg(\pi_{\VDY{k}})
\aand
\deg(\delta)=(0,1)=\deg(\pi_{\VDY{k}}^*)\]
for any $1\leqslant k\leqslant n$.

The algebra $\RDYUA{B}{n}$ inherits this bigrading and is concentrated
in bidegrees $(N,N)$, since a degree $(p,q)$ morphism with source
$\VCDY{1}\ten\cdots\ten\VCDY{n}$ is easily seen to map to $[1]^{\otimes
(q-p)}\otimes\VCDY{1}\ten\cdots\ten\VCDY{n}$. For any $a,b\in\IN$, the
corresponding $\IN$--grading determined by mapping $(1,0),(0,1)$ to
$a,b$ respectively yields the same graded completion $\CRDYUA{B}
{n}$ of $\RDYUA{B}{n}$, so long as $a+b>0$. For definiteness, we
set $a=0$ and $b=1$. The morphisms $\{i_{BB'}\}_{B'\subseteq B}$ naturally
extends to the graded completions and induce on the algebras $\CRDYUA{B}
{n}$, $B\subseteq D$, a lax diagrammatic algebra structure $\CDUA{D}{n}$, which 
extends $\DUA{D}{n}$.  

%----------------------------
\subsection{Invariants}\label{ss:inv}
%----------------------------

For any pair of subdiagrams $B'\subseteq B$, denote by $\CRDYUA {B,B'}
{n}\subseteq\CRDYUA{B}{n}$ the subalgebra of elements which commute
with the diagonal action and coaction of $[\b_{B'}]=([1],\dmp{B'})$ on
$\VDY{1}\otimes\cdots\otimes\VDY{n}$. Note that, by Lemma \ref{le:diag trick}, 
$\CRDYUA{B,B'}{n}$ commutes with the diagonal action of $\CRDYUA{B'}{}$
on $\VDY{1}\otimes\cdots\otimes\VDY{n}$, which is given by 
\[ \CRDYUA{B'}{}\ni x\longrightarrow x_{1,2,\ldots,n}=\Delta^{n-1}_1\circ\cdots\circ\Delta^2_1\circ\Delta^1_1(x)\]

%--------------------------------------------------
\subsection{Associators}\label{ss:assoc}
%--------------------------------------------------

Fix $B\subseteq D$. Define the $r$--matrix $r=r_{\VDY{1},\VDY{2}}\in\CRDYUA
{B}{2}$ as the composition
\[r_{\VDY{1},\VDY{2}}=
\pi_{\VDY{1}}\otimes\id_{\VDY{2}}\circ\,(1\,2)\circ \id_{\VDY{1}}\otimes\pi^*_{\VDY{2}}\]
%\in\pEnd{\MDY{B}{2}}{\VDY{1}\ten\VDY{2}}\] 
(resp. $r^{21}_{\VDY{1},\VDY{2}}=\id_{\VDY{1}}\ten\pi_{\VDY{2}}\circ (1\,2)\circ \pi^*
_{\VDY{1}}\ten\id_{\VDY{2}}$), and set $\Omega=r^{12}+r^{21}$.
An invertible, invariant element $\Phi\in\CRDYUA{B,B}{3}$ is called an \emph{associator}
if the following relations are satisfied (in $\CRDYUA{B}{4}$ and $\CRDYUA{B}{3}$ respectively).
\begin{itemize}
\item {\bf Pentagon relation}
\[
\Phi_{1,2,34}\Phi_{12,3,4}=\Phi_{2,3,4}\Phi_{1,23,4}\Phi_{1,2,3}
\]
\item {\bf Hexagon relations}
\begin{align*}
&e^{\Omega_{12,3}/2}=\Phi_{3,1,2}e^{\Omega_{13}/2}\Phi_{1,3,2}^{-1}e^{\Omega_{23}/
2}\Phi_{1,2,3}\\
&e^{\Omega_{1,23}/2}=\Phi_{2,3,1}^{-1}e^{\Omega_{13}/2}\Phi_{2,1,3}e^{\Omega_{12}/
2}\Phi_{1,2,3}^{-1}
\end{align*}
\item {\bf Duality}
\[
\Phi_{3,2,1}=\Phi_{1,2,3}^{-1}
\]
\item {\bf $2$--jet}
\[
\Phi=1+\frac{1}{24}[\Omega_{12},\Omega_{23}]\qquad\mod(\CRDYUA{B}{3})_{\geqslant 3}
\]
\end{itemize}

%-----------------------------------------------------------------------------------
\subsection{Braided pre--Coxeter structures on $\CDUA{D}{\bullet}$}\label{ss:prop-weak-cox}
%-----------------------------------------------------------------------------------

\begin{definition}
	A \emph{braided pre--Coxeter structure} $(\Phi_B,J_{\F},\DCPA
	{\F}{\G},\redasso{\F}{\F'})$ on $\CDUA{D}{\bullet}$ consists of the following data. 
	\begin{enumerate}\itemsep0.25cm
		\item % associators
		{\bf Associators.} 
		For any $B\subseteq D$, an associator $\Phi_B\in\CRDYUA{B,B}{3}$.
		We set $R_{B}=\exp(\Omega_B/2)\in\CRDYUA{B,B}{2}$.
		\item % relative twists
		{\bf Relative twists.} For any $B^\prime\subseteq B$ and \mns $\F\in\Mns
		{B,B^\prime}$, an invertible element $J_{\F}\in\CRDYUA{B,B'}{2}$ such that
		$(J_{\F})_0=1$ and $\varepsilon_2^1(J_{\F})=1=\varepsilon_2^2(J_{\F})$,
		where $\varepsilon_2^1,\varepsilon_2^2:\CRDYUA{B}{2}\to\CRDYUA{B}{}$
		are the degeneration homomorphisms, and satisfying the following properties.
		\vspace{0.25cm}
		\begin{itemize}\itemsep0.25cm
		\item {\bf Compatibility with associators.}
		The relative twist equation holds,
		\begin{equation}\label{eq:rel twist}
		J_{\F, 1,23}\cdot J_{\F, 23}\cdot \Phi_{B'} = \Phi_{B}\cdot J_{\F, 12,3}\cdot J_{\F, 12}
		\end{equation}
		\item {\bf Normalisation.} For any $B\subseteq D$, $J_B=1$.\footnote
		{Here $B$ is identified with the unique element in $\Mns{B,B}$.}
		\end{itemize}
		\item % DCP associators
		{\bf \DCP associators.}
		For any $B'\subseteq B$ and $\F,\G\in\Mns{B,B'}$, an invertible element
		$\DCPA{\G}{\F}\in\CRDYUA{B,B'}{}$
		%henceforth referred to as a \emph{De Concini--Procesi associator},
		such that $(\DCPA{\G}{\F})_0=1$, $\varepsilon(\DCPA{\G}{\F})=1$, and satisfying the
		following properties.
		\vspace{0.25cm}
		\begin{itemize}\itemsep0.25cm
			\item {\bf Compatibility with $J$.} For any $\F,\G\in\Mns{B,B'}$,  
			\[J_{\G}=\gaugetris{\DCPA{\G}{\F}}{J_{\F}}\]
			\item {\bf Horizontal factorisation.} For any $\F,\G,\H\in\Mns{B,B'}$,
			\[\DCPA{\H}{\F}=\DCPA{\H}{\G}\cdot\DCPA{\G}{\F}\]
			In particular, $\DCPA{\F}{\F}=1$ and $\DCPA{\F}{\G}=\DCPA{\G}{\F}^{-1}$.
		\end{itemize}
		\item % vertical joins
		{\bf Vertical joins.}
		For any $B''\subseteq B'\subseteq B$, $\F\in\Mns{B,B'}$, and $\F'\in\Mns{B',B''}$, an
		invertible element $\redasso{\F}{\F'}\in\CRDYUA{B,B''}{}$
		%, henceforth referred to as a \emph{vertical join},
		such that $(\redasso{\F}{\F'})_0=1$, $\varepsilon(\redasso{\F}{\F'})=1$, 
		and satisfying the following properties.
		\vspace{0.25cm}
		\begin{itemize}\itemsep0.25cm
		\item {\bf Compatibility with $J$ (\emph{vertical $J$--factorisation}).}
		\[J_{\F'\cup\F}=\gaugetris{\redasso{\F}{\F'}}{J_{\F}\cdot J_{\F'}}\]
		\item {\bf Compatibility with $\DCPA{}{}$ (\emph{vertical $\DCPA{}{}$--factorisation}).} 
		For any 
		$\F,\G\in\Mns{B,B'}$ and $\F',\G'\in\Mns{B',B''}$,
		\[\DCPA{(\G\cup\G')}{(\F\cup\F')}\cdot\redasso{\F}{\F'}=
		\redasso{\G}{\G'}\cdot\DCPA{\G}{\F}\cdot\DCPA{\G'}{\F'}\]
		\item {\bf Associativity.} For any $B'''\subseteq B''\subseteq B'\subseteq B$, 
		$\F\in\Mns{B,B'}$, $\F'\in\Mns{B',B''}$, and $\F''\in\Mns{B'',B'''}$,
		\[\redasso{\F'\cup\F}{\F''}\cdot\redasso{\F}{\F'}=\redasso{\F}{\F''\cup\F'}\cdot\redasso{\F'}{\F''}\]
		\item {\bf Normalisation.} For any $\F\in\Mns{B,B'}$,
		\[\redasso{\F}{B'}=1=\redasso{B}{\F}\]
		\end{itemize}
	\end{enumerate}
Consistently with the diagrammatic algebra structure on $\CRDYUA{}{\bullet}$, specifically 
the fact that $\CRDYUA{B_1\sqcup B_2}{n}\cong\CRDYUA{B_1}{n}\ten\CRDYUA
{B_2}{n}$ (Remark \ref{re:inj D}), we further require that $J$, $\DCPA{}{}$ and
$\redasso{}{}$ satisfy the following property.
\begin{enumerate}
\item[]
\begin{itemize}
\item {\bf Orthogonal factorisation.}
If $B_1''\subseteq B_1^\prime\subseteq B_1\perp B_2\supseteq B'_2\supseteq B''_2$,
$(\F_1,\F_2)\in\Mns{B_1\sqcup B_2,B'_1\sqcup B'_2}$, %and 
$(\F_1',\F_2')\in\Mns{B'_1\sqcup B'_2,B''_1\sqcup B''_2}$,
\begin{eqnarray*}
\Phi_{B_1\sqcup B_2}&=&\Phi_{B_1}\cdot\Phi_{B_2}\\
J_{(\F_1,\F_2)}&=&J_{\F_1}\cdot J_{\F_2}\\
\DCPA{(\G_1,\G_2)}{(\F_1,\F_2)}&=&\DCPA{\G_1}{\F_1}\cdot\DCPA{\G_2}{\F_2}	\\
\redasso{(\F_1, \F_2)}{(\F'_1,\F'_2)}&=&\redasso{\F_1}{\F'_1}\cdot\redasso{\F_2}{\F'_2}
\end{eqnarray*}
\end{itemize}
\end{enumerate}
\end{definition}

\noindent
Note that $R_{B_1\sqcup B_2}=R_{B_1}\cdot R_{B_2}$. Moreover, since $\CRDYUA{B_1}
{n}$ and $\CRDYUA{B_2}{n}$ commute, the order of the products in the above identities is
irrelevant. The following is a direct consequence of the orthogonal factorisation and normalisation
of $J$, $\DCPA{}{}$, and $\redasso{}{}$.

\begin{lemma}\label{lem:orthogonality}\hfill
\begin{enumerate}
\item For any $B'\subseteq B\perp B''$, and $\F\in\Mns{B,B'}$, $J_{(\F, B'')}=J_{\F}$.
\item For any $B'\subseteq B\perp B''$, and $\F,\G\in\Mns{B,B'}$, $\DCPA{(\F, B'')}{(\G,B'')}=\DCPA{\F}{\G}$.
\item For any $B'_1\subseteq B_1\perp B_2\supseteq B'_2$, $\F_1\in\Mns{B_1,B'_1}$,
and $\F_2\in\Mns{B_2,B'_2}$, $\redasso{(\F_1,B_2)}{(B'_1,\F_2)}=1=\redasso{(B_1,\F_2)}{(\F_1,B'_2)}$.
\end{enumerate}
\end{lemma}

%%%%%%%%%%

%--------------------------------------------------------------------------------------------------
\subsection{Twisting of braided pre--Coxeter structures on $\CDUA{D}{\bullet}$}\label{ss:twist-weak}
%--------------------------------------------------------------------------------------------------

\begin{definition}\hfill
\begin{enumerate}\itemsep0.25cm
\item A {\it twist} in $\CDUA{D}{\bullet}$ is a pair $\bbT=(u,\twF{})$ where
\begin{enumerate}\itemsep0.25cm
			\item $u=\{u_{\F}\}$ is a collection of invertible elements in $\CRDYUA{B,B'}{}$, indexed
			by a maximal nested set $\F\in\Mns{B,B'}$, which satisfy $\varepsilon(u_{\F})=1$ and
			orthogonal factorisation, \ie for any $B'_1\subseteq B_1\perp B_2\supseteq B'_2$, and
			$(\F_1,\F_2)$ in $\Mns{B_1,B_1'}\times\Mns{B_2,B_2'}=\Mns{B_1\sqcup B_2, B'_1\sqcup B'_2}$,
			\[u_{(\F_1,\F_2)}=u_{\F_1}\cdot u_{\F_2}=u_{\F_2}\cdot u_{\F_1}\]
			\item $\twF{}=\{\twF{B}\}$ is a collection of invertible elements of $\CRDYUA{B,B}{2}$,
			indexed by subdiagrams $B\subseteq D$, which satisfy $\varepsilon_2^1
			(\twF{B})=1=\varepsilon_2^2(\twF{B})$, are symmetric,\ie $(\twF{B})_{21}=\twF{B}$,
			\footnote{
				There is a natural action of $\SS_n$ on $\RDYUA{B}{n}$ 
				given by permutations of $\VDY{1}\ten\cdots\ten\VDY{n}$ (cf. \cite[Sec. 7.2]{ATL2}),
				which is a propic version of the action of $\SS_n$ on $U\gb^{\ten n}$.
			} 
			$d(\twF{B})_1=0$, and such that $\twF{B_1\sqcup B_2}=\twF{B_1}\cdot\twF{B_2}$.
		\end{enumerate}
		\item The {\it twisting} of a braided pre--Coxeter structure
		$\pCox{}=(\Phi_B,J_{\F},\DCPA{\F}{\G})$ %on $\CRDYUA{D}{\bullet}$ 
		by a twist $\bbT=(u,\twF{})$  is the braided pre--Coxeter structure 
		\[\pCox{\bbT}=((\Phi_B)_{F_B},(J_{\F})_{(u,\twF{})},(\DCPA{\F}{\G})_{u}, (\redasso{\F}{\F'})_{u})\]
		given by
		\begin{eqnarray*}
			(\Phi_B)_{\twF{}}			&=&\twistp{\twF{B}}{\Phi_B}\\
			(J_{\F})_{(u,\twF{})}			&=&\gaugebis{u_{\F}}{\twF{B}^{-1}\cdot J_{\F}\cdot\twF{B'}}\\
			(\DCPA{\F}{\G})_{u}		&=&u_{\F}^{-1}\cdot \DCPA{\F}{\G}\cdot u_{\G}\\
			(\redasso{\F}{\F'})_{u}	&=&u_{\F'\cup\F}^{-1}\cdot \redasso{\F}{\F'}\cdot u_{\F'}
			\cdot u_{\F}
		\end{eqnarray*}
	\end{enumerate}
\end{definition}

\noindent
\remark The twisting of a braided pre--Coxeter structure does not affect the
$R$--matrix $R_B=\exp(\Omega_B/2)$ (cf. \cite[Sec. 13.2]{ATL2}).

%---------------------------------------------
\subsection{Gauging of twists transformation}\label{ss:gauge-twist}
%----------------------------------------------

\begin{definition}\hfill
	\begin{enumerate}
		\item A {\it gauge} is a collection $a=\{a_B\}$ of
		invertible elements $a_B\in\CRDYUA{B,B}{}$ indexed 
		by subdiagrams $B\subseteq D$ and satisfying
		\[
		a_{B_1\sqcup B_2}=a_{B_1}\cdot a_{B_2}
		\]
		\item The {\it gauging} of a twist $\bbT=(u,\twF{})$ by $a$ is the
		twist $\bbT_a=(u_a,\twF{a})$ given by
		\begin{align*}
		(u_{\F})_a&=a_{B'}\cdot u_\F\cdot a_{B}^{-1}\\
		(\twF{a})_B	&=\gaugetris{a_B}{\twF{B}}
		\end{align*}
	\end{enumerate}
\end{definition}

\noindent
\begin{remark}
	It is easy to see that if $(u,F)$ is a twist, and $a$
	a gauge, the twist of a braided pre--Coxeter structure on
	$\CRDYUA{}{\bullet}$ by $(u,F)$ is the same as that by $(u_a,F_a)$.
\end{remark}

%%%%%%%%%%

%-----------------------------------------------------------------
\subsection{Deformation Drinfeld--Yetter modules}\label{ss:defDY}
%-----------------------------------------------------------------
Let $\b$ be a split diagrammatic Lie bialgebra and $\gb$ its %diagrammatic
Drinfeld double. We explained in \ref{ss:univ-alg} that $\DUA{D}{n}$ is a
universal analogue of the diagrammatic algebra $U\gb^{\otimes n}$. In a similar
vein, we now show that the completion $\CDUA{D}{n}$ introduced in \ref
{ss:grading} is a universal analogue of the trivially deformed diagrammatic algebra
$\hext{U\gb^{\otimes n}}$.

Let for this purpose $\c$ be a Lie bialgebra and $\defDY{\c}$ the category
of \DYt $\c$--modules in topologically free $\sfk{\fml}$--modules. Scaling
the coaction on $V\in\defDY{\c}$ by $\hbar$ yields an isomorphism between
$\defDY{\c}$ and the category $\hDrY{\hextsup{\c}}{\adm}$ of \DYt modules
over the Lie bialgebra $\hextsup{\c}=(\hext{\c},[\cdot,\cdot],\hbar\delta)$, whose
coaction is divisible by $\hbar$. We denote by $\hDYA{\c}{n}$ the algebra
of endomorphisms of the $n$--fold tensor power of the forgetful functor
$\ff_\c:\defDY{\c}\to{\khvect}$. $\hDYA{\c}{n}$ identifies canonically with
the analogous completion defined for $\hDrY{\hextsub{\c}}{\adm}$. 

In the case of the diagrammatic Lie bialgebra $\b$, the realisation functors 
\[\G_{(\b^\hbar_B;V_1,\dots, V_n)}:\MDY{B}{n}\longrightarrow{\vect_{\hext{\sfk}}}\]
corresponding to $V_1,\ldots,V_n\in\hDrY{\b^\hbar_B}{\adm}\cong \hDrY
{\b_B}{\hbar}$ induce a homomorphism $\wh{\rho}_{\b}^{n}:\DUA{D}{n}\to\hDYA{\b}{n}$
which naturally extends to $\CDUA{D}{n}$.\footnote{Note
	that $\defDY{\c}$ can also be identified with the category of \DYt modules
	over the Lie bialgebra $\hextsub{\c}=(\hext{\c},\hbar[\cdot,\cdot],\delta)$
	whose action is divisible by $\hbar$. The corresponding realisation
	functors for $\b_{B,\hbar}$ yield the same homomorphism 
	$\wh{\rho}_{\b}^{n}:\CDUA{D}{n}\to\hDYA{\b}{n}$.}
In particular,
\[\wh{\rho}_{\b_B}^{1}(\pi_{\VDY{1}}\circ\pi^*_{\VDY{1}})=\hbar\sum_i b_ib^i
\aand
\wh{\rho}_{\b_B}^{2}(r_{\VDY{1},\VDY{2}})=\hbar\sum_i b_i\otimes b^i\]
where $\{b_i\},\{b^i\}$ are dual bases of $\b_B$ and $\b_B^*$.\footnote
{Note that the realisation functors corresponding to the tuples $(\b^\hbar_{B};V_1,
	\ldots,V_n)$ and $(\b_B;V_1,\ldots,V_n)$, where $V_1,\ldots,V_n\in\hDrY
	{\b^\hbar_B}{\adm}\cong \hDrY{\b_B}{\hbar}$ do not lead to the same
	homomorphism $\DUA{D}{n}\to\hDYA{\b}{n}$ because $\b^\hbar_{B}$
	is not isomorphic to $\b_B\fml$ as Lie bialgebras.}
Note also that if $B'\subseteq B$, the definition of the subalgebra of $[\b_{B'}]
$--invariants in $\CMUA{B}{n}$ (\S \ref{ss:assoc}) implies that $\CMUA{B,B'}
{n}$ is mapped by $\wh{\rho}_{\b_B}^{n}$ to elements of $\hDYA{\b_B}{n}$
commuting with the diagonal (co)action of $\b_{B'}$.

%-----------------------------------------------------------------------------------
\subsection{From universal algebras to Drinfeld--Yetter modules}\label{ss:univ-weak-cox}
%-----------------------------------------------------------------------------------

We shall make use of the following standard construction due to Drinfeld.
Let $\b$ be a diagrammatic Lie bialgebra, $B\subseteq D$, and $\Phi_B
\in\CRDYUA{B}{3}$ an associator. Then, $\hDrY{\b_B}{\Phi}$ is the braided
monoidal category with the same objects as $\defDY{\b_B}$,
and commutativity and associativity constraints given respectively by
\[\beta_{\b_B}=(1\,2)\circ \wh{\rho}_{\b_B}^{2}(e^{\Omega_B/2})
\aand
\Phi_{\b_B}=\wh{\rho}_{\b_B}^3(\Phi_B)\]

\begin{proposition}\label{pr:transfer}\label{pr:univ pre Cox}
	Let $\b$ be a diagrammatic Lie bialgebra. 
	\begin{enumerate}
		\item A braided pre--Coxeter structure $\pCox{}$ on $\CDUA{D}{\bullet}$ 
		gives rise to a braided pre--Coxeter category $\Cox{\b}{\pCox{}}$ with
		\begin{itemize}
		\item diagrammatic categories $\Cox{\b, B}{\pCox{}}=\hDrY{\b_B}{\Phi_B}$
		\item restriction functors $F_\F^{\pCox{}}:\hDrY{\b_B}{\Phi_B}\to\hDrY{\b_{B'}}{\Phi_{B'}}$
		of the form $\left(\Res_{\b_{B'}\b_B},J_\F^{\pCox{}}\right)$ for some tensor structure
		$J_\F^{\pCox{}}$ on $\Res_{\b_{B'}\b_B}$
		\end{itemize}
		Moreover, $\Cox{\b}{\pCox{}}$ is a deformation of $\Cox{\b}{}$ 
		(cf. \ref{ss:diag-DYLBA}).
		\item A twist $\bbT$ in $\CDUA{D}{\bullet}$ gives rise to a $1$--isomorphism
		$\bbT_\b:\Cox{\b}{\pCox{}}\to\Cox{\b}{\pCox{\bbT}}$.
		%, where $\pCox{\bbT}$ denotes the twisted braided pre--Coxeter structure. 
		\item A gauge $\bfg$ in $\CDUA{D}{\bullet}$ gives rise to a $2$--isomorphism
		$\bfg_{\b}:\bbT_{\b}\Rightarrow(\bbT_{\bfg})_{\b}$. %, where $\bbT_{\bfg}$ denotes the gauged twist. 
	\end{enumerate}
\end{proposition}

\begin{pf}

	$(1)$ 
	Consider the following data.
	\begin{itemize}
		\item {\em Diagrammatic categories.} For any $B\subseteq D$, set $\Cox{\b,B}{\pCox{}}\coloneqq\hDrY{\b_B}{\Phi_B}$.
		\item {\em Restriction functors.} For any $B'\subseteq B$ and $\F\in\Mns{B,B'}$,
		the action of  $J_{\F}^{\pCox{}}\coloneqq\wh{\rho}_{\b_B}^{2}(J_{\F})$ defines a 
		linear automorphism of $V\ten W$, for any $V,W\in\Cox{\b,B}{\pCox{}}$. 
		By the properties of $J_{\F}$, this defines a tensor structure on the standard restriction 
		functor $\Res_{\b_{B'},\b_{B}}$. Then, we set $F_{\F}^{\pCox{}}\coloneqq
		\left(\Res_{\b_{B'},\b_{B}}, J_{\F}^{\pCox{}}\right)$, where the tensor structure
		is given by the natural isomorphism
		\[(J_{\F}^{\pCox{}})_{V,W}:\Res_{\b_{B'},\b_{B}}(V)\ten\Res_{\b_{B'},\b_{B}}(W)\to
		\Res_{\b_{B'},\b_{B}}(V\ten W)\]
		\item {\em De Concini--Procesi associators.} For any $B'\subseteq B$, and $\F,\G\in\Mns{B,B'}$, the action of $\Upsilon_{\G\F}^{\pCox{}}\coloneqq\wh{\rho}_{\b_B}(\DCPA{\G}{\F})$ defines
		a linear automorphism of $V\in\Cox{\b,B}{\pCox{}}$. By the properties of $\DCPA{\G}{\F}$,
		this defines an isomorphism of tensor functors $F_{\F}^{\pCox{}}\Rightarrow F_{\G}^{\pCox{}}$.
		\item {\em Vertical joins.} For any $B''\subseteq B'\subseteq B$, $\F''\in\Mns{B',B''}$, $\F'\in\Mns{B,B'}$, 
		let ${\redasso{\F'}{\F''}}^{\pCox{}}: F_{\F''}^{\pCox{}}\circ F_{\F'}^{\pCox{}}\Rightarrow F_{\F'\cup\F''}^{\pCox{}}$
		be the tensor isomorphism defined by $\wh{\rho}_{\b_B}(\redasso{\F'}{\F''})$, together with the equality
		$\Res_{\b_{B''},\b_{B'}}\circ \Res_{\b_{B'},\b_{B}}=\Res_{\b_{B''},\b_{B}}$.
		
	\end{itemize}
    These satisfy the conditions of Proposition \ref{prop:alpha-reduction}, so that
	$\Cox{\b}{\pCox{}}=(\Cox{\b,B}{\pCox{}}, F_{\F}^{\pCox{}}$,
	 $J_{\F}^{\pCox{}}$, $\DCPA{\F}{\G}^{\pCox{}},{\redasso{\F'}{\F''}}^{\pCox{}})$ 
	 is a  braided pre--Coxeter category. 
	
	$(2)$ Let $\bbT=(u,F)$ be a twist in $\CDUA{D}{\bullet}$ and $\pCox{\bbT}$ the twisted
	braided pre--Coxeter structure (cf. \ref{ss:twist-weak}). Define a $1$--isomorphism 
	$\bbT_{\b}=(H_{B}^{\bbT}, \gamma_{\F}^{\bbT}): \Cox{\b}{\pCox{}}\to \Cox{\b}{\pCox{\bbT}}$ as follows.
	\begin{itemize}
		\item For any $B\subseteq D$, we denote by $H^\bbT_B$ the identity functor on $\Cox{\b,B}{\pCox{}}$ 
		endowed with the tensor structure $\wh{\rho}_B^{2}(F_{B})$. It follows from Definition \ref{ss:twist-weak}
		that $H^\bbT_B$ is a braided tensor equivalence $\Cox{\b,B}{\pCox{}}\to\Cox{\b,B}{\pCox{\bbT}}$. 
		\item For any $B'\subseteq B\subseteq D$ and $\F\in\Mns{B,B'}$, 
		we denote by $\gamma_{\F}^{\bbT}$ the natural isomorphism $F^{\pCox{\bbT}}_{\F}\circ H^{\bbT}_B\Rightarrow 
		H^{\bbT}_{B'}\circ F^{\pCox{}}_{\F}$ induced by $\wh{\rho}_{BB'}(u_{\F})$. %Therefore, by definition of $u$, 
		$\gamma_{\F}^{\bbT}$ is a well--defined isomorphism of tensor functors satisfying the vertical 
		factorisation property.
	\end{itemize}
	
	$(3)$ Finally, let $\bfg$ be a gauge in $\CDUA{\dsg}{\bullet}$ and $\bbT_{\bfg}$ the gauged twist
	(cf. \ref{ss:gauge-twist}). Define a $2$--isomorphism ${\bfg}_{\b}:\bbT_{\b}\Rightarrow(\bbT_{\bfg})_{\b}$ 
	as follows. For any $B\subseteq D$, denote by $v^{\bfg}_B$ the isomorphism of braided tensor functors
	$H^{\bbT}_{B}\Rightarrow H^{\bbT_{\bfg}}_{B}$ given by $\wh{\rho}_B({\bfg}_B)$. It follows from the definition of $\bfg$ 
	that $\gamma^{\bbT_{\bfg}}_{\F}\circ v^{\bfg}_{B}=v^{\bfg}_{B'}\circ\gamma^{\bbT}_{\F}$.
\end{pf}

\begin{definition}\label{de:univ pre Cox}
Let $\b$ be a diagrammatic Lie bialgebra, and $\{\Phi_B\}_{B\subseteq D}$ a
collection of associators. A braided pre--Coxeter category with diagrammatic
categories $\hDrY{\b_B}{\Phi_B}$ is called {\it universal} if it is of the form $\Cox
{\b}{\pCox{}}$, for some braided pre--Coxeter structure $\pCox{}$ on $\CDUA
{D}{\bullet}$.
\end{definition}

%%%%%%%%%%

\subsection{Coherence and minimal data}\label{ss:coherence}
%-----------------------------------------------------

Let $\pCox{}=(\Phi_B, J_{\F}, \DCPA{\F}{\G}, \redasso{\F}{\F'})$
be a braided pre--Coxeter structure on $\CDUA{D}{\bullet}$. We
show in this section that $\pCox{}$ is determined by its vertical joins,
together with a minimal collection of associators, relative twists,
\DCP associators, vertical joins. We shall need two preliminary results.

\subsubsection{}
Let $B'\subseteq B$, $\F\in\Mns{B,B'}$ and ${\bf C}: B'=B_0\subsetneq
B_1\subsetneq \cdots \subsetneq B_{\ell}=B$ a maximal chain from
$B$ to $B'$ corresponding to $\F$ (cf. \ref{ss:chain-ns}). For any $1\leq k\leq\ell$, denote
by $\F_k\in\Mns{B_k, B_0}$ the restriction of $\F$ to $B_k$, and note
that $\F_k=\F_{k-1}\cup\E_k$, where $\E_k$ is the unique element in
$\Mns{B_k, B_{k-1}}$.

\begin{lemma}\label{le:a and b}
Define $\elemasso{\bf C}\in\CRDYUA{B,B'}{}$ by
\begin{equation}\label{eq:bF}
\elemasso{\bf C}=
\redasso{\E_\ell}{\F_{\ell-1}}\cdot\redasso{\E_{\ell-1}}{\F_{\ell-2}}\cdots\redasso{\E_2}{\F_1}
\end{equation}
\begin{enumerate}
\item $\elemasso{\bf C}$ is independent of the choice of ${\bf C}$,
and will be denoted $\elemasso{\F}$.
\item For any $B''\subseteq B'\subseteq B$, $\F'\in\Mns{B,B'}$ and
$\F''\in\Mns{B', B''}$,
\[\redasso{\F'}{\F''}=\elemasso{\F'\cup\F''}\cdot\elemasso{\F'}^{-1}\cdot\elemasso{\F''}^{-1}\]
\end{enumerate}
\end{lemma} 
\begin{pf}
(1) Lemma \ref{lem:orthogonality} (3) implies that $\elemasso{}$ is constant
on the connected components of the graph ${\bf G}_{B,B'}$ (cf. \ref{ss:chain-ns}).
(2) Let ${\bf C}:B''=B_0\subsetneq B_1\subsetneq \cdots \subsetneq B_{\ell}=B$
be a maximal chain such that $B'=B_p$, for some $1\leq p\leq\ell-1$, and the
restriction of ${\bf C}$ to a chain from $B''$ to $B'$ (resp. $B'$ to $B$) corresponds
to $\F''$ (resp. $\F'$). Note that, with respect to the notation established above, we have
\[
\F''=\F_p=\E_{p}\cup\cdots\cup\E_1\aand\F'=\E_{\ell}\cup\E_{\ell-1}\cup\cdots\cup\E_{p+1}
\]
For $1<k\leqslant\ell-p$, we set $\F'_k\coloneqq\E_{p+k}\cup\cdots\cup\E_{p+1}$.
By definition,
\[\begin{split}
\elemasso{\F'\cup\F''}
&=
\redasso{\E_\ell}{\F_{\ell-1}}\cdots\redasso{\E_{p+2}}{\F_{p+1}}\cdot\redasso{\E_{p+1}}{\F_p}
\cdot
\redasso{\E_p}{\F_{p-1}}\cdots\redasso{\E_{2}}{\E_1}\\
&=\redasso{\E_\ell}{\F_{\ell-1}}\cdots\redasso{\E_{p+2}}{\F_{p+1}}\cdot\redasso{\E_{p+1}}{\F_p}
\cdot\elemasso{\F''}
\end{split}
\]
Note that, by the associativity of the vertical joins, 
\[
\redasso{\E_{p+2}}{\F_{p+1}}\cdot\redasso{\E_{p+1}}{\F_p}=\redasso{\E_{p+2}\cup\E_{p+1}}{\F_{p}}\cdot\redasso{\E_{p+2}}{\E_{p+1}}
\]
More in general, for any $1<k<\ell-p$, we have $\F_{p+k}=\F'_k\cup\F_p$ and
\[
\redasso{\E_{p+k+1}}{\F_{p+k}}\cdot\redasso{\F'_k}{\F_p}=
\redasso{\F'_{k+1}}{\F_p}\cdot\redasso{\E_{p+k+1}}{\F'_k}
\]
Therefore, we get
\[\begin{split}
\elemasso{\F'\cup\F''}
&=
\redasso{\E_\ell}{\F_{\ell-1}}\cdots\redasso{\E_{p+2}}{\F_{p+1}}\cdot\redasso{\E_{p+1}}{\F_p}
\cdot\elemasso{\F''}\\
&=
\redasso{\F'_{\ell-p}}{\F_{p}}\cdot\redasso{\E_{\ell}}{\F'_{\ell-p-1}}\cdots\redasso{\E_{p+2}}{\E_{p+1}}
\cdot\elemasso{\F''}\\
&=
\redasso{\F'}{\F''}\cdot\elemasso{\F'}\cdot\elemasso{\F''}
\end{split}\]
where the second identity follows from iterated applications of the associativity of the vertical joins.
\end{pf}

\subsubsection{}\label{sss:dcp-elem}
Let now $B'\subseteq B$, and $\F,\G\in\Mns{B,B'}$. Assume there is a chain of inclusions
$B'=B_0\subseteq B_1\subseteq B_2\subseteq B_3=B$, and $\F_k,\G_k\in\Mns{B_k,B_
{k-1}}$, $1\leq k\leq 3$, such that
\begin{gather*}
\F=\F_1\cup\F_2\cup\F_3
\qquad
\G=\G_1\cup\G_2\cup\G_3\\
\F_1=\G_1\aand
\F_3=\G_3
\end{gather*}
so that $\F,\G$ only differ in the choice of an element in $\Mns{B_2,B_1}$.

\begin{lemma}
The following holds
\[\elemasso{\G}^{-1}\cdot\DCPA{\G}{\F}\cdot\elemasso{\F}=
\elemasso{\G_2}^{-1}
\cdot\DCPA{\G_2}{\F_2}\cdot
\elemasso{\F_2}\]
\end{lemma}
\begin{pf}
The compatibility of the associators $\Upsilon$ with the vertical joins yields
\[\begin{split}
\DCPA{\G}{\F}
&=
\redasso{\G_3}{\G_2\cup\G_1}\cdot\redasso{\G_2}{\G_1}
\cdot
\DCPA{\G_1}{\F_1}\cdot\DCPA{\G_2}{\F_2}\cdot\DCPA{\G_3}{\F_3}
\cdot
\left(\redasso{\F_3}{\F_2\cup\F_1}\cdot\redasso{\F_2}{\F_1}\right)^{-1}\\
&=
\elemasso{\G}\cdot\elemasso{\G_3}^{-1}\cdot\elemasso{\G_2}^{-1}\cdot\elemasso{\G_1}^{-1}
\cdot\DCPA{\G_2}{\F_2}\cdot
\elemasso{\F_1}\cdot\elemasso{\F_2}\cdot\elemasso{\F_3}\cdot\elemasso{\F}^{-1}\\
&=
\elemasso{\G}\cdot\elemasso{\G_2}^{-1}
\cdot\DCPA{\G_2}{\F_2}\cdot
\elemasso{\F_2}\cdot\elemasso{\F}^{-1}
\end{split}\]
where the second identity follows from Lemma \ref{le:a and b}, $\F_1=
\G_1$, and $\F_3=\G_3$, and the third from the invariance of $\DCPA
{\G_2}{\F_2}$ under $[\b_{B_1}]$ and that of $\elemasso{\G_3}$ under
$[\b_{B_2}]$.
\end{pf}

\noindent
\remark
Recall that if $B'\subseteq B$ and $\F,\G\in\Mns{B,B'}$, there is a sequence
$\F=\H_1,\ldots,\H_m=\G$ in $\Mns{B,B'}$ such that $\H_k$ and $\H_{k-1}$
differ by one element \cite[Prop. 3.26]{vtl-4}. We term such a sequence an
{\it elementary sequence}. Moreover, if $\F,\G\in\Mns{B,B'}$ differ by one element, there
are a unique $\ol{B}\in\F\cap\G$, vertices $i\neq j\in\ol{B}$, and
maximal nested sets $\ol{\F},\ol{\G}\in\Mns{\ol{B}, \ol{B}\setminus\{i,j\}}$
such that \[\F=\H'\cup\ol{\F}\cup\H'' \aand \G=\H'\cup\ol{\G}\cup\H''\]
for some $\H'\in\Mns{B,\wt{B}}$, $\H''\in\Mns{\wt{B}\setminus\{i,j\},B'}$,
where 
\[\wt{B}=\ol{B}\cup\bigcup_{\substack{B''\in\cc{B'}\\ B''\not\subseteq\ol{B}}}B''\] 
Then, it follows from the result above and Lemma \ref{lem:orthogonality} (2), that
\[\elemasso{\G}^{-1}\cdot\DCPA{\G}{\F}\cdot\elemasso{\F}=
\elemasso{\ol{\G}}^{-1}
\cdot\DCPA{\ol{\G}}{\ol{\F}}\cdot
\elemasso{\ol{F}}\]

\subsubsection{}
We show below that $\pCox{}$ is determined by the elements $\elemasso{\H}$,
where $\H$ is any maximal nested set, and $\left(\Phi_B,J_{B,B'},\DCPA{\F}{\G}\right)$, 
where $B$ is connected, $B'\subsetneq B$ is a 1--step maximal chain with $B$ 
connected, and $\F,\G$ are maximal 2--step chains of the form $B''\subsetneq B'_1\subsetneq B$
and $B''\subsetneq B'_2\subsetneq B$ respectively, with $B$ connected. 

\begin{proposition}
Let $\pCox{}=(\Phi_B, J_{\F}, \DCPA{\F}{\G}, \redasso{\F}{\F'})$ be a braided
pre--Coxeter structure on $\CDUA{D}{\bullet}$. Then,
\begin{enumerate}
% associator
\item For any $B\subseteq D$,
\[\Phi_B=\displaystyle\prod\Phi_{B_i}\]
where the product is over the connected components of $B$.
 % twist
\item For any $B'\subseteq B$, $\F\in\Mns{B,B'}$ %and ${\bf C}\in p^{-1}(\F)$, 
\[ J_{\F}=
\gaugetris{\elemasso{\F}}{J_{B_\ell,B_{\ell-1}}\cdots J_{B_1,B_0}}\]
where $B'=B_0\subsetneq\cdots\subsetneq B_\ell=B$
is a maximal chain corresponding to $\F$.
% DCP associator
\item For any $B'\subseteq B$, and elementary sequence $\H_1,\ldots,\H_m$
in $\Mns{B,B'}$, 
\[\DCPA{\H_m}{\H_1}=
\elemasso{\H_m}\cdot
\left(\elemasso{\ol{\H}_m}^{-1}\cdot\DCPA{\ol{\H}_m}{\ol{\H}_{m-1}}\cdot\elemasso{\ol{\H}_{m-1}}\right)
\cdots
\left(\elemasso{\ol{\H}_2}^{-1}\cdot\DCPA{\ol{\H}_2}{\ol{\H}_1}\cdot\elemasso{\ol{\H}_1}\right)
\elemasso{\H_1}^{-1}\]
\end{enumerate} 
\end{proposition}
\begin{pf}
(1) is the orthogonal factorisation property of the associators $\Phi$.

(2) For $k=1,\ldots,\ell$, let $\F_k$ be the restriction of $\F$ to $B_k$ and $\E_k$
the unique element in $\Mns{B_k,B_{k-1}}$, so that $\F_k=\E_k\cup\F_{k-1}$, and
$\F=\F_\ell$. Then,
\[\begin{split}
J_{\F}
&=
\gaugetris{\redasso{\E_\ell}{\F_{\ell-1}}}{J_{\E_\ell}\cdot J_{\F_{\ell-1}}}\\
&=
\gaugetris{\redasso{\E_\ell}{\F_{\ell-1}}}{J_{\E_\ell}\cdot \gaugetris{\redasso{\E_{\ell-1}}{\F_{\ell-2}}}{J_{\E_{\ell-1}}\cdot J_{\F_{\ell-2}}}}\\
&=
\gaugetris
{\redasso{\E_\ell}{\F_{\ell-1}}\cdot \redasso{\E_{\ell-1}}{\F_{\ell-2}}}
{J_{\E_\ell}\cdot J_{\E_{\ell-1}}\cdot J_{\F_{\ell-2}}}\\
&\phantom{||}\vdots\\
&=
\gaugetris
{\redasso{\E_\ell}{\F_{\ell-1}}\cdots \redasso{\E_2}{\E_1}}
{J_{\E_\ell}\cdots J_{\E_1}}\\
&=
\gaugetris
{\elemasso{\F}}
{J_{\E_\ell}\cdots J_{\E_1}}
\end{split}\]
where the second identity follows from the invariance of $J_{\E_\ell}$ under $[\b_
{B_{\ell-1}}]$.

(3) follows from Lemma \ref{sss:dcp-elem} and the subsequent remark.
\end{pf}

\subsection{Strict pre--Coxeter structures}\label{ss:strict-pre-Cox-prop}
%-----------------------------------------------------------------

By Proposition \ref{ss:univ-weak-cox}, a braided pre--Coxeter structure $\pCox{}=(\Phi_B,
J_{\F},\DCPA{\F}{\G}, \redasso{\F}{\F'})$ on $\CDUA{D}{\bullet}$ gives rise to a braided
pre--Coxeter category $\Cox{\b}{\pCox{}}$. The following conditions ensure that $\Cox{\b}
{\pCox{}}$ is $\DCPA{}{}$--\emph{strict} or $\redasso{}{}$--\emph{strict} (cf. \ref
{ss:dcp-str-cox-cat} and \ref{ss:a-str-cox-cat}).

We say that
\begin{itemize}
	\item $\pCox{}$ is $\DCPA{}{}$--\emph{strict} if $\DCPA{\F}{\G}=1$ for any $\F,\G\in\Mns{B,B'}$
	\item $\pCox{}$ is $\redasso{}{}$--\emph{strict} if $\redasso{\F}{\F'}=1$ for any $\F\in\Mns{B,B'}$ and $\F'\in\Mns{B',B''}$
	\footnote{In \cite{ATL2} we only consider $\redasso{}{}$--{strict} braided pre--Coxeter structures and
	for simplicity refer to them as braided pre--Coxeter structures. Note also that such a structure
	is essentially a quasi--Coxeter quasitriangular quasi--Hopf algebra
	structure on the diagrammatic algebra $\CDUA{D}{\bullet}$, as defined in \cite{vtl-4}.}
\end{itemize}
The following result shows that we can always restrict to either of these cases.

\begin{proposition}\label{prop:strictness-univ-cox}
	Let $\pCox{}$ be a braided pre--Coxeter structure on $\CDUA{D}{\bullet}$.
	\begin{enumerate}
		\item $\pCox{}$ is twist equivalent to a $\DCPA{}{}$--strict braided pre--Coxeter structure.
		\item $\pCox{}$ is canonically twist equivalent to an $\redasso{}{}$--strict braided pre--Coxeter structure.
	\end{enumerate}
\end{proposition}

\begin{pf}
	$(1)$ The trivialisation of the associators $\DCPA{\F}{\G}$ follows as in Proposition 
	\ref{ss:dcp-str-cox-cat}, and can be thought of as a universal lift of the fact that every 
	pre--Coxeter category is equivalent to a $\DCPA{}{}$--strict one. Equivalently, it is enough 
	to observe that, for any choice of maximal nested sets $\E=\{\E(B,B')\}_{B'\subseteq B}$, 
	$\bbT_{\E}=(u, F)$ with 
	\[u_{\F}=\DCPA{\E(B,B')}{\F}\aand F_B=1_B\]
	is a twist in $\CDUA{D}{\bullet}$, and $\pCox{\bbT_{\E}}$ is a $\DCPA{}{}$--strict 
	braided pre--Coxeter structure.
	
	$(2)$ The trivialisation of the vertical joins $\redasso{\F}{\F'}$ follows as in Proposition 
	\ref{ss:a-str-cox-cat}. Indeed, the result of Lemma \ref{lem:orthogonality} (3) implies that 
	the propic analogues of the diagrams \eqref{eq:rhombus-cat} are trivial in $\CDUA{D}{\bullet}$.
	Equivalently, it is enough to observe that $\bbT=(u,F)$ with
	\[u_{\F}=\elemasso{\F}^{-1}\aand F_B=1_B\] is a twist 
	in $\CDUA{D}{\bullet}$, and $\pCox{\bbT}$ is an $\redasso{}{}$--strict braided pre--Coxeter 
	structure (cf. Proposition \ref{ss:coherence} (4)).
\end{pf}

\noindent\remark
It is easy to see that Proposition \ref{prop:strictness-univ-cox} cannot 
be used to obtain a braided pre--Coxeter structure on $\CDUA{D}{\bullet}$
which is both $\DCPA{}{}$--strict and $\redasso{}{}$--strict.

\section{An equivalence of braided pre--Coxeter categories}\label{s:qcstructure}
%============================================

In this section, we rely on the results of \cite{ATL1-1} to prove the existence
of a braided pre--Coxeter structure $\pCox{}$ on $\CDUA{D}{\bullet}$. We
then show that, for any diagrammatic Lie bialgebra $\b$, the braided pre--Coxeter
category $\Cox{\b}{\pCox{}}$ determined by $\pCox{}$ and $\b$ is equivalent
to that of admissible \DYt modules over the \nEK quantisation $\Q(\b)$ of
$\b$.

\subsection{Factorisable associators}\label{ss:fact ass}
%-----------------------------------------------

Let $\LBA$ be the $\PROP$ describing Lie bialgebras, and $\CRDYUA{\mathsf
{LBA}}{\bullet}$ the corresponding universal algebra.\footnote{Note that $\LBA$
(resp. $\CRDYUA{\mathsf{LBA}}{\bullet}$) coincides with the $\PROP$ (resp. universal
algebra) $\LBA_D$ (resp. $\CRDYUA{D}{\bullet}$) for a diagram $D$ consisting
of a single vertex.}
Let $\LBAo$ be the $\PROP$ describing a Lie bialgebra $\ubb{}$, which
decomposes as the direct sum $\ubb{}=\ubb{1}\oplus\ubb{2}$ of two Lie
bialgebras, and $\CRDYUA{\mathsf{o}}{\bullet}$ the corresponding universal algebra.
Equivalently, $\LBAo$ is the $\PROP$ generated by a Lie bialgebra object
$\ubb{}$, together with two Lie bialgebra idempotents $\theta_1,\theta_2\in
\End(\ubb{})$ satisfying $\theta_1\cdot\theta_2=0=\theta_2\cdot\theta_1$
and $\theta_1+\theta_2=\id_{\ubb{}}$. It therefore coincides with the
$\PROP$ $\LBA_D$ for the diagram $\vertice{1}\;\;\vertice{2}$ consisting
of two disconnected vertices.

Let $\Phi\in\CRDYUA{\mathsf{LBA}}{3}$ be an associator, and $\Phi_{\ubb{}},
\Phi_{\ubb{1}},\Phi_{\ubb{2}}\in\CRDYUA{\mathsf{o}}{3}$ its images under
the homomorphisms $\CRDYUA{\mathsf{LBA}}{\bullet}\to\CRDYUA{\mathsf{o}}{\bullet}$
corresponding to the Lie bialgebras $\ubb{},\ubb{1}$ and $\ubb{2}$
respectively. $\Phi$ is said to be {\it factorisable} if the following holds in 
$\CRDYUA{\mathsf{o}}{3}$
\footnote{The order of the factors is irrelevant, since the images
of $\CRDYUA{\vertice{1}}{n}$ and $\CRDYUA{\vertice{2}}{n}$ commute
in $\CRDYUA{\vertice{1}\;\vertice{2}}{n}=\CRDYUA{\mathsf{o}}{n}$ by
Proposition \ref{pr:univ D alg}.}
\[\Phi_{\ubb{}}=\Phi_{\ubb{1}}\cdot\Phi_{\ubb{2}}\]
This is the case for example if $\Phi$ is a {\it Lie associator}, that is the
exponential of a Lie series in $\Omega_{12}$ and $\Omega_{23}$.

\subsection{A pre--Coxeter structure on $\CDUA{D}{\bullet}$}\label{ss:main-thm-1}
%-------------------------------------------------------------------------------

Let now $D$ be a fixed diagram. By construction, the generating object in
$\mLBA{\dsg}$ is a split diagrammatic Lie bialgebra. In particular, for any
$B\subseteq D$, the subobject $[\b_B]=([1],\dmp{B})$ is a Lie bialgebra
in $\mLBA{\dsg}$. This induces a functor $\mLBA{}\to\mLBA{\dsg}$
which factors through $\mLBA{B}$, and a homomorphism $\wh{\rho}
^n_B:\CRDYUA{\mathsf{LBA}}{n}\to\CRDYUA{B}{n}$.

\begin{theorem}\label{th:assoc pre Cox}
For any factorisable associator $\Phi$ in $\CRDYUA{\mathsf{LBA}}{3}$,
there is a $\DCPA{}{}$--strict braided pre--Coxeter structure
$\pCox{\Phi}^\gstr$ on $\CDUA{D}{\bullet}$ which is trivial modulo
$\wh{\mathfrak{U}}^{\bullet}_{\geqslant 1}$, and such that $\Phi_B
=\wh{\rho}_B^3(\Phi)$ for any $B\subseteq D$.
\end{theorem}

\noindent 
The proof of Theorem \ref{th:assoc pre Cox} is given in \ref{ss:proof-thm-1}.
It relies on our earlier results in \cite{ATL1-1}, which are reviewed in \ref
{ss:rel-twist}--\ref{ss:ATL main}.\\

\noindent\remarks\hfill

\begin{itemize}
\item Theorem \ref{th:assoc pre Cox} and Proposition \ref{prop:strictness-univ-cox} 
imply the existence of an $\redasso{}{}$--strict braided pre--Coxeter
structure $\pCox{\Phi}^\astr$ on $\CDUA{D}{\bullet}$ with associators
$\Phi_B=\wh{\rho}_B^3(\Phi)$, which is canonically
twist equivalent to $\pCox{\Phi}^\gstr$.
\item As mentioned in \ref{ss:dcp-str-cox-cat}, and proved in \cite{ATL3},
the monodromy of the Casimir connection of a \KMA is encoded by an
$\redasso{}{}$--strict pre--Coxeter structure, which is more naturally
compared with $\pCox{\Phi}^\astr$.
\end{itemize}

\begin{corollary}
Let $\Phi\in\CRDYUA{\mathsf{LBA}}{3}$ be a factorisable associator.
Then, for any split diagrammatic Lie bialgebra $\b$, there is a $\DCPA{}
{}$--strict (resp. $\sfa$--strict) braided pre--Coxeter category $\Cox{\b}
{\Phi,\gstr}$ (resp. $\Cox{\b}{\Phi,\astr}$) with
\begin{itemize}
	\item diagrammatic categories
	$(\Cox{\b}{\Phi,\gstr})_B=\hDrY{\b_B}{\Phi_B}=(\Cox{\b}{\Phi,\astr})_B$
	\item restriction functors $\hDrY{\b_B}{\Phi_B}\to\hDrY{\b_{B'}}{\Phi_{B'}}$
	of the form $\left(\Res_{\b_{B'}\b_B},J_\F\right)$ for some tensor structure
	$J_\F$ on $\Res_{\b_{B'}\b_B}$
\end{itemize}
Moreover, $\Cox{\b}{\Phi,\gstr}$ and $\Cox{\b}{\Phi,\astr}$ are canonically 
equivalent braided pre--Coxeter categories.
\end{corollary}
\begin{pf}
This follows by applying Proposition \ref{pr:univ pre Cox} to the braided
pre--Coxeter structures $\pCox{\Phi}^\gstr$, $\pCox{\Phi}^\astr$, and
setting
\[\Cox{\b}{\Phi,\gstr}\coloneqq\Cox{\b}{\pCox{\Phi}^\gstr}
\aand
\Cox{\b}{\Phi,\astr} \coloneqq\Cox{\b}{\pCox{\Phi}^\astr}\]
By the remark above, $\Cox{\b}{\Phi,\gstr}$ and $\Cox{\b}{\Phi,\astr}$
are canonically equivalent. % braided pre--Coxeter categories.
\end{pf}

%-------------------------------------------
\subsection{A relative fiber functor}\label{ss:rel-twist}
%-------------------------------------------

Let $\sLBA(\sfk)$ be the category whose objects are Lie bialgebras, and
morphisms are split embeddings (cf. \eqref{eq:sLBA hom}). Fix an associator
$\Phi$ in $\CRDYUA{\mathsf{LBA}}{3}$. In \cite{ATL1-1}, we construct a
$2$--functor
\[\mathsf{DY}^{\Phi}:\sLBA(\sfk)\longrightarrow\mathsf{Cat}_{\sfK}^{\ten}\]
which assigns
\begin{itemize}\itemsep.25cm
\item to any Lie bialgebra $\b$, the monoidal category $\hDrY{\b}{\Phi}$ of deformation
Drinfeld--Yetter $\b$--modules with associativity constraint $\Phi_{\b}=\wh{\rho}_{\b}^3(\Phi)$
\item to any split embedding $\a\hookrightarrow\b$, a monoidal structure
$J\resped$ on the restriction functor $\Res_{\a,\b}:\hDrY{\b}{\Phi}\to
\hDrY{\a}{\Phi}$
\item to any chain of split embeddings (a \emph{split triple})
$\a\hookrightarrow\b\hookrightarrow\c$, 
an isomorphism of monoidal functors
\[u_{\a,\b,\c}:(\Res_{\a,\b},J_{\a,\b})\circ (\Res_{\b,\c},J_{\b,\c})\longrightarrow
(\Res_{\a,\c},J_{\a,\c})\]
in such a way that, for any chain $\a\hookrightarrow\b\hookrightarrow\c
\hookrightarrow\d$, one has
\begin{equation}\label{eq:u assoc}
u_{\a,\b,\d}\circ u_{\b,\c,\d}=u_{\a,\c,\d}\circ u_{\a,\b,\c}
\end{equation}
as isomorphisms
\[ (\Res_{\a,\b},J_{\a,\b})\circ (\Res_{\b,\c},J_{\b,\c})\circ(\Res_{\c,\d},J_{\c,\d})
 \longrightarrow (\Res_{\a,\d},J_{\a,\d})\]
\end{itemize}
Moreover, $J_{\a,\a}$, $u_{\a,\a,\b}$, and $u_{\a,\b,\b}$ are trivial and, if
$\Phi$ is factorisable, then 
\[J_{\a_1\oplus\a_2,\b_1\oplus\b_2}=J_{\a_1,\b_1}\cdot J_{\a_2,\b_2}
\aand
u_{\a_1\oplus\a_2,\b_1\oplus\b_2,\c_1\oplus\c_2}=
u_{\a_1,\b_1,\c_1}\cdot u_{\a_2,\b_2,\c_2}\]
\remark When $\a=0$, $J_{\a,\b}$ is gauge equivalent to the monoidal
structure on the forgetful functor $\hDrY{\b}{\Phi}\to\vect_K$ constructed
by \nEK \cite{ek-1}.

%---------------------------------------------------------------
\subsection{Functoriality of the \nEK equivalence}\label{ss:comp-twist-quant}
%---------------------------------------------------------------

In \cite{ek-6}, Etingof and Kazhdan define an equivalence of braided monoidal
categories $\EKeq\b:\hDrY{\b}{\Phi}\to\aDrY{\Q(\b)}$, where $\b$ is a Lie bialgebra
and $\Q(\b)$ its \nEK quantisation. 
We prove in \cite{ATL1-1} that the equivalence $\EKeq\b$ is functorial with respect
to split embeddings. Specifically, let $\sQ:\sLBA(\sfk)\to\sQUE(\sfK)$ be the \nEK
quantisation functor between the categories of split Lie bialgebras and split QUEs.
We show that there is an isomorphism of $2$--functors 
\[
\xymatrix@R=.2cm@C=.5cm{
	\sLBA(\sfk) \ar[ddrr]_{\mathsf{DY}^{\Phi}} \ar[rrrr]^{\sQ} && && \sQUE(\sfk) \ar[ddll]^{\aDrY{}}\ar@{=>}[dlll]\\
	&&&&\\
	&& \mathsf{Cat}_\sfK^{\ten} &&
}
\]
which assigns to a Lie bialgebra $\b$ the equivalence $\EKeq\b$.
In particular,
\begin{itemize}
	\item For any split embedding $\a\hookrightarrow\b$, there is a natural
	isomorphism $v\resped$ making the following diagram commute
	\begin{equation}\label{eq:nat transf}
	\xymatrix@C=2cm{
		\hDrY{\b}{\Phi} \ar[r]^{\EKeq\b} \ar[d]_{(\Res_{\a,\b}, J\resped)}& 
		\aDrY{\Uhb} \ar[d]^{(\Res_{\Uha,\Uhb}, \id)}
		\ar@{<=}[dl]_{v\resped}
		\\
		\hDrY{\a}{\Phi} \ar[r]_{\EKeq\a} & 
		\aDrY{\Uha}
	}
	\end{equation}
where $(\Res_{\a,\b}, J\resped)$ is the monoidal functor from \ref{ss:rel-twist}, and
the functor $\Res_{\Uha,\Uhb}$ is induced by the split embedding $\Uha\hookrightarrow
\Uhb$.
	\item For any chain of split embeddings $\a\hookrightarrow\b\hookrightarrow\c$,
	the following prism is commutative
	\begin{equation}\label{eq:cylinder}
	\xy
	%LBA
	(0,20)*+{\hDrY{\c}{\Phi}}="C";
	(15,0)*+{\hDrY{\b}{\Phi}}="B";
	(0,-20)*+{\hDrY{\a}{\Phi}}="A";
	%HA
	(25,30)*+{\aDrY{\EK\c}}="HC";
	(40,10)*+{\aDrY{\EK\b}}="HB";
	(25,-10)*+{\aDrY{\EK\a}}="HA";
	%ARROWS
	{\ar|(.4){\Res_{\a,\b}}"B";"A"};
	{\ar|{\Res_{\b,\c}}"C";"B"};
	{\ar@/_1pc/ "C";"A"_{\Res_{\a,\c}}};
	{\ar|{\Res_{\Uha,\Uhb}}"HB";"HA"};
	{\ar|{\Res_{\Uhb,\Uhc}}"HC";"HB"};
	{\ar|(.4){\Res_{\Uha,\Uhc}}@{..>}@/_1pc/"HC";"HA"};
	{\ar "C";"HC"^{\EKeq{\c}}};
	{\ar "B";"HB"^{\EKeq{\b}}};
	{\ar "A";"HA"_{\EKeq{\a}}};
	{\ar@{=>}_{u_{\a,\b,\c}}"B";(-3,0)};
	\endxy
	\end{equation}
	where $u_{\a,\b,\c}$ is the isomorphism from \ref{ss:rel-twist}, 
	the back 2--face is the identity, and the
	lateral 2--faces are the isomorphisms $v_{\a,\c},v_{\b,\c},
	v\resped$.\footnote{To alleviate the notation, tensor structures
		are suppressed from the diagram \eqref{eq:cylinder}.}
\end{itemize}

\noindent\remarks\hfill
\begin{itemize}
\item The natural isomorphism $v_{\a,\b}$ is not trivial in general.
Indeed, the strict commutativity of \eqref{eq:nat transf}, even as a
diagram of non--monoidal functors, would contradict Prop.~3.2 in
\cite{vtl-4} (see \cite[Sec.~1.10]{A}). Namely, let $\DYA{\b}{\hbar}$ and  
$\DYA{\Q(\b)}{}$ be the algebras of endomorphisms of the 
forgetful functors $\defDY{\b}\to{\khvect}$ and 
$\aDrY{\Q(\b)}\to{\khvect}$, respectively. 
The Etingof--Kazhdan equivalence $H_{\b}:\defDY{\b}\to\aDrY{\Q(\b)}$ 
intertwines the forgetful functors and
gives rise to an isomorphism of algebras
$\DYA{\b}{\hbar}\to\DYA{\Q(\b)}{}$. 
Through the classical and quantum restriction functors, 
we get canonical inclusions $\DYA{\a}{\hbar}\hookrightarrow\DYA{\b}{\hbar}$ 
and $\DYA{\Q(\a)}{}\hookrightarrow\DYA{\Q(\b)}{}$. Therefore,
the strict commutativity of \eqref{eq:nat transf} is equivalent to the 
commutativity of the diagram
\[
\xymatrix@R=0.4cm@C=0.5cm{
\DYA{\b}{\hbar} \ar[r] & \DYA{\Q(\b)}{}\\
\DYA{\a}{\hbar} \ar@{^(->}[u]\ar[r] & \DYA{\Q(\a)}{} \ar@{^(->}[u]
}
\]
In the case of a semisimple Lie algebra $\g$, $\Q(\g)$ is isomorphic to 
the Drinfeld--Jimbo quantum groups $\DJ{\g}$ as a diagrammatic QUE
(cf. \ref{s:qkm} and Proposition \ref{ss:iso-ek-dj}).
Thus, we obtain a diagrammatic isomorphism $\hext{U\g}\to\Q(\g)\simeq\DJ{\g}$,
which contradicts \cite[Pro. 3.2]{vtl-4}.

\item The natural transformation $u_{\a,\b,\c}$ described in \ref
{ss:rel-twist} is in fact \emph{defined} so as to make \eqref{eq:cylinder}
commutative. Namely, since $H_{\a}$ is invertible, $v_{\a,\b}$
induces a natural isomorphism $w_{\a,\b}$
\begin{equation*}
\xymatrix@R=.5cm@C=2cm{
	\hDrY{\b}{\Phi} \ar[r]^{\EKeq\b} \ar[d]_{(\Res_{\a,\b}, J\resped)}& 
	\aDrY{\Uhb} \ar[d]^{(\Res_{\Uha,\Uhb}, \id)}
	\ar@{<=}[dl]_{w\resped}
	\\
	\hDrY{\a}{\Phi} \ar@{<-}[r]_{\EKeq\a^{-1}} & 
	\aDrY{\Uha}
}
\end{equation*}
The natural transformation $u_{\a,\b,\c}$ is then defined as 
\[u_{\a,\b,\c}=w_{\a,\c}^{-1}\circ w_{\a,\b}\circ w_{\b,\c}\]
In particular, this makes the associativity \eqref{eq:u assoc}
of $u$ manifest. Finally, one observes that $w_{\a,\a}$ is trivial
and $w_{\a_1\oplus\a_2,\b_1\oplus\b_2}=w_{\a_1,\b_1}\cdot w_{\a_2,\b_2}$,
so that the normalisation and factorisation properties of $u$ follow.
\end{itemize}

%----------------------------------------------
\subsection{Auxiliary $\PROP$s}\label{ss:aux-prop}
%----------------------------------------------

The constructions described in \ref{ss:rel-twist}--\ref{ss:comp-twist-quant}
are universal, in that the relative twist $J_{\a,\b}$, the natural transformation
$u_{\a,\b,\c}$ and their properties are induced by analogous elements and
relations in the universal algebras associated to the following $\PROP$s.

Let $\mLBA{\mathsf{sp}}$ be the $\PROP$ generated by a Lie bialgebra
object $([1],\mu,\delta)$, together with a Lie bialgebra idempotent $\dmp
{}:[1]\to[1]$. We denote by $\RDYUA{\mathsf{sp}}{\bullet}$ the corresponding
universal algebras. A split embedding of Lie bialgebras $(i,p):\a\to\b$ is
equivalent to a realisation functor $\G_{\b}:\mLBA{\mathsf{sp}}\to\vect_
{\sfk}$ given by
\[\G_{\b}[1]=\b\aand\G_{\b}\dmp{}=i\circ p\]
It therefore gives rise to a map $\rho^{\bullet}_{\a,\b}:\RDYUA{\mathsf{sp}}{\bullet}
\to\DYA{\b}{\bullet}$. We denote the Lie bialgebra objects $[1], \dmp{}[1]$ by $[\b],
[\a]$, respectively.

Let $\mLBA{\mathsf{st}}$ be the $\PROP$ generated by a Lie bialgebra
object $([1],\mu,\delta)$ with idempotents $\dmp{},\dmp{}':[1]\to[1]$ such that 
$\dmp{}\dmp{}'=\dmp{}'=\dmp{}'\dmp{}$. We denote by $\RDYUA{\mathsf{st}}
{\bullet}$ the corresponding universal algebras. A split triple of Lie bialgebras
$(i,p)\circ(i',p'):\a\to\b\to\c$ is equivalent to a realisation functor $\G_{\c}:
\mLBA{\mathsf{st}}\to\vect_{\sfk}$ given by
\[\G_{\c}[1]=
\c,\qquad\G_{\c}\dmp{}=i\circ p\aand\G_{\c}\dmp{}'=i\circ i'\circ p'\circ p\]
It therefore gives rise to a map $\rho^{\bullet}_{\a,\b,\c}:\RDYUA{\mathsf{st}}
{\bullet}\to\DYA{\c}{\bullet}$. We denote the Lie 
bialgebra objects $[1], \dmp{}[1], \dmp{}'[1]$ by $[\c], [\b], [\a]$, respectively.
The $\PROP$ $\mLBA{\mathsf{sq}}$ and its universal algebra 
$\RDYUA{\mathsf{sq}}{\bullet}$, corresponding to split quadruples, are 
similarly defined.

Let $\mLBA{\mathsf{osp}}$ (resp. $\RDYUA{\mathsf
{osp}}{\bullet}$) be the $\PROP$ (resp. its universal algebra) consisting of a split
pair $\pair{\uaa{}}{\ubb{}}$ which decomposes as the direct sum of two split pairs
$\pair{\uaa{1}}{\ubb{1}}$ and $\pair{\uaa{2}}{\ubb{2}}$. 
The $\PROP$s $\mLBA{\mathsf{ost}}$, $\mLBA{\mathsf{osq}}$ and their universal 
algebras $\RDYUA{\mathsf{ost}}{\bullet}$, $\RDYUA{\mathsf{osq}}{\bullet}$, 
corresponding, respectively, to a split triple and a split quadruple with a direct sum 
decomposition, are similarly defined.

%----------------------------------------------------------
\subsection{Universal relative twists and joins}\label{ss:univ-rel}
%----------------------------------------------------------

Let $\Phi\in\CRDYUA{\mathsf{LBA}}{3}$ be an associator. An  element
$J\in\CRDYUA{\mathsf{sp}}{2}$ is a {\it relative twist} if it is such that $
(J)_0=1$, $\varepsilon_2^1(J)=1=\varepsilon_2^2(J)$, it commutes with
the diagonal action and coaction of $\ua$, and satisfies the relative twist
equation with respect to $\Phi$
\[J_{1,23}\cdot J_{23}\cdot \Phi_{\ua} = \Phi_{\ub}\cdot J_{12,3}\cdot J_{12}\]
$J$ is said to be
\begin{itemize}
% normalisation
\item {\it normalised} if $J_{[\a],[\a]}=1$, where $J_{[\a],[\a]}$ is the
image of $J$ under the map $\CRDYUA{\mathsf{sp}}{\bullet}\to\CRDYUA
{\mathsf{LBA}}{\bullet}$, corresponding to the split pair $([1],[1])$ in
$\LBA$
% factorisation
\item {\it factorisable} if $\Phi$ is a factorisable associator, and
\[J_{\uaa{1}\oplus\uaa{2},\ubb{1}\oplus\ubb{2}}=
J_{\uaa{1},\ubb{1}}\cdot
J_{\uaa{2},\ubb{2}}\]
where $J_{\uaa{1},\ubb{1}}$, $J_{\uaa{2},\ubb{2}}$,
$J_{\uaa{1}\oplus\uaa{2},\ubb{1}\oplus\ubb{2}}$ 
are the images of $J$ under the maps $\CRDYUA{\mathsf{sp}}
{\bullet}\to\CRDYUA{\mathsf{osp}}{\bullet}$ induced  by the corresponding
split pairs in $\mLBA{\mathsf{osp}}$.
\end{itemize}

Let $J$ be a relative twist, and denote by $J_{[\a],[\b]}$, $J_{[\b],[\c]}$,
$J_{[\a],[\c]}$ the images of $J$ under the homomorphisms $\CRDYUA
{\mathsf{sp}}{\bullet}\to\CRDYUA{\mathsf{st}}{\bullet}$ induced by the
corresponding to split pairs in $\mLBA{\mathsf{st}}$. An element $u\in
\CRDYUA{\mathsf{st}}{}$ is a {\it vertical join} if is such that $(u)_0=1$,
$\varepsilon(u)=1$, it commutes with the action and coaction of $\ua$,
and satisfies
\[J_{[\a],[\c]}=u_{12}\cdot J_{[\b],[\c]}\cdot J_{[\a],[\b]}\cdot u_1^{-1}\cdot u_2^{-1}\]
$u$ is said to be
\begin{itemize}
\item
{\it normalised} if $u_{[\a],[\a],[\b]}=1=u_{[\a],[\b],[\b]}$, 
where $u_{[\a],[\a],[\b]}$ and $u_{[\a],[\b],[\b]}$ are the images of $u$
under the homomorphism $\CRDYUA{\mathsf{st}}{\bullet}\to\CRDYUA
{\mathsf{sp}}{\bullet}$, induced by the split triples $(\ua,\ua,\ub)$ and
$(\ua,\ub,\ub)$ in $\mLBA{\mathsf{sp}}$
\item 
{\it associative} if
\[u_{[\a],[\b],[\d]}\cdot u_{[\b],[\c],[\d]}=u_{[\a],[\c],[\d]}\cdot u_{[\a],[\b],[\c]}\]
where $u_{[\a],[\b],[\d]}$, $u_{[\b],[\c],[\d]}$, $u_{[\a],[\c],[\d]}$ and $u_{[\a],[\b],[\c]}$
are the images of $u$ under the homomorphisms $\CRDYUA{\mathsf{st}}
{\bullet}\to\CRDYUA{\mathsf{osq}}{\bullet}$, induced by the corresponding split
triples in $\mLBA{\mathsf{sq}}$
\item {\it factorisable} if $\Phi$ and $J$ are factorisable, and
\[u_{\uaa{1}\oplus\uaa{2},\ubb{1}\oplus\ubb{2},\ucc{1}\oplus\ucc{2}}=
u_{\uaa{1},\ubb{1},\ucc{1}}\cdot u_{\uaa{2},\ubb{2},\ucc{2}}\]
where $u_{\uaa{1}\oplus\uaa{2},\ubb{1}\oplus\ubb{2},\ucc{1}\oplus\ucc{2}}$, 
$u_{\uaa{1},\ubb{1},\ucc{1}}$ and $u_{\uaa{2},\ubb{2},\ucc{2}}$
are the images of $u$ under the homomorphisms $\CRDYUA{\mathsf{st}}
{\bullet}\to\CRDYUA{\mathsf{ost}}{\bullet}$, induced by the corresponding split
triples in $\mLBA{\mathsf{ost}}$
\end{itemize}

\subsection{Existence of a universal relative twist and join}\label{ss:ATL main}
%---------------------------------------------------------------------------

\begin{theorem}\label{thm:J u}
Let $\Phi\in\CRDYUA{\mathsf{LBA}}{3}$ be an associator.
\begin{enumerate}
\item There is a relative twist $J\in\CRDYUA{\mathsf{sp}}{2}$, which
is normalised and such that, for any split pair $\a\hookrightarrow\b$,
$J_{\a,\b}=\wh{\rho}_{\a,\b}^2(J)$.
\item There is a vertical join $u\in\CRDYUA{\mathsf{st}}{}$, which is
normalised, associative and such that, for any split triple $\a\hookrightarrow
\b\hookrightarrow\c$, $u_{\a,\b,\c}=\wh{\rho}_{\a,\b,\c}(u)$.
\end{enumerate}
Moreover, if $\Phi$ is factorisable, then so are $J$ and $u$.
\end{theorem}
\begin{pf}
(1) The existence of a relative twist $J\in\CRDYUA{\mathsf{sp}}{2}$ is
proved in \cite[Prop. 7.7, 8.2.2]{ATL1-1}. By construction, $J$ satisfies
$J_{\a,\b}=\wh{\rho}_{\a,\b}^2(J)$ and, by direct inspection, it is normalised
and factorisable (for the latter property, see also \cite[Prop. 2.25]{ATL1-1}).

(2) We show in \cite[Sec.~6.17]{ATL1-1} that the \nEK equivalence 
$H_\b:\hDrY{\b}{\Phi}\to\aDrY{\Q(\b)}$ is $\PROP$ic. Specifically, let
$\DY{\mathsf{UE_{cP}}}$ (resp. $\DY{\QQUE}$) be the $\PROP$ describing
an admissible \DYt module over a co--Poisson universal enveloping 
algebra (resp. over a QUE). Then, the category $\hDrY{\b}{\Phi}$ (resp. $\aDrY{\Q(\b)}$) 
is equivalent to that of realisation functors from $\DY{\mathsf{UE_{cP}}}$ 
(resp. $\DY{\QQUE}$) to $\vect_{\hext{\sfk}}$. Under these identifications, 
$H_{\b}$ arises as the \emph{pullback} of an isomorphism of $\PROP$s 
$H:\DY{\QQUE}\to\DY{\mathsf{UE_{cP}}}$. Similarly, one shows that the 
natural isomorphism $v_{\a,\b}$ is $\PROP$ic, \ie it is induced by 
\begin{equation*}
\xymatrix@C=2cm{
	\DY{\QQUE} \ar[r] \ar[d] & 
	\DY{\mathsf{UE_{cP}}} \ar[d]
	\ar@{=>}[dl]_{v_{[\a],[\b]}}
	\\
	\DY{\mathsf{QUE, sp}} \ar[r] & \DY{\mathsf{UE_{cP}, sp}}
}
\end{equation*}
where 
\begin{itemize}
\item
$\DY{\mathsf{UE_{cP}, sp}}$ (resp. $\DY{\mathsf{QUE, sp}}$) denote the $\PROP$s 
describing a \DYt module over a split pair of co--Poisson universal 
enveloping algebras $[A_0]\to[B_0]$ (resp. over a split pair of QUEs $[A]\to[B]$);
\item the vertical arrows are the canonical functors mapping the generating objects of 
$\mathsf{UE_{cP}}$ and $\QQUE$ to $[A_0]$ and $[A]$, respectively
\item the horizontal arrows are $\PROP$ic Etingof--Kazhdan equivalences
\end{itemize}
The natural transformation $v_{[\a],[\b]}$ is normalised and factorisable,
\ie $v_{[\a],[\a]}$ is trivial and $v_{[\a_1]\oplus[\a_2],[\b_1]\oplus[\b_2]}=
v_{[\a_1],[\b_1]}\cdot v_{[\a_2],[\b_2]}$. The construction of $u_{[\a],[\b],[\c]}$,
and its normalisation, associativity and factorisability follow as in \ref
{ss:comp-twist-quant}, by considering the $\PROP$ic analogue  of the
diagram \eqref{eq:cylinder}.
\end{pf}

%%%%%%%%%%%%%%%%%%%%%%%%%%%%%%%%%%%%%%

%----------------------------------------------------------------------------------
\subsection{$\DCPA{}{}$--strict braided pre--Coxeter structures}\label{ss:gstr-pre-cox}
%----------------------------------------------------------------------------------

It is useful to observe, in analogy with Proposition \ref{ss:dcp-str-cox-cat}, that a $\DCPA
{}{}$--strict braided pre--Coxeter structure on $\CDUA{D}{\bullet}$ is described by the
datum of
\begin{itemize}
	\item for any $B\subseteq D$, an associator $\Phi_B\in\CRDYUA{B,B}{3}$ 
	\item for any $B'\subseteq B$, a relative twist $J_{B'B}\in\CRDYUA{B,B'}{2}$ satisfying
	\[J_{B'B, 1,23}\cdot J_{B'B, 23}\cdot \Phi_{B'} = \Phi_{B}\cdot J_{B'B, 12,3}\cdot J_{B'B, 12}\]
	together with the normalisation $J_{BB}=1$
	\item for any $B''\subseteq B'\subseteq B$, a vertical join $\redasso{}{B''B'B}\in\CRDYUA{B,B''}{}$
	satisfying
	\[J_{B''B}=\gaugetris{\redasso{}{B''B'B}}{J_{B''B'}\cdot J_{B'B}}\]
	together with the associativity
	\[\oalpho{B}{B''}{B'''}\cdot\oalpho{B}{B'}{B''}=\oalpho{B}{B''}{B'''}\cdot\oalpho{B'}{B''}{B'''}\]
	for any $B'''\subseteq B''\subseteq B'\subseteq B$, and the normalisation $\redasso{}{B'B'B}=1=\redasso{}{B'BB}$
\end{itemize}
Moreover, for any $B''_1\subseteq B'_1\subseteq B_1\perp B_2\supseteq B_2'\supseteq B_2''$,
the following holds
\begin{eqnarray*}
\Phi_{B_1\sqcup B_2}&=&\Phi_{B_1}\cdot\Phi_{B_2}\\
J_{B'_1\sqcup B'_2, B_1\sqcup B_2}&=&J_{B'_1 B_1}\cdot J_{B'_2B_2} \\
\redasso{}{B''_1\sqcup B''_2,B'_1\sqcup B'_2, B_1\sqcup B_2}&=&\redasso{}{B''_1B'_1B_1}\cdot\redasso{}{B''_2B'_2B_2}
\end{eqnarray*}

%----------------------------------------------------------------
\subsection{Proof of Theorem \ref{ss:main-thm-1}}\label{ss:proof-thm-1}
%----------------------------------------------------------------

%Let $\Phi\in\CRDYUA{\mathsf{LBA}}{3}$ be a factorisable associator.
We now construct
a $\DCPA{}{}$--strict braided pre--Coxeter structure $\pCox{\Phi}
^\gstr=(\Phi_B, J_{B'B}, \oalpho{B}{B'}{B''})$ in $\CDUA{D}{\bullet}$.

\subsubsection*{Associators}
%------------------------------------

We use the notation from Section \ref{ss:main-thm-1}. For any $B\subseteq D$,
set $\Phi_B=\wh{\rho}_B^3(\Phi)\in\CRDYUA{B}{3}$. Since $\Phi$ is a factorisable
associator, $\Phi_{B_1\sqcup B_2}=\Phi_{B_1}\cdot\Phi_{B_2}$.

\subsubsection*{Relative twists}\label{sss:twist}
%----------------------------

For any $B'\subseteq B\subseteq D$, the Lie bialgebra objects $[\b_{B}]$ and
$[\b_{B'}]$ a split pair in $\mLBA{\dsg}$. This induces a functor $\G_{[\b_{B'}],
[\b_{B}]}:\mLBA{\mathsf{sp}}\to\mLBA{\dsg}$, and a homomorphism $\wh{\rho}
^n_{B'B}:\CRDYUA{\mathsf{sp}}{n}\to\CRDYUA{B}{n}$. Set $J_{B'B}=\wh{\rho}
_{B'B}^2(J)\in\CRDYUA{B,B'}{2}$. By Theorem \ref{ss:univ-rel}, the relative
twists $J_{B'B}$ satisfy the required properties of normalisation and orthogonal
factorisation.

\subsubsection*{Vertical joins}\label{sss:vert-join}

For any chain of subdiagrams $B''\subseteq B'\subseteq B\subseteq D$, the
Lie bialgebra objects $[\b_{B}]$, $[\b_{B'}]$, and $[\b_{B''}]$ induce a functor
$\G_{[\b_{B''}],[\b_{B'}],[\b_{B}]}:\mLBA{\mathsf{st}}\to\mLBA{\dsg}$, and a homomorphism
$\wh{\rho}^n_{B'',B',B}:\CRDYUA{\mathsf{st}}{n}\to\CRDYUA{B}{n}$.
Then, we set $\oalpho{B}{B'}{B''}=\wh{\rho}_{B'', B', B}(u)\in\CRDYUA{B,B''}{}$.
By Theorem \ref{ss:univ-rel}, the vertical joins $\oalpho{B}{B'}{B''}$
satisfy the required properties of associativity, normalisation, and orthogonal
factorisation.\hfill$\qed$

%---------------------------------------------------------------------------------
\subsection{An equivalence of braided pre--Coxeter categories}\label{ss:main-thm-2}
%---------------------------------------------------------------------------------

We now show that the braided pre--Coxeter structures associated to a diagrammatic
Lie bialgebra $\b$ and to its Etingof--Kazhdan quantisation $\Q(\b)$ are equivalent.
\begin{theorem}
Let $\b$ be a split diagrammatic Lie bialgebra. For any factorisable associator
$\Phi\in\CRDYUA{\mathsf{LBA}}{3}$, there is an equivalence of braided pre--Coxeter
categories 
\[
\Hcox{\b}:\Cox{\b}{\Phi,\gstr}\longrightarrow\Cox{\Q(\b)}{\adm}
\]
where $\Cox{\b}{\Phi,\gstr}$ and $\Cox{\Q(\b)}{\adm}$ are defined
in \ref{ss:main-thm-1} and \ref{ss:quantum-DY}, respectively, and whose 
diagrammatic equivalences are given by the Etingof--Kazhdan functors
$H_{\b_B}:\hDrY{\b_B}{\Phi_B}\to\aDrY{\Q(\b_B)}$, $B\subseteq D$.
\end{theorem}

\begin{pf}
By definition, an equivalence 
$\Hcox{\b}:\Cox{\b}{\Phi,\gstr}\longrightarrow\Cox{\Q(\b)}{\adm}$ of braided pre--Coxeter categories is the datum of 
\begin{itemize}
\item For any $B\subseteq D$, an equivalence of braided monoidal categories 
$H_{B}:\hDrY{\b_B}{\Phi_B}\longrightarrow\aDrY{\Q(\b_B)}$
\item For any $B'\subseteq B$, a natural transformation of monoidal functors
\[
\xymatrix@C=2cm{
\hDrY{\b_B}{\Phi_B} \ar[r]^{H_{B}} \ar[d]_{(\Res_{\b_{B'}\b_B}, J_{B'B})}& 
\aDrY{\Q(\b_B)} \ar[d]^{(\Res_{\Q(\b_{B'}),\Q(\b_{B})}, \id)}
\ar@{<=}[dl]_{\gamma_{B'B}}
\\
\hDrY{\b_{B'}}{\Phi_{B'}} \ar[r]_{H_{B'}} & 
\aDrY{\Q(\b_{B'})}
}
\]
\end{itemize} 
satisfying the properties  \ref{ss:cox}. Then, it is enough to set
set $H_B=\EKeq{\b_B}$ and $\gamma_{B'B}=v_{\b_{B'},\b_{B}}$.
The required properties are easily verified and the result follows.
\end{pf}

%%%%%%%%%%%%%%%%%
\section{Kac--Moody algebras}\label{se:kmas}
%%%%%%%%%%%%%%%%%

Let $\sfk$ be a field of characteristic zero, $\bfI$ a finite set, and
$\sfA=(a_{ij})_{i,j\in\bfI}$ a fixed $|\bfI|\times|\bfI|$ matrix with entries
in $\sfk$. We review in this section the definition and basic properties
of the \KMA associated to $\sfA$. Our treatment is a little more
general than \cite{K}, in that we consider realisations of $\sfA$
whose dimension is not assumed to be minimal. Such realisations
will be used in Section \ref{se:diagr kmas} to endow a \KMA and
its Borel subalgebras with a diagrammatic structure.

\subsection{Realisations} 
%-------------------------------

Departing slightly from the terminology in \cite{K}, we define a
{\it realisation} of $\sfA$ to be a triple $\real$, where\footnote
{In \cite{K}, $V$ is additionally required to be of dimension $2
|\bfI|-\rank(\sfA)$.}
\begin{itemize}
\item $V$ is a \fd vector space over $\sfk$
\item $\Pi=\{\alpha_i\}_{i\in\bfI}$ is a linearly independent subset of $V^*$
\item $\Pi^\vee=\{\hcor{i}\}_{i\in\bfI}$ is a linearly independent subset of $V$
\item $\alpha_i(\hcor{j})=a_{ji}$ for any $i,j\in\bfI$ 
\end{itemize}
Given a realisation $\real$, we denote by $V'\subset V$ the $|\bfI|
$--dimensional subspace spanned by $\Pi^\vee$, and by $\Pi^\perp
\subset V$ the $|\bfI|$--codimensional subspace given by the annihilator
of $\Pi$.

\begin{lemma}\label{le:min dim}
Let $(V,\Pi,\Pi^\vee)$ be a realisation of $\sfA$. Then
\begin{enumerate}
\item $\dim V\geqslant 2|\bfI|-\rank(\sfA)$.
\item $\Pi^\perp\subset V'$ if, and only if $V$ is of minimal dimension $2|\bfI|-\rank(\sfA)$.
\end{enumerate}
\end{lemma}
\begin{pf}
% dimension lower bound
(1) Let $\<\Pi\>\subset V^*$ and $\<\Pi^\vee\>\subset V$ be the
subspaces spanned by $\Pi$ and $\Pi^\vee$. Restriction to $\<
\Pi^\vee\>$ gives rise to a surjection $V^*\to\<\Pi^\vee\>^*$ which
maps $\<\Pi\>$ to a subspace $V^*_\sfA$ of dimension $\rank
(\sfA)$. Thus,
\[\dim V-|\bfI|=
\dim\left(V^*/\<\Pi\>\right)\geq
\dim\left(\<\Pi^\vee\>^*/V^*_\sfA\right)=
|\bfI|-r\]
% Pi^\per
(2) $\Pi^\perp$ is of dimension $\dim V-|\bfI|$, while $\Pi^\perp\cap
V_1$ is of dimension $|\bfI|-\rank(\sfA)$.
\end{pf}

\subsection{Subrealisations}
%-----------------------------------

If $\real$ is a realisation of $\sfA$, a {\it subrealisation} of $V$
is a subspace $U\subseteq V$ such that $\Piv\subset U$ and
the restriction of the linear forms $\{\alpha_i\}_{i\in\bfI}$ to $U$
are linearly independent, so that $(U,\left.\Pi\right|_U,\Piv)$ is
a realisation of $\sfA$.

If $(U,\Pi,\Piv)$ is a realisation of $\sfA$, and $U^0$ a \fd vector
space, then $(V=U\oplus U^0,\Pi,\Piv)$ is a realisation of $\sfA$,
$U$ a subrealisation and $U^0$ a {\it null subspace} that is a 
subspace of $V$ contained in $\Pi^\perp$.

\begin{lemma}\label{le:min subrel}
If $(V,\Pi,\Piv)$ is a realisation of $\sfA$, there is a subrepresentation
$U\subseteq V$ of minimal dimension $2|\bfI|-\rank(\sfA)$ and a
null subspace $U^0\subseteq V$ such that $V$ is equal to the
realisation $U\oplus U^0$.
\end{lemma}
\begin{pf}
Note first that $U\subseteq V$ is a subrepresentation iff $V'\subseteq
U$, and $U^\perp\cap\<\Pi\>=0$ or equivalently $U+\Pi^\perp=V$.
Let now $q:V\to V/V'$ be the quotient map. Since $\Pi^\perp\cap V'$
is of dimension $|\bfI|-\rank(\sfA)$, $q(\Pi^\perp)=\Pi^\perp/\Pi^\perp\cap V'$ is
of dimension $\dim V-(2|\bfI|-\rank(\sfA))$. Thus, if $\ol{U}\subset V/V'$ is a
complementary subspace to $q(\Pi^\perp)$, then $U=q^{-1}(\ol{U})$
is a subrepresentation of $V$ of dimension $2|\bfI|-\rank(\sfA)$.
Note also that $U\cap\Pi^\perp=V'\cap\Pi^\perp$ since the \rhs
is contained in the \lhs and their dimensions agree. Let now
$U^0\subset V$ be a complementary subspace to $V'\cap\Pi^\perp$
in $\Pi^\perp$. Then $U^0$ is a null subspace of $V$ such that
$U\oplus U^0=V$.
\end{pf}

\subsection{Morphisms of realisations}
%------------------------------------------------

A {\it morphism} $(V_1,\Pi_1,\Pi_1^\vee)\to(V_2,\Pi_2,\Pi_2^\vee)$
of realisations is a linear map $T:V_1\to V_2$ such that $T\hcor{1,i}
=\hcor{2,i}$ and $T^t \alpha_{2,i}=\alpha_{1,i}$ for any $i\in\bfI$. We
denote the set of such morphisms by $\HomA(V_1,V_2)$.

\begin{proposition}\label{pr:min real}\hfill
\begin{enumerate}
\item Let $T\in\HomA(V_1,V_2)$ be a morphism of realisations.
\begin{enumerate}
\item If $V_1$ is of minimal dimension, $T$ is injective.
\item If $V_2$ is of minimal dimension, $T$ is surjective.
\end{enumerate}
\item Given two realisations $\{(V_i,\Pi_i,\Pi^\vee_i)\}_{i=1,2}$ of $\sfA$,
the set $\HomA(V_1,V_2)$ is non--empty. Moreover, the map
\[\Hom_\sfk(V_1/V_1',\Pi_2^\perp)\times\HomA(V_1,V_2) 
\to
\HomA(V_1,V_2)\]
defined by $(\delta,T)\to T+\delta$ gives $\HomA(V_1,V_2)$ the
structure of a torsor over the abelian group $\Hom_\sfk(V_1/V_1',\Pi
_2^\perp)$.
\item There is, up to (non--unique) isomorphism, a unique realisation
of $\sfA$ of minimal dimension $2|\bfI|-\rank(\sfA)$.
\end{enumerate}
\end{proposition}
\begin{pf}
% To and from minimal dimension
(1a) Since $\alpha_{2,i}\circ T=\alpha_{i,1}$ for any $i\in\bfI$, $\Ker
(T)\subset\Pi_1^\perp\subseteq V_1'$, where the last inclusion holds
by (2) of Lemma \ref{le:min dim}. Since the restriction of $T$ to $V_1'$ is injective, it follows that
so is $T$. (1b) follows from (1a) since $T^t:(V_2^*,\Pi_2^\vee,\Pi_2)
\to(V_1^*,\Pi_1^\vee,\Pi_1)$ is a morphism of realisations of $\sfA^t$. 

% Torsor claim
(2) The second part of the claim is clear, once the non--emptyness
of $\Hom_{\sfA}(V_1,V_2)$ is proved. A linear map $T\in\Hom_\sfk
(V_1,V_2)$ satisfies $T\hcor{i,1}=\hcor{i,2}$ for any $i\in\bfI$ iff in
a(ny) decomposition $V_1=V_1'\oplus V_1''$, $T$ has the block
form $T=\begin{pmatrix}\imath&*\end{pmatrix}$, where $\imath$ is
the map $V_1'=V_2'\hookrightarrow V_2$, $\hcor{1,i}\to\hcor{2,i}$.
Similarly, given a decomposition $V_2=\Pi_2^\perp\oplus \wt{V}_2$,
let $p$ be the map $V_1\to\<\Pi_1\>^*=\<\Pi_2\>^*=V_2/\Pi_2^\perp
\cong\wt{V}_2$ given by assigning to $v_1\in V_1$ the unique $v_2
\in\wt{V}_2$ such that $\alpha_{2,i}(v_2)=\alpha_{1,i}(v_1)$ for any
$i\in\bfI$. Then, $\alpha_{2,i}\circ T=\alpha_{i,1}$ holds for any $i\in
\bfI$ iff $T$ has the block form $T=\begin{pmatrix}*\\p\end{pmatrix}$.
Combining, we see that $T$ is a morphism of realisations iff it has
the form
\[ T=\begin{pmatrix}\imath_{\Pi_2^\vee}&*\\ \imath_{\wt{V}_2}=p_{V'_1}&p_{V'_2}\end{pmatrix}\]
where the equality $\imath_{\wt{V}_2}=p_{V'_1}$ follows because $\alpha_{2,i}
(\hcor{2,j})=a_{ij}=\alpha_{1,i}(\hcor{1,j})$. In particular, $\HomA
(V_1,V_2)$ is non--empty.

(3) It is easy to see that there is a realisation of $\sfA$ of minimal
dimension. Its uniqueness then follows from (2) and (1).
\end{pf}

Abusing language slightly, we shall refer to a realisation of $\sfA$
of minimal dimension $2|\bfI|-\rank(\sfA)$ as {\it the} realisation of
$\sfA$, and denote the underlying vector space by $\h$.

\subsection{Invariant forms}
%----------------------------------

Recall that $\sfA$ is symmetrisable if there is an invertible diagonal
matrix $\sfD=\Diag(d_i)_{i\in\bfI}$ such that $\sfD\sfA$ is symmetric,
that is such that $d_ia_{ij}=d_ja_{ji}$ for any $i,j\in\bfI$.

If $\sfA$ is symmetrisable, an {\it invariant form} on a realisation
$\real$ is a non--degenerate, symmetric bilinear form $\iip{\cdot}{\cdot}$ on
$V$ such that $\iip{\hcor{i}}{\cdot}=d_i^{-1}\alpha_i$.

\begin{proposition} Assume that $\sfA$ is symmetrisable. Then
\begin{enumerate}
\item If $V$ is a realisation of minimal dimension, then any
symmetric bilinear form on $V$ such that $\iip{\hcor{i}}{-}=d_i^{-1}
\alpha_i$ is non--degenerate, and therefore an invariant
form.
\item Any realisation $(V,\Pi,\Piv)$ of $\sfA$ possesses an
invariant form.
\end{enumerate}
\end{proposition}
\begin{pf}
(1) If $v\in V$ is such that $\iip{v}{\cdot}=0$, then $v\in\Pi^\perp
\subset V'$, where the last inclusion follows by part (2)
of Lemma \ref{pr:min real}. The result then follows
from the fact the map $\nu:V'\to V^*$ given by $\hcor{i}
\to d_i^{-1}\alpha_i=\iip{\hcor{i}}{\cdot}$ is an injection.

(2) By Lemma \ref{le:min subrel}, there is a subrepresentation
$U\subseteq V$ of minimal dimension, and a null subspace
$U^0\subset V$ such that $V=U\oplus U^0$. By (1), $U$
admits an invariant form $\iip{\cdot}{\cdot}$. If $\iip{\cdot}{\cdot}^0$ is a non--degenerate
symmetric bilinear form on $U^0$, $\iip{\cdot}{\cdot}\oplus\iip{\cdot}{\cdot}^0$ is
an invariant form on $V$.
\end{pf}

%-------------------------------------------
\subsection{Kac--Moody algebras}\label{ss:km}
%-------------------------------------------

Let $\real$ be a realisation of $\sfA$, and $\wtg=\wtg(V)$
the Lie algebra generated by $V$ and elements $\{e_i, f_i\}
_{i\in\bfI}$, with relations $[h,h']=0$ for any $h,h'\in V$, and
\[[h,e_i]=\alpha_i(h)e_i
\qquad
[h,f_i]=-\alpha_i(h)f_i
\qquad
[e_i,f_j]=\drc{ij}\hcor{i}\]
The Lie algebra $\wtg$ is graded by the root lattice $\sfQ=
\bigoplus_i\IZ\alpha_i\subset V^*$, that is $\wtg=\bigoplus
_{\alpha\in\sfQ}\wtg_\alpha$, where $\wtg_\alpha=\{X\in
\wtg|[h,X]=\alpha(h)X,\,h\in{V}\}$ is finite--dimensional. In
fact, if $\sfQ_+=\bigoplus_{i\in\bfI}\IZ_{\geqslant0}{\alpha}
_i$, then $\wtg$ has the triangular decomposition
\[\wtg=\wt{\n}_-\oplus \wtg_0 \oplus\wt{\n}_+\]
where $\wt{\n}_\pm=\bigoplus_{\alpha\in\sfQ_+\setminus\{0\}}
\wtg_{\pm\alpha}$, and $\wtg_0=V$.

The \KMA corresponding to $\real$ is the quotient $\g=\g(V)
=\wtg/\wtI$, where $\wtI$ is the sum of all (graded) ideals in
$\wtg$ having trivial intersection with $\wtg_0$. $\g$ inherits
the $\sfQ$--grading and triangular decomposition of $\wtg$, 
and $\g_0=V$.\footnote{
If $\GCM{A}$ is a symmetrisable generalised Cartan matrix
(\ie $a_{ii}=2$, $a_{ij}\in\IZ_{\leqslant 0}$, $i\neq j$, and $a_{ij}=0$ implies
$a_{ji}=0$), the ideal $\wtI$ is generated by the Serre relations 
$\mathsf{ad}(e_i)^{1-a_{ij}}(e_j)=0=\mathsf{ad}(f_i)^{1-a_{ij}}
(f_j)$ for any $i\neq j$ \cite{GK}. Note that our terminology
differs slightly from the one given in \cite{K} where $\g(\sfA)$
is called a Kac--Moody algebra only if $\sfA$ is a generalised
Cartan matrix.} 

\begin{lemma}\label{le:funct rel}
Let $T\in\Hom_\sfA(V_1,V_2)$ be a morphism of realisations
of $\sfA$. Then
\begin{enumerate}
\item The assignments $v_1\to T(v_1)$, $e_i\to e_i$, $f_i\to f_i$
extend uniquely to a Lie algebra homomorphism $\g(T):\g(V_1)
\to\g(V_2)$.
\item $\g(T)$ is homogeneous \wrt the $\sfQ$--grading. Its
restriction to
\[V_1=\g(V_1)_0\to\g(V_2)_0=V_2\]
is equal to $T$, and its restriction to $\g(V_1)_\alpha\to\g(V_2)_\alpha$ is
an isomorphism for any $\alpha\in\sfQ\setminus\{0\}$.
\item If $T_1:V_1\to V_2$ and $T_2:V_2\to V_3$ are morphisms
of realisations, then
\[\g(T_2\circ T_1)=\g(T_2)\circ\g(T_1)\aand \g(\id_{V_1})=\id_{\g(V_1)}\]
\end{enumerate}
\end{lemma}
\begin{pf}
(1) The given assignments clearly uniquely determine a Lie algebra
homomorphism $\wtg(T):\wtg(V_1)\to\wtg(V_2)$. If $I_1\subset\wtg_1$
is an ideal, then $\wtg(T)(I_1)$ is stable under the adjoint action of
$V_2$ since the latter factors through $V_2/\Pi_2^\perp\cong\<\Pi_1
\>^*=\<\Pi_2\>^*\cong V_1/\Pi_1^\perp$. Since $\wtg(T)(I_1)$ is also 
stable under the adjoint action of
$e_i=\wtg(T)(e_i)$ and $f_i=\wtg(T)(f_i)$, it is an ideal in $\wtg_2$
and $\wtg(T)$ descends to $\wtg(V_1)/\wtI_1\to \wtg(V_2)/
\wtI_2$. 

(2) The homogeneity of $\g(T)$ is clear, as is the fact that the
restriction of $\g(T)$ to $V_1\to V_2$ is equal to $T$. $\g(T)$
is surjective in degrees $\alpha\neq 0$ since $\wtg(T)$ is. If
$K\subset\g(T_1)$ is the kernel of $\g(T)$, then $K=\bigoplus
_{\alpha\in\sfQ} K_\alpha$, where $K_\alpha=K\cap\wtg(V_1)
_\alpha$. It is easy to check that $K^\times=
\bigoplus_{\alpha\in\sfQ\setminus 0} K_\alpha$ is an ideal
in $\wtg(V_1)$ with trivial intersection with $V$ hence it is
equal to zero.

(3) is clear.
\end{pf}

Let $\LQ$ be the category of $\sfQ$--graded Lie algebras
$\g$ over $\sfk$ which are generated by $\g_0$ and elements
$e_i\in\g_{\alpha_i}$ and $f_i\in\g_{-\alpha_i}$, $i\in\bfI$, with morphisms
$\g_1\to\g_2$ which are homogeneous \wrt $\sfQ$ and map
$e^1_i,f^1_i$ to $e^2_i,f^2_i$. By Lemma \ref{le:funct rel}, $\g(-)$ is
a faithful functor from the category of realisations of $\sfA$
to $\LQ$. It is easy to see that $\g(-)$ is also full.

\subsection{The derived subalgebra $\g(V)'$}\label{ss:derived}
%---------------------------------------------------------

Lemma \ref{le:funct rel} implies in particular that the derived
subalgebras $\g(V_1)'$ and $\g(V_2)'$ corresponding to any
two realisations of $\sfA$ are canonically isomorphic. Indeed,
as vector spaces, each $\g(V_i)$ is easily seen to be $\n_-
\oplus V_i'\oplus\n_+$, and any morphism $T\in\Hom_\sfA
(V_1,V_2)$ restricts to the canonical identification $V_1'=V_2'$.

Moreover, the derived subalgebra $\g(V)'$ admits a presentation
similar to that of $\g(V)$. Namely, let $\wtg'$ the Lie algebra
generated by elements $\{e_i,f_i,\hcor{i}\}$ with relations
\[[\hcor{j},e_i]=a_{ji}e_i
\qquad
[\hcor{j},f_i]=-a_{ji}f_i
\qquad
[e_i,f_j]=\drc{ij}\hcor{i}\]
$\wtg'$ is graded by $\sfQ$, with $\wtg'_0=\h'$, where the latter
is the $|\bfI|$--dimensional span of $\{\hcor{i}\}_{i,\in\bfI}$. The
quotient of $\wtg'$ by the sum $\wtI'$ of its graded ideals with
trivial intersection with $\wtg'_0$ is easily seen to be canonically
isomorphic to $\g(V)'$.

%---------------------------------------------------------------
\subsection{Symmetrisable Kac--Moody algebras}\label{ss:bil on g}
%---------------------------------------------------------------

Assume  that $\sfA$ is symmetrisable, and fix an invertible diagonal
matrix $\sfD=\Diag(d_i)$ such that $\sfD\sfA$ is symmetric. Let $\real$
be a realisation of $\sfA$ endowed with an invariant form $\iip{\cdot}{\cdot}$. Then,
$\iip{\cdot}{\cdot}$ uniquely extends to a symmetric, invariant, non--degenerate bilinear form
on $\g=\g(V)$, which satisfies $\iip{e_i}{f_j}=\delta_{ij}d_i^{-1}$ \cite[Thm. 2.2]{K}.

Recall that $\g$ has a standard $\IZ$--grading with finite--dimensional
homogeneous components, given by $\deg(f_i)=1=-\deg(e_i)$ and $\deg(V)=0$.
Set $\b_{\pm}=V\oplus\bigoplus_{\alpha\in\sfR_+}\g_{\pm\alpha}\subset\g$.  Then, $\b_{\pm}$
are $\mp\IN$--graded Lie algebras with finite--dimensional components.
Moreover, the bilinear form induces a canonical isomorphisms $\b^{\star}_{\pm}\simeq\b_{\mp}$,
where $\b^{\star}_{\pm}$ is the restricted dual of $\b_{\pm}$, and is equal to
\[
\b^{\star}_{\pm}\coloneqq V^*\oplus\bigoplus_{\alpha\in\sfR_+}\g_{\pm\alpha}^*
\]

These identifications allows to determine on $\b_{\pm}$, and therefore
on $\g$, a natural structure of Lie bialgebra compatible with the grading.

More precisely, consider the Lie algebra $\gtwo=\g\oplus V$, and
endow it with the inner product $\iip{\cdot}{\cdot}\oplus-\left\iip{\cdot}{\cdot}
\right|_{V\times V}$. Let $\pi_0:\g\to\g_0=V$ be the projection, and
$\btwo_\pm\subset\gtwo$ the subalgebra
\[\btwo_\pm=\{(X,v)\in\b_\pm\oplus V|\,\pi(X)=\pm v\}\]
Note that the projection $\gtwo\to\g$ onto the first component restricts
to an isomorphism $\btwo_\pm\to\b_\pm$ with inverse $\b_\pm\ni X\to
(X,\pm\pi_0(X))\in\btwo_\pm$.

Then, the following is easily seen to hold (cf. \cite[Ex. 3.2]{drin-2},
\cite[Prop. 2.1]{ek-6}).

\begin{proposition}\label{pr:double g}\hfill
\begin{enumerate}
% Manin triple
\item $(\gtwo, \btwo_-, \btwo_+)$ is a restricted Manin triple.
In particular, $\btwo_\mp$ and $\gtwo$ are Lie bialgebras,
with cobracket $\delta_{\btwo_\mp}=[\cdot,\cdot]_{\btwo_\pm}^t$
and $\delta_{\gtwo}=\delta_{\btwo_-}-\delta_{\btwo_+}$.
% quotienting
\item The central subalgebra $0\oplus V\subset\gtwo$ is a coideal,
so that the projection $\gtwo\to\g$ induces a Lie bialgebra structure
on $\g$ and $\b_\mp$. 
% formulas
\item The Lie bialgebra structure on $\g$ is given by
\[\left.\delta\right|_V=0
\qquad\quad
\delta(e_i)=d_i \hcor{i}\wedge e_i
\qquad\quad
\delta(f_i)=d_i \hcor{i}\wedge f_i
\]
\end{enumerate}
\end{proposition}

%%%%%%%%%%%%%%%%%%%%
\section{Diagrammatic \KM algebras}  \label{se:diagr kmas}
%%%%%%%%%%%%%%%%%%%%

As pointed out in \ref{ss:ex-simple-LA-diag}, a complex semisimple
Lie algebra $\g$ and its positive Borel subalgebra are diagrammatic
Lie bialgebras \wrt the Dynkin diagram of $\g$. The extension of this
result to an arbitrary \KMA requires the introduction of {\it extended}
\KM algebras which correspond to non--minimal realisations of the
underlying matrices. These realisations are defined in this section,
together with a natural braided Coxeter structure on integrable \DYt
modules over the corresponding Borel subalgebras.

\subsection{}
%--------------

Fix an $|\bfI|\times|\bfI|$ matrix $\sfA$ with entries in $\sfk$,
and let $D$ be the diagram having $\bfI$ as its vertex set
and an edge between $i\neq j$ unless $a_{ij}=a_{ji}=0$.
For any $B\subseteq D$, let $\sfA_B$ be the $|B|\times|B|$
matrix $(a_{ij})_{i,j\in B}$, $\g(\sfA_B)$ the \KMA corresponding
to its minimal realisation, and $\h(\sfA_B)$ its Cartan subalgebra.

As pointed out in \ref{ss:derived}, the derived subalgebra
$\g(\sfA)$ is generated by $\{e_i,f_i,\hcor{i}\}_{i\in D}$. It
possesses a diagrammatic structure over $D$ which is
given by associating to any subdiagram $B\subseteq D$
the derived algebra $\g(\sfA_B)'$, and to each inclusion
$B'\subseteq B$ the morphism $\imath'_{BB'}:\g(\sfA_{B'})'
\to\g(\sfA_B)'$ mapping $e_i^{B'},f_i^{B'},{\hcor{i}}^{B'}$
to $e_i^B,f_i^B,{\hcor{i}}^B$, $i\in B'$. This is a diagrammatic
structure since, if $i\perp j$, $e_i$ (resp.  $f_i$) commutes
with $e_j$ (resp. $f_j$) \cite[Lemma 1.6]{K}.

We say that $\g(\sfA)$ is {\it \cdg diagrammatic} if it is endowed
with a diagrammatic structure such that $\g_B=\g(\sfA_B)$ for
any $B\subseteq D$, and the following diagram commutes
for any $B'\subseteq B$
\[\xymatrix{
\g(\sfA_{B'})\ar[r]^{\imath_{BB'}}&\g(\sfA_{B})\\
\g(\sfA_{B'})'\ar[u]\ar[r]_{\imath_{BB'}'}&\g(\sfA_{B})'\ar[u]
}\]
where the vertial arrows are the natural inclusions.

For any $B\subseteq D$, set $\Pi_B=\{\alpha_i\}_{i\in B}$,
$\Piv_B=\{\hcor{i}\}_{i\in B}$, and let $\<\Pi_B\>
\subset\h(\sfA)^*$ and $\h_B'=\<\Piv_B\>\subset\h(\sfA)$ the
subspaces they span respectively.

\begin{proposition}\label{pr:diagr h}\hfill
\begin{enumerate}
\item If $\g(\sfA)$ is \cdg diagrammatic, each morphism $\imath_{BB'}:
\g(\sfA_{B'})\to\g(\sfA_B)$, $B'\subseteq B$, is an embedding.
\item $\g(\sfA)$ is \cdg diagrammatic iff, for any $B\subseteq D$,
there is a subspace $\h_B\subseteq\h(\sfA)$ such that $\left
(\h_B,\left.\Pi_B\right|_{\h_B},\Piv_B\right)$ is a minimal
realisation of $\sfA_B$, that is
\begin{enumerate}
\item $\h'_B\subseteq\h_B$
\item $\<\Pi_B\>\cap\h_B^\perp=0$
\item $\dim\h_B=2|B|-\rank(\sfA_B)$
\end{enumerate}
and, for any $B,B'\subseteq D$
\begin{enumerate}
\item[(d)] $\h_{B'}\subseteq\h_B$ if $B'\subseteq B$
\item[(e)] $\h_B\subseteq\Pi_{B'}^\perp$ and $\h_{B'}\subseteq\Pi_B^\perp$ if $B\perp B'$
\end{enumerate}
\end{enumerate}
\end{proposition}
\begin{pf}
(1) It suffices to show that the restriction $\imath^\h_{BB'}$ of $\imath_{BB'}$
to a map $\h(\sfA_{B'})\to\h(\sfA_B)$ is injective for any $B'\subseteq B$. Applying $\imath
_{BB'}$ to the relation $[h,e_i^{B'}]=\alpha_i^{B'}(h)e_i^{B'}$ shows that $\alpha_i^B\circ
\imath^\h_{BB'}=\alpha_i^{B'}$ for any $i\in B'$. It follows that $\Ker\imath^\h_
{BB'}$ is contained in $\Pi_{B'}^\perp\subseteq\h(\sfA_{B'})'$, where the inclusion
holds by Lemma \ref{le:min dim}. Since the restriction of $\imath^\h_{BB'}$
to $\h(\sfA_{B'})'$ is injective by assumption, the conclusion follows.

(2) Assume that $\g(\sfA)$ is diagrammatic, and set $\h_B=\imath_{DB}(\h(\sfA_B))$.
Since $\imath_{DB}(\hcor{i})$ $=\hcor{i}$ and $\alpha_i^D\circ\left.\imath_{DB}\right|
_{\h_B}=\alpha_i^B$ for any $i\in B$, $\h_B$ contains $\h_B'$ and the restrictions
of the linear forms $\alpha_i^D$ to $\h_B$ are linearly independent. Moreover,
$\h_B$ has the claimed dimension since $\imath_{DB}$ is injective by (1). The
remaining properties are clear.

Conversely, assume given subspaces $\h_B$ satisfying the above properties.
For any $B$, the triple $(\h_B,\left.\Pi_B\right|_{\h_B},\Pi_B^\vee)$ is a minimal
realisation of $\sfA_B$, which determines a morphism of realisations $\imath_
{DB}^\h:\h(\sfA_B)\to\h$ with image $\h_B$. Since, for any $B'\subseteq B$
the image of $\imath_{DB'}^\h$ is contained in the image of $\imath_{DB}^\h$, there
is a uniquely defined morphism of realisations of $\sfA_{B'}$ such that $\imath_{BB'}^\h:\h_{B'}\to\h_B$ such that $\imath_{DB}
^\h\circ\imath_{BB'}^\h=\imath_{DB'}^\h$. Let now
$B''\subseteq B'\subseteq B$. We wish to show that $\imath_{BB'}^\h\circ\imath
_{B'B''}^\h=\imath_{BB''}^\h$. It suffices to show that this holds after composition
with $\imath_{DB}^\h$ since the latter is injective. However,
\[\imath_{DB}^\h\circ\imath_{BB'}^\h\circ\imath_{B'B''}^\h=
\imath_{DB'}^\h\circ\imath_{B'B''}^\h=
\imath_{DB''}^\h=
\imath_{DB}^\h\circ\imath_{BB''}^\h\]
The morphisms of realisations $\imath_{BB'}^\h$ canonically induce Lie algebra
homomorphisms $\imath_{BB'}:\g(\sfA_{B'})\to\g(\sfA_B)$ which give rise to a
Cartan diagrammatic structure on $\g(\sfA_B)$.
\end{pf}

In \ref{ss:corank 2}--\ref{ss:suff and counter} we give sufficient conditions
for $\g(\sfA)$ to be \cdg diagrammatic, together with a %number of
counterexample which show that $\g(\sfA)$ is not \cdg diagrammatic
in general.

\subsection{}\label{ss:corank 2}
%--------------

\begin{lemma}\label{le:corank 2}
If $\det(\sfA_B)\neq0$ for any $B\subset D$ with $|D\setminus B|
\geqslant 2$, then  $\g(\sfA)$ is \cdg diagrammatic.
\end{lemma}
\begin{pf}
We rely on part (2) of Proposition \ref{pr:diagr h}. For any $B$ such that
$|D\setminus B|\geqslant 2$, set $\h_B=\h_B'$. If $|D\setminus B|=1$,
Lemma \ref{le:min subrel} implies that $\h(\sfA)$ contains a subspace
$\h_B$ such that $(\h_B,\left.\Pi_B\right|_{\h_B},\Piv_B)$ is a minimal
realisation of $\sfA_B$. If $B$ is perpendicular to the single vertex $i$
in $D\setminus B$, we require additionally that $\h_B$ be chosen in
$\Ker(\alpha_i)$. Finally, if $B=D$, set $\h_B=\h(\sfA)$. It is easy to
see that the subspaces $\h_B$ satisfy the conditions of Proposition
\ref{pr:diagr h} except possibly the orthogonality condition (d) when
$B$ is such that $|D\setminus B|=1$. If $i$ is the single vertex in $D
\setminus B$ and $a_{ii}\neq 0$, then (d) holds with $B'=i$ by
construction. If $a_{ii}=0$ then, by assumption, $\sfA$ must be
the diagonal matrix $\Diag(*,0)$, and $\g(\sfA)$ is readily seen
to be diagrammatic in this case.
\end{pf}

\noindent\remark
The converse of Lemma \ref{le:corank 2} does not hold. 
Indeed, let $\sfA$ be the zero matrix, which for $n\geqslant 3$ does
not satisfy the above condition. Its minimal realisation can be taken
to be the 2$|\bfI|$--dimensional vector space $\h$ with basis
$\{\hcor{i}\}_{i\in\bfI}\cup\{\partial_i\}_{i\in\bfI}$, and $\{\alpha_i\}_{i\in
\bfI}\subset\h^*$ the last $|\bfI|$ elements of the corresponding dual
basis, so that $\alpha_i(\hcor{j})=0$ and $\alpha_i(\partial_j)=\delta
_{ij}$ for any $i,j\in\bfI$. The corresponding \KMA $\g(\sfA)$ is \cdg
diagrammatic, with $\g_B$ the Lie subalgebra of $\g(\sfA)$
generated by $\{e_i,f_i,\hcor{i},\partial_i\}_{i\in B}$, $B\subseteq D$.

\subsection{}\label{ss:suff and counter}
%--------------

Assume in this paragraph that $\sfk=\IQ$, and that $\sfA$ is such
that $a_{ij}\leqslant 0$ for $i\neq j$ and that $a_{ij}=0\Leftrightarrow a_{ji}
=0$. Recall that if $\sfA$ is indecomposable, it is called {\it finite} if $\rank(\sfA)
=|\bfI|$, {\it affine} if $\rank(\sfA)=|\bfI|-1$, and {\it indefinite} otherwise.
$\sfA$ is {\it hyperbolic} if it is indefinite, but the irreducible components
of any $\sfA_B$, with $B\subsetneq D$, are all of finite or affine type.
In $\sfA$ is finite or affine, then any submatrix $\sfA_B$, with $B
\subsetneq D$ decomposes into a direct sum of matrices of finite
type \cite[Chap. 4]{K}.

If $\sfA$ is a direct sum of indecomposable matrices $\sfA_1\oplus
\cdots\oplus\sfA_m$. Then $\g(\sfA)\cong\g(\sfA_1)\oplus\cdots\oplus
\g(\sfA_m)$ is \cdg diagrammatic iff each $\g(\sfA_i)$ is.

\begin{proposition}
Assume that $\sfA$ is indecomposable. Then
\begin{enumerate}
\item $\g(\sfA)$ is \cdg diagrammatic if $\sfA$ is of finite, affine or
hyperbolic type.
\item $\g(\sfA)$ is not \cdg diagrammatic in general.
\end{enumerate}
\end{proposition}

\begin{pf}
$(1)$ is an immediate consequence of Lemma \ref{le:corank 2}.
To prove $(2)$, we consider the following example.
Let $\sfA$ be the generalised Cartan matrix
\[\sfA=\left[
\begin{array}{rrrr}
2&-1&0&0\\
-1&2&-2&0\\
0&-2&2&-1\\
0&0&-1&2
\end{array}
\right]\]
Note that $\sfA_B$ is of full rank if $|B|=3$, so that $\h_B=\h'_B$
for any such $B$. Then, $\dim\h_{23}=3$, while $\h_{123}\cap\h_{234}=
\h'_{123}\cap\h'_{234}=\h'_{23}$ is of dimension 2. Therefore the condition 
$\h_{23}\subseteq\h_{123}\cap\h_{234}$ cannot be satisfied.
\end{pf}

\subsection{The canonical realisation}\label{ss:diag-KM}
%------------------------------------------------

To remedy the fact that $\g(\sfA)$ is not diagrammatic in general,
we follow a suggestion of P. Etingof, and give in \ref{ss:ext-KM}
a modified definition of $\g(\sfA)$ along the lines of \cite{FZ}. The
corresponding Cartan subalgebra is given by the following
(non--minimal) realisation of $\sfA$.

Let $(\olh,\olPi,\olPiv)$ be the realisation given by $\olh\cong\sfk^
{2|\bfI|}$ with basis $\{\hcor{i}\}_{i\in\bfI}\cup\{\cow{i}\}_{i\in\bfI}$, $\olPiv
=\{\hcor{i}\}_{i\in\bfI}$ and $\olPi=\{\alpha_i\}_{i\in\bfI}$, where $\alpha_i$
is given by
\[\alpha_i(\hcor{j})=a_{ji}
\aand
\alpha_i(\cow{j})=\delta_{ij}\]
We refer to $(\olh,\olPi,\olPiv)$ as the {\it canonical} realisation of
$\sfA$, and denote by $\Lambdav\subset\olh$ the $|\bfI|$--dimensional
subspace spanned by $\{\cow{i}\}_{i\in\bfI}$. 

\begin{proposition}\label{pr:homs ext}
Let $(V,\Pi,\Piv)$ be a realisation of $\sfA$.
\hfill
\begin{enumerate}
% \olh\to V
\item If $p\in\HomA(\olh,V)$, then $p(\Lambdav)\subset V$ is
a complementary subspace to $\Pi^\perp$. Moreover, the map
\[\HomA(\olh,V)\to \{\wtV\subseteq V|\,\Pi^\perp\oplus\wtV=V\},
\quad\quad
p\to p(\Lambdav)\]
is a bijection.
% \h\to\olh
\item If $\imath\in\HomA(V,\olh)$, then $\imath^{-1}(\Lambdav)
\subset V$ is a complementary subspace to $V'$. Moreover,
the map
\[\HomA(V,\olh)\to \{V''\subseteq V|\,V'\oplus V''=V\},
\quad\quad
\imath\to \imath^{-1}(\Lambdav)\]
is a bijection.
% compatibility
\item If $\imath\in\HomA(V,\olh)$ and $p\in\HomA(\olh,V)$
correspond to the subspaces $V'',\wtV\subset V$ respectively,
then $p\circ\imath=\id_V$ if, and only if, $V''\subset\wtV$.
\end{enumerate}
\end{proposition}
\begin{pf}
(1) Since $p$ is a morphism, $\Ker(p)\subset p^{-1}(\Pi^\perp)
\subseteq\olPi^\perp$. It follows in particular that $p(\Lambdav)
\subset V$ is an $|\bfI|$--dimensional subspace with trivial
intersection with $\Pi^\perp$ since $\Lambdav\cap\olPi^\perp
=0$. Let now $\wtV\subset V$ be a complementary subspace
to $\Pi^\perp$. Then, $\wtV\cong\Pi^*=\olPi^*\cong\Lambdav$
so there is a unique map $q:\Lambdav\to\wtV$ such that
$\alpha_i\circ q=\alpha_i$, and therefore a unique morphism
of realisations $p=\id_{\olh'}\oplus q:\olh\to V$ such that
$p(\Lambdav)=\wtV$.\\
(2) $\imath^{-1}(\Lambdav)$ has trivial intersection with $V'$
since $\imath(V')\subseteq\olh'$. Moreover, $V=V'+\imath^{-1}(V)$.
Indeed, let $\imath',\imath''$ be the components of $\imath$
corresponding to the decomposition $\olh=\olh'\oplus\Lambdav$.
Then, for any $v\in V$,
\[\imath(v)=\imath'(v)+\imath''(v)=
\imath(\left.\imath\right|_{V'}^{-1}\circ\imath'(v))+\imath''(v)\]
so that $v-\left.\imath\right|_{V'}^{-1}\circ\imath'(v)\in\imath^
{-1}(\Lambdav)$. Finally, note that the restriction of $\imath$
to $\imath^{-1}(\Lambdav)$ is necessarily given by
$\imath(v)=\sum_i \alpha_i\circ\imath(v)\cow{i}=\sum_i
\alpha_i(v)\cow{i}$, so that $\imath$ is uniquely determined
by the subspace $\imath^{-1}(\Lambdav)$. Conversely,
given a decomposition $V=V'\oplus V''$, then $\imath=
\imath'\oplus\imath'':V\to\olh$, where $\imath'$ is the
canonical identification $V'\to\olh'$, and $\imath'':V''
\to\Lambdav$ is given by $v\to\sum_i\alpha_i(v)\cow{i}$
is easily seen to be the unique morphism of realisations
such that $V''=\imath^{-1}(\Lambdav)$.\\
(3) If $p\circ\imath=\id_V$, then $V''=p\circ\imath(V'')
\subseteq p(\Lambdav)=\wtV$ since $V''=\imath^{-1}
(\Lambdav)$. To prove the converse, it suffices to show
that the restriction of $p\circ\imath$ to $V''$ is the identity.
This follows from the fact that
a) for any $v''\in V''$,  $\imath(v'')$ is the unique $\lambda
^\vee\in\Lambdav$ such that $\alpha_i(v'')=\alpha_i
(\lambda^\vee)$ for any $i\in\bfI$, b) for any $\lambda
^\vee\in\Lambdav$, $p(\lambda)$ is the unique element
$\wt{v}\in\wtV$ such that $\alpha_i(\lambda^\vee)=
\alpha_i(\wt{v})$ for any $i\in\bfI$ and c) $V''\subseteq
\wtV$.
\end{pf}

%---------------------------------------------
\subsection{Extended \KM algebras}\label{ss:ext-KM}
%---------------------------------------------

We denote by $\olg=\olg(\sfA)$ the {\it extended} \KMA
corresponding to $\sfA$, that is the Lie algebra associated
to the canonical realisation of $\sfA$. In particular, $\olg$ is
generated by $\{e_i,f_i,\hcor{i},\cow{i}\}_{i\in\bfI}$, with
relations $[\hcor{i},\hcor{j}]=0$, $[\fcw_i,\fcw_j]=0$, $[\hcor{i},\fcw_j]=0$, 
\[[\hcor{i},e_j]=a_{ji}e_j,
\qquad [\hcor{ji},f_j]=-a_{ji}f_j,
\qquad
[\fcw_i,e_j]=\drc{ij}e_j ,
\qquad [\fcw_i,f_j]=-\drc{ij}f_j,
\]
and $[e_i,f_j]=\drc{ij}h_i$, for any $i,j\in\bfI$.
Unlike $\g(\sfA)$, $\olg(\sfA)$ always possesses a diagrammatic
structure over the Dynkin diagram $D$ of $\sfA$.
\begin{proposition}
The extended Kac--Moody algebra $\olg$ is a diagrammatic
Lie algebra, with diagrammatic Lie subalgebras $\olg_B\coloneqq
\langle e_i,f_i,\hcor{i},\fcw_i\;|\; i\in B\rangle
=\olg(\sfA_B)$, $B\subseteq D$.
\end{proposition}

%-----------------------------------------------------------------
\subsection{Relation between $\g$ and $\olg$}\label{ss:h''}
%-----------------------------------------------------------------

The following shows that $\olg$ is non--canonically a split central
extension of $\g$, with a $\rank(\sfA)$--dimensional kernel. Let
$\LQ$ be the category of $\sfQ$--graded Lie algebras defined in
\ref{ss:km}.

\begin{proposition}\hfill
\begin{enumerate}
\item Any $p\in\Hom_{\LQ}(\olg,\g)$ is surjective, and $\Ker(p)$
is a $\rank(\sfA)$--dimensional subspace of $\olPi^\perp=\zee
(\olg)$ which is complementary to $\olPi^\perp\cap\olh'$.
\item There is a bijection between $\Hom_{\LQ}(\olg,\g)$ and the
set of  subspaces $\wth\subset\h$ which are complementary to 
$\Pi^\perp$, given by mapping $p:\olg\to\g$ to $p(\Lambdav)$.
\item Any $i\in\Hom_{\LQ}(\g,\olg)$ is injective.
\item There is a bijection between $\Hom_{\LQ}(\g,\olg)$ and the
set of  subspaces $\wth\subset\h$ which are complementary to 
$\h'$, given by mapping $i:\g\to\olg$ to $i^{-1}(\Lambdav)$.
\item If $p\in\Hom_{\LQ}(\olg,\g)$ and $i\in\Hom_{\LQ}(\g,\olg)$
correspond to the subspaces $\wth$ and $\h''\subset\h$ respectively,
then $p\circ i=\id_\g$ if, and only if $\wth\subset\h''$.
\end{enumerate}
\end{proposition}
\begin{pf}
(1) By \ref{ss:km}, $p$ is of the form $\g(p_0)$ for a unique
$p_0\in\Hom_A(\olh,\h)$. $p$ is surjective by part (2) of Lemma
\ref{le:funct rel} and part (1b) of Proposition \ref{pr:min real}. 
Moreover, $\Ker(p)=\Ker(p_0)$ is a $\rank(\sfA)$ dimensional
subspace of $\olPi^\perp$ since $\alpha_i\circ p_0=\alpha_i$.
Since $p_0$ is injective on $\olh'$, $\Ker(p_0)\cap(\olPi^\perp
\cap\olh')=0$ and it follows that the two spaces are in direct
sum since their dimensions add up to $|\bfI|=\dim\olPi^\perp$.

(3) The injectivity of $i$ follows from \ref{ss:km} and part (1a)
of Proposition \ref{pr:min real}. 

(2), (4) and (5) Follow from \ref{ss:km} and Proposition \ref
{pr:homs ext}.
\end{pf}

\subsection{Split diagrammatic structure}\label{ss:split-diag-olb}
%----------------------------------------------------

Assume henceforth that $\sfA$ is symmetrisable. Fix $\sfD=\Diag
(d_i)$ such that $\sfD\sfA$ is symmetric, and an invariant form
$\iip{\cdot}{\cdot}$ on $\ol{\h}$. Then, by Proposition \ref{ss:bil
on g}, there is a standard Lie bialgebra structure on $\olg=\olg
(\sfA)$ given by
\[\left.
\delta(\hcor{i})=0=\delta(\fcw_{i})
\qquad\quad
\delta(e_i)=d_i \hcor{i}\wedge e_i
\qquad\quad
\delta(f_i)=d_i \hcor{i}\wedge f_i
\right.\]
It follows as in \ref{ss:ext-KM} that $\olg$ 
is a diagrammatic Lie bialgebra with Lie subbialgebras
$\olg_B=\langle e_i,f_i,\hcor{i}, \fcw_i\;|\; i\in B\rangle$, 
$B\subseteq D$.\\  

As in the finite--dimensional case described in Example 
\ref{ss:ex-simple-LA-diag}, the diagrammatic structure on 
$\olg$ determines a split diagrammatic one on $\olb_\pm$. 
For any $B\subseteq D$, let $\olb_{\pm,B}=\olb_\pm\cap\olg_B$ be the
Lie subbialgebras generated by $\{\hcor{i},\fcw_i,e_i\}_{i\in B}$ and $
\{\hcor{i},\fcw_i,f_i\}_{i\in B}$ respectively. If $B'\subseteq B$, 
let $i_{\pm,BB'}:\olb_{\pm,B'}\to\olb_{\pm,B}$ be the standard embedding, 
and regard $p_{\pm,B'B}=i_{\mp,BB'}^t$ as a map $\olb_{\pm,B}\to\olb_ {\pm,B'}$
via the identifications $\olb_{\mp,C}^\star\cong\olb_{\pm,C}$
given by the inner product, where as usual $\olb^{\star}_{\pm}$ 
is the restricted dual of $\olb_{\pm}$, and is equal to
\[
\olb^{\star}_{\pm}\coloneqq \olh^*\oplus\bigoplus_{\alpha\in\sfR_+}\olg_{\pm\alpha}^*
\] 
Then, $\ker(p_{\pm, B'B})$ is a Lie subalgebra in $\olb_{\pm, B}$, 
and therefore $\{p_{\pm,B'B}\}$ give the required splitting of the Lie
bialgebra $\olb_\pm$ (cf.~\ref{ss:diag-sLBA}). The
splitting can also be explicitly described as follows. Set $\ekm{\n}
_{B,\pm}=\bigoplus_{\alpha\in\sfR_{B,+}}\olg_{\pm\alpha}
\subset\olb_{B,\pm}$.

\begin{lemma}
The projection $p_{\pm,B'B}:\olb_{\pm,B}\to\olb_{\pm,B'}$ corresponds
to the splitting
\[\ekm{\n}_{B,\pm}=
\ekm{\n}_{B',\pm}\oplus\ekm{\n}_{B'B,\pm} 
\qquad\mbox{where}\qquad
\ekm{\n}_{B'B,\pm}=\bigoplus_{\alpha\in\sfR_{B,+}\setminus\sfR_{B',+}}\olg_{\pm\alpha}\]
together with the orthogonal splitting
\[
\olh_B=\olh_{B'}\oplus\olh_{B'B} 
\qquad\mbox{where}\qquad
 \olh_{B'B}=\bigoplus_{j\in B\setminus B'}\sfk\cdot\fcw_j\oplus\sfk\cdot\omega^{\vee}_{B'B, j}
\]
and $\omega^{\vee}_{B'B,j}$ is given by $\hcor{j}-\sum_{i\in B'}\alpha_i(\hcor{j})\fcw_i$.
In particular, $ \olh_{B'B}\subset\bigcap_{i\in B'}\Ker(\alpha_i)$.
\end{lemma}
\begin{pf}
It is enough to observe that for any $i\in B'$ and $j\in B\setminus B'$, 
\[\iip{\hcor{i}}{\fcw_j}=0=\iip{\fcw_i}{\fcw_j} 
\aand
\iip{\hcor{i}}{\omega^{\vee}_{B'B,j}}=0=\iip{\fcw_i}{\omega^{\vee}_{B'B,j}}\]
\end{pf}

%%%%%%%%%%%%%%%%%%%%%%%%%%%%%%%%%%%%%%%%%%%
% IV. Representation theory for extended KM algebras: category O and DY modules	%
%%%%%%%%%%%%%%%%%%%%%%%%%%%%%%%%%%%%%%%%%%%

%-----------------------------------------------
\subsection{The category $\O_{\olg}$}\label{ss:cat-O-lba}
%-----------------------------------------------

A $\olg$--module $V$ is in category $\O_{\olg}$ if the following holds.
\begin{enumerate}
	\item[($\O1$)] $V=\bigoplus_{\lambda\in\olh^*}V_{\lambda}$, where
	$V_{\lambda}=\{v\in V|\,h\,v=\lambda(h)v,\,h\in\olh\}$\\[.05ex]
	\item[($\O2$)]  $\dim V_{\lambda}<\infty$ for any $\lambda\in\sfP(V)
	=\{\lambda\in\olh^*|\,V_{\lambda}\neq0\}$\\[.05ex]
	\item[($\O3$)]  $\sfP(V)\subseteq D(\lambda_1)\cup\cdots\cup D
	(\lambda_m)$, for some $\lambda_1,\ldots,\lambda_m\in\olh^*$
\end{enumerate}
where $D(\lambda)=\{\mu\in\olh^*\;|\;\mu\leqslant\lambda\}$, with
$\mu\leqslant\lambda$ iff $\lambda-\mu\in\sfQ_+$.
The category $\O_{\olg}$ has a natural symmetric tensor structure
inherited from $\Rep\olg$.

We observed in \ref{ss:bil on g} that the restricted Drinfeld
double of the negative Borel subalgebra $\olb_-$ of $\olg$ is isomorphic
to the trivial central extension $\olgtwo=\olg\oplus\zolh$ of $\olg$ by $\zolh
=\olh$. It follows by \ref{ss:pre-DY}--\ref{ss:drinf-db-rep} that the category
of \DYt modules over $\olb_-$ is equivalent to the category $\E_{\olgtwo}
$ of $\olgtwo$--modules, where $\olgtwo=\olg\oplus\zolh$, which carry a
locally finite action of $\olbtwo_+\subset\olgtwo$. This implies the following.

\begin{proposition}
\hfill
\begin{enumerate}
\item Category $\O_{\olg}$ is isomorphic to the full tensor subcategory
of $\E_{\olgtwo}$ consisting of those modules carrying a trivial action
of $\zolh$ and satisfying, as a module over $\olh\subset\olg\subset
\olgtwo$, the conditions $(\O1)$--$(\O3)$ above. 
\item Under the equivalence $\E_{\olgtwo}\simeq\DrY{\olb_-}$, $\O_
{\olg}$ is isomorphic to the full tensor subcategory of $\DrY{\olb_-}$
consisting of those modules $V$ such that the action $\rho_V$ and 
coaction $\rho^*_V$ of $\olh$ on $V$\footnote{The (co)action of $\olh$
is defined by restricting that of $\olb_-$ as in \ref{ss:diag-DYLBA},
since the inclusion $i_0:\olh\to\olb_-$ is a split embedding with left
inverse $p_0$: $\rho_V=\pi_V\circ i_0\otimes\id_V$, $\rho^*_V=p_
0\otimes\id_V\circ\pi^*_V$.}
coincide under $\iip{\cdot}{\cdot}_{\olh}$, \ie
\begin{equation}\label{eq:DY-ss-0}
\rho_V=
\iip{\cdot}{\cdot}_{\olh}\ten\id_V\circ\id_{\olh}\ten\rho_V^*
\end{equation}
as maps  $\olh\otimes V\to V$ and, as a module over $\olh\subset\olb_-$,
$V$ satisfies the conditions $(\O1)$--$(\O3)$ above.
\end{enumerate}
\end{proposition}

\subsection{Pre--Coxeter structures and category $\O_\infty$}\label{ss:cat-O-infty}
%-------------------------------------------------------------------------------

Condition ($\O2$)  on the finite--dimensionality of weight spaces in
\ref{ss:cat-O-lba} is not stable under restriction from $\olg=\olg_D$
to $\olg_B$ if $B\subsetneq D$, which makes category $\O_{\olg}$
unsuitable to the axiomatic framework of braided pre--Coxeter structures. We
therefore omit it, and denote  by $\O_{\infty,\olg}$ the category of
$\olg$--modules satisfying conditions ($\O1$) and ($\O3$). Proposition
\ref{ss:cat-O-lba} shows that $\O_{\infty,\olg}$ is a full subcategory
of $\DrY{\olb_-}$. Moreover, the universal braided pre--Coxeter
structure on $\{\defDY{\olb_{B,-}}\}_{B\subseteq D}$ restricts to
one on $\{\O_{\infty,\olg_B}^{\hbar}\}_{B\subseteq D}$. 

\subsection{Braid group actions}\label{ss:qC-LBA}
%------------------------------------------------------

Assume now that $\sfA$ is a symmetrisable generalised Cartan matrix,
let $W$ be the corresponding Weyl group with set of simple reflections
$\{s_i\}_{i\in\bfI}$, and set $\ulm=(m_{ij})$, where $m_{ij}$ is the order
of $s_is_j$ in $W$. 

Let $\M_{\olg}^{\sint}$ be the category of integrable $\olg$--modules,
\ie $\olh$--semisimple
modules endowed with a locally nilpotent action of the elements $\{e_i,f_i\}_{i\in\bfI}$.
For any $i\in D$, let $\wt{s}_i\in\sfEnd{\M_{\olg}^{\sint}\to\vect}$ be the triple exponential
\[\wt{s}_i=\exp(e_i)\cdot\exp(-f_i)\cdot\exp(e_i)\]
It is well--known (cf. \cite{tits}) that these satisfy the generalised braid relations
\eqref{eq:gen-braid}.

Let $\hDrY{\olb_-}{\sint }$ be the category of integrable Drinfeld--Yetter $\olb_-
$--modules in $\DrY{\olb_-}$, \ie $\olh$--diagonalisable, endowed with a locally
nilpotent action of the elements $\{f_i\}_{i\in D}\subseteq\olb_-$, and satisfying
\eqref{eq:DY-ss-0}, so as to give rise to integrable modules over $\olg$.
In particular, the triple exponential $\wt{s}_i$ acts on the objects in $\hDrY{\olb_-}{\sint }$ and
the subcategory of integrable modules in $\O_{\infty,\olg}$, denoted $\O^{\sint}_{\infty,\olg}$, 
is isomorphic to a braided tensor subcategory of $\hDrY{\olb_-}{\sint }$.
The following is straightforward.

\begin{proposition}
There is a canonical $(\oalpho{}{}{},\DCPA{}{})$--strict symmetric Coxeter category 
$\Cox{\olb_-}{\sint }$ of type $(D,\ulm)$, defined as follows
\begin{itemize}
\item For any $B\subseteq D$, $\Cox{\olb_-,B}{\sint }$ is the symmetric
monoidal category $\hDrY{\olb_{B,-}}{\sint }$
\item For any $B'\subseteq B$, the functor $F_{B'B}:\Cox{\olb_-,B}{\sint }\to\Cox{\olb_-,B'}{\sint }$
is the restriction $\Res_{\bb{B'},\bb{B}}:\hDrY{\olb_{B,-}}{\sint }\to
\hDrY{\olb_{B',-}}{\sint }$
\item for any $i\in D$, $S_i=\wt{s}_i$
\end{itemize}
There is a natural symmetric Coxeter category $\OCox{\infty, \olg}{\sint}$
obtained from $\Cox{\olb_-}{\sint }$ by restriction to the subcategories
$\O^{\sint}_{\infty,\olg_B}$, $B\subseteq D$. 
\end{proposition}
\begin{pf}
It is enough to observe that $\wt{s}_i$ is group--like and therefore satisfies the
coproduct identity \eqref{eq:coprod-id}, which for the symmetric category $\hDrY{\olb_{i}}{\sint }$ reduces
precisely to the condition $\Delta(\wt{s}_i)=\wt{s}_i\ten\wt{s}_i$.
\end{pf}

%---------------------------------------------------------------------------------------------
\subsection{Universal braided Coxeter structures on Kac--Moody algebras}\label{ss:univ-cox-str}
%---------------------------------------------------------------------------------------------

Let $\intDY{\olb_{B,-}}$ be the category of integrable deformation \DYt
$\olb_{B,-}$--modules. Recall that $\hDYA{B}{n}$ and $\hDYA{B,0}{n}$
denote the algebras of endomorphisms of the forgetful functors $\ff_B^
{\boxtimes n}:(\defDY{\olb_{B,-}})^{n}\to{\khvect}$ and $\ff_{B,0}^{\boxtimes
n}:(\defDY{\olb_{B,-}})^{ n}\to\defDY{\olh}{}$, respectively. For any $X\in
\hDYA{B}{n}$, we denote by $\sfp(X)$ the induced endomorphism of the
forgetful functor $(\intDY{\olb_{B,-}})^{n}\to{\khvect}$.

\begin{definition}
	A braided Coxeter structure of type $(D,\ulm)$ with diagrammatic
	categories $\{\intDY{\olb_{B,-}}\}_{B\subseteq D}$ is \emph{universal}
	if the underlying braided pre--Coxeter structure is (cf. \ref{pr:univ pre Cox}),
	and its local monodromies have the form
	\begin{equation*}\label{eq:locmon}
	S_i=\wt{s}_i\cdot\sfp(\ul{S}_i)
	\end{equation*}
	where $\ul{S}_i\in\hDYA{\vrtx{i},0}{1}$, $\ul{S}_i=1\mod\hbar$, and
	$\wt{s}_i=\exp(e_i)\exp(-f_i)\exp(e_i)$.
\end{definition}

\noindent\remark 
Since $\defDY{\ekm{\b}_i}{}\simeq\Rep\hext{U\olgtwo_i}$ with $\g_i=\mathfrak
{sl}_2^{\alpha_i}$, we have $\hDYA{\vrtx{i}}{n}=\hext{(U\olgtwo_i)^{\ten n}}$. 
In particular, $\sfp(\ul{S}_i)$ is an element in $\hext{(U\olg_i)^{\olh_i}}$.

%%%%%%%%%%%%%%%%%%%%%%%%%%%%%%%%%%%%%%%%
%	QUANTUM KAC-MOODY ALGEBRAS AND TRANSFER THEOREM		%
%%%%%%%%%%%%%%%%%%%%%%%%%%%%%%%%%%%%%%%%

\section{Quantum Kac--Moody algebras}\label{s:qkm}
%---------------------------------------------------------------------

We show in this section that integrable, category $\O_{\infty}$ representations 
of a quantised extended \KMA $U_\hbar\olg$ give rise to a braided Coxeter
category, with local monodromies given by Lusztig's quantum Weyl group
operators. Using the fact that $U_\hbar\olg$ is isomorphic to the \nEK
quantisation of $\olg$ \cite{ek-6}, together with % as a diagrammatic QUE,
the results of Section \ref{s:qcstructure}, we then transport this structure
to integrable, category $\O_{\infty}$ representations of $\olg$.

%%%%%%%%% DRINFELD-JIMBO

\subsection{The extended Drinfeld--Jimbo quantum group}\label{ss:DJ-qg}
%---------------------------------------------------------------------------

Throughout this section, $\GCM{A}=\{a_{ij}\}_{i,j\in\bfI}$ denotes a fixed,
symmetrisable generalised Cartan matrix, \ie $a_{ii}=2$, $a_{ij}\in\IZ_
{\leqslant0}$ if $i\neq j$, and there is a non--singular diagonal matrix
$\GCM{D}$ such that $\GCM{B}=\GCM{DA}$ is symmetric (in particular,
$a_{ij}=0$ if and only if $a_{ji}=0$). The matrix $\GCM{D}$ is determined
uniquely by requiring that $d_i\in\IZ_+$ and $\mathsf{gcd}\{d_i\}=1$.

Let $\olg=\olg(\sfA)$ be the corresponding extended Kac--Moody algebra
with the standard diagrammatic Lie bialgebra structure described in \ref
{ss:split-diag-olb}, and set $q_i=\exp(\hbar/2\cdot d_i)$, $i\in\bfI$. The
following is a straightforward generalisation to extended \KMAs of the
Drinfeld--Jimbo quantum group $\DJ{\g}$ \cite[Example 6.2]{drin-2},\cite{j}.

\begin{definition}
The Drinfeld--Jimbo quantum group of $\olg$ is the unital associative
algebra $\DJ{\olg}$ over $\hext{\sfk}$ topologically generated by $\olh$
and the elements $\{E_i, F_i\}_{ i\in\bfI}$, with relations 
\begin{align*}
	[h,h']&=0&
	[h, E_i]&=\alpha_i(h)E_i&
	[h, F_i]&=-\alpha_i(h)F_i\\
	&&
	[E_i, F_i]&=\frac{q_i^{h_i}-q_i^{-h_i}}{q_i-q_i^{-1}}&
	&
\end{align*}
for any $h,h'\in\olh$, $i\in\bfI$, where $h_i=\cor{i}$, and
\[\sum_{m=0}^{1-a_{ij}}(-1)^mX_i^{(1-a_{ij}-m)}X_jX_i^{(m)}=0\]
for $X=E,F$, $i\neq j\in\bfI$, where $X_i^{(r)}=X_i^r/[r]_{q_i}!$

$\DJ{\olg}$ is a Hopf algebra, with counit $\varepsilon(h)=\varepsilon(E_i)=
\varepsilon(F_i)=0$, coproduct
\[\Delta(h)=
h\ten1+1\ten h\qquad
\Delta(E_i)=E_i\ten q_i^{h_i}+1\ten E_i\qquad
\Delta(F_i)=F_i\ten 1+q_i^{-h_i}\ten F_i\]
and antipode $S(h)=-h$, $S(E_i)=-E_iq_i^{-h_i}$, and $S(F_i)=-q_i^{h_i}F_i$, for any
$h\in\olh$ and $i\in\bfI$.
\end{definition}

The following result is well--known for $\DJ{\g}$ (cf. \cite[Sec. 13]{drin-2} and
\cite[Sec. 8.3]{CP}). It readily extends to $\DJ{\olg}$ through the isomorphism
of Hopf algebras $\DJ{\olg}\simeq\DJ{\g}\ten\DJ{\c}$, where $\DJ{\c}=\hext{S\c}$,
which quantises the decomposition $\olg\simeq\g\oplus\c$ (cf.\;\ref{ss:h''}).

\begin{proposition}\cite{drin-2,CP}\hfill
	\begin{enumerate}
		\item
		The Hopf algebra $\DJ{\olg}$ is a quantisation of the Lie bialgebra $\olg$. 
		\item 
		Let $\DJ{\ekm{\b}_\mp}\subset\DJ{\olg}$ be the Hopf subalgebra topologically
		generated by $\h$ and $\{F_i\}_{i\in\bfI}$ (resp. $\h$ and $\{E_i\}_{i\in\bfI}$). 
		Then, $\DJ{\ekm{\b}_\mp}$ is a quantisation of the Lie bialgebra 
		$\ekm{\b}_\mp$, and there is a unique non--degenerate Hopf pairing 
		$\iip{\cdot}{\cdot}_{\D}:\DJ{\ekm{\b}_-}\ten\DJ{\ekm{\b}_+}\to\sfk(\negthinspace(\hbar)\negthinspace)$, 
		defined on the generators by
		\[
		\begin{array}{ccccccc}
		\iip{1}{1}_{\D} = 1 &&&
		\iip{h}{h'}_{\D} =\displaystyle{\frac{1}{\hbar}\iip{h}{h'}} &&&
		\iip{F_i}{E_j}_{\D}=\displaystyle{\frac{\delta_{ij}}{q-q^{-1}}}
		\end{array}
		\]
		and zero otherwise. 
		
		\item The Hopf pairing $\iip{\cdot}{\cdot}_{\D}$ induces an isomorphism of
		finitely $\IN$--graded QUEs $\DJ{\ekm{\b}_-}\simeq(\DJ{\ekm{\b}_+})^{\star}$,
		where the latter is the restricted QUE dual (cf. \ref{ss:dualityQUE}). This
		gives rise to an isomorphism of QUE $\DJ{\olg}\simeq\resqD{\DJ{\ekm{\b}_-}}/
		(\h\simeq\h^*)$. In particular, $\DJ{\olg}$ is a quasitriangular Hopf algebra, with 
		$R$--matrix
		\begin{equation}\label{eq:R-mx}
		\ol{R}=q^{\sum_{i}u_i\ten u^i}\cdot\sum_pX_p\ten X^p,
		\end{equation}
		where $\{u_i\},\{u^i\}\subset\olh$ are dual bases \wrt $\iip{\cdot}{\cdot}$, and
		$\{X_p\}\subset\DJ{\ekm{\n}_-}$, $\{X^p\}\subset\DJ{\ekm{\n}_+}$ are dual
		bases \wrt $\iip{\cdot}{\cdot}_{\D}$.
	\end{enumerate}
\end{proposition}

%%%%%%%%%%%%%%%%%%%%%%%%%%%%%%%%%%%%

\subsection{Diagrammatic structures on $\DJ{\olg}$}\label{ss:diag-on-qG}
%------------------------------------------------------------------

The quantum group $\DJ{\olg}$ is canonically endowed with the structure
of diagrammatic Hopf algebra, with subalgebras $\DJ{\olg_B}=\langle
\hcor{i},\fcw_i, E_i, F_i\rangle_{i\in B}$, $B\subseteq D$. 

As in the classical case (cf. \ref{ss:split-diag-olb}), the diagrammatic structure 
of $\DJ{\olg}$ induces a split diagrammatic one on $\DJ{\olb_\pm}$.
Namely, for any $B\subseteq D$, let $\DJ{\olb}_{\pm,B}=\DJ{\olb}_\pm\cap\DJ{\olg_B}$ 
be the Hopf sualgebras topologically generated by $\{\hcor{i},\fcw_i,e_i\}_{i\in B}$ and 
$\{\hcor{i},\fcw_i,f_i\}_{i\in B}$ respectively. For $B'\subseteq B$, 
let $i_{\pm,BB',\hbar}:\DJ{\olb}_{\pm,B'}\to\DJ{\olb}_{\pm,B}$ be the standard 
embedding, and regard $p_{\pm,B'B,\hbar}=i_{\mp,BB',\hbar}^t$ as a map 
$\DJ{\olb}_{\pm,B}\to\DJ{\olb}_ {\pm,B'}$ via the identifications
$\DJ{\olb}_{\mp,C}^\star\cong\DJ{\olb}_{\pm,C}$ given by the inner 
product $\iip{\cdot}{\cdot}_{\D}$. The map $(D\DJ{\olb}_-)\res\to\DJ{\olg}$ 
from Proposition \ref{ss:DJ-qg} (3) is then a morphism of diagrammatic
Hopf algebras.

%--------------------------------
% QCQTQBA ON Uhg
%--------------------------------

\subsection{Coxeter structures on quantum groups}\label{ss:qC-on-qG}
%------------------------------------------------------------------

Let $W$ be the Weyl group of $\olg$, $\{s_i\}_{i\in\bfI}$ its generators,
and set $\ulm=(m_{ij})$, where $m_{ij}$ is the order of $s_is_j$ in $W$.
Thus, for any $B\subseteq D$, the generalised braid group $\BBm$ is
the Tits braid group of the standard parabolic subgroup of $W$
generated by $\{s_i\}_{i\in B}$. 

Let $\hDrY{\DJ{\ekm{\b}_{B,-}}}{\adm}$ be the braided monoidal category 
of admissible Drinfeld--Yetter $\DJ{\ekm{\b}_{B,-}}$--modules. As in \ref
{ss:cat-O-lba}, denote by $\hDrY{\DJ{\ekm{\b}_{B,-}}}{\adm,\sint }$ the 
full subcategory of $\olh_B$--diagonalisable, integrable Drinfeld--Yetter
$\DJ{\ekm{\b}_{B,-}}$--modules $\V$ such that the action and coaction
of $\olh$ on $\V$ coincide under $\iip{\cdot}{\cdot}_{\olh}$, that is satisfy
\begin{equation*}
\rho_{\V}=
\iip{\cdot}{\cdot}_{\olh}\ten\id_{\V}\circ\id_{\olh}\ten\rho_{\V}^*
\end{equation*}
so as to give rise to integrable modules over $\DJ{\ekm{\g}_B}$. 
\begin{proposition}
There is a canonical $(\oalpho{}{}{},\DCPA{}{})$--strict braided Coxeter
category $\Cox{\DJ{\olb_-}}{\adm,\sint }$ of type $(D, \ulm)$, with
\begin{itemize}
\item diagrammatic categories $\hDrY{\DJ\ekm{\b}_{B,-}}{\adm,\sint }$, $B\subseteq D$
\item standard restriction functors $\hDrY{\DJ\ekm{\b}_{B,-}}{\adm,\sint}\to\hDrY{\DJ\ekm
{\b}_{B',-}}{\adm,\sint}$ determined by the split diagrammatic structure of $\DJ{\olb_-}$
\item local monodromies given by Lusztig's quantum Weyl group operators $S_i^{\hbar}$
\end{itemize}
\end{proposition}
\begin{pf}
	The $(\oalpho{}{}{},\DCPA{}{})$--strict braided pre--Coxeter structure on $\Cox{\DJ{\olb_-}}
	{\adm,\sint }$ is defined in \ref{ss:DYdiag}. For the Coxeter structure, we proceed as in
	\ref{ss:qC-LBA}. Denote by $\M^{\sint}_{\DJ\olg}$ the category
	of integrable $\DJ\olg$--modules. Following \cite{l}, the quantum Weyl group operator 
	of $\DJ{\olg}$ corresponding to $i\in\bfI$ is the element $S_i^{\hbar}\in\sfEnd{\M_{\DJ\olg}^{\sint}\to\vect_{\sfK}}$ acting
	on $\V\in\M_{\DJ\olg}^{\sint}$ as
	\[
	S_i^{\hbar}(v)=\sum_{\substack{a,b,c\in\IZ\\ a-b+c=-\lambda(\alpha^{\vee}_i)}}
	(-1)^bq_i^{\frac{h_i^2}{4}+b-ac}E_i^{(a)}F_i^{(b)}E_i^{(c)}v
	\]
	where $v\in\V_{\lambda}$ for $\lambda\in\h^*$. The quantum
	Weyl group operators $S^{\hbar}_i$ satisfy the braid relations \eqref{eq:gen-braid}, together with the
	coproduct identity
	\[
	\Delta(S_i^{\hbar})=R_i^{21}\cdot(S_i^{\hbar}\ten S_i^{\hbar})
	\]
	Each $S^{\hbar}_i$, acts on any $\V_i\in\hDrY{\DJ{\ekm{\b}_{\{i\},-}}}{\adm,\sint }$ 
	and they complete the $(\oalpho{}{}{},\DCPA{}{})$--strict braided Coxeter structure on $\Cox{\DJ{\olb_-}}{\adm,\sint }$.
\end{pf}

%%%%%%%%%%%%%%%%%%%%%%%%%%%%%%%%%%%%%%

\subsection{\nEK quantisation} % of Kac--Moody algebras}
\label{ss:qKM}
%-----------------------------------------------------------------------

Let $\Q(\olg)$ (resp. $\Q(\ekm{\b}_{\pm})$) be the \nEK quantisation
of the extended \KMA $\olg$ (resp. the Borel subalgebras $\ekm{\b}
_{\pm}\subset\olg$).
\begin{proposition}\hfill
\begin{enumerate}
\item $\Q(\olg)$ is a diagrammatic QUE, with subalgebras %$\Q(\olg)_B=
$\Q(\olg_B)$, $B\subseteq D$.
\item $\Q(\ekm{\b}_{\pm})$ is a split diagrammatic QUE, with
subalgebras %$\Q(\ekm{\b}_{\pm})_B=
$\Q(\ekm{\b}_{B,\pm})$. 
\item The quantised embeddings $\Q(\ekm{\b}_{B,-})\to\Q(\olg_B)$,
$B\subseteq D$, give rise to a morphism of diagrammatic QUEs 
$\Q(\ekm{\b}_{-})\to\Q(\olg)$.
\item
The following data defines an $(\oalpho{}{}{},\DCPA{}{})$--strict braided
pre--Coxeter category $\Cox{\Q(\olb_-)}{\adm,\sint}$
\begin{itemize}
	\item For any $B\subseteq D$, $\Cox{\Q(\ekm{\b}_-), B}{\adm,\sint}$ is the braided monoidal category $\hDrY{\Q(\b_{B,-})}{\adm,\sint}$
	\item For any $B'\subseteq B$, the functor
	$F_{B'B}:\Cox{\Q(\ekm{\b}_-), B}{\adm,\sint}\to\Cox{\Q(\ekm{\b}_-), B'}{\adm,\sint}$
	is the restriction functor $\Res_{\Q(\ekm{\b}_{B',-}),\Q(\ekm{\b}_{B,-})}$.
\end{itemize}
\item The braided pre--Coxeter category $\Cox{\Q(\olb_-)}{\adm,\sint}$ is a 
deformation of $\Cox{\olb_-}{\sint}$.
\end{enumerate}
\end{proposition}

\begin{pf}
(1) and (2) follow from the compatibility of the quantisation functor
with the diagrammatic and split diagrammatic structures of $\olg$ 
and $\olb_{\pm}$, respectively (Corollary \ref{ss:ek-quantisation}).
(3) follows from the functoriality of $\Q$, and the canonical morphism
of diagrammatic Lie bialgebras $\olb_{-}\to\olg$ (Proposition ~\ref
{ss:split-diag-olb} (4)). (4) is given by Corollary \ref{ss:quantum-DY}.
(5) is clear.
\end{pf}

%%%%%%%%%%%%%%%%%%%%%%%%%%%%%%%%%%%%%%%%

\subsection{Quantum double construction of $\Q(\olg)$}\label{ss:ek-qKM-qD}
%-----------------------------------------------------------------------

By \cite{EG}, the \nEK quantisation functor $\Q$ is compatible with taking
duals and doubles. This is used in \cite{ek-6} to show that $\Q(\g)$ is a
quotient of the quantum double of $\Q(\b_-)$, and that it is isomorphic
to the quantum group $\DJ{\g}$. The argument is easily adapted to the
extended \KMA $\olg$, since the latter is a %non--canonical 
central extension of the former $\g$ (cf. \ref{ss:h''}).
Specifically, by Proposition \ref{ss:bil on g}, $\olg$ is isomorphic to the quotient of 
the Drinfeld double of $\olb_-$ by the ideal generated by the identification 
of $\phi:\olh\to\olh^*$, \ie $\olg\simeq\resDb{\olb_-}/(\olh\simeq\olh^*)$. Since 
$\Q$ is compatible with doubling operations, there is an isomorphism $\Q(\resDb{\olb_-})
\simeq\resqD{\Q(\olb_-)}$, which is the identity on $\olh\oplus\olh^*$. This
yields an isomorphism of Hopf algebras
\[
\slantfrac{\resqD{\Q(\olb_-)}}{\olh\simeq\olh^*}\simeq\Q(\olg).
\]
which shows, in particular, that $\Q(\olg)$ is quasitriangular. Finally, one
proves the following

\begin{theorem}\cite{ek-6}\label{thm:ek-iso}
	\begin{enumerate}
	\item There is a (non--canonical) isomorphism of 
	QUEs $\varphi^{\olb_-}:\DJ{\olb_-}\to\Q(\olb_-)$, which is
	the identity on $\olh$.
	\item By the quantum double construction of $\Q(\olg)$ and $\DJ{\olg}$
	(cf. Proposition \ref{ss:DJ-qg} (3)), $\varphi^{\olb_-}$ induces an isomorphism 
	of quasitriangular QUEs $\varphi^{\olg}:\DJ{\olg}\to\Q(\olg)$.
	\end{enumerate}
\end{theorem}

%--------------------------------------------------------------------------------------------------
\subsection{Diagrammatic isomorphism between $\Q(\olg)$ and $\DJ{\olg}$}\label{ss:iso-ek-dj}
%--------------------------------------------------------------------------------------------------

We now show that the isomorphism between $\Q(\olg)$ and $\DJ{\olg}$ 
can be chosen so as to preserve the diagrammatic structures.

\begin{proposition}
	\hfill
	\begin{enumerate}
		\item
		There is an isomorphism of split diagrammatic QUEs
		$\psi^{\olb_-}:\DJ{{\olb_-}}\to\Q({\olb_-})$, which is the
		identity on $\olh$.
		\item 
		By the quantum double construction, $\psi^{\olb_-}$ induces
		an isomorphism of diagrammatic QUEs $\psi^{\olg}:\DJ{{\olg}}\to\Q({\olg})$.
	\end{enumerate}
\end{proposition}

\begin{pf}
	(1) For any $j\in D$, use Theorem \ref{ss:ek-qKM-qD} (1) to choose
	an isomorphism of Hopf algebras $\psi^{\olb_-}_j\coloneqq\varphi^
	{\olb_{\{j\},-}}:\DJ{\ekm{\b}_{\{j\},-}}\to\Q(\ekm{\b}_{\{j\},-})$. Then,
	for any $B\subseteq D$, we get an isomorphism of Hopf algebras
	$\psi^{\olb_-}_B:\DJ{\ekm{\b}_B}\to\Q(\ekm{\b}_B)$ by 
	\[
	\psi_B^{\olb_-}(F_j)\coloneqq\Q(i^{\olb_-}_j)\circ\psi^{\olb_-}_j(F_j)
	\]
	where $j\in B$. The collection $\psi^{\olb_-}=\{\psi_B^{\olb_-}\}_{B\subseteq D}$ 
	gives an isomorphism of split diagrammatic Hopf algebras. (2) is clear.
\end{pf}

%%%%%%%%%%%%%%%%%%%%%%%%%%%%%%%%%%%%

%--------------------
% FINAL SECTION CONCLUSION
%--------------------
%---------------------------------------------------------------------------
\subsection{An equivalence of braided Coxeter categories}\label{ss:main-thm-3}
%---------------------------------------------------------------------------

We now prove the main result of this paper. We show that the Coxeter
structure on integrable \DYt modules for $\DJ{\olb}$, which accounts for
the \qW operators of $\DJ{\olg}$, can be transferred to a Coxeter structure
on integrable \DYt modules for $\olb$, with standard restriction functors.

\begin{theorem}
	Let $\Phi\in\CRDYUA{\mathsf{LBA}}{3}$ be a factorisable associator.
	\begin{enumerate}[leftmargin=3em]
	\item
	There is an equivalence of braided pre--Coxeter categories
	\[\Hcox{\olb_-}:\Cox{\olb_-}{\Phi,\astr,\sint}\longrightarrow\Cox{\DJ{\olb_-}}{\adm,\sint}\]
	where
	\begin{enumerate}\itemsep0.25cm
% DJ side
		\item $\Cox{\DJ{\olb_-}}{\adm,\sint}$ is the $(\sfa,\DCPA{}{})$--strict
		structure defined in \ref{ss:qC-on-qG}, with
		\begin{itemize}
		\item diagrammatic categories $\hDrY{\DJ{\olb_{B,-}}}{\adm,\sint}$
		\item standard monoidal restriction functors
		\[\Res_{\DJ{\olb_{B',-}}, \DJ{\olb_{B,-}}}:\hDrY{\DJ{\olb_{B,-}}}{\adm,\sint}\to\hDrY{\DJ{\olb_{B',-}}}{\adm,\sint}\]
		\end{itemize}
% classical side
		\item $\Cox{\olb_-}{\Phi,\astr,\sint}$ is $\sfa$--strict and universal (see \ref{de:univ pre Cox}), with
		\begin{itemize}
			\item diagrammatic categories $\hDrY{\olb_{B,-}}{\Phi_B,\sint}$ 
			\item restriction functors of the form
			% $\hDrY{\olb_{B,-}}{\Phi_B,\sint}\to\hDrY{\olb_{B',-}}{\Phi_B,\sint}$
			$(\Res_{\olb_{B',-}\olb_{B,-}}, J_\F):\hDrY{\olb_{B,-}}{\Phi_B,\sint}\to\hDrY{\olb_{B',-}}{\Phi_B,\sint}$,
			for some monoidal structure $J_\F$
		\end{itemize}
% form of the equivalence
		\item the equivalence $\Hcox{\olb_-}$ is given the composition 
		\begin{equation}\label{eq:equiv-comp}
		\Cox{\ekm{\b}_-}{\Phi,\astr,\sint}\to\Cox{\ekm{\b}_-}{\Phi,\gstr,\sint}\to\Cox{\Q(\ekm{\b}_-)}{\adm,\sint }\to\Cox{\DJ\ekm{\b}_-}{\adm,\sint }
		\end{equation}
		where the first equivalence is given by Corollary \ref{ss:main-thm-1},
		the second one by the transfer Theorem \ref{ss:main-thm-2}, and the third one by the isomorphism 
		of diagrammatic QUEs $\Q(\ekm{\b}_-)\simeq\DJ\ekm{\b}_-$ (Prop. \ref{ss:iso-ek-dj}).
	\end{enumerate}
	\item There is a unique braided Coxeter structure on $\Cox{\olb_-}{\Phi,\astr,\sint}$,
	which extends the pre--Coxeter structure, and is such that 
		\[\Hcox{\olb_-}:\Cox{\olb_-}{\Phi,\astr,\sint}\to\Cox{\DJ{\olb_-}}{\adm,\sint}\]
	is an equivalence of Coxeter categories with respect to the Coxeter structure on 
	$\Cox{\DJ{\olb_-}}{\adm,\sint}$ arising from the quantum Weyl group operators of $\DJ{\olg}$ (cf. \ref{ss:qC-on-qG}).
	Moreover, the braided Coxeter structure on $\Cox{\olb_-}{\Phi,\astr,\sint}$ is universal in the sense of Definition 
	\ref{ss:univ-cox-str}.
% equivalence of braid group representations
\item The representations of the generalised braid groups $\BBm$, $B
\subseteq D$, arising from the quantum Weyl group operators of
$\DJ{\olg}$ and the Coxeter category $\Cox{\olb_-}{\Phi,\astr,\sint}$ 
are equivalent. 
	\end{enumerate}
\end{theorem}

\begin{pf}
(1) Let $\Cox{\olb_-}{\Phi,\astr}$ be the universal $\sfa$--strict braided
pre--Coxeter category associated to $\olb_-$ by Corollary \ref{ss:main-thm-1}. 
By Theorem \ref{ss:main-thm-2}, there is
an equivalence of braided pre--Coxeter categories 
$\Cox{\olb_-}{\Phi,\astr}\to\Cox{\Q(\olb_-)}{\adm}$.
The isomorphism of split diagrammatic Hopf algebras 
$\Q(\olb_-)\simeq\DJ{\olb_-}$ constructed in Proposition \ref{ss:iso-ek-dj}, 
then allows to extend it to an equivalence 
\[\wt{\mathbb{H}}_{\olb_-}:\Cox{\olb_-}{\Phi,\astr}\to
\Cox{\DJ{\olb_-}}{\adm}\]
Let $\Cox{\olb_-}{\Phi,\astr,\sint}$ be the braided pre--Coxeter subcategory of 
$\Cox{\olb_-}{\Phi,\astr}$, with underlying diagrammatic categories 
$\hDrY{\ekm{\b}_{B,-}}{\Phi_B,\sint}$, $B\subseteq D$.
Since the Etingof--Kazhdan functors preserve integrable modules
\cite[Prop. 6.5]{A}, the restriction of $\wt{\mathbb{H}}_{\olb_-}$ give rise to an equivalence
of braided pre--Coxeter categories $\Hcox{\olb_-}:\Cox{\olb_-}{\Phi,\astr,\sint}\to
\Cox{\DJ{\olb_-}}{\adm,\sint}$.

(2) By \ref{ss:qC-on-qG}, the \qW operators of $\DJ{\olg}$ define a Coxeter
structure on $\Cox{\DJ{\olb_-}}{\adm,\sint}$. The requirement that $\Hcox
{\olb_-}$ be an equivalence of braided Coxeter  categories therefore uniquely
determines a Coxeter structure on $\Cox{\olb_-}{\Phi,\astr,\sint}$. Namely,
let $\Psi_i$ denote the pullback on the algebra of endomorphisms of the
forgetful functor along
\[
H_{i}:\hDrY{\ekm{\b}_{\{i\},-}}{\Phi_{\{i\}},\sint}\to\hDrY{\Q(\ekm{\b}_{\{i\},-})}{\adm,\sint }\to\hDrY{\DJ\ekm{\b}_{\{i\},-}}{\adm,\sint }
\]
Then, the operators $\Psi_i(S_i^{\hbar})$ extend the braided pre--Coxeter
structure of $\Cox{\olb_-}{\Phi,\astr,\sint}$ to a braided Coxeter  structure.
It is then clear by \ref{ss:qC-on-qG} that $\Psi_i(S_i^{\hbar})$ satisfies the
conditions of Definition \ref{ss:univ-cox-str}, and therefore that this structure
is universal.

(3) By construction, the action of the generalised braid groups $\BBm$
on $\V\in\hDrY{\DJ{\ekm{\b}_{B,-}}}{\adm,\sint }$ arising from the 
Coxeter category $\Cox{\DJ{\olb_-}}{\adm,\sint }$ coincides with 
the action of the quantum Weyl group operators of $\DJ{\olg_B}$
(cf. \ref{ss:qC-on-qG}). The result then follows from (2) and Proposition
\ref{ss:braid-gp-act} (2).
\end{pf}

\subsection{Coxeter structures and category $\O_{\infty}$}\label{ss:Cox-catO-infty}
%---------------------------------------------------------------------------

Fix $B\subseteq D$. Recall that a $\DJ\olg_B$--module $\V$ is in category $\O^{\int}_{\DJ
\olg_B}$ if it is topologically free over $\hext{\sfk}$, integrable, and
satisfies the conditions $(\O1)-(\O3)$ of \ref{ss:cat-O-lba}. % as an $\olh$--module.
Let $\O^{\sint}_{\infty,\DJ\olg_B}$ be the category of $\DJ\olg_B$--modules
satisfying conditions $(\O1)$ and $(\O3)$, but not necessarily
the finite--dimensionality of weight spaces. The realisation
of $\DJ\olg_B$ as a quotient of the quantum double
of $\DJ\ekm{\b}_{B,-}$ (\ref{ss:DJ-qg}) gives rise to a full embedding
$\O^{\sint}_{\infty,\DJ\olg_B}\subset\hDrY{\DJ\ekm{\b}_{B,-}}{\adm,\sint }$. 

Since the \nEK functor $\hDrY{\olb_{B,-}}{\Phi_B}\to\aDrY{\Q(\olb_{B,-})}$
is the identity on $\olh_B$--modules, the equivalence \eqref{eq:equiv-comp}
preserves the categories $\O^{\Phi_B,\sint}_{\infty,\olg_B}\subset\hDrY{\ekm
{\b}_{B,-}}{\Phi_B,\sint}$ and $\O^{\sint}_{\infty,\DJ\olg_B}\subset\hDrY
{\DJ\ekm{\b}_{B,-}}{\adm,\sint }$. This yields the following 

\begin{theorem}
Let $\Phi\in\CRDYUA{\mathsf{LBA}}{3}$ be factorisable associator.
Then, there is an equivalence of braided Coxeter categories
\[
\Hcox{\olg}:\OCox{\infty,\olg}{\Phi,\sint}\to\OCox{\infty,\DJ{\olg}}{\sint}
\]
where $\OCox{\infty,\olg}{\Phi,\sint}$ (resp. $\OCox{\infty,\DJ{\olg}}{\sint}$)
is the braided Coxeter category obtained from $\Cox{\olb_-}{\Phi,\astr,\sint}$ 
(resp. $\Cox{\DJ{\olb_-}}{\adm,\sint }$) by restriction to integrable, 
category $\O_{\infty}$ representations.
\end{theorem}

\subsection{Coxeter structures, Levi subalgebras and category $\O$}\label{ss:Levi}
%----------------------------------------------------------------------------------------

As mentioned in \ref{ss:diag-KM} and \ref{ss:cat-O-infty}, the reason
behind the introduction of extended \KM algebras and of category $
\O_{\infty}$ is the construction of a diagrammatic structure endowed
with well--defined restriction functors. 

There is, however, a weaker notion of diagrammatic structure which
leads to an analogue of Theorem \ref{ss:Cox-catO-infty} expressed
solely in terms of standard \KMAs and category $\O$ representations.
Indeed, the facts that minimal realisations of \KMAs do not give rise
to a diagrammatic structure (Prop. \ref{ss:suff and counter}), and
that category $\O$ representations are not stable under restriction,
due to the requirement on the finite--dimensionality of weight spaces,
can both be overcome by considering instead the Levi subalgebras
$\ll_B=\langle e_i,f_i,\h\rangle_{i\in B}$ of a given \KMA $\g$.

The collection $\{\ll_B\}$ does not, however, define a diagrammatic
structure on $\g$, since it does not satisfy the orthogonality condition
$[\ll_B,\ll_{B'}]=0$ for $B\perp B'$. As mentioned in \ref{ss:Cox D-alg},
this condition is convenient in the construction of $\PROP$ic structures,
but not required by the axioms of a Coxeter category.
It is in fact possible to adapt the definition of universal pre--Coxeter
structure, and consequently Sections \ref{s:diag-prop}--\ref{s:univ pC},
by removing the orthogonal factorisation axiom in Definition \ref{ss:prop-weak-cox}.
In this new setting, Proposition \ref{prop:strictness-univ-cox} does not
hold, \ie a non--orthogonal structure cannot be $\sfa$--strictified in
general. With this exception, all other results from Section \ref{s:qcstructure}
can be adapted, and applied to the case of  the Levi subalgebras
$\ll_B$.

As observed in \ref{ss:main-thm-3} and \ref{ss:Cox-catO-infty}, for any
$B\subseteq D$, the \nEK equivalence $\hDrY{\b_{B,-}}{\Phi_B}\to\aDrY
{\Q(\b_{B,-})}$ preserves integrable modules, and is the identity on
$\h$--modules. It therefore restricts to an equivalence of braided
monoidal categories
\[H_{\ll_B}:\O^{\Phi_B,\sint}_{\ll_B}\to\O_{\DJ\ll_B}^{\sint}\]
Together with the fact that the universal constructions described
in Section \ref{s:qcstructure} are easily seen to yield twists, associators
and joins which are invariant under $\h$, this yields the following
analogue of Theorem \ref{ss:Cox-catO-infty}.

\begin{theorem}\label{th:latest}
Let $\Phi\in\CRDYUA{\mathsf{LBA}}{3}$ be an associator. Then,
there is an equivalence of braided Coxeter categories
\[
\Hcox{\g}:\OCox{\g}{\Phi,\sint}\to\OCox{\DJ{\g}}{\sint}
\]
where $\OCox{\g}{\Phi,\sint}$ (resp. $\OCox{\DJ{\g}}{\sint}$)
is the braided Coxeter category obtained from $\Cox{\b_-}{\Phi,\gstr,\sint}$ 
(resp. $\Cox{\DJ{\b_-}}{\adm,\sint}$) by restriction to the categories 
$\O_{{\ll}_{B}}^{\Phi,\sint}$ (resp. $\O_{\DJ{\ll}_{B}}^{\sint}$).\footnote
{Note that in order to have an action of the \qW operators $S_i^\hbar$,
which do not commute with the action of $\h$, the diagrammatic
categories $(\OCox{\g}{\Phi,\sint})_\emptyset$ and $(\OCox{\DJ\g}
{\sint})_\emptyset$ have to be taken to be $\vect_{\hext{\sfk}}$,
rather than category $\O$ for $\ll_\emptyset=\h$.
Note also that the example
in Proposition \ref{ss:suff and counter} shows that the minimal
realisation of \KMAs does not lead to a diagrammatic structure,
even if the orthogonality requirement is omitted. It is therefore
not possible in general to formulate an analogue of Theorem 
\ref{th:latest} involving minimal realisations, rather than Levi
subalgebras.}%, and category $\O_\infty$.}

\end{theorem}

%%%%%%%%%%%%%%%%%%%%%%%%%%%%%%%%%%%%%%%
%%%%%%%%%%%%%%%%%%%%%%%%%%%%%%%%%%%%%%%
\appendix %%%%%%%%%%%%%%%%%%%%%%%%%%%%%%%%%%
%%%%%%%%%%%%%%%%%%%%%%%%%%%%%%%%%%%%%%%
%%%%%%%%%%%%%%%%%%%%%%%%%%%%%%%%%%%%%%%

\section{Graphical calculus for Coxeter objects}\label{app-graph-calc}

We describe below the axioms of Coxeter objects in a $2$--category
$\XX$ in terms of graphical calculus.

\subsection{Graphical notation}
In the following, we use the graphical notation of string diagrams to describe relations 
between $2$--morphisms in a $2$--category $\mathfrak{X}$
(\eg \cite{joyal-street,lauda}). We represent objects, $1$--morphisms, and
$2$--morphisms with two dimensional, one dimensional and zero dimensional cells, respectively.
Let $X,Y\in\mathfrak{X}$, $F,G\in\mathfrak{X}^{(1)}(X,Y)$
and $\alpha\in\mathfrak{X}^{(2)}(F,G)$. Then, we represent $\alpha$ as 
\begin{center}
\begin{tikzpicture}[scale=.3]
	\coordinate [label=above:G]  (FA) at (0,2);
  	\coordinate [label=below:F] (FB) at (0,-2);
	\node (L) at (0,0) [circle, draw] {$\alpha$};
	\draw (FA) -- (L);
	\draw (FB) -- (L);
	\draw (-2,0) node {$Y$}; 
	\draw (2,0) node {$X$}; 
	\node at (-6,0) {$\rightsquigarrow$};
	\node at (-12.5,0){
	$\xymatrix@C=1.5cm{Y \ar@{<-}@/_10pt/[r]_{F}="A" \ar@{<-}@/^10pt/[r]^{G}="B" & X \ar@{=>} "A";"B"_{\alpha}}$
	};
\end{tikzpicture}
\end{center}
where the diagram on the \rhs is read from bottom to top, and from right to left.
Similarly, a $2$--morphism $\alpha:F\circ G\to H$ will be represented as follows:
\begin{center}
\begin{tikzpicture}[scale=.3]
	\coordinate [label=below:$F$]  (FA) at (-3,-2);
  	\coordinate [label=below:$G$] (FB) at (3,-2);
	\coordinate [label=above:$H$] (FC) at (0,2);
	\node  (L) at (0,0) [circle, draw] {$\alpha$};
	\draw  (FA) -- (L);
	\draw (FB) -- (L);
	\draw  (FC) -- (L);
	\draw (-2,1) node {$Z$}; 
	\draw (2,1) node {$X$}; 
	\draw (0,-2) node {$Y$}; 	
	\node at (-6,0) {$\rightsquigarrow$};
	\node at (-17.5,0){
	$\xymatrix@R=0.5cm@C=1.5cm{Z \ar@{<-}@/_5pt/[dr]_{F}="A" \ar@{<-}@/^10pt/[rr]^{H}="B" &  & X\\
	&Y \ar@{=>} "B"_{\alpha} \ar@{<-}@/_5pt/[ur]_{G}="A"&}$
	};
\end{tikzpicture}
\end{center}
and more generally we represent $\alpha: F_n\circ\dots\circ F_1\Rightarrow G_m\circ\dots\circ G_1$ as
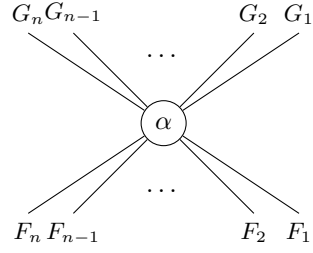
\begin{center}
\begin{tikzpicture}[scale=.3]
	\coordinate [label=below:\small $F_n$]  (FA) at (-6,-4);
  	\coordinate [label=below:\small$F_{n-1}$] (FB) at (-4,-4);
	\coordinate [label=below:\small$F_2$]  (FC) at (4,-4);
	\coordinate [label=below:\small$F_1$]  (FD) at (6,-4);
	\coordinate [label=above:\small$G_m\;\;\;$]  (GA) at (-6,4);
  	\coordinate [label=above:\small$G_{m-1}$] (GB) at (-4,4);
	\coordinate [label=above:\small$G_2$]  (GC) at (4,4);
	\coordinate [label=above:\small$G_1$]  (GD) at (6,4);
	\node  (L) at (0,0) [circle, draw] {$\alpha$};
	\draw  (FA) -- (L);
	\draw  (FB) -- (L);
	\draw  (FC) -- (L);
	\draw  (FD) -- (L);
	\draw  (GA) -- (L);
	\draw (GB) -- (L);
	\draw  (GC) -- (L);
	\draw  (GD) -- (L);
	\node at (0,3) {$\dots$};
	\node at (0,-3) {$\dots$};
\end{tikzpicture}
\end{center}
When no confusion is possible, we omit the labels and identify the $1$--morphisms 
with the color of the string, and the $2$--morphism with the underlying diagram.
\subsection{Coxeter objects (cf.\;\ref{ss:cox})}
A Coxeter object in a $2$--category $\mathfrak{X}$ is the datum of
\begin{itemize}
%objects
\item for any $B\subseteq D$, an object $X_{B}$
%one morphisms
\item for any $\F\in\Mns{B,B'}$, a $1$--morphism $F_{\F}: X_{B}\to X_{B'}$ 
which we represent as the identity $2$--morphisms $\id_{F_{\F}}$
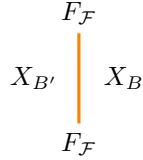
\[
\begin{tikzpicture}[scale=.3]
	\coordinate [label=above:$F_{\F}$]  (FA) at (0,4);
  	\coordinate [label=below:$F_{\F}$] (FB) at (0,0);
	\draw [very thick, orange] (FA) -- (FB);
	\draw (-2,2) node {$X_{B'}$}; 
	\draw (2,2) node {$X_{B}$}; 
\end{tikzpicture}
\]
%two morphisms
\item for any $\F'\in\Mns{B,B'}, \F''\in\Mns{B',B''}$ and $\F=\F'\cup\F''$,
a $2$--morphism
\[
\xymatrix@C=1.5cm{
F_{\F''}\circ F_{\F'} \ar@<3pt>[r]^-{\redasso{\F'}{\F''}} 
& F_{\F} \ar@<2pt>[l]^-{(\redasso{\F'}{\F''})^{-1}}
}
\]
represented as
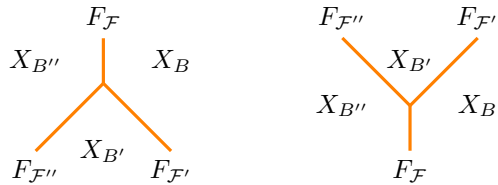
\[
\begin{tikzpicture}[scale=.3]
	\coordinate [label=below:$F_{\F''}$]  (FA) at (-3,-1);
  	\coordinate [label=below:$F_{\F'}$] (FB) at (3,-1);
	\coordinate [label=above:$F_{\F}$] (FC) at (0,4);
	\coordinate  (L) at (0,2);
	\draw [very thick, orange] (FA) -- (L);
	\draw [very thick, orange] (FB) -- (L);
	\draw [very thick, orange] (FC) -- (L);
	\draw (-3,3) node {$X_{B''}$}; 
	\draw (3,3) node {$X_{B}$}; 
	\draw (0,-1) node {$X_{B'}$}; 
\end{tikzpicture}
\qquad\qquad
\begin{tikzpicture}[scale=.3]
	\coordinate [label=above:$F_{\F''}$]  (FA) at (-3,5);
  	\coordinate [label=above:$F_{\F'}$] (FB) at (3,5);
	\coordinate [label=below:$F_{\F}$] (FC) at (0,0);
	\coordinate  (L) at (0,2);
	\draw [very thick, orange] (FA) -- (L);
	\draw [very thick, orange] (FB) -- (L);
	\draw [very thick, orange] (FC) -- (L);
	\draw (-3,2) node {$X_{B''}$}; 
	\draw (3,2) node {$X_{B}$}; 
	\draw (0,4) node {$X_{B'}$}; 
\end{tikzpicture}
\]
\item for any $\F,\G\in\Mns{B,B'}$ a pair of $2$--morphisms
\[
\xymatrix@C=1.5cm{
F_{\F} \ar@<3pt>[r]^-{\DCPA{\G}{\F}} 
& F_{\G} \ar@<3pt>[l]^-{\DCPA{\G}{\F}^{-1}}
}
\]
represented as \emph{fake crossings}
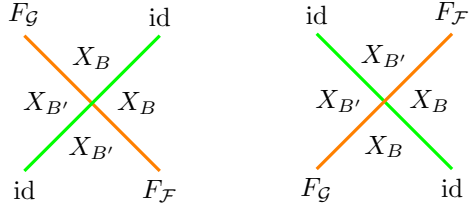
\[
\begin{tikzpicture}[scale=.3]
	\coordinate [label=below:$F_{\F}$]  (FA) at (6,0);
  	\coordinate [label=below:$\id$] (GA) at (0,0);
	\coordinate [label=above:$F_{\G}$] (FB) at (0,6);
	\coordinate [label=above:$\id$] (GB) at (6,6);
	\draw [very thick, orange] (FA) -- (FB);
	\draw [very thick, green] (GA) -- (GB);
	\draw (1,3) node {$X_{B'}$}; 
	\draw (5,3) node {$X_{B}$}; 
	\draw (3,1) node {$X_{B'}$}; 
	\draw (3,5) node {$X_{B}$}; 
\end{tikzpicture}
\qquad\qquad
\begin{tikzpicture}[scale=.3]
	\coordinate  [label=below:$\id$] (FA) at (6,0);
  	\coordinate  [label=below:$F_{\G}$] (GA) at (0,0);
	\coordinate  [label=above:$\id$] (FB) at (0,6);
	\coordinate  [label=above:$F_{\F}$](GB) at (6,6);
	\draw [very thick, green] (FA) -- (FB);
	\draw [very thick, orange] (GA) -- (GB);
	\draw (1,3) node {$X_{B'}$}; 
	\draw (5,3) node {$X_{B}$}; 
	\draw (3,1) node {$X_{B}$}; 
	\draw (3,5) node {$X_{B'}$}; 
\end{tikzpicture}
\]
\item for any $i\in D$, an invertible $1$--morphism $S_i:F_{\{\emptyset,i\}}\to F_{\{\emptyset,i\}}$
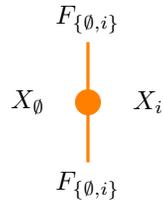
\[
\begin{tikzpicture}[scale=.4]
	\coordinate [label=above: $F_{\{\emptyset,i\}}$] (FA) at (0,4);
	\coordinate [label=below: $F_{\{\emptyset,i\}}$] (FB) at (0,0);
	\draw [very thick, orange] (FA) -- (FB);
	\draw (-2,2) node {$X_{\emptyset}$}; 
	\draw (2,2) node {$X_{i}$}; 
	\node at (0,2) [circle,draw=orange,fill=orange] {};
\end{tikzpicture}
\]
\end{itemize}
satisfying the following relations. To alleviate the notation, the labels of objects and $1$--morphisms
are omitted unless necessary.
\begin{itemize}
\item {\bf Invertibility.}
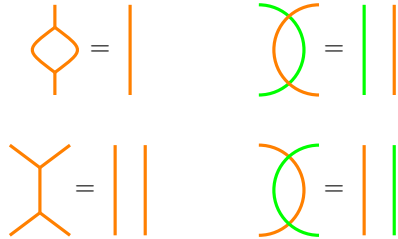
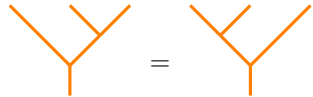
\[
\begin{tabular}{cc}
\begin{tikzpicture}[scale=.2]
	\coordinate (P1) at (0,6);
	\coordinate (P2) at (0,4.5);
	\coordinate (P3) at (0,1.5);
	\coordinate (P4) at (0,0);
	\coordinate (C1) at (-2,3);
	\coordinate (C2) at (2,3);
	\draw [very thick, orange] (P1) -- (P2);
	\draw [very thick, orange] (P2) .. controls (C1) .. (P3);
	\draw [very thick, orange] (P2) .. controls (C2) .. (P3);
	\draw [very thick, orange] (P3) -- (P4);
	\node at (3,3) {$=$};
	\coordinate (P5) at (5,6);
	\coordinate (P6) at (5,0);
	\draw [very thick, orange] (P5) -- (P6);
\end{tikzpicture}
\hspace{1cm}\vspace{0.5cm}
&
\begin{tikzpicture}[scale=.2]
	\coordinate (P1) at (-2,6);
	\coordinate (P2) at (2,6);
	\coordinate (P3) at (-2,0);
	\coordinate (P4) at (2,0);
	\coordinate (C1) at (-2,3);
	\coordinate (C2) at (2,3);
	\draw [very thick, green] (P3) arc (-90:90:3);
	\draw [very thick, orange] (P2) arc (90:270:3);
	\node at (3,3) {$=$};
	\coordinate (P5) at (5,6);
	\coordinate (P6) at (5,0);
	\draw [very thick, green] (P5) -- (P6);
	\coordinate (P7) at (7,6);
	\coordinate (P8) at (7,0);
	\draw [very thick, orange] (P7) -- (P8);
\end{tikzpicture}
\\
\begin{tikzpicture}[scale=.2]
	\coordinate (P1) at (-2,6);
	\coordinate (P2) at (2,6);
	\coordinate (P3) at (0,4.5);
	\coordinate (P4) at (0,1.5);
	\coordinate (P5) at (-2,0);
	\coordinate (P6) at (2,0);
	\draw [very thick, orange] (P1) -- (P3);
	\draw [very thick, orange] (P2) -- (P3);
	\draw [very thick, orange] (P4) -- (P3);
	\draw [very thick, orange] (P5) -- (P4);
	\draw [very thick, orange] (P6) -- (P4);
	\node at (3,3) {$=$};
	\coordinate (P9) at (5,6);
	\coordinate (P10) at (5,0);
	\draw [very thick, orange] (P9) -- (P10);
	\coordinate (P7) at (7,6);
	\coordinate (P8) at (7,0);
	\draw [very thick, orange] (P7) -- (P8);
\end{tikzpicture}
\hspace{1cm}
&
\begin{tikzpicture}[scale=.2]
	\coordinate (P1) at (-2,6);
	\coordinate (P2) at (2,6);
	\coordinate (P3) at (-2,0);
	\coordinate (P4) at (2,0);
	\coordinate (C1) at (-2,3);
	\coordinate (C2) at (2,3);
	\draw [very thick, orange] (P3) arc (-90:90:3);
	\draw [very thick, green] (P2) arc (90:270:3);
	\node at (3,3) {$=$};
	\coordinate (P5) at (5,6);
	\coordinate (P6) at (5,0);
	\draw [very thick, orange] (P5) -- (P6);
	\coordinate (P7) at (7,6);
	\coordinate (P8) at (7,0);
	\draw [very thick, green] (P7) -- (P8);
\end{tikzpicture}
\end{tabular}
\]
\item {\bf Associativity.}
\[
\begin{tabular}{cc}
\begin{tikzpicture}[scale=.2]
	\coordinate (P1) at (-4,8);
	\coordinate (P2) at (0,8);
	\coordinate (P3) at (4,8);
	\coordinate (P4) at (2,6);
	\draw [very thick, orange] (P2) -- (P4);
	\draw [very thick, orange] (P3) -- (P4);
	\coordinate (P5) at (0,4);
	\draw [very thick, orange] (P1) -- (P5);
	\draw [very thick, orange] (P4) -- (P5);
	\coordinate (P6) at (0,2);
	\draw [very thick, orange] (P5) -- (P6);
	\node at (6,4) {$=$};
	\coordinate (R1) at (8,8);
	\coordinate (R2) at (12,8);
	\coordinate (R3) at (16,8);
	\coordinate (R4) at (10,6);
	\draw [very thick, orange] (R2) -- (R4);
	\draw [very thick, orange] (R1) -- (R4);
	\coordinate (R5) at (12,4);
	\draw [very thick, orange] (R3) -- (R5);
	\draw [very thick, orange] (R4) -- (R5);
	\coordinate (R6) at (12,2);
	\draw [very thick, orange] (R5) -- (R6);
\end{tikzpicture}
\end{tabular}
\]
\item {\bf Vertical and horizontal factorisation.}
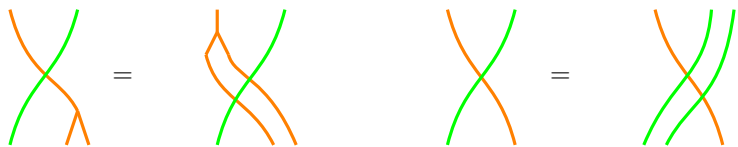
\[
\begin{tabular}{cc}
\begin{tikzpicture}[scale=.15]
	\coordinate (FA) at (6,0);
  	\coordinate (GA) at (0,0);
	\coordinate (FB) at (0,12);
	\coordinate (GB) at (6,12);
	\coordinate (CL) at (1.5,6);
	\coordinate (CR) at (4.5,6);
	\coordinate (G1) at (0,3);
	\coordinate (G2) at (-3,0);
	\coordinate (G3) at (3,0);
	\coordinate (G4) at (12,9);
	\coordinate (G5) at (9,6);
	\coordinate (G6) at (15,6);
	\coordinate (O1) at (6,3);
	\coordinate (O2) at (5,0);
	\coordinate (O3) at (7,0);
	\coordinate (O4) at (0,10);
	\coordinate (O5) at (-1,8);
	\coordinate (O6) at (1,8);
	\coordinate (COL) at (0,4);
	\coordinate (COR) at (3,4);
	\draw [very thick, orange] (FB) .. controls (CL) and (CR) .. (O1);
	\draw [very thick, orange] (O1) -- (O2);
	\draw [very thick, orange] (O1) -- (O3);
	\draw [very thick, green] (GA) .. controls (CL) and (CR) .. (GB);
	\draw (10,6) node {$=$};
\end{tikzpicture}
\hspace{0.6cm}
\begin{tikzpicture}[scale=.15]
	\coordinate (FA) at (6,0);
  	\coordinate (GA) at (0,0);
	\coordinate (FB) at (0,12);
	\coordinate (GB) at (6,12);
	\coordinate (CL) at (1.5,6);
	\coordinate (CR) at (4.5,6);
	\coordinate (G1) at (0,3);
	\coordinate (G2) at (-3,0);
	\coordinate (G3) at (3,0);
	\coordinate (G4) at (12,9);
	\coordinate (G5) at (9,6);
	\coordinate (G6) at (15,6);
	\coordinate (O1) at (6,3);
	\coordinate (O2) at (5,0);
	\coordinate (O3) at (7,0);
	\coordinate (O4) at (0,10);
	\coordinate (O5) at (-1,8);
	\coordinate (O6) at (1,8);
	\coordinate (COL) at (0,4);
	\coordinate (COR) at (3,4);
	\draw [very thick, orange] (FB) -- (O4);
	\draw [very thick, orange] (O4) -- (O5);
	\draw [very thick, orange] (O4) -- (O6);
	\draw [very thick, orange] (O5) .. controls (COL) and (COR) .. (O2);
	\draw [very thick, orange] (O6) .. controls (CL) and (CR) .. (O3);
	\draw [very thick, green] (GA) .. controls (CL) and (CR) .. (GB);
\end{tikzpicture}
\hspace{1.5cm}
&
\begin{tikzpicture}[scale=.15]
	\coordinate (FA) at (6,0);
  	\coordinate (GA) at (0,0);
	\coordinate (FB) at (0,12);
	\coordinate (GB) at (6,12);
	\coordinate (CL) at (1.5,6);
	\coordinate (CR) at (4.5,6);
	\coordinate (G1) at (0,0);
	\coordinate (G2) at (-1,0);
	\coordinate (G3) at (1,0);
	\coordinate (G4) at (6,10);
	\coordinate (G5) at (5,8);
	\coordinate (G6) at (7,8);
	\coordinate (O1) at (6,3);
	\coordinate (O2) at (5,0);
	\coordinate (O3) at (7,0);
	\coordinate (O4) at (0,10);
	\coordinate (O5) at (-1,8);
	\coordinate (O6) at (1,8);
	\coordinate (COL) at (0,4);
	\coordinate (COR) at (3,4);
	\draw [very thick, orange] (FB) .. controls (CL) and (CR) .. (FA);
	\draw [very thick, green] (G1) .. controls (CL) and (CR) .. (GB);
	\draw (10,6) node {$=$};
\end{tikzpicture}
\hspace{0.6cm}
\begin{tikzpicture}[scale=.15]
	\coordinate (FA) at (6,0);
  	\coordinate (GA) at (0,0);
	\coordinate (FB) at (0,12);
	\coordinate (GB) at (6,12);
	\coordinate (CL) at (1.5,6);
	\coordinate (CR) at (4.5,6);
	\coordinate (G1) at (0,3);
	\coordinate (G2) at (-1,0);
	\coordinate (G3) at (1,0);
	\coordinate (G4) at (6,10);
	\coordinate (G5) at (5,12);
	\coordinate (G6) at (7,12);
	\coordinate (O1) at (6,3);
	\coordinate (O2) at (5,0);
	\coordinate (O3) at (7,0);
	\coordinate (O4) at (0,10);
	\coordinate (O5) at (-1,8);
	\coordinate (O6) at (1,8);
	\coordinate (COL) at (0,4);
	\coordinate (COR) at (3,4);
	\coordinate (CGL) at (3,4);
	\coordinate (CGR) at (6,4);
	\draw [very thick, orange] (FB) .. controls (CL) and (CR) .. (FA);
	\draw [very thick, green] (G5) .. controls (CR) and (CL) .. (G2);
	\draw [very thick, green] (G3) .. controls (CGL) and (CGR) .. (G6);
\end{tikzpicture}
\end{tabular}
\]
\item {\bf Braid relations.} For any $i,j\in B\subseteq D$, $\F,\G,\H\in\Mns{B}$ such
that $i\neq j$,  $m_{ij}<\infty$, $\{i\}\in\H$, $\{j\}\in\G$, 
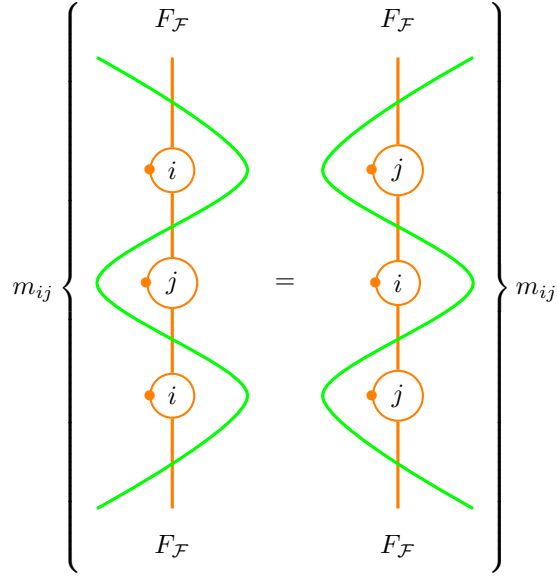
\[
m_{ij}\left\{
\begin{array}{c}
\begin{tikzpicture}[scale=.5]
	\node at (0,-2) {$F_{\F}$};
	\node at (0,12) {$F_{\F}$};
	\node (F1) at (0,2) [shape=circle, draw=orange, thick] {$i$};
	\node at (-0.6,2) {\color{orange}$\bullet$};
	\node (F2) at (0,5) [shape=circle, draw=orange, thick] {$j$};
	\node at (-0.7,5) {\color{orange}$\bullet$};
	\node (F3) at (0,8) [shape=circle, draw=orange, thick] {$i$};
	\node at (-0.6,8) {\color{orange}$\bullet$};
	\coordinate (P2) at (0,-1);
	\coordinate (P3) at (0,11);
	\draw [very thick, orange] (P2) -- (F1.south);
	\draw [very thick, orange] (F1.north) -- (F2.south);
	\draw [very thick, orange] (F2.north) -- (F3.south);
	\draw [very thick, orange] (P3) -- (F3.north);
	\coordinate (V1) at (-2,-1);
	\coordinate (V2) at (2,2);
	\coordinate (V3) at (-2,5);
	\coordinate (V4) at (2,8);
	\coordinate (V5) at (-2,11);	
	%
	%\draw [very thick, green] (V1) -- (V2) -- (V3) -- (V4) -- (V5);
	\draw [very thick, green, xshift=4cm] plot [smooth] coordinates { (V1) (V2) (V3) (V4) (V5)};
	\node (G1) at (6,2) [shape=circle, draw=orange, thick] {$j$};
	\node at (5.3,2) {\color{orange}$\bullet$};
	\node (G2) at (6,5) [shape=circle, draw=orange, thick] {$i$};
	\node at (5.4,5) {\color{orange}$\bullet$};
	\node (G3) at (6,8) [shape=circle, draw=orange, thick] {$j$};
	\node at (5.3,8) {\color{orange}$\bullet$};
	\node at (6,-2) {$F_{\F}$};
	\node at (6,12) {$F_{\F}$};
	\coordinate (Q2) at (6,-1);
	\coordinate (Q3) at (6,11);
	\draw [very thick, orange] (Q2) -- (G1.south);
	\draw [very thick, orange] (G1.north) -- (G2.south);
	\draw [very thick, orange] (G2.north) -- (G3.south);
	\draw [very thick, orange] (Q3) -- (G3.north);
	\coordinate (Z1) at (8,-1);
	\coordinate (Z2) at (4,2);
	\coordinate (Z3) at (8,5);
	\coordinate (Z4) at (4,8);
	\coordinate (Z5) at (8,11);	
	\draw [very thick, green, xshift=4cm] plot [smooth] coordinates { (Z1) (Z2) (Z3) (Z4) (Z5)};
	\node at (3,5) {$=$};
	\end{tikzpicture}
\end{array}
\right\}m_{ij}
\]
\end{itemize}
\subsection{$1$--Morphisms }
Let $X,X'$ be Coxeter objects in $\mathfrak{X}$. We distinguish between them by assigning a different color to
their $2$--cells (specifically, yellow for $X$, gray for $X'$). 
We represent their defining data as
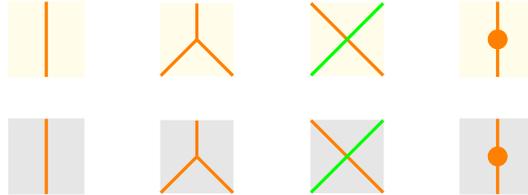
\begin{center}
\begin{tikzpicture}
	\node at (0,0) {
	\begin{tikzpicture}[scale=.25]
	\coordinate (P1) at (0,4);
	\coordinate (P2) at (0,0);
	\draw [very thick, orange] (P1) -- (P2);
	\begin{scope}[on background layer]
    	\draw [yellow!10, fill=yellow!10] (-2,0) rectangle (2,4);
  	\end{scope}
	\end{tikzpicture}
	};
	\node at (2,0) {
	\begin{tikzpicture}[scale=.16]
	\coordinate  (FA) at (-3,-1);
  	\coordinate (FB) at (3,-1);
	\coordinate  (FC) at (0,5);
	\coordinate  (L) at (0,2);
	\draw [very thick, orange] (FA) -- (L);
	\draw [very thick, orange] (FB) -- (L);
	\draw [very thick, orange] (FC) -- (L);
	\begin{scope}[on background layer]
    	\draw [yellow!10, fill=yellow!10] (-3,-1) rectangle (3,5);
  	\end{scope}
	\end{tikzpicture}
	};
	\node at (4,0) {
	\begin{tikzpicture}[scale=.16]
	\coordinate (FA) at (6,0);
  	\coordinate (GA) at (0,0);
	\coordinate (FB) at (0,6);
	\coordinate (GB) at (6,6);
	\draw [very thick, orange] (FA) -- (FB);
	\draw [very thick, green] (GA) -- (GB);
	\begin{scope}[on background layer]
    	\draw [yellow!10, fill=yellow!10] (0,0) rectangle (6,6);
  	\end{scope}
	\end{tikzpicture}
	};
	\node at (6,0) {
	\begin{tikzpicture}[scale=.25]
	\coordinate (FA) at (0,4);
	\coordinate (FB) at (0,0);
	\draw [very thick, orange] (FA) -- (FB);
	\draw [orange, fill=orange] (0,2) circle (0.5);
	\begin{scope}[on background layer]
    	\draw [yellow!10, fill=yellow!10] (-2,0) rectangle (2,4);
  	\end{scope}
	\end{tikzpicture}
	};
\end{tikzpicture}
\end{center}
\begin{center}
\begin{tikzpicture}
	\node at (0,0) {
	\begin{tikzpicture}[scale=.25]
	\coordinate (P1) at (0,4);
	\coordinate (P2) at (0,0);
	\draw [very thick, orange] (P1) -- (P2);
	\begin{scope}[on background layer]
    	\draw [black!10, fill=black!10] (-2,0) rectangle (2,4);
  	\end{scope}
	\end{tikzpicture}
	};
	\node at (2,0) {
	\begin{tikzpicture}[scale=.16]
	\coordinate  (FA) at (-3,-1);
  	\coordinate (FB) at (3,-1);
	\coordinate  (FC) at (0,5);
	\coordinate  (L) at (0,2);
	\draw [very thick, orange] (FA) -- (L);
	\draw [very thick, orange] (FB) -- (L);
	\draw [very thick, orange] (FC) -- (L);
	\begin{scope}[on background layer]
    	\draw [black!10, fill=black!10] (-3,-1) rectangle (3,5);
  	\end{scope}
	\end{tikzpicture}
	};
	\node at (4,0) {
	\begin{tikzpicture}[scale=.16]
	\coordinate (FA) at (6,0);
  	\coordinate (GA) at (0,0);
	\coordinate (FB) at (0,6);
	\coordinate (GB) at (6,6);
	\draw [very thick, orange] (FA) -- (FB);
	\draw [very thick, green] (GA) -- (GB);
	\begin{scope}[on background layer]
    	\draw [black!10, fill=black!10] (0,0) rectangle (6,6);
  	\end{scope}
	\end{tikzpicture}
	};
	\node at (6,0) {
	\begin{tikzpicture}[scale=.25]
	\coordinate (FA) at (0,4);
	\coordinate (FB) at (0,0);
	\draw [very thick, orange] (FA) -- (FB);
	\draw [orange, fill=orange] (0,2) circle (0.5);
	\begin{scope}[on background layer]
    	\draw [black!10, fill=black!10] (-2,0) rectangle (2,4);
  	\end{scope}
	\end{tikzpicture}
	};
\end{tikzpicture}
\end{center}
Then a $1$--morphisms of Coxeter objects $H:X\Rightarrow X'$ is the datum of
\begin{itemize}
\item for any $B\subseteq D$, a $1$--morphism $H_{B}:X_{B}\to X'_{B}$
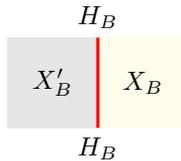
\begin{center}
\begin{tikzpicture}[scale=.3]
	\coordinate [label=above:$H_{B}$]  (FA) at (0,4);
  	\coordinate [label=below:$H_{B}$] (FB) at (0,0);
	\draw [very thick, red] (FA) -- (FB);
	\begin{scope}[on background layer]
	\draw [black!10, fill=black!10] (0,0) rectangle (-4,4);
	\draw [yellow!10, fill=yellow!10] (0,0) rectangle (4,4);
	\end{scope}
	\draw (-2,2) node {$X'_{B}$}; 
	\draw (2,2) node {$X_{B}$}; 
\end{tikzpicture}
\end{center}
\item for any $\F\in\Mns{B,B'}$ a pair of $2$--morphisms
\[
\xymatrix@C=1.5cm{
F'_{\F}\circ H_{B}\ar@<3pt>[r]^-{\gamma_{\F}} 
& H_{B'}\circ F_{\F}  \ar@<3pt>[l]^-{\gamma^{-1}_{\F}}
}
\]
represented as
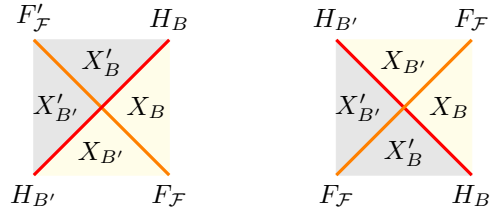
\[
\begin{tikzpicture}[scale=.3]
	\coordinate [label=below:$F_{\F}$]  (FA) at (6,0);
  	\coordinate [label=below:$H_{B'}$] (GA) at (0,0);
	\coordinate [label=above:$F'_{\F}$] (FB) at (0,6);
	\coordinate [label=above:$H_{B}$] (GB) at (6,6);
	\draw [very thick, orange] (FA) -- (FB);
	\draw [very thick, red] (GA) -- (GB);
	\draw (1,3) node {$X'_{B'}$}; 
	\draw (5,3) node {$X_{B}$}; 
	\draw (3,1) node {$X_{B'}$}; 
	\draw (3,5) node {$X'_{B}$}; 
	\begin{scope}[on background layer]
	\draw [black!10, fill=black!10] (GA) -- (GB) -- (FB) -- cycle;
	\draw [yellow!10, fill=yellow!10] (GA) -- (GB) -- (FA) -- cycle;
	\end{scope}
\end{tikzpicture}
\qquad\qquad
\begin{tikzpicture}[scale=.3]
	\coordinate  [label=below:$H_{B}$] (FA) at (6,0);
  	\coordinate  [label=below:$F_{\F}$] (GA) at (0,0);
	\coordinate  [label=above:$H_{B'}$] (FB) at (0,6);
	\coordinate  [label=above:$F_{\F}$](GB) at (6,6);
	\draw [very thick, red] (FA) -- (FB);
	\draw [very thick, orange] (GA) -- (GB);
	\draw (1,3) node {$X'_{B'}$}; 
	\draw (5,3) node {$X_{B}$}; 
	\draw (3,1) node {$X'_{B}$}; 
	\draw (3,5) node {$X_{B'}$}; 
	\begin{scope}[on background layer]
	\draw [yellow!10, fill=yellow!10] (FA) -- (GB) -- (FB) -- cycle;
	\draw [black!10, fill=black!10] (FA) -- (GA) -- (FB) -- cycle;
	\end{scope}
\end{tikzpicture}
\]
\end{itemize}
satisfying the following relations
\begin{itemize}
\item {\bf Invertibility.}
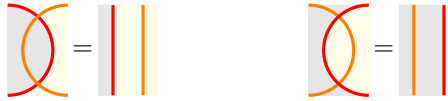
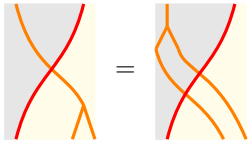
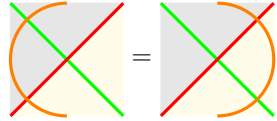
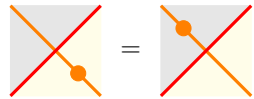
\begin{center}
\begin{tikzpicture}
	\node at (0,0) {
	\begin{tikzpicture}[scale=.2]
	\coordinate (P1) at (-2,6);
	\coordinate (P2) at (2,6);
	\coordinate (P3) at (-2,0);
	\coordinate (P4) at (2,0);
	\coordinate (C1) at (-2,3);
	\coordinate (C2) at (2,3);
	\draw [very thick, red] (P3) arc (-90:90:3);
	\draw [very thick, orange] (P2) arc (90:270:3);
	\node at (3,3) {$=$};
	\coordinate (P5) at (5,6);
	\coordinate (P6) at (5,0);
	\draw [very thick, red] (P5) -- (P6);
	\coordinate (P7) at (7,6);
	\coordinate (P8) at (7,0);
	\draw [very thick, orange] (P7) -- (P8);
	\begin{scope}[on background layer]
	\draw [yellow!10, fill=yellow!10] (P1) -- (P2) -- (P4) -- (P3) arc (-90:90:3) -- cycle;
	\draw [black!10, fill=black!10] (P3) arc (-90:90:3) -- cycle;
	\draw [black!10, fill=black!10] (4,6) -- (P5) -- (P6) -- (4,0);
	\draw [yellow!10, fill=yellow!10] (8,6) -- (P5) -- (P6) -- (8,0);
	\end{scope}
	\end{tikzpicture}
	};
	\node at (4,0) {
	\begin{tikzpicture}[scale=.2]
	\coordinate (P1) at (-2,6);
	\coordinate (P2) at (2,6);
	\coordinate (P3) at (-2,0);
	\coordinate (P4) at (2,0);
	\coordinate (C1) at (-2,3);
	\coordinate (C2) at (2,3);
	\draw [very thick, orange] (P3) arc (-90:90:3);
	\draw [very thick, red] (P2) arc (90:270:3);
	\node at (3,3) {$=$};
	\coordinate (P5) at (5,6);
	\coordinate (P6) at (5,0);
	\draw [very thick, orange] (P5) -- (P6);
	\coordinate (P7) at (7,6);
	\coordinate (P8) at (7,0);
	\draw [very thick, red] (P7) -- (P8);
	\begin{scope}[on background layer]
	\draw [black!10, fill=black!10]  (P1) -- (P2) arc (90:270:3) -- (P3) -- cycle;
	\draw [yellow!10, fill=yellow!10] (P2) arc (90:270:3) -- cycle;
	\draw [black!10, fill=black!10] (4,6) -- (P7) -- (P8) -- (4,0);
	\draw [yellow!10, fill=yellow!10] (8,6) -- (P7) -- (P8) -- (8,0);
	\end{scope}
	\end{tikzpicture}
	};
\end{tikzpicture}
\end{center}
\item {\bf Vertical factorization.}
\begin{center}
\begin{tikzpicture}
	\node at (0,0) {
	\begin{tikzpicture}[scale=.15]
	\coordinate (FA) at (6,0);
  	\coordinate (GA) at (0,0);
	\coordinate (FB) at (0,12);
	\coordinate (GB) at (6,12);
	\coordinate (CL) at (1.5,6);
	\coordinate (CR) at (4.5,6);
	\coordinate (G1) at (0,3);
	\coordinate (G2) at (-3,0);
	\coordinate (G3) at (3,0);
	\coordinate (G4) at (12,9);
	\coordinate (G5) at (9,6);
	\coordinate (G6) at (15,6);
	\coordinate (O1) at (6,3);
	\coordinate (O2) at (5,0);
	\coordinate (O3) at (7,0);
	\coordinate (O4) at (0,10);
	\coordinate (O5) at (-1,8);
	\coordinate (O6) at (1,8);
	\coordinate (COL) at (0,4);
	\coordinate (COR) at (3,4);
	\draw [very thick, orange] (FB) .. controls (CL) and (CR) .. (O1);
	\draw [very thick, orange] (O1) -- (O2);
	\draw [very thick, orange] (O1) -- (O3);
	\draw [very thick, red] (GA) .. controls (CL) and (CR) .. (GB);
	\begin{scope}[on background layer]
	\draw [black!10, fill=black!10]  (-1,12) -- (-1,0) -- (GA) .. controls (CL) and (CR) .. (GB) -- cycle;
	\draw [yellow!10, fill=yellow!10]  (7,0) -- (GA) .. controls (CL) and (CR) .. (GB) -- (7,12) -- cycle;
	\end{scope}
	\end{tikzpicture}
	};
	\draw (1,0) node {$=$};
	\node at (2,0) {
	\begin{tikzpicture}[scale=.15]
	\coordinate (FA) at (6,0);
  	\coordinate (GA) at (0,0);
	\coordinate (FB) at (0,12);
	\coordinate (GB) at (6,12);
	\coordinate (CL) at (1.5,6);
	\coordinate (CR) at (4.5,6);
	\coordinate (G1) at (0,3);
	\coordinate (G2) at (-3,0);
	\coordinate (G3) at (3,0);
	\coordinate (G4) at (12,9);
	\coordinate (G5) at (9,6);
	\coordinate (G6) at (15,6);
	\coordinate (O1) at (6,3);
	\coordinate (O2) at (5,0);
	\coordinate (O3) at (7,0);
	\coordinate (O4) at (0,10);
	\coordinate (O5) at (-1,8);
	\coordinate (O6) at (1,8);
	\coordinate (COL) at (0,4);
	\coordinate (COR) at (3,4);
	\draw [very thick, orange] (FB) -- (O4);
	\draw [very thick, orange] (O4) -- (O5);
	\draw [very thick, orange] (O4) -- (O6);
	\draw [very thick, orange] (O5) .. controls (COL) and (COR) .. (O2);
	\draw [very thick, orange] (O6) .. controls (CL) and (CR) .. (O3);
	\draw [very thick, red] (GA) .. controls (CL) and (CR) .. (GB);
	\begin{scope}[on background layer]
	\draw [black!10, fill=black!10]  (-1,12) -- (-1,0) -- (GA) .. controls (CL) and (CR) .. (GB) -- cycle;
	\draw [yellow!10, fill=yellow!10]  (7,0) -- (GA) .. controls (CL) and (CR) .. (GB) -- (7,12) -- cycle;
	\end{scope}
	\end{tikzpicture}
	};
\end{tikzpicture}
\end{center}
\item {\bf Preserving associators.}
\footnote{The crossings 
\begin{center}
\begin{tikzpicture}
	\node at (6,0) {
	\begin{tikzpicture}[scale=.2]
	\coordinate  (FA) at (6,0);
  	\coordinate (GA) at (0,0);
	\coordinate  (FB) at (0,6);
	\coordinate  (GB) at (6,6);
	\draw [very thick, green] (FA) -- (FB);
	\draw [very thick, red] (GA) -- (GB);
	\begin{scope}[on background layer]
	\draw [black!10, fill=black!10] (GA) -- (GB) -- (FB) -- cycle;
	\draw [yellow!10, fill=yellow!10] (GA) -- (GB) -- (FA) -- cycle;
	\end{scope}
	\end{tikzpicture}
	};
	\node at (8,0) {
	\begin{tikzpicture}[scale=.2]
	\coordinate  (FA) at (6,0);
  	\coordinate (GA) at (0,0);
	\coordinate  (FB) at (0,6);
	\coordinate  (GB) at (6,6);
	\draw [very thick, red] (FA) -- (FB);
	\draw [very thick, green] (GA) -- (GB);
	\begin{scope}[on background layer]
	\draw [black!10, fill=black!10] (GA) -- (GB) -- (FB) -- cycle;
	\draw [yellow!10, fill=yellow!10] (GA) -- (GB) -- (FA) -- cycle;
	\end{scope}
	\end{tikzpicture}
	};
\end{tikzpicture}
\end{center}
represent the identity on $H_B$.
}
\begin{center}
\begin{tikzpicture}
	\node at (-1,0) {
	\begin{tikzpicture}[scale=.25]
	\coordinate (P1) at (-3,3);
	\coordinate (P2) at (0,3);
	\coordinate (P3) at (3,3);
	\coordinate (P4) at (-3,-3);
	\coordinate (P5) at (0,-3);
	\coordinate (P6) at (3,-3);
	\draw [very thick, green] (P1) -- (P6);
	\draw [very thick, red] (P3) -- (P4);
	\draw [very thick, orange] (P2) arc (90:270:3);
	\begin{scope}[on background layer]
	\draw [black!10, fill=black!10] (P3) -- (P4) -- (P1) -- cycle;
	\draw [yellow!10, fill=yellow!10] (P3) -- (P6) -- (P4) -- cycle;
	\end{scope}
	\end{tikzpicture}
	};
	\node at (0,0) {$=$};
	\node at (1,0) {
	\begin{tikzpicture}[scale=.25]
	\coordinate (P1) at (-3,3);
	\coordinate (P2) at (0,3);
	\coordinate (P3) at (3,3);
	\coordinate (P4) at (-3,-3);
	\coordinate (P5) at (0,-3);
	\coordinate (P6) at (3,-3);
	\draw [very thick, green] (P1) -- (P6);
	\draw [very thick, red] (P3) -- (P4);
	\draw [very thick, orange] (P5) arc (-90:90:3);
	\begin{scope}[on background layer]
	\draw [black!10, fill=black!10] (P3) -- (P4) -- (P1) -- cycle;
	\draw [yellow!10, fill=yellow!10] (P3) -- (P6) -- (P4) -- cycle;
	\end{scope}
	\end{tikzpicture}
	};
\end{tikzpicture}
\end{center}
\item {\bf Preserving local monodromies.}
\begin{center}
\begin{tikzpicture}
	\node at (-1,0) {
	\begin{tikzpicture}[scale=.2]
	\coordinate  (FA) at (6,0);
  	\coordinate (GA) at (0,0);
	\coordinate (FB) at (0,6);
	\coordinate  (GB) at (6,6);
	\draw [very thick, orange] (FA) -- (FB);
	\draw [very thick, red] (GA) -- (GB);
	\begin{scope}[on background layer]
	\draw [black!10, fill=black!10] (GA) -- (GB) -- (FB) -- cycle;
	\draw [yellow!10, fill=yellow!10] (GA) -- (GB) -- (FA) -- cycle;
	\draw [orange, fill=orange] (4.5,1.5) circle (0.5);
	\end{scope}
	\end{tikzpicture}
	};
	\node at (0,0) {$=$};
	\node at (1,0) {
	\begin{tikzpicture}[scale=.2]
	\coordinate  (FA) at (6,0);
  	\coordinate (GA) at (0,0);
	\coordinate (FB) at (0,6);
	\coordinate  (GB) at (6,6);
	\draw [very thick, orange] (FA) -- (FB);
	\draw [very thick, red] (GA) -- (GB);
	\begin{scope}[on background layer]
	\draw [black!10, fill=black!10] (GA) -- (GB) -- (FB) -- cycle;
	\draw [yellow!10, fill=yellow!10] (GA) -- (GB) -- (FA) -- cycle;
	\draw [orange, fill=orange] (1.5,4.5) circle (0.5);
	\end{scope}
	\end{tikzpicture}
	};
\end{tikzpicture}
\end{center}
\end{itemize}
\subsection{$2$--Morphisms}
Let $H,H':X\rightarrow X'$ be two $1$--morphisms,
\begin{center}
\begin{tikzpicture}
	\node at (0,0) {
	\begin{tikzpicture}[scale=.3]
	\coordinate (FA) at (0,4);
  	\coordinate (FB) at (0,0);
	\draw [very thick, red] (FA) -- (FB);
	\begin{scope}[on background layer]
	\draw [black!10, fill=black!10] (0,0) rectangle (-2,4);
	\draw [yellow!10, fill=yellow!10] (0,0) rectangle (2,4);
	\end{scope}
	\end{tikzpicture}
	};
	\node at (2,0) {
	\begin{tikzpicture}[scale=.2]
	\coordinate  (FA) at (6,0);
  	\coordinate (GA) at (0,0);
	\coordinate  (FB) at (0,6);
	\coordinate  (GB) at (6,6);
	\draw [very thick, orange] (FA) -- (FB);
	\draw [very thick, red] (GA) -- (GB);
	\begin{scope}[on background layer]
	\draw [black!10, fill=black!10] (GA) -- (GB) -- (FB) -- cycle;
	\draw [yellow!10, fill=yellow!10] (GA) -- (GB) -- (FA) -- cycle;
	\end{scope}
	\end{tikzpicture}
	};
	\node at (4,0) {
	\begin{tikzpicture}[scale=.2]
	\coordinate  (FA) at (6,0);
  	\coordinate (GA) at (0,0);
	\coordinate  (FB) at (0,6);
	\coordinate  (GB) at (6,6);
	\draw [very thick, red] (FA) -- (FB);
	\draw [very thick, orange] (GA) -- (GB);
	\begin{scope}[on background layer]
	\draw [black!10, fill=black!10] (GA) -- (GB) -- (FB) -- cycle;
	\draw [yellow!10, fill=yellow!10] (GA) -- (GB) -- (FA) -- cycle;
	\end{scope}
	\end{tikzpicture}
	};
\end{tikzpicture}
\end{center}
\begin{center}
\begin{tikzpicture}
	\node at (0,0) {
	\begin{tikzpicture}[scale=.3]
	\coordinate (FA) at (0,4);
  	\coordinate (FB) at (0,0);
	\draw [very thick, blue] (FA) -- (FB);
	\begin{scope}[on background layer]
	\draw [black!10, fill=black!10] (0,0) rectangle (-2,4);
	\draw [yellow!10, fill=yellow!10] (0,0) rectangle (2,4);
	\end{scope}
	\end{tikzpicture}
	};
	\node at (2,0) {
	\begin{tikzpicture}[scale=.2]
	\coordinate  (FA) at (6,0);
  	\coordinate (GA) at (0,0);
	\coordinate  (FB) at (0,6);
	\coordinate  (GB) at (6,6);
	\draw [very thick, orange] (FA) -- (FB);
	\draw [very thick, blue] (GA) -- (GB);
	\begin{scope}[on background layer]
	\draw [black!10, fill=black!10] (GA) -- (GB) -- (FB) -- cycle;
	\draw [yellow!10, fill=yellow!10] (GA) -- (GB) -- (FA) -- cycle;
	\end{scope}
	\end{tikzpicture}
	};
	\node at (4,0) {
	\begin{tikzpicture}[scale=.2]
	\coordinate  (FA) at (6,0);
  	\coordinate (GA) at (0,0);
	\coordinate  (FB) at (0,6);
	\coordinate  (GB) at (6,6);
	\draw [very thick, blue] (FA) -- (FB);
	\draw [very thick, orange] (GA) -- (GB);
	\begin{scope}[on background layer]
	\draw [black!10, fill=black!10] (GA) -- (GB) -- (FB) -- cycle;
	\draw [yellow!10, fill=yellow!10] (GA) -- (GB) -- (FA) -- cycle;
	\end{scope}
	\end{tikzpicture}
	};
\end{tikzpicture}
\end{center}
A $2$--morphism $u:H\Rightarrow H'$ is the datum, for any $B\subseteq D$, of an invertible $2$--morphism 
$u_{B}:H_{B}\Rightarrow H'_{B}$
\begin{center}
\begin{tikzpicture}[scale=.3]
	\coordinate (FA) at (0,4);
  	\coordinate (FB) at (0,0);
	\coordinate (L) at (0,2);
	\draw [very thick, blue] (FA) -- (L);
	\draw [very thick, red] (FB) -- (L);
	\draw [purple, fill=violet] (L) circle (0.3);
	\begin{scope}[on background layer]
	\draw [black!10, fill=black!10] (0,0) rectangle (-2,4);
	\draw [yellow!10, fill=yellow!10] (0,0) rectangle (2,4);
	\end{scope}
	\end{tikzpicture}
\end{center}
satisfying 
\begin{center}
\begin{tikzpicture}
	\node at (0,0) {
	\begin{tikzpicture}[scale=.2]
	\coordinate  (FA) at (6,0);
  	\coordinate (GA) at (0,0);
	\coordinate  (FB) at (0,6);
	\coordinate  (GB) at (6,6);
	\coordinate (U) at (1.5,1.5);
	\draw [very thick, orange] (FA) -- (FB);
	\draw [very thick, blue] (GA) -- (U);
	\draw [very thick, red] (GB) -- (U);
	\draw [violet, fill=violet] (U) circle (0.3);
	\begin{scope}[on background layer]
	\draw [black!10, fill=black!10] (GA) -- (GB) -- (FB) -- cycle;
	\draw [yellow!10, fill=yellow!10] (GA) -- (GB) -- (FA) -- cycle;
	\end{scope}
	\end{tikzpicture}
	};
	\node at (1,0) {$=$};
	\node at (2,0) {
	\begin{tikzpicture}[scale=.2]
	\coordinate  (FA) at (6,0);
  	\coordinate (GA) at (0,0);
	\coordinate  (FB) at (0,6);
	\coordinate  (GB) at (6,6);
	%\coordinate (U) at (1.5,1.5);
	%\coordinate (U) at (4.5,1.5);
	\coordinate (U) at (4.5,4.5);
	%\coordinate (U) at (1.5,4.5);
	\draw [very thick, orange] (FA) -- (FB);
	\draw [very thick, blue] (GA) -- (U);
	\draw [very thick, red] (GB) -- (U);
	\draw [violet, fill=violet] (U) circle (0.3);
	\begin{scope}[on background layer]
	\draw [black!10, fill=black!10] (GA) -- (GB) -- (FB) -- cycle;
	\draw [yellow!10, fill=yellow!10] (GA) -- (GB) -- (FA) -- cycle;
	\end{scope}
	\end{tikzpicture}
	};
	\node at (5,0) {
	\begin{tikzpicture}[scale=.2]
	\coordinate  (GA) at (6,0);
  	\coordinate (FA) at (0,0);
	\coordinate  (GB) at (0,6);
	\coordinate  (FB) at (6,6);
	%\coordinate (U) at (1.5,1.5);
	\coordinate (U) at (4.5,1.5);
	%\coordinate (U) at (4.5,4.5);
	%\coordinate (U) at (1.5,4.5);
	\draw [very thick, orange] (FA) -- (FB);
	\draw [very thick, blue] (GA) -- (U);
	\draw [very thick, red] (GB) -- (U);
	\draw [violet, fill=violet] (U) circle (0.3);
	\begin{scope}[on background layer]
	\draw [yellow!10, fill=yellow!10]  (GA) -- (GB) -- (FB) -- cycle;
	\draw [black!10, fill=black!10] (GA) -- (GB) -- (FA) -- cycle;
	\end{scope}
	\end{tikzpicture}
	};
	\node at (6,0) {$=$};
	\node at (7,0) {
	\begin{tikzpicture}[scale=.2]
	\coordinate  (GA) at (6,0);
  	\coordinate (FA) at (0,0);
	\coordinate  (GB) at (0,6);
	\coordinate  (FB) at (6,6);
	\coordinate (U) at (1.5,4.5);
	\draw [very thick, orange] (FA) -- (FB);
	\draw [very thick, blue] (GA) -- (U);
	\draw [very thick, red] (GB) -- (U);
	\draw [violet, fill=violet] (U) circle (0.3);
	\begin{scope}[on background layer]
	\draw[yellow!10, fill=yellow!10]  (GA) -- (GB) -- (FB) -- cycle;
	\draw [black!10, fill=black!10] (GA) -- (GB) -- (FA) -- cycle;
	\end{scope}
	\end{tikzpicture}
	};
\end{tikzpicture}
\end{center}

%---------------------------------------------
% BIBLIOGRAPHY
%--------------------------------------------

\end{document}